\ifx\optionkeymacros\undefined\else\endinput\fi

\catcode`\Œ=\active\defŒ{{\aa}}       
\catcode`\º=\active\defº{\int}        
\catcode`\=\active\def{\c c}        
\catcode`\¶=\active\def¶{\partial}    
\catcode`\Ä=\active\defÄ{\oint}       
\catcode`\Æ=\active\defÆ{\triangle}   
\catcode`\Â=\active\defÂ{\neg}        
\catcode`\µ=\active\defµ{\mu}         
\catcode`\¿=\active\def¿{{\o}}        
\catcode`\¹=\active\def¹{\pi}         
\catcode`\Ï=\active\defÏ{{\oe}}       
\catcode`\§=\active\def§{{\ss}}       
\catcode`\ =\active\def {\dagger}     
\catcode`\Ã=\active\defÃ{\sqrt}       
\catcode`\·=\active\def·{\Sigma}      
\catcode`\Å=\active\defÅ{\approx}     
\catcode`\½=\active\def½{\Omega}      
\catcode`\£=\active\def£{{\it\$}}     
\catcode`\°=\active\def°{\infty}      
\catcode`\¤=\active\def¤{{\S}}        
\catcode`\¦=\active\def¦{{\P}}        
\catcode`\¥=\active\def¥{\bullet}     
\catcode`\»=\active\def»{\leavevmode\raise.585ex\hbox{\b a}}      
\catcode`\¼=\active\def¼{\leavevmode\raise.6ex\hbox{\b o}}        
\catcode`\­=\active\def­{\not=}       
\catcode`\²=\active\def²{\leq}        
\catcode`\³=\active\def³{\geq}        
\catcode`\Ö=\active\defÖ{\div}        
\catcode`\É=\active\defÉ{{\dots}}     
\catcode`\¾=\active\def¾{{\ae}}       
\catcode`\Ç=\active\defÇ{\ll}         
\catcode`\Ò=\active\defÒ{``}          
\catcode`\Á=\active\defÁ{!`}          
\catcode`\¢=\active\def¢{\rlap/c}     
\catcode`\Ô=\active\defÔ{`}           
\catcode`\Õ=\active\defÕ{'}           


\catcode`\=\active\def{{\AA}}       
\catcode`\'=\active\def'{\c C}        
\catcode`\¯=\active\def¯{{\O}}        
\catcode`\¸=\active\def¸{\Pi}         
\catcode`\Î=\active\defÎ{{\OE}}       
\catcode`\®=\active\def®{{\AE}}       
\catcode`\×=\active\def×{\diamond}    
\catcode`\¡=\active\def¡{\accent'27}  
\catcode`\Ó=\active\defÓ{''}          
\catcode`\±=\active\def±{\pm}         
\catcode`\È=\active\defÈ{\gg}         
\catcode`\À=\active\defÀ{?`}          
\catcode`\Ð=\active\defÐ{--}          
\catcode`\Ñ=\active\defÑ{---}         


\catcode`\Š=\active\defŠ{\"a}        
\catcode`\'=\active\def'{\"e}        
\catcode`\•=\active\def•{\"{\i}}     
\catcode`\š=\active\defš{\"o}        
\catcode`\Ÿ=\active\defŸ{\"u}        
\catcode`\Ø=\active\defØ{\"y}        
\catcode`\€=\active\def€{\"A}        
\catcode`\…=\active\def…{\"O}        
\catcode`\†=\active\def†{\"U}        
\catcode`\‡=\active\def‡{\'a}        
\catcode`\Ž=\active\defŽ{\'e}        
\catcode`\'=\active\def'{\'{\i}}     
\catcode`\—=\active\def—{\'o}        
\catcode`\œ=\active\defœ{\'u}        
\catcode`\ƒ=\active\defƒ{\'E}        
\catcode`\ˆ=\active\defˆ{\`a}        
\catcode`\=\active\def{\`e}        
\catcode`\"=\active\def"{\`{\i}}     
\catcode`\˜=\active\def˜{\`o}        
\catcode`\=\active\def{\`u}        
\catcode`\Ë=\active\defË{\`A}        
\catcode`\‹=\active\def‹{\~a}        
\catcode`\–=\active\def–{\~n}        
\catcode`\›=\active\def›{\~o}        
\catcode`\Ì=\active\defÌ{\~A}        
\catcode`\"=\active\def"{\~N}        
\catcode`\Í=\active\defÍ{\~O}        
\catcode`\‰=\active\def‰{\^a}        
\catcode`\=\active\def{\^e}        
\catcode`\"=\active\def"{\^{\i}}     
\catcode`\™=\active\def™{\^o}        
\catcode`\ž=\active\defž{\^u}        

\let\optionkeymacros\null

%
%
%
%
%
\catcode`\@=11
\def\W@{\immediate\write\sixt@@n}
\def\spaces@{\space\space\space\space\space}
\def\spaces@@{\spaces@\spaces@\spaces@}

\voffset=-8mm
\hoffset=-1mm
\hsize=16cm
\vsize=25.8cm

\newdimen\ex@
\newif\iftop@
\newif\ifbot@

\def\thfrac#1#2#3#4#5{\ifmmode{\vbox{\hbox{$#4$}\kern#2\p@}
            \above#1\ex@\vbox{\kern#3\p@\hbox{$#5$}}}\else\relax\fi}
\def\thunderline#1#2#3{\ifmmode\vtop{\ialign{##\crcr
     $\hfil\displaystyle{#3}\hfil$\crcr
     \noalign{\kern#2\p@\nointerlineskip}
     \leaders\hrule height#1\ex@\hfill\crcr}}
     \else$\setbox\z@\hbox{#3} \dp\z@=0pt \thunderline{#1}{#2}{\box\z@}$\fi}
\def\thoverline#1#2#3{\ifmmode\vbox{\ialign{##\crcr
     \leaders\hrule height
     #1\ex@\hfill\crcr\noalign{\kern#2\p@\nointerlineskip}
     $\hfil\displaystyle{#3}\hfil$\crcr}}
     \else$\thoverline{#1}{#2}{\hbox{#3}}$\fi}
\def\thsemiboxed#1#2#3{\ifmmode\setbox\z@\hbox{$\displaystyle{#3}$}
     \hbox{\lower#1\ex@\hbox{\lower#2\p@\hbox{\lower\dp\z@
     \hbox{\vbox{\hbox{\vbox{\box\z@\kern#2\p@}%
     \kern#2\p@\vrule width#1\ex@}\hrule height#1\ex@}}}}}
     \else$\setbox\z@\hbox{#3}\dp\z@=0pt\thsemiboxed{#1}{#2}{\box\z@}$\fi}
\def\thboxed#1#2#3{\ifmmode\setbox\z@\hbox{$\displaystyle{#3}$}
     \hbox{\lower#1\ex@ \hbox{\lower#2\p@\hbox{\lower\dp\z@
     \hbox{\vbox{\hrule height#1\ex@\hbox{\vrule width#1\ex@
     \kern#2\p@\vbox{\kern#2\p@\box\z@\kern#2\p@}%
     \kern#2\p@\vrule width#1\ex@}\hrule height#1\ex@}}}}}
     \else$\setbox\z@\hbox{#3}\dp\z@=0pt\thboxed{#1}{#2}{\box\z@}$\fi}
\def\thsquared#1#2#3{\ifmmode\setbox\z@\hbox{$\displaystyle{#3}$}
       \vcenter{\hrule height#1\ex@\hbox{\vrule width#1\ex@
        \kern#2\p@\vbox{\kern#2\p@\box\z@\kern#2\p@}%
        \kern#2\p@\vrule width#1\ex@}\hrule height#1\ex@}
      \else$\setbox\z@\hbox{#3}\dp\z@=0pt\thsquared{#1}{#2}{\box\z@}$\fi}

\def\NoBlackBoxes{\global\overfullrule\z@}
\def\BlackBoxes{\global\overfullrule5\p@}
\def\,{\relax\ifmmode\mskip\thinmuskip\relax\else\kern.16667em\fi}
\def\!{\relax\ifmmode\mskip-\thinmuskip\relax\else\kern-.16667em\fi}
\def\>{\relax\ifmmode\mskip\medmuskip\relax\else\kern.222222em\fi}
\def\;{\relax\ifmmode\mskip\thickmuskip\relax\else\kern.27777em\fi}
\def\negmedspace{\relax\ifmmode\mskip-\medmuskip\relax\else\kern-.222222em\fi}
\def\negthickspace{\relax\ifmmode\mskip-\thickmuskip\relax\else
 \kern-.27777em\fi}
\def\snug{\unskip\kern-\mathsurround}
\let\thinspace\,
\let\negthinspace\!
\let\medspace\>
\let\thickspace\;

\def\topsmash{\top@true\bot@false\smash@}
\def\botsmash{\top@false\bot@true\smash@}
\def\smash{\top@true\bot@true\smash@}
\def\smash@{\relax\ifmmode\def\next{\mathpalette\mathsm@sh}\else
 \let\next\makesm@sh\fi\next}
\def\finsm@sh{\iftop@\ht\z@\z@\fi\ifbot@\dp\z@\z@\fi\box\z@}

\newfam\msafam  
\newfam\msbfam
\newfam\eufmfam  
\newfam\sffam  
\newfam\bifam

\newdimen\t@ille

\def\taille#1{\t@ille=#1pt
\p@=.1\t@ille
\ex@=.02\t@ille
\font\tenrm=cmr10 at \t@ille                  
\font\sevenrm=cmr7 at .7\t@ille
\font\fiverm=cmr5 at .5\t@ille
\font\teni=cmmi10 at \t@ille                  
\font\seveni=cmmi7 at .7\t@ille
\font\fivei=cmmi5 at .5\t@ille
\font\tensy=cmsy10 at \t@ille                 
\font\sevensy=cmsy7 at .7\t@ille
\font\fivesy=cmsy5 at .5\t@ille
\font\tenex=cmex10 at \t@ille                 
\font\tenbf=cmbx10 at \t@ille  								   			 
\font\sevenbf=cmbx7 at .7\t@ille
\font\fivebf=cmbx5 at .5\t@ille
\font\tentt=cmtt10 at \t@ille                 
\font\tenbi=cmbxti10 at \t@ille
\font\tensl=cmsl10 at \t@ille                 
\font\tenit=cmti10 at \t@ille                 
\font\tensf=cmss10 at \t@ille                 
\font\tenmsa=msam10 at \t@ille
\font\sevenmsa=msam7 at .7\t@ille
\font\fivemsa=msam5 at .5\t@ille
\font\tenmsb=msbm10 at \t@ille                
\font\sevenmsb=msbm7 at .7\t@ille
\font\fivemsb=msbm5 at .5\t@ille
\font\teneufm=eufm10 at \t@ille
\font\seveneufm=eufm7 at .7\t@ille
\font\fiveeufm=eufm5 at .5\t@ille
\textfont0=\tenrm \scriptfont0=\sevenrm \scriptscriptfont0=\fiverm
\textfont1=\teni  \scriptfont1=\seveni  \scriptscriptfont1=\fivei
\textfont2=\tensy \scriptfont2=\sevensy \scriptscriptfont2=\fivesy
\textfont3=\tenex \scriptfont3=\tenex   \scriptscriptfont3=\tenex
\textfont\itfam=\tenit
\textfont\slfam=\tensl
\textfont\bffam=\tenbf   \scriptfont\bffam=\sevenbf
\scriptscriptfont\bffam=\fivebf
\textfont\ttfam=\tentt
\textfont\msafam=\tenmsa \scriptfont\msafam=\sevenmsa
\scriptscriptfont\msafam=\fivemsa
\textfont\msbfam=\tenmsb \scriptfont\msbfam=\sevenmsb
\scriptscriptfont\msbfam=\fivemsb
\textfont\eufmfam=\teneufm \scriptfont\eufmfam=\seveneufm
\scriptscriptfont\eufmfam=\fiveeufm
\textfont\sffam=\tensf
\textfont\bifam=\tenbi
\abovedisplayskip=1.2\t@ille plus .3\t@ille minus .9\t@ille
\abovedisplayshortskip=0pt plus .3\t@ille
\belowdisplayskip=1.2\t@ille plus .3\t@ille minus .9\t@ille
\belowdisplayshortskip=.7\t@ille plus .3\t@ille minus .4\t@ille
\jot=.3\t@ille
\smallskipamount=.3\t@ille plus .1\t@ille minus .1\t@ille
\medskipamount=.6\t@ille plus .2\t@ille minus .2\t@ille
\bigskipamount=1.2\t@ille plus .4\t@ille minus .4\t@ille
\normalbaselineskip=1.2\t@ille
\normallineskip=.1\t@ille
\normallineskiplimit=0pt
\normalbaselines
\tenrm}
\taille{10}

\let\teneuf=\teneufm

\let\euffam=\eufmfam

\def\euf{\fam\euffam\teneuf}
 
\def\bi{\fam\bifam\tenbi}

\def\hexnumber@#1{\ifcase#1 0\or1\or2\or3\or4\or5\or6\or7\or8\or9\or
 A\or B\or C\or D\or E\or F\fi}
\edef\bffam@{\hexnumber@\bffam}
\edef\msa@{\hexnumber@\msafam}
\edef\msb@{\hexnumber@\msbfam}

\mathchardef\arr@bas="0040

\mathchardef\Gamma="0000
\mathchardef\Delta="0001
\mathchardef\Theta="0002
\mathchardef\Lambda="0003
\mathchardef\Xi="0004
\mathchardef\Pi="0005
\mathchardef\Sigma="0006
\mathchardef\Upsilon="0007
\mathchardef\Phi="0008
\mathchardef\Psi="0009
\mathchardef\Omega="000A
\mathchardef\varGamma="0100
\mathchardef\varDelta="0101
\mathchardef\varTheta="0102
\mathchardef\varLambda="0103
\mathchardef\varXi="0104
\mathchardef\varPi="0105
\mathchardef\varSigma="0106
\mathchardef\varUpsilon="0107
\mathchardef\varPhi="0108
\mathchardef\varPsi="0109
\mathchardef\varOmega="010A
\mathchardef\boldGamma="0\bffam@00
\mathchardef\boldDelta="0\bffam@01
\mathchardef\boldTheta="0\bffam@02
\mathchardef\boldLambda="0\bffam@03
\mathchardef\boldXi="0\bffam@04
\mathchardef\boldPi="0\bffam@05
\mathchardef\boldSigma="0\bffam@06
\mathchardef\boldUpsilon="0\bffam@07
\mathchardef\boldPhi="0\bffam@08
\mathchardef\boldPsi="0\bffam@09
\mathchardef\boldOmega="0\bffam@0A

\mathchardef\boxdot="2\msa@00
\mathchardef\boxplus="2\msa@01
\mathchardef\boxtimes="2\msa@02
\mathchardef\square="0\msa@03
\mathchardef\blacksquare="0\msa@04
\mathchardef\centerdot="2\msa@05
\mathchardef\lozenge="0\msa@06
\mathchardef\blacklozenge="0\msa@07
\mathchardef\circlearrowright="3\msa@08
\mathchardef\circlearrowleft="3\msa@09
\mathchardef\rightleftharpoons="3\msa@0A
\mathchardef\leftrightharpoons="3\msa@0B
\mathchardef\boxminus="2\msa@0C
\mathchardef\Vdash="3\msa@0D
\mathchardef\Vvdash="3\msa@0E
\mathchardef\vDash="3\msa@0F
\mathchardef\twoheadrightarrow="3\msa@10
\mathchardef\twoheadleftarrow="3\msa@11
\mathchardef\leftleftarrows="3\msa@12
\mathchardef\rightrightarrows="3\msa@13
\mathchardef\upuparrows="3\msa@14
\mathchardef\downdownarrows="3\msa@15
\mathchardef\upharpoonright="3\msa@16

\mathchardef\downharpoonright="3\msa@17
\mathchardef\upharpoonleft="3\msa@18
\mathchardef\downharpoonleft="3\msa@19
\mathchardef\rightarrowtail="3\msa@1A
\mathchardef\leftarrowtail="3\msa@1B
\mathchardef\leftrightarrows="3\msa@1C
\mathchardef\rightleftarrows="3\msa@1D
\mathchardef\Lsh="3\msa@1E
\mathchardef\Rsh="3\msa@1F
\mathchardef\rightsquigarrow="3\msa@20
\mathchardef\leftrightsquigarrow="3\msa@21
\mathchardef\looparrowleft="3\msa@22
\mathchardef\looparrowright="3\msa@23
\mathchardef\circeq="3\msa@24
\mathchardef\succsim="3\msa@25
\mathchardef\gtrsim="3\msa@26
\mathchardef\gtrapprox="3\msa@27
\mathchardef\multimap="3\msa@28
\mathchardef\therefore="3\msa@29
\mathchardef\because="3\msa@2A
\mathchardef\doteqdot="3\msa@2B

\mathchardef\triangleq="3\msa@2C
\mathchardef\precsim="3\msa@2D
\mathchardef\lesssim="3\msa@2E
\mathchardef\lessapprox="3\msa@2F
\mathchardef\eqslantless="3\msa@30
\mathchardef\eqslantgtr="3\msa@31
\mathchardef\curlyeqprec="3\msa@32
\mathchardef\curlyeqsucc="3\msa@33
\mathchardef\preccurlyeq="3\msa@34
\mathchardef\leqq="3\msa@35
\mathchardef\leqslant="3\msa@36
\mathchardef\lessgtr="3\msa@37
\mathchardef\backprime="0\msa@38
\mathchardef\risingdotseq="3\msa@3A
\mathchardef\fallingdotseq="3\msa@3B
\mathchardef\succcurlyeq="3\msa@3C
\mathchardef\geqq="3\msa@3D
\mathchardef\geqslant="3\msa@3E
\mathchardef\gtrless="3\msa@3F
\mathchardef\sqsubset="3\msa@40
\mathchardef\sqsupset="3\msa@41
\mathchardef\vartriangleright="3\msa@42
\mathchardef\vartriangleleft ="3\msa@43
\mathchardef\trianglerighteq="3\msa@44
\mathchardef\trianglelefteq="3\msa@45
\mathchardef\bigstar="0\msa@46
\mathchardef\between="3\msa@47
\mathchardef\blacktriangledown="0\msa@48
\mathchardef\blacktriangleright="3\msa@49
\mathchardef\blacktriangleleft="3\msa@4A
\mathchardef\vartriangle="0\msa@4D
\mathchardef\blacktriangle="0\msa@4E
\mathchardef\triangledown="0\msa@4F
\mathchardef\eqcirc="3\msa@50
\mathchardef\lesseqgtr="3\msa@51
\mathchardef\gtreqless="3\msa@52
\mathchardef\lesseqqgtr="3\msa@53
\mathchardef\gtreqqless="3\msa@54
\mathchardef\Rrightarrow="3\msa@56
\mathchardef\Lleftarrow="3\msa@57
\mathchardef\veebar="2\msa@59
\mathchardef\barwedge="2\msa@5A
\mathchardef\doublebarwedge="2\msa@5B
\mathchardef\angle="0\msa@5C
\mathchardef\measuredangle="0\msa@5D
\mathchardef\sphericalangle="0\msa@5E
\mathchardef\varpropto="3\msa@5F
\mathchardef\smallsmile="3\msa@60
\mathchardef\smallfrown="3\msa@61
\mathchardef\Subset="3\msa@62
\mathchardef\Supset="3\msa@63
\mathchardef\Cup="2\msa@64

\mathchardef\Cap="2\msa@65

\mathchardef\curlywedge="2\msa@66
\mathchardef\curlyvee="2\msa@67
\mathchardef\leftthreetimes="2\msa@68
\mathchardef\rightthreetimes="2\msa@69
\mathchardef\subseteqq="3\msa@6A
\mathchardef\supseteqq="3\msa@6B
\mathchardef\bumpeq="3\msa@6C
\mathchardef\Bumpeq="3\msa@6D
\mathchardef\lll="3\msa@6E

\mathchardef\ggg="3\msa@6F

\mathchardef\circledS="0\msa@73
\mathchardef\pitchfork="3\msa@74
\mathchardef\dotplus="2\msa@75
\mathchardef\backsim="3\msa@76
\mathchardef\backsimeq="3\msa@77
\mathchardef\complement="0\msa@7B
\mathchardef\intercal="2\msa@7C
\mathchardef\circledcirc="2\msa@7D
\mathchardef\circledast="2\msa@7E
\mathchardef\circleddash="2\msa@7F
\def\ulcorner{\delimiter"4\msa@70\msa@70 }
\def\urcorner{\delimiter"5\msa@71\msa@71 }
\def\llcorner{\delimiter"4\msa@78\msa@78 }
\def\lrcorner{\delimiter"5\msa@79\msa@79 }

\mathchardef\y@n="0\msa@55
\mathchardef\ch@ckmark="0\msa@58
\mathchardef\c@rcledR="0\msa@72
\mathchardef\m@ltese="0\msa@7A
\def\yen{\relax\ifmmode \y@n\else $\m@th\y@n$ \fi}
\def\checkmark{\relax\ifmmode \ch@ckmark\else $\m@th\ch@ckmark$ \fi}
\def\circledR{\relax\ifmmode \c@rcledR\else $\m@th\c@rcledR$ \fi}
\def\maltese{\relax\ifmmode \m@ltese\else $\m@th\m@ltese$ \fi}

\mathchardef\lvertneqq="3\msb@00
\mathchardef\gvertneqq="3\msb@01
\mathchardef\nleq="3\msb@02
\mathchardef\ngeq="3\msb@03
\mathchardef\nless="3\msb@04
\mathchardef\ngtr="3\msb@05
\mathchardef\nprec="3\msb@06
\mathchardef\nsucc="3\msb@07
\mathchardef\lneqq="3\msb@08
\mathchardef\gneqq="3\msb@09
\mathchardef\nleqslant="3\msb@0A
\mathchardef\ngeqslant="3\msb@0B
\mathchardef\lneq="3\msb@0C
\mathchardef\gneq="3\msb@0D
\mathchardef\npreceq="3\msb@0E
\mathchardef\nsucceq="3\msb@0F
\mathchardef\precnsim="3\msb@10
\mathchardef\succnsim="3\msb@11
\mathchardef\lnsim="3\msb@12
\mathchardef\gnsim="3\msb@13
\mathchardef\nleqq="3\msb@14
\mathchardef\ngeqq="3\msb@15
\mathchardef\precneqq="3\msb@16
\mathchardef\succneqq="3\msb@17
\mathchardef\precnapprox="3\msb@18
\mathchardef\succnapprox="3\msb@19
\mathchardef\lnapprox="3\msb@1A
\mathchardef\gnapprox="3\msb@1B
\mathchardef\nsim="3\msb@1C
\mathchardef\napprox="3\msb@1D
\mathchardef\ncong="3\msb@1D
\def\napprox{\not\approx}
\mathchardef\varsubsetneq="3\msb@20
\mathchardef\varsupsetneq="3\msb@21
\mathchardef\nsubseteqq="3\msb@22
\mathchardef\nsupseteqq="3\msb@23
\mathchardef\subsetneqq="3\msb@24
\mathchardef\supsetneqq="3\msb@25
\mathchardef\varsubsetneqq="3\msb@26
\mathchardef\varsupsetneqq="3\msb@27
\mathchardef\subsetneq="3\msb@28
\mathchardef\supsetneq="3\msb@29
\mathchardef\nsubseteq="3\msb@2A
\mathchardef\nsupseteq="3\msb@2B
\mathchardef\nparallel="3\msb@2C
\mathchardef\nmid="3\msb@2D
\mathchardef\nshortmid="3\msb@2E
\mathchardef\nshortparallel="3\msb@2F
\mathchardef\nvdash="3\msb@30
\mathchardef\nVdash="3\msb@31
\mathchardef\nvDash="3\msb@32
\mathchardef\nVDash="3\msb@33
\mathchardef\ntrianglerighteq="3\msb@34
\mathchardef\ntrianglelefteq="3\msb@35
\mathchardef\ntriangleleft="3\msb@36
\mathchardef\ntriangleright="3\msb@37
\mathchardef\nleftarrow="3\msb@38
\mathchardef\nrightarrow="3\msb@39
\mathchardef\nLeftarrow="3\msb@3A
\mathchardef\nRightarrow="3\msb@3B
\mathchardef\nLeftrightarrow="3\msb@3C
\mathchardef\nleftrightarrow="3\msb@3D
\mathchardef\divideontimes="2\msb@3E
\mathchardef\varnothing="0\msb@3F
\mathchardef\nexists="0\msb@40
\mathchardef\mho="0\msb@66
\mathchardef\eth="0\msb@67
\mathchardef\eqsim="3\msb@68
\mathchardef\beth="0\msb@69
\mathchardef\gimel="0\msb@6A
\mathchardef\daleth="0\msb@6B
\mathchardef\lessdot="3\msb@6C
\mathchardef\gtrdot="3\msb@6D
\mathchardef\ltimes="2\msb@6E
\mathchardef\rtimes="2\msb@6F
\mathchardef\shortmid="3\msb@70
\mathchardef\shortparallel="3\msb@71
\mathchardef\smallsetminus="2\msb@72
\mathchardef\thicksim="3\msb@73
\mathchardef\thickapprox="3\msb@74
\mathchardef\approxeq="3\msb@75
\mathchardef\succapprox="3\msb@76
\mathchardef\precapprox="3\msb@77
\mathchardef\curvearrowleft="3\msb@78
\mathchardef\curvearrowright="3\msb@79
\mathchardef\digamma="0\msb@7A
\mathchardef\varkappa="0\msb@7B
\mathchardef\hslash="0\msb@7D
\mathchardef\hbar="0\msb@7E
\mathchardef\backepsilon="3\msb@7F

\mathchardef\bbk="0\msb@7C
\mathchardef\bbA="0\msb@41
\mathchardef\bbB="0\msb@42
\mathchardef\bbC="0\msb@43
\mathchardef\bbD="0\msb@44
\mathchardef\bbE="0\msb@45
\mathchardef\bbF="0\msb@46
\mathchardef\bbG="0\msb@47
\mathchardef\bbH="0\msb@48
\mathchardef\bbI="0\msb@49
\mathchardef\bbJ="0\msb@4A
\mathchardef\bbK="0\msb@4B
\mathchardef\bbL="0\msb@4C
\mathchardef\bbM="0\msb@4D
\mathchardef\bbN="0\msb@4E
\mathchardef\bbO="0\msb@4F
\mathchardef\bbP="0\msb@50
\mathchardef\bbQ="0\msb@51
\mathchardef\bbR="0\msb@52
\mathchardef\bbS="0\msb@53
\mathchardef\bbT="0\msb@54
\mathchardef\bbU="0\msb@55
\mathchardef\bbV="0\msb@56
\mathchardef\bbW="0\msb@57
\mathchardef\bbX="0\msb@58
\mathchardef\bbY="0\msb@59
\mathchardef\bbZ="0\msb@5A

\def\k{\relax\ifmmode\bbk\else $\bbk$\fi}
\def\F{\relax\ifmmode\bbF\else $\bbF$\fi}
\def\N{\relax\ifmmode\bbN\else $\bbN$\fi}
\def\Z{\relax\ifmmode\bbZ\else $\bbZ$\fi}
\def\Q{\relax\ifmmode\bbQ\else $\bbQ$\fi}
\def\R{\relax\ifmmode\bbR\else $\bbR$\fi}
\def\C{\relax\ifmmode\bbC\else $\bbC$\fi}
\def\T{\relax\ifmmode\bbT\else $\bbT$\fi}
\def\arrabas{\relax\ifmmode\arr@bas\else $\m@th\arr@bas$\fi}

\let\B\=
\let\D\.

\edef\msafam@{\hexnumber@\msafam}
\mathchardef\dabar@"0\msafam@39
\def\dashrightarrow{\mathrel{\dabar@\dabar@\mathchar"0\msafam@4B}}
\def\dashleftarrow{\mathrel{\mathchar"0\msafam@4C\dabar@\dabar@}}
\edef\msbfam@{\hexnumber@\msbfam}

\def\widehat#1{\setbox\z@\hbox{$\m@th#1$}%
 \ifdim\wd\z@>\tw@ em\mathaccent"0\msbfam@5B{#1}%
 \else\mathaccent"0362{#1}\fi}
\def\widetilde#1{\setbox\z@\hbox{$\m@th#1$}%
 \ifdim\wd\z@>\tw@ em\mathaccent"0\msbfam@5D{#1}%
 \else\mathaccent"0365{#1}\fi}


\def\newcodes@{\catcode`\\=12 \catcode`\{=12 \catcode`\}=12 \catcode`\#=12
 \catcode`\%=12\relax}
\def\oldcodes@{\catcode`\\=0 \catcode`\{=1 \catcode`\}=2 \catcode`\#=6
 \catcode`\%=14\relax}
\def\comment{\newcodes@\endlinechar=10 \comment@}
{\lccode`\!=`\\
\lowercase{\gdef\comment@#1^^J{\comment@@#1!endcomment\comment@@@}%
\gdef\comment@@#1!endcomment{\futurelet\next\comment@@@}%
\gdef\comment@@@#1\comment@@@{\ifx\next\comment@@@\let
\next=\comment@\else\def\next{\oldcodes@\endlinechar=`\^^M\relax}%
 \fi\next}}}

\font\truetenrm=cmr10
\footline={\hss\truetenrm\folio\hss}

\catcode`\@=12



 \ifx\MYUNDEFINED\BoxedEPSF
   \let\temp\relax
 \else
   \message{}
   \message{ !!! BoxedEPS %
         or BoxedArt macros already defined !!!}
   \let\temp\endinput
 \fi
  \temp
 
 \chardef\EPSFCatAt\the\catcode`\@
 \catcode`\@=11

 \chardef\C@tColon\the\catcode`\:
 \chardef\C@tSemicolon\the\catcode`\;
 \chardef\C@tQmark\the\catcode`\?
 \chardef\C@tEmark\the\catcode`\!
 \chardef\C@tDqt\the\catcode`\"

 \def\PunctOther@{\catcode`\:=12
   \catcode`\;=12 \catcode`\?=12 \catcode`\!=12 \catcode`\"=12}
 \PunctOther@

 \let\wlog@ld\wlog 
 \def\wlog#1{\relax} 

 \newif\ifIN@
 \newdimen\XShift@ \newdimen\YShift@ 
 \newtoks\Realtoks
 
  %
 \newdimen\Wd@ \newdimen\Ht@
 \newdimen\Wd@@ \newdimen\Ht@@
 \newdimen\TT@
 \newdimen\LT@
 \newdimen\BT@
 \newdimen\RT@
 \newdimen\XSlide@ \newdimen\YSlide@ 
 \newdimen\TheScale  
 \newdimen\FigScale  
 \newdimen\ForcedDim@@

 \newtoks\EPSFDirectorytoks@
 \newtoks\EPSFNametoks@
 \newtoks\BdBoxtoks@
 \newtoks\LLXtoks@  
 \newtoks\LLYtoks@

 \newif\ifNotIn@
 \newif\ifForcedDim@
 \newif\ifForceOn@
 \newif\ifForcedHeight@
 \newif\ifPSOrigin

 \newread\EPSFile@ 
 
  \def\ms@g{\immediate\write16}

 \newif\ifIN@\def\IN@{\expandafter\INN@\expandafter}
  \long\def\INN@0#1@#2@{\long\def\NI@##1#1##2##3\ENDNI@
    {\ifx\m@rker##2\IN@false\else\IN@true\fi}%
     \expandafter\NI@#2@@#1\m@rker\ENDNI@}
  \def\m@rker{\m@@rker}

  \newtoks\Initialtoks@  \newtoks\Terminaltoks@
  \def\SPLIT@{\expandafter\SPLITT@\expandafter}
  \def\SPLITT@0#1@#2@{\def\TTILPS@##1#1##2@{%
     \Initialtoks@{##1}\Terminaltoks@{##2}}\expandafter\TTILPS@#2@}


  \newtoks\Trimtoks@

 \def\ForeTrim@{\expandafter\ForeTrim@@\expandafter}
 \def\ForePrim@0 #1@{\Trimtoks@{#1}}
 \def\ForeTrim@@0#1@{\IN@0\m@rker. @\m@rker.#1@%
     \ifIN@\ForePrim@0#1@%
     \else\Trimtoks@\expandafter{#1}\fi}

  \def\Trim@0#1@{%
      \ForeTrim@0#1@%
      \IN@0 @\the\Trimtoks@ @%
        \ifIN@ 
             \SPLIT@0 @\the\Trimtoks@ @\Trimtoks@\Initialtoks@
             \IN@0\the\Terminaltoks@ @ @%
                 \ifIN@
                 \else \Trimtoks@ {FigNameWithSpace}%
                 \fi
        \fi
      }


   \newtoks\pt@ks
   \def \getpt@ks 0.0#1@{\pt@ks{#1}}
   \dimen0=0pt\relax\expandafter\getpt@ks\the\dimen0@

  \newtoks\Realtoks
  \def\Real#1{%
    \dimen2=#1%
      \SPLIT@0\the\pt@ks @\the\dimen2@
       \Realtoks=\Initialtoks@
            }

   \newdimen\Product
   \def\Mult#1#2{%
     \dimen4=#1\relax
     \dimen6=#2%
     \Real{\dimen4}%
     \Product=\the\Realtoks\dimen6%
        }

 \newdimen\Inverse
 \newdimen\hmxdim@ \hmxdim@=8192pt
 \def\Invert#1{%
  \Inverse=\hmxdim@
  \dimen0=#1%
  \divide\Inverse \dimen0%
  \multiply\Inverse 8}

   \def\Rescale#1#2#3{
              \divide #1 by 100\relax
              \dimen2=#3\divide\dimen2 by 100 \Invert{\dimen2}%
              \Mult{#1}{#2}%
              \Mult\Product\Inverse 
              #1=\Product}

  \def\Scale#1{\dimen0=\TheScale %
      \divide #1 by  1280 
      \divide \dimen0 by 5120 %
      \multiply#1 by \dimen0 
      \divide#1 by 10   
     }
 

 \newbox\scrunchbox

 \def\Scrunched#1{{\setbox\scrunchbox\hbox{#1}%
   \wd\scrunchbox=0pt
   \ht\scrunchbox=0pt
   \dp\scrunchbox=0pt
   \box\scrunchbox}}

 \def\Shifted@#1{%
   \vbox {\kern-\YShift@
       \hbox {\kern\XShift@\hbox{#1}\kern-\XShift@}%
           \kern\YShift@}}


 \def\cBoxedEPSF#1{{\leavevmode 
   \ReadNameAndScale@{#1}%
   \SetEPSFSpec@
   \ReadEPSFile@ \ReadBdB@x  
     \TrimFigDims@ 
     \CalculateFigScale@  
     \ScaleFigDims@
     \SetInkShift@
   \hbox{$\mathsurround=0pt\relax
         \vcenter{\hbox{%
             \FrameSpider{\hskip-.4pt\vrule}%
             \vbox to \Ht@{\offinterlineskip\parindent=\z@%
                \FrameSpider{\vskip-.4pt\hrule}\vfil 
                \hbox to \Wd@{\hfil}%
                \vfil
                \InkShift@{\EPSFSpecial{\EPSFSpec@}{\FigSc@leReal}}%
             \FrameSpider{\hrule\vskip-.4pt}}%
         \FrameSpider{\vrule\hskip-.4pt}}}%
     $}%
    \CleanRegisters@ 
    \ms@g{ *** Box composed for the %
         EPSF file \the\EPSFNametoks@}%
    }}
 
 \def\tBoxedEPSF#1{\setbox4\hbox{\cBoxedEPSF{#1}}%
     \setbox4\hbox{\raise -\ht4 \hbox{\box4}}%
     \box4
      }

 \def\bBoxedEPSF#1{\setbox4\hbox{\cBoxedEPSF{#1}}%
     \setbox4\hbox{\raise \dp4 \hbox{\box4}}%
     \box4
      }

  \let\BoxedEPSF\cBoxedEPSF

   %

   %
  \def\gLinefigure[#1scaled#2]_#3{%
        \BoxedEPSF{#3 scaled #2}}
    
   %

  \def\EPSFxsize{\afterassignment\ForceW@\ForcedDim@@}
      \def\ForceW@{\ForcedDim@true\ForcedHeight@false}
  
  \def\EPSFysize{\afterassignment\ForceH@\ForcedDim@@}
      \def\ForceH@{\ForcedDim@true\ForcedHeight@true}

  \def\EmulateRokicki{%
       \let\epsfbox\bBoxedEPSF \let\epsffile\bBoxedEPSF
       \let\epsfxsize\EPSFxsize \let\epsfysize\EPSFysize} 
 
  %
 \def\ReadNameAndScale@#1{\IN@0 scaled@#1@
   \ifIN@\ReadNameAndScale@@0#1@%
   \else \ReadNameAndScale@@0#1 scaled\DefaultMilScale @%
   \fi}
  
 \def\ReadNameAndScale@@0#1scaled#2@{
    \let\OldBackslash@\\%
    \def\\{\OtherB@ckslash}%
    \edef\temp@{#1}%
    \Trim@0\temp@ @%
    \EPSFNametoks@\expandafter{\the\Trimtoks@ }%
    \FigScale=#2 pt%
    \let\\\OldBackslash@
    }
 
 \def\SetDefaultEPSFScale#1{%
      \global\def\DefaultMilScale{#1}}

 \SetDefaultEPSFScale{1000}

  %
 \def \SetBogusBbox@{%
     \global\BdBoxtoks@{ BoundingBox:0 0 100 100 }%
     \global\def\BdBoxLine@{ BoundingBox:0 0 100 100 }%
     \ms@g{ !!! Will use placeholder !!!}%
     }

 {\catcode`\%=12\gdef\P@S@{

 \def\ReadEPSFile@{
     \openin\EPSFile@\EPSFSpec@
     \relax  
  \ifeof\EPSFile@
     \ms@g{}%
     \ms@g{ !!! EPS FILE \the\EPSFDirectorytoks@
       \the\EPSFNametoks@\space WAS NOT FOUND !!!}%
     \SetBogusBbox@
  \else
   \begingroup
   \catcode`\%=12\catcode`\:=12\catcode`\!=12
   \catcode`\G=14\catcode`\\=14\relax
   \global\read\EPSFile@ to \BdBoxLine@
   \IN@0\P@S@ @\BdBoxLine@ @%
   \ifIN@ 
     \NotIn@true
     \loop   
       \ifeof\EPSFile@\NotIn@false 
         \ms@g{}%
         \ms@g{ !!! BoundingBox NOT FOUND IN %
            \the\EPSFDirectorytoks@\the\EPSFNametoks@\space!!! }%
         \SetBogusBbox@
       \else\global\read\EPSFile@ to \BdBoxLine@
       \fi
       \global\BdBoxtoks@\expandafter{\BdBoxLine@}%
       \IN@0BoundingBox:@\the\BdBoxtoks@ @%
       \ifIN@\NotIn@false\fi%
     \ifNotIn@\repeat
   \else
         \ms@g{}%
         \ms@g{ !!! \the\EPSFNametoks@\space not PS!\space !!!}%
         \SetBogusBbox@
   \fi
  \endgroup\relax
  \fi
  \closein\EPSFile@ 
   }

  \def\ReadBdB@x{
   \expandafter\ReadBdB@x@\the\BdBoxtoks@ @}
  
  \def\ReadBdB@x@#1BoundingBox:#2@{
    \ForeTrim@0#2@%
    \IN@0atend@\the\Trimtoks@ @%
       \ifIN@\Trimtoks@={0 0 100 100 }%
         \ms@g{}%
         \ms@g{ !!! BoundingBox not found in %
         \the\EPSFDirectorytoks@\the\EPSFNametoks@\space !!!}%
         \ms@g{ !!! It must not be at end of EPSF !!!}%
         \ms@g{ !!! Will use placeholder !!!}%
       \fi
    \expandafter\ReadBdB@x@@\the\Trimtoks@ @%
   }
    
  \def\ReadBdB@x@@#1 #2 #3 #4@{
      \Wd@=#3bp\advance\Wd@ by -#1bp%
      \Ht@=#4bp\advance\Ht@ by-#2bp%
       \Wd@@=\Wd@ \Ht@@=\Ht@ 
       \LLXtoks@={#1}\LLYtoks@={#2}
      \ifPSOrigin\XShift@=-#1bp\YShift@=-#2bp\fi 
     }

   %
   \def\G@bbl@#1{}
   \bgroup
     \global\edef\OtherB@ckslash{\expandafter\G@bbl@\string\\}
   \egroup

  \def\SetEPSFDirectory{
           \bgroup\PunctOther@\relax
           \let\\\OtherB@ckslash
           \SetEPSFDirectory@}

 \def\SetEPSFDirectory@#1{
    \edef\temp@{#1}%
    \Trim@0\temp@ @
    \global\toks1\expandafter{\the\Trimtoks@ }\relax
    \egroup
    \EPSFDirectorytoks@=\toks1
    }

 \def\SetEPSFSpec@{%
     \bgroup
     \let\\=\OtherB@ckslash
     \global\edef\EPSFSpec@{%
        \the\EPSFDirectorytoks@\the\EPSFNametoks@}%
     \global\edef\EPSFSpec@{\EPSFSpec@}%
     \egroup}

  %
 \def\TrimTop#1{\advance\TT@ by #1}
 \def\TrimLeft#1{\advance\LT@ by #1}
 \def\TrimBottom#1{\advance\BT@ by #1}
 \def\TrimRight#1{\advance\RT@ by #1}

 \def\TrimBoundingBox#1{%
   \TrimTop{#1}%
   \TrimLeft{#1}%
   \TrimBottom{#1}%
   \TrimRight{#1}%
       }

 \def\TrimFigDims@{%
    \advance\Wd@ by -\LT@ 
    \advance\Wd@ by -\RT@ \RT@=\z@
    \advance\Ht@ by -\TT@ \TT@=\z@
    \advance\Ht@ by -\BT@ 
    }

  %
  \def\ForceWidth#1{\ForcedDim@true
       \ForcedDim@@#1\ForcedHeight@false}
  
  \def\ForceHeight#1{\ForcedDim@true
       \ForcedDim@@=#1\ForcedHeight@true}

  \def\ForceOn{\ForceOn@true}
  \def\ForceOff{\ForceOn@false\ForcedDim@false}
  
  \def\CalculateFigScale@{%
     \ifForcedDim@\FigScale=1000pt
           \ifForcedHeight@
                \Rescale\FigScale\ForcedDim@@\Ht@
           \else
                \Rescale\FigScale\ForcedDim@@\Wd@
           \fi
     \fi
     \Real{\FigScale}%
     \edef\FigSc@leReal{\the\Realtoks}%
     }
   
  \def\ScaleFigDims@{\TheScale=\FigScale
      \ifForcedDim@
           \ifForcedHeight@ \Ht@=\ForcedDim@@  \Scale\Wd@
           \else \Wd@=\ForcedDim@@ \Scale\Ht@
           \fi
      \else \Scale\Wd@\Scale\Ht@        
      \fi
      \ifForceOn@\relax\else\global\ForcedDim@false\fi
      \Scale\LT@\Scale\BT@  
      \Scale\XShift@\Scale\YShift@
      }
      
 \def\HideReservedBoxes{\global\def\FrameSpider##1{\null}}
 \def\ShowReservedBoxes{\global\def\FrameSpider##1{##1}}
 \let\HideDisplacementBoxes\HideReservedBoxes  
 \let\ShowDisplacementBoxes\ShowReservedBoxes
 \let\HideFigureFrames\HideReservedBoxes
 \let\ShowFigureFrames\ShowReservedBoxes
  \ShowDisplacementBoxes
 
 \def\hSlide#1{\advance\XSlide@ by #1}
 \def\vSlide#1{\advance\YSlide@ by #1}
 
  \def\SetInkShift@{%
            \advance\XShift@ by -\LT@
            \advance\XShift@ by \XSlide@
            \advance\YShift@ by -\BT@
            \advance\YShift@ by -\YSlide@
             }
  \def\InkShift@#1{\Shifted@{\Scrunched{#1}}}
 
   %
  \def\CleanRegisters@{%
      \globaldefs=1\relax
        \XShift@=\z@\YShift@=\z@\XSlide@=\z@\YSlide@=\z@
        \TT@=\z@\LT@=\z@\BT@=\z@\RT@=\z@
      \globaldefs=0\relax}

 
 \def\SetTexturesEPSFSpecial{\PSOriginfalse
  \gdef\EPSFSpecial##1##2{\relax
    \edef\specialthis{##2}%
    \SPLIT@0.@\specialthis.@\relax
    \special{illustration ##1 scaled
                        \the\Initialtoks@}}}
 
  \def\SetUnixCoopEPSFSpecial{\PSOrigintrue 
   \gdef\EPSFSpecial##1##2{%
      \dimen4=##2pt
      \divide\dimen4 by 1000\relax
      \Real{\dimen4}
      \edef\Aux@{\the\Realtoks}%
      \includegraphics{##1\space}}}

  \def\SetBechtolsheimEPSFSpecial@{
   \PSOrigintrue
   \special{\DriverTag@ Include0 "psfig.pro"}%
   \gdef\EPSFSpecial##1##2{%
      \dimen4=##2pt 
      \divide\dimen4 by 1000\relax
      \Real{\dimen4} 
      \edef\Aux@{\the\Realtoks}
      \special{\DriverTag@ Literal "10 10 0 0 10 10 startTexFig
           \the\mag\space 1000 div 3.25 neg mul 
           \the\mag\space 1000 div .25 neg mul translate 
           \the\mag\space 1000 div \Aux@\space mul 
           \the\mag\space 1000 div \Aux@\space mul scale "}%
      \special{\DriverTag@ Include1 "##1"}%
      \special{\DriverTag@ Literal "endTexFig "}%
        }}

  \def\SetBechtolsheimEPSFSpecial@{
   \PSOrigintrue
   \special{\DriverTag@ Include0 "psfig.pro"}%
   \gdef\EPSFSpecial##1##2{%
      \dimen4=##2pt 
      \divide\dimen4 by 1000\relax
      \Real{\dimen4} 
      \edef\Aux@{\the\Realtoks}
      \special{\DriverTag@ Literal "10 10 0 0 10 10 startTexFig
           \the\mag\space 1000 div 
           dup 3.25 neg mul 2 index .25 neg mul translate 
           \Aux@\space mul dup scale "}%
      \special{\DriverTag@ Include1 "##1"}%
      \special{\DriverTag@ Literal "endTexFig "}%
        }}

  \def\SetBechtolsheimDVITPSEPSFSpecial{\def\DriverTag@{dvitps: }%
      \SetBechtolsheimEPSFSpecial@}

  \def\SetBechtolsheimDVI2PSEPSFSSpecial{\def\DriverTag@{DVI2PS: }%
      \SetBechtolsheimEPSFSpecial@}

  \def\SetLisEPSFSpecial{\PSOrigintrue 
   \gdef\EPSFSpecial##1##2{%
      \dimen4=##2pt
      \divide\dimen4 by 1000\relax
      \Real{\dimen4}
      \edef\Aux@{\the\Realtoks}%
      \special{pstext="10 10 0 0 10 10 startTexFig\space
           \the\mag\space 1000 div \Aux@\space mul 
           \the\mag\space 1000 div \Aux@\space mul scale"}%
      \includegraphics{##1}%
      \special{pstext=endTexFig}%
        }}

  \def\SetRokickiEPSFSpecial{\PSOrigintrue 
   \gdef\EPSFSpecial##1##2{%
      \dimen4=##2pt
      \divide\dimen4 by 10\relax
      \Real{\dimen4}
      \edef\Aux@{\the\Realtoks}%
      \includegraphics{##1}}}

  \def\SetInlineRokickiEPSFSpecial{\PSOrigintrue 
   \gdef\EPSFSpecial##1##2{%
      \dimen4=##2pt
      \divide\dimen4 by 1000\relax
      \Real{\dimen4}
      \edef\Aux@{\the\Realtoks}%
      \special{ps::[begin] 10 10 0 0 10 10 startTexFig\space
           \the\mag\space 1000 div \Aux@\space mul 
           \the\mag\space 1000 div \Aux@\space mul scale}%
      \special{ps: plotfile ##1}%
      \special{ps::[end] endTexFig}%
        }}

 \def\SetOzTeXEPSFSpecial{\PSOrigintrue
 \gdef\EPSFSpecial##1##2{%
 \dimen4=##2pt
 \divide\dimen4 by 1000\relax
 \Real{\dimen4}
 \edef\Aux@{\the\Realtoks}
 \special{epsf=\string"##1\string"\space scale=\Aux@}%
 }} 

  \def\SetPSprintEPSFSpecial{\PSOriginFALSE 
   \gdef\EPSFSpecial##1##2{
     \special{##1\space 
       ##2 1000 div \the\mag\space 1000 div mul
       ##2 1000 div \the\mag\space 1000 div mul scale
       \the\LLXtoks@\space neg \the\LLYtoks@\space neg translate
       }}}

 \def\SetArborEPSFSpecial{\PSOriginfalse 
   \gdef\EPSFSpecial##1##2{%
     \edef\specialthis{##2}%
     \SPLIT@0.@\specialthis.@\relax 
     \special{ps: epsfile ##1\space \the\Initialtoks@}}}

 \def\SetClarkEPSFSpecial{\PSOriginfalse 
   \gdef\EPSFSpecial##1##2{%
     \Rescale {\Wd@@}{##2pt}{1000pt}%
     \Rescale {\Ht@@}{##2pt}{1000pt}%
     \special{dvitops: import 
           ##1\space\the\Wd@@\space\the\Ht@@}}}

  \let\SetDVIPSONEEPSFSpecial\SetUnixCoopEPSFSpecial
  \let\SetDVIPSoneEPSFSpecial\SetUnixCoopEPSFSpecial

  \def\SetBeebeEPSFSpecial{
   \PSOriginfalse%
   \gdef\EPSFSpecial##1##2{\relax
    \special{language "PS",
      literal "##2 1000 div ##2 1000 div scale",
      position = "bottom left",
      include "##1"}}}
  \let\SetDVIALWEPSFSpecial\SetBeebeEPSFSpecial

  \def\SetNorthlakeEPSFSpecial{\PSOrigintrue
   \gdef\EPSFSpecial##1##2{%
     \edef\specialthis{##2}%
     \SPLIT@0.@\specialthis.@\relax 
     \special{insert ##1,magnification=\the\Initialtoks@}}}

 \def\SetStandardEPSFSpecial{%
   \gdef\EPSFSpecial##1##2{%
     \ms@g{}
     \ms@g{%
       !!! Sorry! There is still no standard for \string%
       \special\space EPSF integration !!!}%
     \ms@g{%
      --- So you will have to identify your driver using a command}%
     \ms@g{%
      --- of the form \string\Set...EPSFSpecial, in order to get}%
     \ms@g{%
      --- your graphics to print.  See BoxedEPS.doc.}%
     \ms@g{}
     \gdef\EPSFSpecial####1####2{}
     }}

  \SetStandardEPSFSpecial 
 
 \let\wlog\wlog@ld 

 \catcode`\:=\C@tColon
 \catcode`\;=\C@tSemicolon
 \catcode`\?=\C@tQmark
 \catcode`\!=\C@tEmark
 \catcode`\"=\C@tDqt

 \catcode`\@=\EPSFCatAt


 %
 %
 %
 %
 %

\SetTexturesEPSFSpecial
\HideDisplacementBoxes
\def\english{\language=0\lccode`\'=0}
\def\french{\language=1\lccode`\'=`\'}
\french

\font\bi=msbm10 at 16pt
\font\doub=msbm10 at 10pt
\font\para=cmbx12 at 18pt
\font\titre=cmbx12 at 16pt
\font\soustitre=cmbx12 at 13pt
\font\bib=cmmi12 at 13pt
\font\papa=cmsy10 at 16pt
\font\toto=msam10 at 20pt
\font\tot=msbm10 at 10pt
\font\tut=cmr12 at 16pt
\font\couri=Courier at 10pt
\font\couria=Courier at 7pt
\font\gdcm=cmr12 at 20pt
\NoBlackBoxes
\newfam\regfam  




\def\OO#1{{\reg O}\ifmmode _{#1}\else $_{#1}$\fi}
\def\id{{\rm id}} 
\let \lra=\longrightarrow
\def\petitsaut{\smallbreak}
\def\moyensaut{\medbreak}
\def\grandsaut{\bigbreak}
\def\hra{\lhook\joinrel\lra}
\def\strot{\mathchoice{\vrule width-1pt height 8.5pt depth 3.5pt}%
   {\vrule width-1pt height 8.5pt depth 3.5pt}{\vrule width-0.7pt height 5.95pt depth 2.45pt}%
   {\vrule width-0.5pt height 4.25pt depth 1.75pt}}

\def\retr#1,{\leavevmode\kern-5mm{\bf #1}}
\def\prop {\retr Proposition {\the\chapnomb}.{\the\parnomb}.{\the\nomb},\global\advance\nomb by 1\par\nobreak\smallskip}
\def\th {\retr ThŽorme {\the\chapnomb}.{\the\parnomb}.{\the\nomb},\global\advance\nomb by 1\par\nobreak\smallskip}
\def\defi {\retr DŽfinition {\the\chapnomb}.{\the\parnomb}.{\the\nomb},\global\advance\nomb by 1\par\nobreak\smallskip}
\def\defis {\retr DŽfinitions {\the\chapnomb}.{\the\parnomb}.{\the\nomb},\global\advance\nomb by 1\par\nobreak\smallskip}
\def\coro {\retr Corollaire {\the\chapnomb}.{\the\parnomb}.{\the\nomb},\global\advance\nomb by 1\par\nobreak\smallskip}
\def\lem {\retr Lemme {\the\chapnomb}.{\the\parnomb}.{\the\nomb},\global\advance\nomb by 1\par\nobreak\smallskip}

\def\propa {\retr Proposition A{\the\chapnomb}.{\the\nomb},\global\advance\nomb by 1\par\nobreak\smallskip}
\def\tha {\retr ThŽorme A{\the\chapnomb}.{\the\nomb},\global\advance\nomb by 1\par\nobreak\smallskip}
\def\defia {\retr DŽfinition A{\the\chapnomb}.{\the\nomb},\global\advance\nomb by 1\par\nobreak\smallskip}
\def\defisa {\retr DŽfinitions A{\the\chapnomb}.{\the\nomb},\global\advance\nomb by 1\par\nobreak\smallskip}
\def\coroa {\retr Corollaire A{\the\chapnomb}.{\the\nomb},\global\advance\nomb by 1\par\nobreak\smallskip}
\def\lema {\retr Lemme A{\the\chapnomb}.{\the\nomb},\global\advance\nomb by 1\par\nobreak\smallskip}

\def\rems{\retr Remarques :,\par\nobreak\smallskip}
\def\rem{\retr Remarque :,\par\nobreak\smallskip}
\def\exs{\retr Exemples :,\par\nobreak\smallskip}
\def\ex{\retr Exemple :,\par\nobreak\smallskip}
\def\theo {\retr ThŽorme {\the\nomb},\global\advance\nomb by 1\par\nobreak\smallskip}
\def\defin {\retr DŽfinition {\the\nomb},\global\advance\nomb by 1\par\nobreak\smallskip}
\def\lemm {\retr Lemme {\the\nomb},\global\advance\nomb by 1\par\nobreak\smallskip}
\def\corol {\retr Corollaire {\the\nomb},\global\advance\nomb by 1\par\nobreak\smallskip}

\def\demo{\retr DŽmonstration :,\par\nobreak\smallskip}
\def\dst{\displaystyle}
\def\findemol{\hbox{\di ,}}
\def\findemo{\hbox{\di \char'144}}
\def\findemon{\hfill\hbox{\special{illustration pio}}\grandsaut}
\def\cucuplus{\mathrel{\lower .9pt\rlap{$\square$}+}}
\def\caplus{\boxplus}
\def\<{\mathopen<}
\def\>{\mathclose>}
\def\P{\hbox{\doub P}}
\def\rom#1{\uppercase\expandafter{\romannumeral #1}} 
\long\def\art#1{{\parindent0pt\item{#1}}\hangindent=7mm\hangafter=-20}
\long\def\artart#1{{\parindent0pt\item{#1}}\hangindent=12mm\hangafter=-20}
\def\tr#1{\,{\vphantom{#1}}^t\!#1}
\def\momo#1\over#2{\mathrel{\mathop{#2}\limits^{\lower5pt\hbox{$\scriptstyle #1$}}}}
\def\isom{\momo\sim\over\lra}

\smallskipamount4pt plus1pt minus1pt
\medskipamount8pt plus2pt minus2pt
\bigskipamount16pt plus4pt minus4pt
\baselineskip16pt

\catcode`\Ê=\active
\defÊ{papa g‰teaux!}

\def\dfleche#1#2{\smash{\mathop{\hbox to 9mm{\rightarrowfill}}\limits^{#1}_{#2}}}
\def\bfleche#1#2{\llap{$\vcenter{\hbox{$\scriptstyle#1$}}$}\left\downarrow\vbox to 4.5mm{}\right.
  \rlap{$\vcenter{\hbox{$\!\scriptstyle#2$}}$}}
\def\gfleche#1#2{\smash{\mathop{\hbox to 9mm{\leftarrowfill}}\limits^{\ #1}_{\ #2}}}
\def\hfleche#1#2{\llap{$\vcenter{\hbox{$\scriptstyle#1$}}$}\left\uparrow\vbox to 4.5mm{}\right.
  \rlap{$\vcenter{\hbox{$\!\scriptstyle#2$}}$}} 
 
\def\iso{\buildrel \sim \over{\vphantom{.}\smash\longrightarrow}}

\def\diagram#1{{\def\normalbaselines{\baselineskip18pt \lineskip3pt \lineskiplimit3pt} \matrix{#1}}}

\def\carre#1,#2,#3,#4;#5,#6,#7,#8.{\diagram{#1&\dfleche{#5}{}&#2\cr
  \bfleche{#6}{}&&\bfleche{}{#7}\cr
  #3&\dfleche{#8}{}&#4\cr}}
\def\cartesien#1,#2,#3;#4,#5,#6,#7.{\carre{#2\times_#3 #1},#1,#2,#3;
  #4,#5,#6,#7.}
\def\invcarre#1,#2,#3,#4;#5,#6,#7,#8.{\diagram{#1&\gfleche{#5}{}&#2\cr
  \hfleche{#6}{}&&\hfleche{}{#7}\cr
  #3&\gfleche{#8}{}&#4\cr}}

\def\rest#1{\lower 2pt \hbox{$|$}_{#1}}
\def\grtilde#1{\,\widetilde{\!#1}}
\def\Spec{\mathop{\rm Spec}}
\def\Proj{\mathop{\rm Proj}}
\def\discr{\mathop{\rm discr}\nolimits}
\def\Mod{\mathop{\euf Mod}}
\def\Qco{\mathop{\euf Qco}}
\def\Coh{\mathop{\euf Coh}}
\def\hom{\mathop{\rm Hom}\nolimits}
\def\pp{{{\euf p}\strot}}
\def\aa{{\euf a}}
\def\qq{{{\euf P}\strot}}
\def\rr{{\euf r}}
\def\bb{{\euf b}}
\def\ot{\widetilde{O}}
\def\pmod#1{\allowbreak\mkern5mu({\rm mod}\,\,#1)}
\def\tra{\mathop{\rm Tr}\nolimits}
\def\res{\mathop{\rm Res}\nolimits}
\def\fl#1{\buildrel #1\over {\hbox{\toto  \char 032}}}
\def\fmin#1{\buildrel #1\over {\hbox{\raise 3pt\hbox{\tot p}\hskip-2pt\toto\char 032}}}

\def\zeta{{\mathchar"0110\strot}}

\def\longsur{\relbar\joinrel\twoheadrightarrow}

\def\longinj{\lhook\joinrel\longrightarrow}
\def\findemo{{\unskip\nobreak\hfil\penalty50\hskip1em
             \hbox{}\nobreak\hfil\hbox{\di \char'144} \parfillskip=0pt
             \finalhyphendemerits=0 \par}}

 \input miniltx

\ifx\pdfoutput\undefined
  \def\Gin@driver{dvips.def} 
\else
  \def\Gin@driver{pdftex.def} 
\fi
 
\input graphicx.sty
\resetatcatcode

\baselineskip 11pt

\font\nom=cmr12 at 15pt
\font\tut=cmr12 at 13pt
\font\nic=cmr14 at 24pt

 \pageno=0
\line{GYMNASE DE BEAULIEU\hfill}


\kern4pt
\hrule

\vskip5cm
\centerline{\nic La conjecture de Catalan}

\vskip2cm\centerline{\tut RacontŽe ˆ un ami qui a le temps}
\vskip3cm\centerline{par}
\vskip.5cm\centerline{\nom Maurice Mischler}

\vskip1cm\centerline{ˆ partir d'un sŽminaire guidŽ par}
\vskip1cm\centerline{\nom Jacques BoŽchat }

\vskip3cm
\centerline{Mai 2005}

\vfill\eject

\long\def\art#1{{\parindent0pt\item{#1}}\hangindent=7mm\hangafter=-20}
\long\def\artart#1{{\parindent0pt\item{#1}}\hangindent=12mm\hangafter=-20}
\font\para=cmbx12 at 18pt
\def\O{\hbox{$\cal O$}}
\def\U{\hbox{$\cal U$}}
\def\m{\hbox{\rs m\!}}
\def\dst{\displaystyle}
\font\doub=msbm10 at 10pt

\def\qed{\hfill$\square$}

\newcount\chapnomb \chapnomb=1
\newcount\parnomb \parnomb=1
\pageno =-1

\parindent0pt

\centerline {\para Introduction }
\bigskip

Toute personne qui lira ce texte, n'est peut-tre
pas {\it a priori} un ami, mais on espre sincrement qu'elle le deviendra au fil de sa
lecture.  

Tout commence pour nous en juin 2002,
lorsqu'un collgue apprend ˆ l'un de nous que Preda Mih$\breve{\rm a}$ilescu avait rŽussi ˆ montrer cette conjecture. Nous avons donc dŽcidŽ avec enthousiasme, de faire un sŽminaire sur la
preuve de ce rŽsultat. Nous avons trouvŽ des notes de Youri Bilu [Bil] et de RenŽ Schoof [Sch] sur
internet et avons utilisŽ le livre de Paulo Ribenboim [Rib] sur le sujet ainsi qu'un article de
Mih$\breve{\rm a}$ilescu [Mih]. Nous avons voulu tre le plus Òself-contained" possible. Alors, nous
avons dŽmontrŽ le plus possible de rŽsultats, qui semblaient tre bien connus de leurs auteurs, mais
pas (ou mal) par nous. On aurait pu appeler ce texte ÒLa conjecture de Catalan pour les nuls", mais
on a vite remarquŽ que la difficultŽ de la preuve allait croissante. S'il est possible de lire les 
quatre premiers chapitres avac un bagage minimal, il faut au moins avoir suivi un cours
de deuxime cycle de thŽorie des nombres pour comprendre les chapitres suivants. Vous remarquerez
d'ailleurs que le style de ce texte est assez lŽger au dŽbut et qu'il deviendra de plus en plus
austre au fur et ˆ mesure qu'on avancera dans les difficultŽs.

Ont participŽ de manire plus ou moins suivie ˆ ce sŽminaire Henri Joris, Emmanuel
Preissmann, StŽphane Materna, Michel-StŽphane Dupertuis et Vincent Brayer.

Cette histoire  commence en 1844, avec Eugne Catalan. Il a
posŽ dans un journal trs lu (le journal de Crelle) la question
suivante~:

$Ç${\sl Je vous prie, Monsieur, de bien vouloir Žnoncer,
dans votre recueil, le thŽorme suivant, que je crois vrai,
bien que je n'aie pas encore rŽussi ˆ le dŽmontrer compltement, d'autres seront
peut-tre plus heureux~: Deux nombres entiers consŽcutifs, autres que 8 et
9, ne peuvent tre des puissances exactes; autrement dit~:
{\it l'Žquation $x^m-y^n = 1$, dans laquelle les inconnues sont
entires et positives, n'admet qu'une seule solution.}}$È$

Evidemment, il a dž s'y prendre ainsi~: il a remarquŽ que $3^2-2^3=1$ et a essayŽ avec
d'autres puissances, sans succs; il a peut-tre prouvŽ ce rŽsultat pour
des valeurs de $m$ et $n$ particulires. Alors il a Žmis ce qu'on appelle une Òconjecture"
c'est-ˆ-dire un rŽsultat mathŽmatique que l'on croit vrai, mais qu'on ne sait pas prouver.

Voilˆ, voilˆ... des conjecture comme a, il y en a des centaines, voire des
milliers. Mais il faut reconna"tre que celle-ci est particulirement simple~: on prend
deux nombres qui sont des puissances de nombres entiers, alors ils ne sont consŽcutifs que lorsque
ces nombres sont 8 et 9. Evidemment, quand on parle de Òpuissance", cela veut dire comme le dit
Catalan de Òpuissance exacte", c'est-ˆ-dire que l'exposant sera plus grand que 1, car sinon, on
aurait des solutions comme $34^1-33^1=1$ ou alors $2^5-31^1=1$ ou encore $2^1-144^0=1$.

Regardons ces calculs~: 

$$1024=257^2-255^2=130^2-126^2=68^2-60^2=40^2-24^2=32^2-0^2.$$

C'est le cas $n=5$ de l'Žquation 

$$2^{2n}=(2^{2n-k-1}+2^{k-1})^2-(2^{2n-k-1}-2^{k-1})^2\quad\hbox{pour }k=1,\ldots ,n.$$

Cela veut dire que pour tout entier $n$, il existe $m$ ( ici, c'est $2^{2n}$) tel que
l'Žquation $x^y-z^t=m$ possde au moins $n$ solutions diffŽrentes pour des entiers
$x,y,z,t$ plus grand que $1$. Si on remarque cela, c'est pour
donner un bŽmol ˆ la conjecture suivante~: 

\hbox{ÒPour tout entier naturel $m$, il n'existe qu'un nombre fini de solutions ˆ
l'Žquation }$x^y-z^t=m \hbox{ pour autant que } x,y,z,t\hbox{ soient }>1$"

La remarque qui prŽcde montre que ce nombre peut malgrŽ tout tre aussi grand qu'on
veut. Nous allons donc montrer cette conjecture pour $m=1$.

Au lieu de prendre des $x,y$ entiers supŽrieurs ˆ 1, on prendra des entiers non nuls ($x\ne 0$ et
$y\ne 0$).  

On veut donc trouver tous les entiers $(x,y,m,n)$ tels que $x^m-y^n=1$. On peut
dŽjˆ supposer que $n$ et $m$ sont des nombres premiers. en effet~: si par exemple $x^{21}-y^{20}=1$
possde une solution qu'on notera $(x_0,y_0)$, alors l'Žquation $x^7-y^5=1$ possde aussi une
solution~: c'est $(x_0^3,y_0^4)$. 

Donc on se restreint ˆ rŽsoudre l'Žquation $\thboxed 55{x^p-y^q=1}$ avec $p$ et $q$ des nombres
premiers et $x,y\ne 0$. On va montrer que les seules solutions sont $(p,q,x,y)=(2,3,\pm 3,2)$.

D'ores et dŽjˆ, on peut supposer que $p$ et $q$ ne sont pas tous les deux le nombre 2. Parce que sinon
$x^2-y^2=(x-y)(x+y)=1$, et ceci est  absurde car on aurait $x-y=\pm 1$ et $x+y=\pm 1$, ce qui
voudrait dire que $x=0$ ou $y=0$ ce qui est contraire ˆ l'hypothse. 

Voici un petit rŽsumŽ de ce qui va suivre~:

La premire chose (enfin, ce n'est pas la premire, mais presque !) qu'on va montrer, a
nous prendra dŽjˆ 8 pages, c'est que ni $p$, ni $q$ ne sont Žgaux ˆ 2 sauf pour le cas
$x^2-y^3=1$, ce sera les chapitres 1,2 et~3. Le chapitre 4 est consacrŽ aux identitŽs de Cassels,
c'est-ˆ-dire que si $x^p-y^q= 1$, alors $q$ divise $x$ et $p$ divise $y$. Ces quatre premiers
chapitres ne requirent aucune connaissance approfondie en thŽorie des nombres. Le chapitre 5 est
le premier Ònon ŽlŽmentaire" il portera sur le ThŽorme de Stickelberger. Les Chapitres 6 ˆ 9 donnent
les thŽormes de Mih$\breve{\rm a}$ilescu, ils sont au nombre de quatre (un par chapitre). Ils
permettent alors de prouver la conjecture de Catalan en une page au Chapitre 10. Ensuite on donne
deux appendices, l'un sur les anneaux semi-simples, l'autre sur le thŽorme de Thaine. 

Nous avons voulu montrer les choses le plus soigneusement possible sans que le lecteur soit obligŽ
(comme c'est hŽlas trop souvent le cas) de consulter des centaines de sources. Nous n'avons pas pu
nŽanmoins tre totalement autonomes. Les rŽsultats Òclassiques", comme les solutions de
l'Žquation de Pell, les triplets pythagoriciens, la thŽorie de Galois, le thŽorme de Dirichlet sur
les unitŽs d'un corps de nombres, les rŽsultats sur les sŽries $L$ de Dirichlet, le thŽorme de
Hensel ou le thŽorme de Hilbert 90 seront supposŽs connus et sont prouvŽs de manire exhaustive dans
des ouvrages de rŽfŽrence comme [Sier], [Nar], ou [Ser]. En revanche, deux gros thŽormes, dont la
preuve est difficile ˆ comprendre, mme dans les ouvrages de rŽfŽrence sont utilisŽs dans la preuve du
thŽorme de Thaine~: le thŽorme de
$\check{\rm C}$ebotarev et l'existence du corps de Hilbert. Mais nous sommes actuellement en train de
faire un nouveau sŽminaire lˆ-dessus, peut-tre rŽdigerons-nous des notes sur ce sujet~?

Enfin, nous tenons ˆ remercier chaleureusement Preda Mih$\breve {\rm a}$ilescu. Gr‰ce ˆ lui,
aux notes de RenŽ Schoof et ˆ celles Youri Bilu, nous avons passŽ deux ans de pur
bonheur. Puisse ce texte retranscrire le plaisir que nous avons eu ˆ essayer de comprendre
cette merveilleuse preuve.

\vfill\eject

\long\def\art#1{{\parindent0pt\item{#1}}\hangindent=7mm\hangafter=-20}
\long\def\artart#1{{\parindent0pt\item{#1}}\hangindent=12mm\hangafter=-20}
\font\para=cmbx12 at 18pt
\def\O{\hbox{$\cal O$}}
\def\U{\hbox{$\cal U$}}
\def\m{\hbox{\rs m\!}}
\def\dst{\displaystyle}
\font\doub=msbm10 at 10pt
\def\lra{\longrightarrow}
\def\qed{\hfill$\square$}
\def\gfP{\relax\ifmmode\bbP\else $\bbP$\fi}
\def\gP{{\euf P}}
\def\P{{\cal P}}
\def\QQ{{\cal Q}}
\def\Log{{\rm Log}}

\def\ggP{{\bf P}}

\parindent0pt
\pageno=-3
\NoBlackBoxes
\centerline{\para Table des matires}
\vskip 3cm
{\bf Chapitre 1 : Un thŽorme d'Euler (la solution non triviale)}\dotfill 1

\bigskip
{\bf Chapitre 2 : Le thŽorme de Lebesgue (le cas $\taille{10} q=2$)}\dotfill 4

\bigskip
{\bf Chapitre 3 : Le thŽorme de Ko-Chao (le cas $\taille{10} p=2$)}\dotfill 6
\bigskip
{\bf Chapitre 4 :  Les relations de Cassels}\dotfill 9
\bigskip
{\bf Chapitre 5 :  Le thŽorme de Stickelberger}\dotfill 16
\bigskip
{\bf Chapitre 6 :  Premier ThŽorme de Mih$\breve{\bf a}$ilescu}\dotfill 29

\bigskip
{\bf Chapitre 7 : Premiers contacts avec le groupe $\taille {10} H$ et petites valeurs de
$\taille {10}p$ et $\taille {10}q$ }\dotfill 32
\bigskip
{\bf Chapitre 8 : Troisime thŽorme de Mih$\breve{\bf a}$ilescu~: $\taille {10} p<4q^2$ et $\taille {10}
q<4p^2$ }\dotfill 40
\bigskip
{\bf Chapitre 9 : Quatrime thŽorme de Mih$\breve{\bf a}$ilescu~: $\taille {10} p\equiv 1\pmod q$ ou $\taille {10} q\equiv 1\pmod p$ }\dotfill 46
\bigskip
{\bf Chapitre 10 : Preuve de la Conjecture de Catalan }\dotfill 54
\bigskip
{\bf Appendice 1 : Deux mots sur les anneaux semi-simples }\dotfill 55
\bigskip
{\bf Appendice 2 : Le thŽorme de Thaine }\dotfill 58
\bigskip
{\bf Bibliographie }\dotfill 70

\vfill\eject

\long\def\art#1{{\parindent0pt\item{#1}}\hangindent=7mm\hangafter=-20}
\long\def\artart#1{{\parindent0pt\item{#1}}\hangindent=12mm\hangafter=-20}
\font\para=cmbx12 at 18pt
\def\O{\hbox{$\cal O$}}
\def\U{\hbox{$\cal U$}}
\def\m{\hbox{\rs m\!}}
\def\dst{\displaystyle}
\font\doub=msbm10 at 10pt

\def\qed{\hfill$\square$}

\newcount\chapnomb \chapnomb=1
\newcount\parnomb \parnomb=1
\pageno =1

\parindent0pt
\centerline{\para CHAPITRE 1}

\bigskip
{\para Un thŽorme d'Euler (la solution non triviale)}
\bigskip

Dans ce chapitre, on va montrer que la seule solution non triviale de l'Žquation $x^2-y^3=1$ est
$x=\pm3$ et $y=2$. On dira que c'est le cas $p=2$, $q=3$.

On a pris cette preuve dans le livre {\it Elementary theory of numbers} de Sierpinski [Sier]. Mais il
nous semblait opportun de la donner tout de mme.
\bigskip

{\bf Notations}

On considre connu la notion de pgcd$(x,y)$. Pour raccourcir, on Žcrira $(x,y)$ pour
pgcd$(x,y)$. Quand  $x$ divise $y$, on note $x|y$ ou alors $y\equiv 0\pmod x$. De manire
gŽnŽrale, $y\equiv y'\pmod x$ veut dire que $x|y-y'$.

{\bf Lemme 1}

{\sl Soient $a,b,d,n,m$ des nombres entiers.

\art{a)} Si $(a,b)=1$ et $(a,d)=1$ alors $(a,bd)=1$.

\art{b)} Si $a^n|b^n$ alors $a|b$.

\art{c)} Si $(a,b)=1$, alors $(a^n,b^m)=1$.

\art{d)} Supposons que $a,b\geq 1$, $(a,b)=1$ et $ab=c^n$, alors $a=a_1^n$ et $b=b_1^n$ pour
des entiers naturels $a_1,b_1$. Remarque que si la puissance $n$ est impaire, on peut se passer de
l'hypothse que $a,b\geq 1$, mais la conclusion sera que $a_1,b_1$ seront des entiers
Žventuellement nŽgatifs. {\bf Attention, ce rŽsultat sera trs souvent utilisŽ par la suite}

\art{e)}Supposons que $a,b\geq 1$, $(a,b)=1$ et $ab=p^kc^n$, avec $p$ un nombre premier, alors
$a=p^k a_1^n$ et $b=b_1^n$ ou alors $a=a_1^n$ et $b=p^k b_1^n$ pour des entiers naturels
$a_1,b_1$. 

}

{\bf Preuve}

Ce sont des rŽsultats ŽlŽmentaires, tous montrŽs dans [Sier].\qed

\bigskip

{\bf Lemme 2}

{\sl l'Žquation $$x^4+9x^2y^2+27y^4=z^2\eqno {(1)}$$  n'a pas de solution dans les nombres
entiers non nuls.}

{\bf Preuve}

Supposons que $x,y,z\in\N$ soit une solution de $(1)$
avec la valeur de $z$ la plus petite possible. Si $(x,y)=d>1$ alors $x=dx_1$ et $y=dx_2$, en
incorporant cela dans $(1)$, on voit que $d^4|z^2$, $d^2|z$ et $z=d^2z_1$, avec
$x_1,y_1,z_1\in\N$. Divisant (1) par $d^4$, on voit que $x_1,y_1,z_1$ est une solution, ce qui
contredit la minimalitŽ de $z_1$. Donc, on a $(x,y)=1$.

Si $2|x$, $(1)$ implique que $4|27 y^4-z^2$ et donc $2|y$ (car si $y$ impair $y^2\equiv 1\pmod 4$
et l'Žquation $3-z^2\equiv 0\pmod 4$ est impossible) ce qui contredit
$(x,y)=1$.

Donc $x$ est impair. Si $y$ est aussi impair, alors $8|5-z^2$ ce qui est aussi impossible. En
dŽfinitive,

$$\hbox{$x$ est impair et $y$ est pair}\eqno {(i)}$$

Si $3|x$, alors $27| z^2$, donc $9|z$ et $81|27 y^4$ et $3|y$ contrairement ˆ $(x,y)=1$, donc
$(x,3)=1$.

On a aussi $(x,z)=1$, en effet, si $(x,z)=d$, $(1)$ implique que $d|27 y^4$. Or, $(x,3y)=1$, donc
$(d,27y^2)=1$ et alors $d=1$. D'autre part, puisque $x$ est impair, $(1)$ implique que $z$ est
impair.

Posons maintenant $y=2y_1$. L'Žquation $(1)$ s'Žcrit alors (petit calcul)~:

$$27 y_1^4=\left({z+x^2\over 2}+9y_1^2\right )\left({z-x^2\over 2}-9y_1^2\right )=:A\cdot B.$$

Les facteurs $A$ et $B$ sont positifs puisque le premier l'est clairement et que leur produit
l'est aussi. On va montrer que $A$ et $B$ sont premiers entre eux. Soit $d_1$ leur pgcd.
Alors $d_1^2 | 27y_1^4$, donc $d_1^2|81 y_1^4$. Il suit (Lemme 1 b)) que $d_1|9 y_1^2$.
En additionnant les deux facteurs, on trouve $z$, donc $d_1|z$. En les soustrayant, on trouve
$x^2+18y_1^2$, donc $d_1| x^2$, puisque $d_1|9 y_1^2$. Ainsi, $d_1|(x^2,z)=1$, car
$(x,z)=(x^2,z)=1$. Cela veut dire que (Lemme 1 e))~:

$$ A=27 a^4,\quad B=b^4,\quad y_1=ab\eqno{ (ii)}$$

ou alors

$$ A= a^4,\quad B=27 b^4,\quad y_1=ab\quad
\eqno{(iii)}$$

o $a$ et $b$ sont des entiers positifs premiers entre eux. Mais le systme $(ii)$ est impossible,
car cela donnerait (en soustrayant les deux termes) que $x^2+18 a^2b^2=27a^4-b^4$. En regardant
cette ŽgalitŽ modulo 3 et en se souvenant que $(x,3)=1$ (donc $x^2\equiv 1\pmod 3$), on en dŽduit
que $b^4\equiv -1\pmod 3$, ce qui est absurde. Donc, c'est l'Žquation $(iii)$ qui est
possible. A nouveau en soustrayant, on obtient

$$x^2+18a^2b^2=a^4-27 b^4.\eqno{(iv)}$$ 

Gr‰ce ˆ $(i)$, on sait que $x$ est impair. Donc, regardant modulo 2, on trouve que $a$ ou $b$
est pair (le Òou" est ici exclusif, comme rarement en mathŽmatique). Si $a$ Žtait pair, alors
on aurait $a^4=x^2+18a^2b^2+27b^4\equiv 4\pmod 8$, ce qui est impossible. Donc, $a$ est impair
et $b$ est pair, et on trouve en triturant $(iv)$ que 

$$27 b^4=\left ( {a^2+x\over 2}-{9\over 2}b^2\right )\left ( {a^2-x\over 2}-{9\over 2}b^2\right
)=:C\cdot D.$$

Posons

$$d_2=\left (C,D\right
).$$

Nous trouvons que $d_2^2|27 b^4$, donc $d_2| 9b^2$. De plus (en soustrayant les deux termes du
pgcd) $d_2|x$. Donc, $d_2|(9y^2,x)$. Mais comme $(3y,x)=1$, on en dŽduit que $d_2=1$. Si les
nombres $C$ et $D$ Žtaient les deux nŽgatifs, on en dŽduirait que
$a^2<9b^2$ donc $a^4<9a^2b^2<18a^2b^2$. Donc, par $(iv)$, que $0>a^4-18a^2b^2=x^2+27b^4$, ce qui
est impossible. Donc, ces deux nombres sont positifs et premiers entre eux. On en dŽduit gr‰ce au
Lemme 1 e) que 

$${a^2\pm x\over 2}-{9\over 2}b^2=m^4,\quad {a^2\mp x\over 2}-{9\over 2}b^2=27n^4,\quad b=mn$$

Pour des entiers positifs $m$ et $n$. En additionnant tout a, on obtient 

$$a^2=m^4+9m^2n^2+27n^4$$

Mais, souvenons-nous que $a\leq y_1<y<z$. Cela contredit la minimalitŽ de $z$. Cela prouve par
ce qu'on appelle Òdescente infinie" que l'Žquation $x^4+9x^2y^2+27y^4=z^2$ n'a pas de solution
entire.\qed

{\bf Lemme 3}

{\sl l'Žquation $$x^3+y^3=2z^3\eqno {(2)}$$  n'a pas de solution dans les nombres
entiers tels que $x\ne y$ et $z\ne 0$.}

{\bf Preuve}

Bon et bien pour changer, on supposera que $x,y,z$ sont des solutions de $(2)$ avec $x\ne y$ et
$z\ne 0$. On peut dŽjˆ supposer que $(x,y)=1$, car si $(x,y)=d>1$, alors $d^3|2z^3$, donc $d|z$,
ce qui voudrait dire que ${x\over d},{y\over d},{z\over d}$ seraient aussi des solutions.

Puisque $x^3+y^3=2z^3$, alors $x+y$ et $x-y$ sont des nombres pairs. Posons $u={x+y\over 2}$ et
$v={x-y\over 2}$. Comme $(x,y)=1$, $(u,v)=1$, car $x=u+v$ et $y=u-v$. L'Žquation (2) devient
alors $(u+v)^3+(u-v)^3=2z^3$. Cela veut dire que $u(u^2+3v^2)=z^3$ et, puisque $x\ne y$ et $z\ne
0$, $uvz={1\over 4}(x^2-y^2)z\ne 0$. A partir de maintenant, deux cas se profilent.

\art{a)}Si $(u,3)=1$, alors $(u,u^2+3v^2)=1$, puisque $(u,v)=1$. Donc, (Lemme 1 d)) il existe
des entiers $z_1$ et $z_2$, premiers entre eux, tels que $u=z_1^3$ et $u^2+3v^2=z_2^3$. Ainsi,
$z_2^3-z_1^6=3v^2$ et en factorisant, cela donne $(z_2-z_1^2)\left ( (z_2-z_1^2)^2+3z_2z_1^2\right
)=3v^2$. Posons
$t=z_2-z_1^2\geq 0$. Alors $(t,z_1)=1$, puisque $(z_1,z_2)=1$. De plus,
$t(t^2+3tz_1^2+3z_1^4)=3v^2$. Cela implique que $3|t$. Posons $t=3t_1$. Alors on trouve que
$t_1(9t_1^2+9t_1z_1^2+3z_1^4)=v^2$, ce qui veut dire que $3|v$ et $v=3v_1$. Puisque $(u,3)=1$,
$(z_1,3)=1$, donc le nombre $9t_1^2+9t_1z_1^2+3z_1^4$ n'est pas divisible par 9. Mais $9|v^2$,
donc $3|t_1$ et alors $t_1=3t_2$. On trouve alors que $t_2(27t_2^2+9t_2z_1^2+z_1^4)=v_1^2$.
Puisque $(t,z_1)=1$, $(t_2,z_1)=1$ et donc $(t_2,27t_2^2+9t_2z_1^2+z_1^4)=1$. Ainsi, le Lemme 1
d) nous assure l'existence de $b$ et $c$ tels que $t_2=b^2$ et $27b^4+9b^2z_1^2+z_1^4=c^2$. Reste
a voir que $b$ et $|z_1|$ sont non nuls. Si $b=0$, alors $t_2$ et donc $t=0$, ce qui veut dire que
$z_2=z_1^2$. Le fait que $(z_2,z_1)=1$, implique  que $z_1=\pm 1$ et $z_2=1$ et donc que $v=0$
et donc que $x=y$ contrairement ˆ l'hypothse. D'autre part, si $z_1=0$, alors $u=0$ et donc
$z^3=u(u^2+3v^2)=0$ donc, $z=0$ contrairement ˆ l'hypothse. Cela implique que l'Žquation
$x^4+9x^2y^2+27y^4=z^2$ possderait des solutions entires non nulles ce qui contredit le lemme
prŽcŽdent.

\art{b)}Supposons que $3|u$. Puisque $(u,v)=1$, alors $(v,3)=1$. Posons $u=3u_1$, et de
$u(u^2+3v^2)=z^3$, on peut poser $z=3z_1$ et $u_1(3u_1^2+v^2)=3z_1^3$. Puisque $(v,3)=1$, alors
$3|u_1$ et on pose $u_1=3u_2$ et alors $u_2(27u_2^2+v^2)=z_1^3$. On a $(u_2,v)=1$, donc
$(u_2,27u_2^2+v^2)=1$. Il existe alors $a,b$, premiers entre eux, tel que $(a,b)=1$ et $u_2=a^3$ et
$27a^6+v^2=b^3$. Posons $t=b-3a^2$. On a $t^3\equiv b^3\equiv v^2\pmod 3$, donc $(t,3)=1$. On voit
que $v^2=b^3-27 a^6=(b-3a^2)(b^2+3a^2b+9a^4)=(b-3a^2)((b-3a^2)^2+9a^2(b-3a^2)+27
a^4)=t(t^2+9a^2t+27a^4)$. Puisque $(a,b)=1$, alors $(a,t)=1$ et on a vu que $(t,3)=1$, donc
$(t,t^2+9a^2t+27a^4)=1$. Par l'Žternel Lemme 1 d), il existe $a_1$ et $b_1$, premiers entre eux,
tels que $t=a_1^2$ et
$a_1^4+9a^2a_1^2+27a^2=b_1^2$. On va finalement voir que $a_1$ et $a$ sont non nuls. Si $a_1=0$
alors $t=0$ et donc $b=3a^2$ ce qui contredit $(a,b)=1$. Si $a=0$ alors $u=0$ et donc, comme
avant, $z=0$ contrairement ˆ l'hypothse. On conclut comme ˆ la partie a) que l'Žquation
$x^4+9x^2y^2+27y^4=z^2$ possderait des solutions entires non nulles.\qed

\bigskip\goodbreak

{\bf ThŽorme d'Euler}

{\sl L'Žquation $x^2-y^3=1$ n'a pas d'autres solutions que $x=0,\,y=-1$ ou $x=\pm 1,\, y=0$ ou encore
$x=\pm 3,y=2$. }

{\bf Preuve}

On a $y^3=x^2-1=(x-1)(x+1)$. 

Si $x$ est pair, alors $(x-1,x+1)=1$. Il existe donc des entiers
$a$ et $b$ premiers entre eux tels que $x-1=a^3$ et $x+1=b^3$. On a donc
$b^3+(-a)^3=2\cdot 1^3$. Ainsi, par le lemme prŽcŽdent, on a $a=-b$, donc $x-1=a^3=-b^3=-x-1$ ce
qui implique que $x=0$ et donc $y=-1$.

Si $x$ est impair alors $y=2y_1$ est pair et $({x-1\over 2},{x+1\over 2})=1$, et on trouve
$2y_1^3=({x-1\over 2})({x+1\over 2})$. Il existe donc des entiers $a$ et $b$ tels que $x\pm
1=4a^3$ et $x\mp 1=2b^3$. En Žliminant le $x$, on trouve $b^3+(\pm 1)^3=2a^3$. Le lemme
prŽcŽdent nous dit alors que $a=0$ donc $x=\pm 1$ et $y=0$. Si $a\ne 0$, alors $b=\pm 1$. Si
$b=1$, alors $x=3$, donc $y=2$ (attention aux $\pm$ et aux $\mp$). Enfin, si $b=-1$, alors
$x=-3$ et
$y=2$.\qed
\bigskip

\vfill\eject

\centerline{\para CHAPITRE 2}
\bigskip
{\para  Le thŽorme de Lebesgue (le cas $\taille{15} q=2$)}
\bigskip

La seconde chose, qu'on va voir, c'est un thŽorme de Monsieur Victor AmŽdŽe Lebesgue. Les
mathŽma-ticiens connaissent un autre Lebesgue, Henri LŽon, qui a inventŽ une intŽgrale qui porte
son nom, mais celui dont on parle n'est pas le mme. Il vaudrait mieux, pour ce chapitre,
conna"tre ce qu'est un anneau factoriel.

{\bf Lemme 1}

{\sl 
\art{a)}Soit $A$ un anneau factoriel et $x,y$ premiers entre eux (cela veut dire qu'ils n'ont pas de
diviseurs communs autre que les unitŽs), alors si
$xy=z^m$, alors $x=ux_1^m$ et $y=u^{-1}y_1^m$, avec $x_1,y_1\in A$ et $u\in U(A)=$ les inversibles de
$A$ (c'est l'Žquivalent du Lemme 1 d) du chapitre prŽcŽdent).

\art{b)}L'anneau $\Z[i]=\{a+bi\mid i=\sqrt{-1}\hbox{ et }a,b\in \Z\}$, qui est l'anneau des entiers
de Gauss, est factoriel.

}

{\bf Preuve}

\art{a)} Suit de la dŽfinition d'anneau factoriel

\art{b)}[Mar, ex 7 p. 7]

\bigskip

{\bf ThŽorme de Lebesgue}

{\sl Il n'existe pas de solution $(x,y)\in\Z^2$ avec $x,y\ne 0$, tels que $x^m-y^2=1$ avec
$m\in\N$.}

{\bf Preuve}
On peut supposer $m$ impair, car sinon, $x^2-y^2=1$ aurait une
solution et on a vu que c'Žtait impossible. On peut aussi supposer que $x,y\geq 1$ (car $y$ est au
carrŽ et si $x$ est nŽgatif, l'Žquation est trivialement fausse). Supposons par
l'absurde qu'il existe $x$ et $y$ tels que $x^m=y^2+1$ avec $m$ impair. Remarque que $x$ et $y$ sont de
paritŽs diffŽrentes. Car si $x$ est pair, $x^m$ aussi et donc $y^2+1$ aussi, donc $y^2$ est impair et
donc $y$ aussi, et rŽciproquement, si $x$ est impair, $y$ aussi.

Supposons que $y$ soit impair et donc que $x$ soit  pair. Alors lˆ, si on regarde
modulo 4, on remarque que alors $x^m\equiv y^2+1\equiv 2\pmod 4$. Et tu peux remarquer
facilement que pour tout nombre pair et tout $m\geq 2$, on a $x^m\equiv 0\pmod 4$.

On peux donc supposer que $y$ est pair et $x$ impair. Bon, maintenant on va passer
dans l'anneau $\Z[i]$. Alors, l'Žquation de dŽpart $x^m=y^2+1$ s'Žcrit $x^m=(y-i)(y+i)$. On va
d'abord montrer que $(y-i)$ et $(y+i)$ sont premiers entre eux dans $\Z[i]$ (c'est-ˆ-dire que
si un nombre premier de $\Z[i]$ divise $(y+i)$ et $(y-i)$ alors il y a un contradiction).
Supposons donc que $\pi$, premier de $\Z[i]$, divise $(y+i)$ et $(y-i)$. On Žcrit $a|b$ pour $a$
divise $b$ (comme dans $\Z$). Alors $\pi | (y+i)-(y-i)=2i$. Puisque $\pi$ est premier, cela
implique que $\pi|2$ ou $\pi|i$. On sait que $i$ est inversible dans $\Z[i]$ ($i\cdot (-i)=1$).
Donc $\pi$ ne peut pas diviser un inversible, puisqu'il est premier (un premier n'est pas
inversible). Donc $\pi$ divise $2$. Mais lˆ, a ne va pas non plus~: puisque $y$ est pair, on
aurait  que $\pi$ divise 2 qui divise $y$. Mais alors $\pi| (y+i)-y=i$ et on vient de dire que
c'est impossible... Donc $(y+i)$ et $(y-i)$ sont premiers entre eux. Quand un produit de deux
nombres premiers entre eux est un puissance $m$-ime, alors chacun de ces nombres est une
puissance $m$-ime (cf. Lemme 1 a) et b)). Ainsi, il existe $u,v\in\Z$ et $0\leq s\leq 3$
tels que 

$$ (y+i)=(u+iv)^m\cdot i^s\qquad \hbox{et}\qquad (y-i)=(u-iv)^m\cdot (-i)^s$$

En multipliant ces deux termes, on trouve $x^m=(u^2+v^2)^m$, donc $x=u^2+v^2$, ainsi, $u$ et $v$
ne sont pas de mme paritŽ, puisque $x$ est impair. En soustrayant les deux mmes termes, on
obtient 
$$2i=\left ( (u+iv)^m-(u-iv)^m(-1)^s\right )\cdot i^s$$

\art{a)} Supposons que $s=2r$ (=0 ou 2). On en dŽduit ici que $i^s=(-1)^r$ et $(-1)^s=1$. En
dŽveloppant gr‰ce au bin™me de Newton, en Žgalant les parties imaginaires, en divisant par 2 et
en se souvenant que
$m$ est impair, on obtient 

$$1=( -1)^r\left ( \pmatrix{m\cr 1\cr} u^{m-1}v-\pmatrix{m\cr 3\cr} u^{m-3}v^3+\cdots \pm
v^m\right )$$

On en dŽduit que $v$ divise 1, donc $v=\pm 1$ et donc que $u$ est pair, car il n'ont pas la mme
paritŽ.

\art{b)}Si $s=2r+1$, on obtient (en comparant les parties rŽelles cette fois) que 
$$1=(-1)^r\left ( u^m-\pmatrix{m\cr 2\cr} u^{m-2} v^2+\cdots +\pmatrix{m\cr 1\cr}uv^{m-1}\right
)$$

Par le mme raisonnement qu'avant, on en dŽduit cette fois que $u=\pm 1$ et que c'est $v$ qui est
pair. Dans tous les cas, posons $w$ l'ŽlŽment qui de $u$ ou $v$ est pair. Puisque, pour
$k=0,\ldots ,m$, on a  
$\pmatrix{m\cr k\cr}=\pmatrix{m\cr m-k\cr}$ et que $m$ est impair, on trouve~: 

$$1-\pmatrix{m\cr 2\cr}w^2+\pmatrix{m\cr 4\cr}w^4-\cdots \pm m w^{m-1}=\pm 1$$

Comme $w$ est pair et donc que les $w^{2j}$ sont multiples de 4, alors le signe du terme de
droite est $+1$ et il reste en divisant par $w^2$ 

$$\pmatrix{m\cr 2\cr}-\pmatrix{m\cr 4\cr}w^2+\cdots \pm m w^{m-3}=0\eqno{(i)}$$

Posons $v_2(k)$ la valuation 2-adique de $k$, cela veut dire que si
$k=2^t\cdot u$ avec $u$ un nombre impair, alors $v_2(k)=t$. Si on considre l'Žquation $(i)$
modulo 2, on voit que $\pmatrix{m\cr 2\cr}$ est pair. Posons alors $t=v_2(\pmatrix{m\cr
2\cr})\geq 1$. Les autres termes sont de la forme $\pmatrix{m\cr 2k\cr}w^{2k-2}$ avec $k\geq 2$.
Si $k\geq 2$, alors 
$$\pmatrix{m\cr 2k\cr}w^{2k-2}=\pmatrix{m\cr 2\cr}\pmatrix{m-2\cr 2k-2 \cr}{2\over
2k(2k-1)} w^{2k-2}$$

D'autre part, pour $k\geq 2$ , on a $2^{2k-2}>k$, donc $v_2(k)<2k-2$. En utilisant le fait
facile que $v_2(ab)=v_2(a)+v_2(b)$, on en dŽduit que $v_2(\pmatrix{m\cr 2k\cr}w^{2k-2})\geq 
t+(2k-2)-v_2(k)>t+(2k-2)-(2k-2)=t$ (n'oublions pas que $v_2(2k-1)=0$). On en dŽduit que l'Žquation
$(i)$ est impossible. En effet, elle s'Žcrirait $2^t\alpha-2^{t+l}\beta=0$ o $\alpha$ est un
nombre impair, $l$ un nombre positif et $\beta$ un entier non nul. Et on en dŽduirait que
$\alpha=2^l\cdot \beta$, ce qui est absurde. Donc l'Žquation $(i)$ est impossible.\qed

\vfill\eject

\long\def\art#1{{\parindent0pt\item{#1}}\hangindent=7mm\hangafter=-20}
\long\def\artart#1{{\parindent0pt\item{#1}}\hangindent=12mm\hangafter=-20}
\font\para=cmbx12 at 18pt
\def\O{\hbox{$\cal O$}}
\def\U{\hbox{$\cal U$}}
\def\m{\hbox{\rs m\!}}
\def\dst{\displaystyle}
\font\doub=msbm10 at 10pt

\def\qed{\hfill$\square$}

\newcount\chapnomb \chapnomb=1
\newcount\parnomb \parnomb=1
\pageno =6

\parindent0pt
\bigskip
\centerline{\para CHAPITRE 3}

\bigskip

{\para Le thŽorme de Ko-Chao (le cas $\taille{15} p=2$)}
\bigskip
Dans ce chapitre, nous allons prouver que l'Žquation $x^2-y^q=1$ n'a pas de solution si
$q\geq 5$.
\bigskip
{\bf Lemme de Diophante}

{\sl Si $x,y,z$ forment un {\it triplet pythagoricien primitif}, i.e.  sont des entiers
premiers entre eux tels que $x^2+y^2=z^2$, alors il existe des entiers $c$ et $d$, premiers
entre eux et de paritŽ diffŽrentes tels que $x=c^2-d^2$, $y=2bc$ et $z=c^2+d^2$. }

{\bf Preuve}

C'est un rŽsultat connu depuis plus de 3000 ans, cf. par exemple [Sier, ThŽorme 1,
p.38].\qed
\bigskip
{\bf Lemme de Pell}

{\sl Soit $D$ un nombre entier positif qui n'est pas le carrŽ d'un nombre entier. Alors
l'Žquation 
$$x^2-Dy^2=1\eqno{(3)}$$

possde une infinitŽ de solutions entires non nulles. Plus prŽcisŽment, si $x_1,y_1$ est la
solution fondamen-tale (c'est-ˆ-dire la solution de $(3)$ en entiers positifs dont $y$ est le
petit possible (non nul)). Alors les autres solutions positives sont les $x_m,y_m$ tels que 

$$(x_m+y_m\sqrt{D})=(x_1+y_1\sqrt{D})^m,\ m\geq 1\eqno{(4)}$$}

{\bf Preuve}

La dŽmonstration de ce fait est donnŽe dans tout bon livre de thŽorie des
nombres, par exemple [Sier, ThŽorme 15, p.98].\qed
\bigskip

{\bf DŽfinition}

Sous les hypothses du Lemme de Pell, soit $m$ un entier positif. On dit que $m$ a la {\it
propriŽtŽ de Stoermer} si pour tout premier
$p$, $p|y_m$ implique $p|D$. 

{\bf Lemme de Stoermer (1898)}

{\sl Si $m$ satisfait la propriŽtŽ de Stoermer, alors $m=1$. Autrement dit, le seul cas o la
condition de Stoermer peut tre possible est la solution fondamentale.}

{\bf Preuve}

L'Žquation $(4)$ dŽveloppŽe donne

$$\eqalign{x_n&=x_1^n+\pmatrix{n\cr 2\cr} x_1^{n-2}y_1^2D+\cdots\cr y_n&=n
x_1^{n-1}y_1+\pmatrix{n\cr 3\cr} x_1^{n-3}y_1^3D+\cdots .\cr}\eqno{(i)}$$

Donc $y_1$ divise $y_n$ (mais attention, si $n$ est pair, $x_1$ ne divise pas forcŽment $x_n$,
le dernier terme du dŽveloppement de $x_n$, ne contient pas de $x_1$). Posons alors $z_1=1$ et
$y_n=z_ny_1$. Ainsi,
$z_n=n x_1^{n-1}+\pmatrix{n\cr 3\cr} x_1^{n-3}y_1^2D+\cdots$. Posons $A=x_1^2-1=Dy_1^2>1$, donc
$\sqrt{A}=y_1\sqrt{D}$. L'Žquation $(4)$ devient $(x_n+z_n\sqrt{A})=(x_1+\sqrt{A})^n$ et donc, 
$$z_n=nx_1^{n-1}+\pmatrix{n\cr 3\cr}x_1^{n-3}A+\cdots .\eqno{(ii)}$$

Si $q|n$ alors $z_q|z_n$. En effet, si $n=ql$,
$(x_n+z_n\sqrt{A})=(x_1+\sqrt{A})^n=(x_q+z_q\sqrt{A})^l$. D'o $z_q$ divise $z_n$, de la mme
manire que $y_1$ divise $y_n$. 

Fixons $m\geq 1$ ayant la propriŽtŽ de Stoermer et soit $q$ premier tel que $q|m$. Donc,
$z_q|z_m$. Soit $p$ premier tel que $p|z_q$. Alors $p|z_m$, donc $p|y_m$, par dŽfinition de
$z_m$. Ainsi $p|D$ par propriŽtŽ de Stoermer et $p|A=Dy_1^2$. On se souvient $(ii)$ que
$$z_q=qx_1^{q-1}+\pmatrix{q\cr 3\cr}x_1^{q-3}A+\cdots\eqno{(iii)}$$ 
donc $p|q x_1^{q-1}$.
Comme $x_1^2=A+1$ et que $p|A$, alors $p\not \hskip-0.4pt |\  x_1$ et ainsi, $p|q$ et donc
$p=q$. Posons
$z_q=q^r$, $r\geq 0$. Si $r=0$, alors $z_q=1$, ce qui est incompatible avec $(iii)$. Supposons
que $q>3$ et
$r\geq 2$. La relation $(iii)$ nous donne $q^2|qx_1^{q-1}$ (car $q$ divise $A$ et $\pmatrix{q\cr
3\cr}$ et les autres coefficients bin™miaux). Donc $q|x_1$, ce qui est impossible, comme on l'a
vu. Ainsi, si
$q>3$, alors $r=1$, et donc $z_q=q$, mais c'est ˆ nouveau impossible car $q x_1^{q-1}>q$. Donc
$q=2$ ou 3. Supposons $z_2=2^r$, $r\geq 1$. L'Žquation $(iii)$ devient $z_2=2x_1$, donc
$x_1=2^{r-1}$. Comme
$2\not \hskip-0.4pt |\  x_1$, $r=1$, donc $x_1=1$, mais $x_1>1$ !!

RŽsumons-nous~: on a dŽmontrŽ que si $m$ a la propriŽtŽ de Stoermer, alors $m=3^t$. On va
montrer d'abord que si $m=3n$ possde la propriŽtŽ de Stoermer alors $n$ aussi. En
effet, 
$$(x_{3n}+y_{3n}\sqrt{D})=(x_{1}+y_{1}\sqrt{D})^{3n}=(x_{n}+y_{n}\sqrt{D})^3$$
Donc,
$$\eqalignno{x_{3n}&=3x_ny_n^2D&(iv)\cr
y_{3n}&=3x_n^2y_n+y_n^3D=y_n(3x_n^2+y_n^2D)&\cr}$$

Maintenant, si $p|y_n$, alors $p|y_{3n}$, donc $p|D$, ce qui implique que $n$ possde la
propriŽtŽ de Stoermer. 

Reste ˆ voir que 3 n'a pas la propriŽtŽ de Stoermer. Supposons donc que 3 possde la propriŽtŽ
de Stoermer. Supposons aussi que $p|3x_1^2+y_1^2D$. Par $(iv)$ et $n=1$, on a alors
$p|y_3$, donc (Stoermer) $p|D$, donc $p|3x_1^2$. Or, on a vu que si $p|D$, alors $p$ ne divise
pas $x_1$, donc $(p,x_1)=1$. Donc, $p=3$. Ainsi, $3x_1^2+y_1^2D=3^s$ avec $s\geq 2$, car
$3x_1^2+y_1^2D>4$. On obtient alors
$3^s=3x_1^2+y_1^2D=3x_1^2+(x_1^2-1)=4x_1^2-1=(2x_1+1)(2x_1-1)$. On tire le systme

$$\left\{\eqalign{2x_1+1&=3^r\cr 2x_1-1&=3^t\cr}\right. \quad 0\leq t<r,\ r+t\geq 2$$

En soustrayant, on trouve $2=3^r-3^t=3^t(3^{r-t}-1)$, ce qui donne $t=0$ et $r=1$, contredisant
$r+t\geq 2$.\qed
\bigskip
{\bf Petit Lemme}

{\sl Soit $x$ et $y$ des entiers premier entre eux et $n\geq 1$.  Alors

$$(x\pm y,{x^n\pm y^n\over x\pm y})=(x\pm y,n),$$

avec $n$ impair si on considre $x+y$.}

{\bf Preuve}

On a ${x^n-y^n\over x-y}=x^{n-1}+x^{n-2}y+\cdots +y^{n-1}$. Or, pour tout $k\geq 1$,
$x^k=y^k+(x^k-y^k)=y^{k}+(x-y)(\cdots )$. Ainsi, ${x^n-y^n\over x-y}=ny^{n-1}+(x-y)(\cdots)$.
Donc $(x-y,{x^n-y^n\over x-y})=(x-y,ny^{n-1})=(x-y,n)$, car $(x,y)=1$.

De mme, si $n$ est impair, ${x^n+y^n\over x+y}=x^{n-1}-x^{n-2}y+\cdots +y^{n-1}$. Si $k$ est
impair, on a
$x^k=-y^k+(y^k+x^k)=-y^k+(x+y)(\cdots)$. Et,
$x^{2k}=y^{2k}+(x^{2k}-y^{2k})=y^{2k}+(x^2-y^2)(\cdots)=y^{2k}+(x+y)(\cdots)$. Ainsi,
${x^n+y^n\over x+y}=ny^{n-1}+(x+y)(\cdots)$. Donc, ˆ nouveau, $(x+y,{x^n+y^n\over
x+y})=(x+y,ny^{n-1})=(x+y,n)$, car $(x,y)=1$.\qed
\bigskip
{\bf Lemme de Nagell (1921)}

{\sl Soit $x$ et $y$ des entiers strictement positifs et $q\geq 3$ un nombre premier.  Supposons
que
$x^2-y^q=1$ alors $2|y$ et $q|x$}

{\bf Preuve}

Supposons $y$ impair. On a $y^q=x^2-1=(x+1)(x-1)$ et $(x+1,x-1)=1$, donc il existe des entiers
positifs $a$ et $b$, tels que $x+1=a^q$ et $x-1=b^q$, et donc
$2=(a^q-b^q)=(a-b)(a^{q-1}+\cdots +b^{q-1})$. Si $b=0$ alors $y=0$, impossible. De mme, $a>0$.
Et donc
$(a^{q-1}+\cdots +b^{q-1})\geq 3$. Mais c'est impossible, car on devrait avoir
$(a^{q-1}+\cdots +b^{q-1})=1$ ou 2. Il ne reste plus que $x$ impair et $y$ pair, car il est
Žvident qu'ils ne peuvent pas tre les deux impairs. Donc $2|y$.

Supposons par l'absurde que $q\not \hskip-0.4pt |\   x$. On a $x^2=y^q+1=(y+1){y^q+1\over
y+1}$. Par le Petit Lemme, $(y+1,{y^q+1\over y+1})=q$ ou $1$. Si c'est $q$, cela veut dire que
$q|x$. Donc
$(y+1,{y^q+1\over y+1})=1$. Il existe donc $c>1$ et $d>0$ tels que $y+1=c^2$ et ${y^q+1\over
y+1}=d^2$. L'Žquation $x^2-y^q=1$ devient
$$x^2-\underbrace{(c^2-1)}_{:=D}\left [(c^2-1)^{q-1\over 2}\right]^2=1.$$

Ce qui veut dire que $x$, $(c^2-1)^{q-1\over 2}$ est une solution de l'Žquation de Pell
$X^2-DY^2=1$, disons la solution $X_m$ $Y_m$ (on met des majuscules ici pour Žviter les
confusions).  La solution fondamentale est visiblement
$c$,1. D'autre part, si $p|Y_m$, alors $p|D$, donc, par le Lemme de Stoermer $x$,
$(c^2-1)^{q-1\over 2}$ est la solution fondamentale, et alors $c^2-1=1$, ce qui est absurde, car
$\sqrt{2}$, n'est pas entier.\qed 
\bigskip\goodbreak
{\bf ThŽorme de Ko-Chao (1964)}

{\sl Il n'existe pas de solution $(x,y)\in \Z^2$ avec $x,y\ne 0$, tels que $x^2-y^q=1$ et $q>3$
premier.}

{\bf Preuve}

Nous allons donner la preuve de Chein de ce rŽsultat (1976).

Supposons par l'absurde que $x^2-y^q=1$. On peut supposer que $x,y>1$, car $x$ est au carrŽ
et si $y$ est nŽgatif, l'Žquation est trivialement fausse.

On peut aussi supposer (cf. Lemme de Nagell) que $x$ est impair et $y$

est pair. Dans ce cas,
$(x+1,x-1)=2$, et comme $y^q=x^2-1=(x-1)(x+1)$, il existe des entiers strictement positifs $a,b$
premiers entre eux avec $a$ impair tels que $y=2ab$ et
$$\left\{\eqalign{x+1&=2\cdot a^q\cr x-1&=2^{q-1}\cdot b^q\cr}\right .\eqno{ (I) }$$
 ou 
$$\left\{\eqalign{x+1&=2^{q-1}\cdot b^q\cr x-1&=2\cdot a^q.\cr}\right .\eqno{(II)}$$

Nous allons traiter les deux cas ensemble, pour le cas $(I)$, on lit le signe du haut, et pour
le cas $(II)$, on lit le signe du bas. En soustrayant les deux Žquations et en divisant par 2,
on trouve 
$$a^q-2^{q-2}b^q=\pm 1.\eqno{(*)}$$

En particulier on a
$$a^q\geq -1+2^{q-2}b^q>b^q$$

En effet, la dernire inŽgalitŽ est Žquivalente ˆ $(2^{q-2}-1) b^q>1$ qui est vrai puisque
$q>3$. On a ainsi montrŽ que $a>b$ (souvenons-nous bien de cela). Calculons~:
$$(a^q\mp2)^2=(a^2)^q\mp4a^q+4\buildrel (*)\over = (a^2)^q\mp
4(2^{q-2}b^q\pm1)+4=(a^2)^q\mp(2b)^q.$$

Mais, on a $a^q\mp 2={x\pm 1\over 2}\mp 2={x\mp 3\over 2}$, donc $\left ({x\mp
3\over 2}\right )^2=(a^2)^q\mp(2b)^q$, ou encore 
$$\left ({x\mp 3\over 2}\right )^2=(a^2\mp 2b)\cdot {(a^2)^q\mp(2b)^q\over a^2\mp
2b}.\eqno{(**)}$$

Le Petit Lemme s'applique alors, car $(a^2,2b)=1$. Cela veut dire que $(a^2\mp 2b,
{(a^2)^q\mp(2b)^q\over a^2\mp 2b})=1$ ou $q$. Le Lemme de Nagell, nous affirme que $q|x$, donc,
$q\not \hskip-0.4pt |\  {x\mp 3\over 2}$, (car sinon $q|3$). Ainsi, $(a^2\mp 2b)$ et $
{(a^2)^q\mp(2b)^q\over a^2\mp 2b}$ sont premiers entre eux, et sont donc les deux des carrŽs.
En particulier, $a^2\mp 2b=h^2$ et $h|{x\mp3\over 2}$. Puisque $a$ est impair, $h$ aussi et
donc $b$ est pair (en regardant modulo 4). Par suite $(h\cdot a)^2+b^2=a^4\mp2a^2b+b^2=(a^2\mp
b)^2$. On a affaire ˆ un triplet pythagoricien primitif (car $(b,a^2\mp b)=1$) avec $b$ pair.
Donc, par le Lemme de Diophante, il existe $c>d>0$ tels que $ha=c^2-d^2$, $b=2cd$ et $a^2\mp
b=c^2+d^2$. D'o $(c\pm d)^2=(a^2\mp b)\pm b=a^2$, c'est-ˆ-dire $a=c\pm d$.

Si on est dans le cas $(I)$, on a $b-a=2cd-(c+d)=(c-1)(d-1)+(cd-1)>0$. Si on est dans le cas
$(II)$, on a $b-a=2cd-(c-d)=c(2d-1)+d>0$. Dans tous les cas, on a $b>a$, mais on se souvient que
$a>b$. C'est une contradiction et le thŽorme est prouvŽ.\qed

\vfill\eject

\long\def\art#1{{\parindent0pt\item{#1}}\hangindent=7mm\hangafter=-20}
\long\def\artart#1{{\parindent0pt\item{#1}}\hangindent=12mm\hangafter=-20}
\font\para=cmbx12 at 18pt
\def\O{\hbox{$\cal O$}}
\def\U{\hbox{$\cal U$}}
\def\m{\hbox{\rs m\!}}
\def\dst{\displaystyle}
\font\doub=msbm10 at 10pt

\def\qed{\hfill$\square$}

\newcount\chapnomb \chapnomb=1
\newcount\parnomb \parnomb=1
\pageno =9

\parindent0pt
\centerline{\para CHAPITRE 4}

\bigskip

\bigskip
{\para  Les relations de Cassels}
\bigskip

A partir de maintenant, on va Žtudier l'Žquation $x^p-y^q=\pm 1$, avec $p$ et $q$ premiers impairs et
$x,y>1$. Et on va dŽmontrer que  dans ce cas $p|y$ et $q|x$. Mais pour cela, il va falloir travailler assez
ferme.

Commenons par un petit lemme d'analyse.

{\bf Lemme 1}

{\sl Soit $a,b,t$ des nombres rŽels tels que $b>0$, $t>1$ et $a+b^t>0$. On considre $f_{a,b}(t)=(a+b^t)^{1\over
t}$. Alors

$$f'_{a,b}(t)\matrix{\raise-3pt\hbox{$<$}\cr \raise3pt\hbox{$>$}\cr}0\Longleftrightarrow
b^t\log(b^t)\matrix{\raise-3pt\hbox{$<$}\cr
\raise3pt\hbox{$>$}\cr}(a+b^t)\cdot \log(a+b^t).$$

En particulier, si $m>n>1$ sont des entiers et $z>1$, alors

$$\eqalign{(z^n-1)^m&<(z^m-1)^n,\cr (z^m+1)^n&<(z^n+1)^m.\cr}$$

}

{\bf Preuve}

Tout le monde devrait savoir que $\left(f(t)^{g(t)}\right)'=f(t)^{g(t)}\cdot \left(g'(t)\log(f(t))+{g(t)f'(t)\over
f(t)}\right)$. Donc,

$$f'_{a,b}(t)={(a+b^t)^{1\over t}\over t}\cdot\left [{b^t\log(b)\over a+b^t}-{1\over t}\log(a+b^t)\right
].$$ Ainsi,

$$f'_{a,b}(t)\matrix{\raise-3pt\hbox{$<$}\cr \raise3pt\hbox{$>$}\cr}0\Longleftrightarrow {b^t\log(b)\over
a+b^t}\matrix{\raise-3pt\hbox{$<$}\cr \raise3pt\hbox{$>$}\cr}{1\over t}\log(a+b^t)\Longleftrightarrow
b^t\log(b^t)\matrix{\raise-3pt\hbox{$<$}\cr \raise3pt\hbox{$>$}\cr}(a+b^t)\cdot \log(a+b^t).$$

Prenons $a=-1$, $b=z>1$ et $t>1$. Alors $z^t>z^t-1>0$, donc $\log(z^t)>\log(z^t-1)$ et ainsi,
$z^t\log(z^t)>(z^t-1)\log(z^t-1)$. Cela implique que $f'_{-1,z}$ est croissante. Donc, si $m>n>1$, on a
$(z^n-1)^{1\over n}<(z^m-1)^{1\over m}$ ou encore $(z^n-1)^m<(z^m-1)^n$.

Posons maintenant $a=1$, $0<b={1\over z}<1$ et $t>1$. Alors $0<{1\over z^t}<1$. Donc ${1\over z^t}\log({1\over
z^t})<0<(1+{1\over z^t})\log(1+{1\over z^t})$. Cela implique que $f_{1,{1\over z}}$ est dŽcroissante. Donc, si
$m>n>1$, on a $(1+{1\over z^m})^{1\over m}<(1+{1\over z^n})^{1\over n}$ ou encore $(z^m+1)^n<(z^n+1)^m$.\qed

{\bf DŽfinition}

On rappelle que si $p$ est un nombre premier, tout nombre rationnel s'Žcrit $x=p^s\cdot u$ o $u$ est un nombre
rationnel dont $p$ ne divise ni le numŽrateur, ni le dŽnominateur. Dans ce cas, la {\it valuation
$p$-adique} de ce nombre, notŽe $v_(x)$ vaut $s\in\Z\cup\{\infty\}$, car $v_p(0)=\infty$, par dŽfinition.
On montre facilement que $v_p(a\cdot b)=v_p(a)+v_p(b)$ et $v_p(a+b)\geq \min(v_p(a),v_p(b))$ avec ŽgalitŽ
si $v_p(a)\ne v_p(b)$.

{\bf Lemme 2}

{\sl Soit $r,m$ et $n$, des entiers positifs non nuls et $p$ un nombre premier tel que $p\not
\hskip-0.4pt |\   n$. Alors

$$v_p(r!)\leq v_p\left ({m\over n}\left ({m\over n}-1\right )\cdots  \left ({m\over n}-(r-1)\right )\right).$$

}

{\bf Preuve}

Posons $a={m\over n}\left ({m\over n}-1\right )\cdots  \left ({m\over n}-(r-1)\right )$. Souvent, on note
${a\over r!}=\pmatrix{{m\over n}\cr r\cr}$. Si $a=0$, c'est Žvidemment vrai. Supposons donc $a\ne 0$. Posons
$e=v_p(a)$. Il est clair que $e$ est un entier supŽrieur ou Žgal ˆ 0, car $p\not \hskip-0.4pt |\   n$.
Mettant au mme dŽnominateur, on a $v_p(a)=v_p(m(m-n)(m-2n)\cdots (m-(r-1)n))$, car $v_p(xn)=v_p(x)$ pour
tout $x$, ceci car
$v_p(n)=0$. Choisissons $n'\in\N$, tel que
$n\cdot n'\equiv 1\pmod {p^{e+1}}$. En particulier, $p\not \hskip-0.4pt |\   n'$. Donc,
$v_p(a)=v_p(mn'(mn'-nn')\cdots (mn'-(r-1)mn'))$. Posons
$m'=m\cdot n'$. Si $1\leq j\leq (r-1)$, alors $m'-jnn'\equiv m'-j\pmod {p^{e+1}}$. Donc, 
$$m'(m'-nn')\cdots (m'-(r-1)nn')\equiv m'(m'-1)\cdots (m'-(r-1))\pmod {p^{e+1}}$$
Ainsi, $v_p(a)=v_p(m'(m'-1)\cdots (m'-(r-1))+k\cdot p^{e+1})=v_p(m'(m'-1)\cdots (m'-(r-1)))$. Mais $r!|
m'(m'-1)\cdots (m'-(r-1))$, le quotient est
$\pmatrix{m'\cr r\cr}$. Et on obtient $v_p(r!)\leq v_p(a)$.\qed

{\bf Lemme 3}

{\sl Soit $p>q>2$ des nombres premiers et $x,y$ des entiers strictement positifs tels que $x^p-y^q=\pm1$. Alors

$$(x\mp 1)^p\cdot q^{(p-1)q}>(y\pm 1)^q.$$

}

{\bf Preuve}

Il est clair que $x\geq 2$ (car sinon $y=0$ ou $q=1$), donc $x\pm1\buildrel (i) \over \geq {x\over 2}$. De mme,
$y\geq 2$,  donc $x^p=y^p\pm 1\buildrel (ii) \over > {y^q\over 2}$. Et aussi $y\buildrel (iii) \over
>{y\pm1 \over 2}$. Enfin, puisque $p$ et $q$ sont  premiers, et $p>q$, on a $(p-1)(q-1)\geq
(q+1)(q-1)=q^2-q+q-1=\underbrace{q(q-1)}_{\geq 6}+q+1>q+2$. Donc, $(p-1)q=(p-1)(q-1)+p-1>(q+2)+(p-1)=p+q+1$. Donc,
$q^{(p-1)q} \buildrel (iv) \over >2^{p+q+1}$. On obtient alors
$$(x\pm1)^p\buildrel (i) \over \geq \left ({x\over 2}\right )^p \buildrel (ii) \over >{y^q\over 2^{p+1}}\buildrel
(iii)\over > {(y\pm 1)^q\over 2^{q+p+1}}\buildrel (iv) \over > {(y\pm 1)^q\over q^{(p-1)q}}.$$
\qed
\medskip
{\bf Lemme 4}

{\sl Soit $a,b$ des entiers non nuls premiers entre eux et $q$ premier, alors il existe un entier $u$ tel que 

$${a^q- b^q\over a- b}=k(a- b)+qb^{q-1},$$

avec $k=(a-b)^{q-2}+ubq$. Remarquons que si $q=2$, alors $u=0$ et donc $k=1$.}

{\bf Preuve}

 On a $\dst{a^q- b^q\over a- b}={[(a- b)+ b]^q- b^q\over a-b}=(a- b)^{q-1}+ \pmatrix{q\cr
1\cr}b(a-b)^{q-2}+\cdots +\pmatrix{q\cr q-2\cr}b^{q-2}(a-b)+qb^{q-1}=k(a-b)+qb^{q-1}$, Avec
$k=(a-b)^{q-2}+\underbrace{\pmatrix{q\cr 1\cr}b(a-b)^{q-3}+\cdots +\pmatrix{q\cr q-2\cr}b^{q-2}}_{=uqb}$, 

car $q|\pmatrix{q\cr j\cr}$ si $1\leq j\leq q-1$.\qed
\medskip
{\bf Corollaire}

{\sl Soit $a,b$ des entiers non nuls premiers entre eux et $q$ premier, alors $\dst \left ({a^q-b^q\over
a-b},a-b\right )=(q,a-b)=1$ ou $q$. De plus, si ce pgcd vaut $q$ et que $q\geq 3$, alors $q|k$ et $q|\dst
{a^q-b^q\over a-b}$, mais $q^2\not \hskip-0.4pt |\  \dst {a^q-b^q\over a-b}$.

}
{\bf Preuve}

La premire partie a dŽjˆ ŽtŽ montrŽe lors du Petit Lemme du Chapitre prŽcŽdent. Supposons que ce pgcd soit
$q$. Le fait qu'alors $q|k$, vient du fait que $q|(a-b)$.  On peut poser $k=qk'$ et $(a-b)=qs$. Supposons
que $q^2|\dst {a^q-b^q\over a-b}=q(qsk'+b^{q-1})$, par le lemme prŽcŽdent. Ainsi, $q|b^{q-1}$, donc $q|b$,
et puisque $q|a-b$, $q|a$. Mais c'est impossible, car $a$ et $b$ sont premiers entre eux.\qed 

\bigskip

{\bf Lemme 5}

{\sl Soit $a>1$ et $b$ des nombres rŽels positifs. La fonction $f(x)=(1-{a\over b^x})^x$ est croissante si
${2a\over b^x}<1$. En particulier, 
$(1-{2\over 3^x})^x>{1\over 3}$  si $x\geq 2$.

}

{\bf Preuve}

Remarquons tout d'abord que 
$${2a\over b^x}<1\Longleftrightarrow {a\over b^x}<1-{a\over b^x}\eqno{(i)}$$

D'autre part,

$$f'(x)=f(x)\cdot \left (\log (1-{a\over b^x})+{a\log (b)x\over b^x(1-{a\over b^x})}\right )$$

Donc $f'(x)$ est positive si $\dst\log (1-{a\over b^x})>{-a\log (b)x\over b^x(1-{a\over b^x})}$, ou encore si
$ (1-{a\over b^x})\log (1-{a\over b^x})>{a\over b^x}\log({1\over b^x})$. Mais c'est le cas, car 

$${a\over b^x}\log\left ({1\over b^x}\right)<{a\over b^x}\log\left({a\over b^x}\right)<(1-{a\over b^x})\log
(1-{a\over b^x})$$

par $(i)$ et car $x\log(x)$ est croissante. La seconde partie du Lemme rŽsulte du fait que ${4\over 3^x}<1$ si
$x\geq 2$ et $(1-{2\over 9})^2={49\over 81}>{1\over 3}$.\qed

\bigskip
\goodbreak
{\bf Lemme 6}

{\sl Soit $p$ un nombre premier  et $n$ un entier naturel. Le dŽveloppement $p$-adique de $n$ est $n=\sum_{i=0}^k
n_i p^i$ avec $0\leq n_i\leq p-1$ pour $i=0,\ldots , k=\left [{\log(n)\over \log (p)}\right ]$
($\left[x\right ]$ signifie la partie entire de $x$). Posons encore
$S=\sum_{i=0}^k n_i$. Alors on a 
$$v_p(n !)=\sum_{i=1}^\infty \left [{n\over p^i}\right ]={n-S\over p-1}.$$

Remarquons que la somme $\sum_{i=1}^\infty \left [{n\over p^i}\right ]$ s'arrte en fait ˆ $k$.
}

{\bf Preuve}

Les multiples de $p$ infŽrieurs ˆ $n$ sont 

$$p,2p,\ldots , \left [{n\over p}\right ] p.$$

parmi ces multiples, on doit compter dans $v_p(n !)$, deux fois les multiples de $p^2$ qui sont 
$$p^2,2p^2,\ldots , \left [{n\over p^2}\right ] p^2.$$

Et ainsi de suite. La premire ŽgalitŽ est donc prouvŽe. Pour la deuxime, il est clair que  pour tout $j$, $ \left
[{n\over p^j}\right ]=\sum_{i=j}^k n_i p^{i-j}$. Donc

$$\eqalign{\sum_{j=1}^k \left [{n\over p^j}\right ]&=\sum_{j=1}^k\sum_{i=j}^k n_i
p^{i-j}=\sum_{j=1}^k(n_jp^0+n_{j+1}p^1+\cdots +n_k p^{k-j})\cr
&=\sum_{j=1}^k n_j p^0+\sum_{j=2}^kn_j p^1+\sum_{j=3}^k n_j p^2+\cdots +\sum_{j=k}^k n_j p^{k-1}\cr
&=\sum_{j=1}^k n_j \sum_{s=0}^{j-1} p^s=\sum_{j=1}^k n_j {p^j-1\over p-1}={1\over p-1}\cdot \left ( \sum_{j=1}^k
n_j p^j-\sum_{j=1}^k n_j\right )\cr
&={1\over p-1}\cdot (n-S).\cr}$$

\qed

{\bf ThŽorme (Cassels, 1953, 1961)}

{\sl Soient $p,q$ des nombres premiers impairs et $x,y>0$ des entiers naturels tels que $x^p-y^q=\pm 1$, alors
$q|x$ et
$p|y$.}

{\bf Preuve}

On supposera que $p> q$ (le cas $p=q$ Žtant Žvidemment impossible). On va d'abord voir que $q|x$.

Supposons  que $q\not \hskip-0.4pt |\  x$, alors $q\not \hskip-0.4pt |\   y^q\pm 1$. Ainsi, en vertu du
Lemme 4 et de son corollaire, on en dŽduit que $\dst \left ({y^q\pm 1\over y\pm 1},y\pm 1\right )=(q,y\pm
1)=1$. De l'ŽgalitŽ 
$${y^q\pm 1\over y\pm 1}\cdot (y\pm1)=y^q\pm 1=x^p,$$

on dŽduit qu'il existe $b$ entier tel que $y\pm 1=b^p$.

\art{a)} (signe du haut) Si $y+1=b^p$, alors $b\geq 2$. Et donc, $x^p=y^q+1=(b^p-1)^q+1<b^{pq}$ (cette dernire
inŽgalitŽ vient du fait que $X^q-(X-1)^q=X^{q-1}+X^{q-2}(X-1)+\cdots +X(X-1)^{q-2}+(X-1)^{q-1}>1$ si $X\geq 2$ et
$q\geq 3$ et on pose $X=b^p$). Donc, $x<b^q$ et alors $x\leq b^q-1$. Par le Lemme 1, $(b^q-1)^p<(b^p-1)^q$, car
$p>q$. Finalement, $y^q+1=x^p\leq (b^q-1)^p<(b^p-1)^q=y^q$, ce qui est absurde. 

\art{b)} (signe du bas) Si $y-1=b^p$, alors $b\geq 2$ aussi, car $x^p=y^q-1$ et $p>q$, donc $y>x$ et $y\geq 3$. On
a alors  $x^p=(b^p+1)^q-1>b^{pq}$, ce qui veut dire $x>b^q$ ou $x\geq b^q+1$. En utilisant ˆ nouveau le Lemme 1,
on a $y^q-1=x^p\geq (b^q+1)^p>(b^p+1)^q=y^q$, ce qui est ˆ nouveau absurde.

On a montrŽ que $q|x$.

Montrons que $p|y$. Lˆ, ce sera un peu plus long, mais on a le temps hein ? Remarquons un certain nombre de
chose avant de supposer que
$p$ ne divise pas $y$. Tout d'abord, $y^q\geq 8$, car $y\geq 2$ est $q\geq 3$. Ensuite, $q\leq x$, car on vient de
voir que $q|x$. Donc, $q^p|x^p=y^q\pm 1=(y\pm 1){y^q\pm 1\over y\pm 1}$. On a vu au Lemme 4 et ˆ son
corollaire que  $\dst \left ({y^q\pm 1\over y\pm 1},y\pm 1\right )=1$ ou $q$. Si ce pgcd est 1, on a
$q\not \hskip-0.4pt |\  y\pm 1$. Donc,
$q^p| {y^q\pm 1\over y\pm 1}$ et on sait, par le Lemme 4 (en posant $a=\pm 1$ et $b=-y$), que $${y^q\pm
1\over y\pm 1}=k(y\pm 1)+qy^{q-1}.\eqno{(*)}$$
Donc $q|k$. Le mme rŽsultat nous dit que $k=(y\pm 1)^{q-2}- quy$. En rŽduisant modulo
$q$, on en dŽduit que $q|(y\pm 1)$, contradiction. Donc $\dst \left ({y^q\pm 1\over y\pm 1},y\pm 1\right )=q$.
Le corollaire du Lemme 4 nous apprend que $q^2\not \hskip-0.4pt |\   {y^q\pm 1\over y\pm 1}$. Il existe
donc $b,c>0$ tels que $(b,c)=1$ et
$$\left\{\eqalign{y\pm 1&=q^{p-1} b^p\cr {y^q\pm 1\over y\pm 1}&=qc^p\cr }\right. \hbox{ avec $q\not
\hskip-0.4pt |\   c$ et
$x=qbc$} \eqno{(**)}$$
Concernant $c$, on peut remarquer que~:
\art{a)} $c\ne 1$. Sinon, on aurait $y\pm 1=q^{p-1} b^p>q={y^q\pm 1\over y\pm 1}$. Donc $(y\pm 1)^2>y^q\pm 1\geq
y^3\pm1$. C'est impossible si $y\geq 2$, en effet, dans le cas du signe du haut, on aurait
$y^2+2y+1>y^q+1\geq y^3+1$; donc, $y(y-2)(y+1)<0$ ce qui est impossible si $y\geq 2$. Dans le cas du signe du bas,
c'est encore plus clair, on aurait $y^2-2y+1>y^3-1$, donc $(y^2+2)(y-1)<0$, ˆ nouveau impossible.

\art{b)} $c\equiv 1\pmod {q^{p-1}}$. En rŽutilisant $(*)$, on a $qc^p={y^q\pm 1\over y\pm 1}=k(y\pm 1)+q y^{q-1}$.
On a vu il y a 10 lignes que $q|k$, en plus $q^{p-1}|y\pm 1$. En divisant par $q$, on obtient $c^p\equiv
y^{q-1}\pmod{q^{p-1}}$. Puisque $y\equiv \mp 1\pmod{q^{p-1}}$ et puisque $q-1$ est pair, on a $c^p\equiv 1\pmod
{q^{p-1}}$. L'ordre de $c$ modulo $q^{p-1}$ est donc 1 ou $p$. Si c'est $p$, alors $p$ divise l'ordre du groupe
des inversibles modulo $q^{p-1}$ qui vaut $\varphi(q^{p-1})=q^{p-2}(q-1)$. Donc $p|q-1$, c'est impossible, car
$p>q$. Donc $c\equiv 1\pmod {q^{p-1}}$ ou encore $qc\equiv q\pmod{q^p}$ (la fonction $\varphi(\cdot)$ est
dŽfinie au chapitre 5, pp 18-19).

De cette analyse de $c$, on tire facilement (mais vraiment) 
$$x\ne qb\hbox{ et } x\equiv  qb\pmod{q^p}.\eqno{(***)}$$

A partir de maintenant, on va supposer (par l'absurde) que $p\not \hskip-0.6pt |\  y$. On a $\dst y^q=(x\mp
1){x^p\mp 1\over x\mp 1}$, et on a $\left ( x\mp 1, {x^p\mp 1\over x\mp 1}\right )=(x\mp 1,p)=1$ ou $p$. Il
n'est pas possible que ce pgcd vaille $p$, car sinon $p$ diviserait $y$, donc il vaut 1. Il existe ainsi
$a>1$  tel que  $x\mp 1=a^q$. On trouve alors
$$a^{pq}=(x\mp 1)^p\buildrel {\rm Lemme\ 3}\over >{(y\pm 1)^q\over q^{(p-1)q}}\buildrel(**)\over = b^{pq}.$$ 

Donc $a>b$. Par $(***)$, on a $q^p\leq |x-qb|=|a^q\pm 1-qb|\leq a^q+qb\pm 1$. 

Si on avait  $a^q<{1\over 2} q^p$, alors on aurait $qb\pm1 >{1\over 2} q^p$, en particulier on aurait $b\geq
2$, et alors $b^q\geq qb+1$ et donc $a^q>b^q\geq qb+1\geq qb\pm 1>{1\over 2} q^p$, ce qui est contradictoire.
On a donc montrŽ que 

$$a^q\geq {1\over 2} q^p.\eqno{(+)}$$

Par suite $x^p=(a^q\pm 1)^p\geq (a^q-1)^p$ et $y^q=x^p\mp 1=(a^q\pm 1)^p\mp 1\geq (a^q-1)^p$. Remarquons
d'autre part que $(1-{2\over q^p})^p\geq (1-{2\over 3^p})^p\buildrel {\rm Lemme\ 5}\over > {1\over 3}\geq
{1\over q}$. On obtient

$$\min (x^p,y^q)\geq (a^q-1)^p=a^{pq}(1-{1\over a^q})^p\buildrel (+)\over \geq a^{pq}(1-{2\over
q^p})\geq {a^{pq}\over q}.\eqno{(i)}$$

L'ŽgalitŽ $(x^{p\over q}-y){(x^{p\over q})^q-y^q\over x^{p\over q}-y}=x^p-y^q=\pm 1$ entra"ne

$$|x^{p\over q}-y|={1\over |\sum_{i=0}^{q-1} x^{p\cdot i\over q} y^{q-1-i}|}.\eqno{(ii)}$$ 

Pour $i=0,\ldots q-1$, on observe et cela n'est pas de la magie que  

$$x^{p\cdot i\over q}y^{q-1-i}\buildrel {(i)}\over \geq \left ({a^{pq}\over q}\right )^{{i\over q}+{q-1-i\over
q}}=a^{p(q-1)}\cdot{1\over q^{{q-1}\over q}}>a^{p(q-1)}\cdot {1\over q}.$$

Ainsi, l'Žquation $(ii)$ devient, via $q$,

$$|x^{p\over q}-y|<{1\over a^{p(q-1)}}.\eqno{(iii)}$$

Ecrivons 
$$x^{p\over q}=(a^q\pm 1)^{p\over q}=a^p(1\pm {1\over a^q})^{p\over q}=\sum_{r=0}^\infty t_r,\eqno{(iv)}$$

o, par le dŽveloppement de Taylor, 

$$t_r=(\pm 1)^r{{p\over q}({p\over q}-1)\cdots ({p\over q}-r+1)\over r!}\cdot a^{p-rq}\ne 0.\eqno{(v)}$$

A noter que $t_0=a^p$. Soit $l$ un nombre premier diffŽrent de $q$ et $r\geq 1$. On sait, par le Lemme 2 que
$v_l(r!)\leq v_l({p\over q}({p\over q}-1)\cdots ({p\over q}-r+1))$. D'o

$$v_l(t_r)\geq v_l(a^{p-rq})\quad\hbox{c'est aussi vrai si $r=0$.}\eqno{(vi)}$$

Posons $R=\left [{p\over q}\right]+1$ et $\rho=\left [{R\over q-1}\right ]$. En particulier, $Rq>p$. Par le Lemme
6, $v_q(R!)={R-S\over q-1}$ o $S$ est la somme des chiffres du dŽveloppement $q$-adique de $R$. D'o 

$$v_q(R!)<{R\over q-1}\hbox{ et donc }v_q(R!)\leq \rho.\eqno{(vii)}$$

Si $0\leq r\leq R$, on calcule pour $l\ne q$, 

$$v_l(t_r\cdot q^{R+\rho}\cdot a^{Rq-p})\buildrel (vi)\over \geq
v_l(a^{p-rq})+v_l(a^{Rq-p})=v_l(a(R-r))\geq 0.$$

Et aussi,

$$\eqalign{v_q(t_r\cdot q^{R+\rho}\cdot a^{Rq-p})&=\underbrace{v_q({p\over q}({p\over q}-1)\cdots ({p\over
q}-r+1))}_{\geq -r}-v_q(r!)+R+\rho+(Rq-p)v_q(a)\cr
&\geq R-r+(\rho- v_q(R!))+(Rq-p)v_q(a)\geq 0.\cr}$$

On a utilisŽ que $q$ ne divisait pas $a^{p-qr}$, car sinon il diviserait $y$ et comme il divise $x$, cela voudrait
dire que $q$ divise $\pm 1$, ce qui est absurde ! On a aussi utilisŽ $(vii)$ et le fait que $v_q(r!)\leq
v_q(R!)$. Remarquons qu'un nombre rationnel qui a des valuations positives pour tout premier est forcŽment
entier. On a donc prouvŽ que
$t_r\cdot q^{R+\rho}\cdot a^{Rq-p}\in\Z$, ce qui montre que 

$$\eqalign{I&:=a^{Rq-p}\cdot q^{R+\rho}\cdot \left ((y-x^{p\over q})+\sum_{r\geq R+1}t_r\right )\cr
&=a^{Rq-p}\cdot q^{R+\rho}\cdot \left (y-\sum_{r=0}^R t_r\right )\in\Z.\cr}$$

On va passer un bon moment ˆ prouver que $I\ne 0$ et c'est cela qui nous donnera la contradiction.

On Žcrit $I=I_1+I_2+I_3$ o

$$\eqalign{I_1&=a^{Rq-p}\cdot q^{R+\rho}\cdot(y-x^{p\over q})\cr
I_2&=a^{Rq-p}\cdot q^{R+\rho}\cdot t_{R+1}\ne 0\cr
I_3&=a^{Rq-p}\cdot q^{R+\rho}\cdot \sum_{r>R+1} t_r.\cr}$$

On va montrer que $\left |{I_3\over I_2}\right |<{1\over 10}$ et  $\left |{I_1\over I_2}\right |<{1\over 10}$.
Montrons la premire inŽgalitŽ, ce ne sera pas trop dur~:

Soit $r>R$. On a $\left |{p\over q}-r\right |=r-{p\over q}<r+1$. Donc $\left | {t_{r+1}\over t_r}\right |=\left
|{{p\over q}-r\over r+1}\right |\cdot{1\over a^q}<{1\over a^q}\buildrel(+)\over \leq {2\over q^p}$.
\goodbreak
Ainsi,

$$\eqalign{\left |{I_3\over I_2}\right |&=\left |\sum_{r>R+1} {t_{r+1}\over t_{R+1}}\right |\leq
\sum_{r>R+1}\left |{t_{r+1}\over t_{R+1}}\right |\cr
&=\left | {t_{R+2}\over t_{R+1}}\right |+\left | {t_{R+3}\over t_{R+2}}\right |\cdot \left | {t_{R+2}\over
t_{R+1}}\right |+\cdots\cr
&<\left ({2\over q^p}\right)+\left ({2\over q^p}\right)^2+\cdots ={2\over q^p}\cdot{1\over {1-{2\over
q^p}}}={2\over q^p-2}\cr
&\leq {1\over 3^5-2}<{1\over 10}.\cr}$$

En avant pour la seconde inŽgalitŽ. Observons les choses suivantes~: $(R-{p\over q})+({p\over q}+1-R)=1$ et chacun
des termes de cette somme est strictement positif, on en dŽduit que $|{p\over q}-R|\cdot |{p\over q}-R+1|\leq
{1\over 4}$. En effet~: si $x+y=1, x, y>0$, alors $xy\leq {1\over
4}\Leftrightarrow x(1-x)-{1\over 4}\leq 0\Leftrightarrow 4x^2-4x+1=(2x-1)^2\geq 0$ et c'est trivial. Donc,

$$\left |{p\over q}\left ({p\over q}-1\right )\cdots \left ({p\over q}-R\right )\right |\leq R(R-1)\cdots 2\cdot
\left |{p\over q}-R+1\right |\cdot \left |{p\over q}-R\right |\leq  R!\cdot {1\over 4}.\eqno{(viii)}$$
 D'autre part,

$$\left |{p\over q}\left ({p\over q}-1\right )\cdots \left ({p\over q}-R\right )\right |\geq (R-1)(R-2)\cdots
1\cdot \left |{p\over q}-R+1\right |\cdot \left |{p\over q}-R\right |\geq (R-1)!\cdot {1\over q^2},$$
car $R-{p\over q}\geq {1\over q}$ et ${p\over q}-(R-1)\geq {1\over q}$. Donc,

$$|t_{R+1}|=\left |{{p\over q}({p\over q}-1)\cdots ({p\over q}-R)\over (R+1)!}\right |\cdot a^{p-(R+1)q}\geq
{a^{p-(R+1)q}\over q^2 \cdot R\cdot (R+1)}.\eqno{(ix)}$$

Calculons~:

$$\left |{I_1\over I_2}\right |=\left |{y-x^{p\over q}\over t_{R+1}}\right|\buildrel (ix)\over \leq
a^{(R+1)q-p}q^2\cdot R(R+1)\cdot
\left |y-x^{p\over q}\right |\buildrel (iii)\over < {a^{(R+1)q-p}q^2(R+1)^2\over a^{p(q-1)}}={q^2(R+1)^2\over
a^{(p-R-1})q}.$$

Or,

$$p-R-1\geq 2\quad\hbox{ et donc aussi }  R+1\leq p.\eqno{(x)}$$

On va montrer cela (ce ne sera pas super ŽlŽgant, mais enfin..) $p-R-1=p-\left [p\over q\right ]-2\geq 2$ si
$p-\left [p\over q\right ]\geq 4$. Distinguons trois cas~: si $p=5$, et donc $q=3$ alors $p-\left [p\over q\right
]= 4$ et c'est en ordre. Si $p= 7$, alors $q=3$ ou $5$, et $p-\left [p\over q\right
]= 5$ ou 6 et c'est ˆ nouveau bon. Si $p\geq 11$, alors $p\left ({q-1\over q}\right )\geq p^{2\over
3}\Leftrightarrow p^{1\over 3}\geq {q\over q-1}$ et on a $p^{1\over 3}\geq 11^{1\over 3}\geq {3\over 2}$. Donc,
dans notre cas, $p\left ({q-1\over q}\right )\geq p^{2\over 3}\geq 4$. D'o, $p-\left [p\over q\right
]\geq p-{p\over q}=p\left ({q-1\over q}\right )\geq 4$. Donc $(x)$ est montrŽ.

On peut ainsi terminer de majorer $\left |{I_1\over I_2}\right |$~:

$$\left |{I_1\over I_2}\right |\leq {q^2(R+1)^2\over a^{2q}}\buildrel (+)\over <{q^2(R+1)^2\over ({1\over 2}
q^p)^2}\buildrel (x)\over \leq \left ( {2p\over q^{p-1}}\right )^2\leq  \left ( {2p\over 3^{p-1}}\right )^2\leq 
\left ( {2\cdot 5\over 3^{4}}\right )^2\leq {1\over 10}.$$

L'avant-dernire ŽgalitŽ venant du fait que la fonction $\left ({2x\over 3^{x-1}}\right )^2$ dŽcro"t si $x\geq 5$.
En effet,  sa dŽrivŽe vaut ${8x\over 3^{2x-2}} (1-\log(3)x)$.

On en dŽduit que $I$ est diffŽrent de zŽro. Si c'Žtait le cas, on aurait $1={-I_1\over I_2}+{-I_3\over I_2}$ et
donc 

$$1=\left |{-I_1\over I_2}+{-I_3\over I_2}\right |\leq \left |{I_1\over I_2}\right |+\left |{I_3\over I_2}\right
|<{1\over 10}+{1\over 10},$$

ce qui est absurde, donc $I\ne 0$. Or, on a montrŽ que $I$ Žtait entier,  donc  $|I|\geq 1$.

D'autre part,

$$\left |a^{(R+1)q-p}\cdot t_{R+1}\right |=\left |{{p\over q}\left ({p\over q}-1\right)\cdots \left ({p\over
q}-R\right )\over (R+1)!}\right |\buildrel (viii)\over\leq {1\over 4(R+1)}\leq {1\over 4}.\eqno{(xi)} $$

On en dŽduit,

$$I_2=\left |{ q^{R+\rho}\cdot a^{(R+1)q-p}\cdot t_{R+1}\over a^q}\right |\buildrel (xi)\over \leq {
q^{R+\rho}\over 4a^q}\buildrel(+)\over \leq {1\over 2} q^{R+\rho-p}.\eqno{(xii)}$$

Donc,

$$1\leq |I|=|I_2|\cdot\left |1+{I_1\over I_2}+{I_3\over I_2}\right |\buildrel (xii)\over \leq {1\over
2}q^{R+\rho-p}\left (1+{1\over 10}+{1\over 10}\right )<q^{R+\rho-p}.$$

Ce qui prouve que $\thboxed 15{R+\rho-p>0}$ Mais c'est imppossiiiiiible~: on se souvient que $\rho=\left [{R\over
q-1}\right]$, donc

$$R+\rho\leq R\left(1+{1\over q-1}\right)\leq \left({p\over q}+1\right)\left ({q\over q-1}\right)={p+q\over
q-1}<{2p\over q-1}\leq p.$$

Ce qui donne $\thboxed 15{R+\rho-p<0}$, et voilˆ enfin notre contradiction !!! Donc $p|y$.\qed

\bigskip

{\bf Corollaire (les relations de Cassels)}

{\sl Soient $p,q$ des nombre premiers impairs et $x,y>0$ des entiers tels que $x^p-y^q=\pm 1$, alors

\art{a)}Il existe $b,c>0$ tels
que $(b,c)=1$ et
$$x=qbc\quad y\pm 1=q^{p-1} b^p\quad  {y^q\pm 1\over y\pm 1}=qc^p\quad \hbox{ avec $q\not
\hskip-0.4pt |\   c$ }.$$

\art{b)}Il existe $u,v>0$ tels
que $(u,v)=1$ et
$$y=puv\quad x\mp 1=p^{q-1} u^q\quad  {x^q\mp 1\over x\mp 1}=pv^q\quad \hbox{ avec $p\not \hskip-0.4pt
|\   v$ }.$$

}

{\bf Preuve}

La partie a) est la relation $(**)$ du thŽorme prŽcŽdent et la partie b) se dŽmontre de manire identique sachant
que $p|y$.\qed

\vfill\eject

\long\def\art#1{{\parindent0pt\item{#1}}\hangindent=7mm\hangafter=-20}
\long\def\artart#1{{\parindent0pt\item{#1}}\hangindent=12mm\hangafter=-20}
\font\para=cmbx12 at 18pt
\def\O{\hbox{$\cal O$}}
\def\U{\hbox{$\cal U$}}
\def\m{\hbox{\rs m\!}}
\def\dst{\displaystyle}
\font\doub=msbm10 at 10pt
\def\lra{\longrightarrow}
\def\qed{\hfill$\square$}
\def\gfP{\relax\ifmmode\bbP\else $\bbP$\fi}
\def\gP{{\euf P}}
\def\P{{\cal P}}

\def\ggP{{\bf P}}
\newcount\chapnomb \chapnomb=1
\newcount\parnomb \parnomb=1
\pageno =16

\parindent0pt
\bigskip

\centerline{\para CHAPITRE 5}
\bigskip
{\para Le thŽorme de Stickelberger}
\bigskip

Ce chapitre est le plus long de tous. Nous nous sommes basŽs en partie sur le livre de Lemmermeyer [Lem]
pour Žcrire ce chapitre. On va montrer que si
$p$ est un nombre premier, si
$G=\{\sigma_1,\ldots ,\sigma_{p-1}\}$ est le groupe de Galois de l'extension $\Q(\zeta_p)/\Q$, et si
$-\Theta_2=-\sum_{t={p+1\over 2}}^{p-1}\sigma_t^{-1}\in \Z[G]$, alors pour tout idŽal $\euf a$ de
$\Q(\zeta_p)$, l'idŽal ${\euf a}^{-\Theta_2}$ est principal. Ceci sera utilisŽ au Chapitre 6. Nous
dŽfinirons aussi l'idŽal de Stickelberger $I_{st}$ et l'idŽal $I=(1-\iota) I_{st}$ dont nous prouverons
qu'ils sont des $\Z$-modules libres de rang ${p+1\over 2}$ et $p-1\over 2$ respectivement. Nous donnerons
explicitement des gŽnŽrateurs de ces idŽaux. Nous retrouverons ces objets aux Chapitres 7 et 8. Si vous
n'avez pas compris ce qui vient d'tre dit, c'est normal, on n'a encore rien dŽfini. Tout d'abord quelques
rappels sur les caractres.  

\bigskip

{\bf DŽfinitions et rappels}

Soit $G=(G,*,1)$ un groupe abŽlien fini. On appelle {\it caractre} de $G$ tout homomorphisme
$\chi \,  : G\,
\lra \C^*$. On note ${\bf 1}$ le caractre de $G$ qui envoie tout ŽlŽment de $G$ sur $1$.
L'ensemble des caractres est lui-mme un groupe isomorphe ˆ $G$~:
$\chi_1\chi_2(g):=\chi_1(g)\cdot\chi_2(g)$, donc par le thŽorme bien connu de Lagrange, l'ordre
d'un caractre divise l'ordre de $G$. Remarquons que l'on a 

$$\sum_{g\in G}\chi(g)=\cases{0& si $\chi\ne {\bf 1}$\cr |G|& si $\chi= {\bf 1}.$\cr}\eqno{(5)}$$

En effet, c'est clair si $\chi={\bf 1}$. Supposons que $\chi\ne {\bf 1}$, donc il existe $h\in G$ tel que
$\chi(h)\ne 1$. L'application $g\mapsto hg$ est clairement une bijection de $G$ dans lui-mme. Donc, 

$$\sum_{g\in G}\chi(g)=\sum_{g\in G}\chi(hg)=\chi(h)\cdot \sum_{g\in G}\chi(g).$$

Si $\sum_{g\in G}\chi(g)\ne 0$, on en dŽduirait que $\chi(h)=1$, ce qui est contradiction.

\bigskip

On fixera pour un moment $\F=\F_q$ un corps fini ˆ $q$ ŽlŽments. Il est clair que la
{\it caractŽristique de $\F$}, qui est le plus petit entier $n$ tel que $\underbrace{1+1+\cdots +
1}_{n\ \rm fois}=0$, est un nombre premier, disons $p$. Donc, le corps $\F_p=\Z/p\Z$ agit sur $\F$,
faisant de lui un $\F_p$-espace vectoriel de dimension Žvidemment finie. Donc $q=p^f$ o $f$ est
la dimension de $\F$ sur $\F_q$. De plus le sous-ensemble $\{0,1,1+1,\ldots ,
\underbrace{1+1+\cdots + 1}_{p-1\ \rm fois}\}$ est un sous-corps de $\F$ isomorphe ˆ $\F_p$ qu'on
notera par abus encore $\F_p$.

{\bf Lemme 1}

{\sl $\F^*=\F\setminus\{0\}$ est {\it cyclique}, c'est ˆ dire qu'il existe $c\in \F^*$ tel que
$\F^*=\{c,c^2,c^3,\ldots ,c^{q-1}=1\}$}

{\bf Preuve}

C'est un rŽsultat bien connu, cf. [Jac1, Theorem 2.18, p.132].\qed
\bigskip
Sur $\F$ on peut dŽfinir deux sortes de caractres~: des {\it caractres additifs de $\F$}
($G=(F,+,0)$) ou des {\it caractres multiplicatifs de $\F$} ($G=(\F^*,\cdot,1)$). Par convention,
si $\chi$ est un caractre multiplicatif de $\F$, on le prolonge ˆ $\F$ tout entier en posant
$\chi(0)=0$ si $\chi\ne {\bf 1}$ et $\chi(0)=1$, si $\chi={\bf 1}$.

On dŽfinit la {\it trace} de $\F$ sur $\F_p$ comme Žtant l'application $Tr_{\F/\F_p}$ notŽe $Tr$ qui
envoie tout
$t\in \F$ sur $t+t^p+t^{p^2}+\cdots +t^{p^{f-1}}\in \F_p$. Le fait que
$Tr(t)\in \F_p$ n'est pas trivial, mais bien connu (cf. [Sam, ¤2.6]). De plus,
$Tr(t_1+t_2)=Tr(t_1)+Tr(t_2)$ pour tout $t_1,t_2$. 

Soit $m\in\N$. On note $\zeta_m=e^{2i\pi\over m}$ qui est une racine primitive $m$-ime de
l'unitŽ.

On notera $\psi$ le caractre additif sur $\F$, dŽfini par $\psi(t)=\zeta_p^{Tr(t)}$. Soyons
attentifs au fait que puisque le caractre est additif, on a Žvidemment
$\psi(t+s)=\psi(t)\cdot\psi(s)$, c'est redondant, mais, on peut se faire avoir si on n'y prend
pas garde... . D'autre part, par (5) et puisque $\psi\ne {\bf 1}$, on a $\sum_{t\in\F}\psi(t)=0$.
De mme, si $\chi$ est un caractre multiplicatif diffŽrent de ${\bf 1}$, alors
$\sum_{t\in\F^*}\chi(t)=\sum_{t\in\F}\chi(t)=0$. 

\goodbreak
Si $\chi$ est un caractre multiplicatif, on dŽfinit la {\it somme de Gauss}

$$G(\chi)=-\sum_{t\in \F^*}\chi(t)\cdot\psi(t).$$

On observera que $G({\bf 1})=-\sum_{t\in \F^*}\psi(t)=1-\sum_{t\in \F}\psi(t)\buildrel (5)\over=1$. De
manire gŽnŽrale, si $\chi$ est un caractre d'ordre $m$, alors $G(\chi)\in
\Q(\zeta_m,\zeta_p)=\Q(\zeta_{mp})$, car $p$ et $m$ sont premiers entre eux (en effet, $m||\F^*|=q-1$).
\goodbreak
Si $\chi_1$ et $\chi_2$ sont des caractres multiplicatifs de $\F$, alors on dŽfinit la {\it somme de
Jacobi}

$$J(\chi_1,\chi_2)=-\sum_{t\in \F}\chi_1(t)\chi_2(1-t).$$
\goodbreak
\bigskip
{\bf Lemme 2}

{\sl 

\art{a)}Si $\chi_1$ et $\chi_2$ sont des caractres multiplicatifs de $\F$ diffŽrents de ${\bf 1}$
et tels que $\chi_1\chi_2\ne{\bf 1}$, alors

$$G(\chi_1)G(\chi_2)=G(\chi_1\chi_2)\cdot J(\chi_1,\chi_2)$$

\art{b)}Si $\chi$ est un caractre multiplicatif de $\F$ diffŽrent de ${\bf 1}$, alors on a~:

$$G(\chi)G(\chi^{-1})=q\cdot\chi(-1).$$

}

{\bf Preuve}

Remarquons que si $\chi\ne {\bf 1}$, alors $G(\chi)=-\sum_{t\in F}\chi(t)\cdot\psi(t)$, le fait de
prendre $0$ n'ajoute rien par la convention qu'on s'est donnŽe. Calculons donc~:

$$\eqalignno{G(\chi_1)G(\chi_2)&=\sum_{a,b\in\F}\chi_1(a)\chi_2(b)\psi(a+b)\buildrel c=a+b\over =
\sum_{a,c\in\F}\chi_1(a)\chi_2(c-a)\psi(c)\cr
&=\sum_{a\in\F}\chi_1(a)\chi_2(-a)+\sum_{a,c\in\F^*}\chi_1(a)\chi_2(c-a)\psi(c)&(i)\cr}$$

Or, d'une part

$$\sum_{a,c\in\F^*}\chi_1(a)\chi_2(c-a)\psi(c)\buildrel a=c\cdot t\over
=\sum_{t,c\in\F^*}\underbrace{\chi_1(c)\chi_2(c)}_{=\chi_1\chi_2(c)}\psi(c)\chi_1(t)\chi_2(1-t)=G(\chi_1\chi_2)\cdot
J(\chi_1,\chi_2)\eqno{(ii)}$$

et d'autre part, si $\chi_1\ne\chi_2^{-1}$, alors

$$\sum_{a\in\F}\chi_1(a)\chi_2(-a)=\chi_2(-1)\cdot\sum_{a\in\F}\chi_1\chi_2(a)\buildrel (5)\over =0.$$
On a donc prouvŽ la partie a).

Pour la partie b), on commence comme la partie a), et on retrouve les ŽgalitŽs $(i)$ et $(ii)$~:

$$G(\chi)G(\chi^{-1})=\sum_{a\in\F}\chi(a)\chi^{-1}(-a)+G(\chi\chi^{-1})\cdot J(\chi,\chi^{-1})$$

D'abord, il est clair que $\chi^{-1}(-1)=\chi(-1)=\pm 1$, car $\chi^{-1}(-1)\chi(-1)=1=\chi^2(-1)$.
Ainsi,

$$\sum_{a\in\F}\chi(a)\chi^{-1}(-a)=\chi^{-1}(-1)\cdot\sum_{a\in\F^*}\chi(a)\chi^{-1}(a)=\chi(-1)\cdot
\sum_{a\in\F^*}{\bf 1}(a)=\chi(-1)\cdot(q-1).$$

Ensuite, $G(\chi\chi^{-1})=G({\bf 1})=1$, comme on l'a vu en introduction. Finalement,

$$J(\chi,\chi^{-1})=-\sum_{t\in\F}\chi(t)\chi^{-1}(1-t)=-\sum_{t\in\F\setminus\{1\}}\chi\left({t\over
1-t}\right )\buildrel
(*)\over=-\sum_{s\in\F\setminus\{-1\}}\chi(s)=\chi(-1)-\sum_{t\in\F}\chi(s)\buildrel (5)\over =\chi(-1).$$

\art{(*)}Vient du fait que si $t$ parcourt $\F\setminus\{1\}$, alors ${t\over 1-t}$ parcourt
$\F\setminus\{-1\}$. On a donc montrŽ que  $G(\chi)G(\chi^{-1})=(q-1)\chi(-1)+\chi(-1)=q\cdot\chi(-1)$. \qed
\bigskip
\goodbreak
{\bf Lemme 3}

{\sl Soit $\chi\ne{\bf 1}$ un caractre multiplicatif de $\F$. Alors on a

\art{a)}$\overline{G(\chi)}=\chi(-1)G(\chi^{-1})$

\art{b)}$G(\chi)\overline{G(\chi)}=q$

\art{c)}si $\chi$ est d'ordre $m$, alors $G(\chi)^m\in\Q(\zeta_m)$.

}

{\bf Preuve}

\art{a)}Il est clair que puisque $\chi(t)$ et $\psi(t)$ sont des racines de l'unitŽ, on a 
$\chi(t)\overline{\chi(t)}=\psi(t)\overline{\psi(t)}=1$. Donc,
$\overline{\chi(t)}=\chi(t)^{-1}=\chi^{-1}(t)$ et $\overline{\psi(t)}=\psi(t)^{-1}=\psi(-t)$. On a
donc, 

$$\overline{G(\chi)}=-\sum_{t\in \F^*}\overline{\chi(t)}\cdot\overline{\psi(t)}=-\sum_{t\in
\F^*}\chi^{-1}(t)\cdot\psi(-t)=-\sum_{t\in \F^*}\chi^{-1}(-t)\cdot\psi(t)=\chi(-1)G(\chi^{-1}).$$

\art{b)}$G(\chi)\overline{G(\chi)}\buildrel a)\over=G(\chi)\chi(-1)G(\chi^{-1})\buildrel {\rm Lemme\
2\ b)}\over =q\chi(-1)\chi(-1)=q$.

\art{c)}Remarquons que si $\chi_1,\chi_2$ sont d'un ordre qui divise $m$, alors
$J(\chi_1,\chi_2)\in\Q(\zeta_m)$. Cela dŽcoule directement de la dŽfinition de $J$. On a
$G(\chi)G(\chi)\buildrel{\rm Lemme\ 2\ a)}\over =G(\chi^2)J(\chi,\chi)$. Par le mme Lemme,
$G(\chi)^3=G(\chi^3)J(\chi,\chi)J(\chi,\chi^2)$. Et ainsi de suite, on en dŽduit, par rŽcurrence,
que~:

$$G(\chi)^{m-1}=G(\chi^{m-1})J(\chi,\chi)J(\chi,\chi^2)\cdots J(\chi,\chi^{m-2}).$$

En multipliant de part et d'autre par $G(\chi)$, on trouve~:

$$\eqalign{G(\chi)^{m}&=J(\chi,\chi)J(\chi,\chi^2)\cdots
J(\chi,\chi^{m-2})\underbrace{G(\chi^{m-1})}_{=G(\chi^{-1})}G(\chi)\cr
&\buildrel {\rm Lemme\ 2\ b)}\over=q\cdot\chi(-1)J(\chi,\chi)J(\chi,\chi^2)\cdots
J(\chi,\chi^{m-2})\in\Q(\zeta_m),\cr}$$

car on a vu que $J(\chi_1,\chi_2)\in\Q(\zeta_m)$ si $\chi_1,\chi_2$ sont d'un ordre qui divise $m$, et
puis $\chi(-1)=\pm1$.\qed
\bigskip
Maintenant, on va tre obligŽ de considŽrer comme connu un certain nombre de rŽsultats Òclassiques"
de thŽorie algŽbrique des nombres.
\bigskip
{\bf DŽfinitions-ThŽormes ÒRappels sur les corps de nombres et les  corps cyclotomiques"}

{\sl
Soit $K\supset\Q$ un corps, qu'on supposera inclus dans $\C$. On dit que $K$ est {\it un corps de
nombres} s'il est de dimension finie comme $\Q$-espace vectoriel. On note $[K:\Q]=n$ cette dimension.
On peut montrer qu'il existe
$\alpha\in K$ tel que $K=\Q[\alpha]$, c'est-ˆ-dire que $1,\alpha,\ldots ,\alpha^{n-1}$ est une
$\Q$-base de $K$. Chaque ŽlŽment $k$ de $K$ est {\it algŽbrique} sur $\Q$, c'est-ˆ-dire qu'il
existe un polyn™me, dŽpendant de $k$, $f_k(x)=a_m x^m+\cdots +a_1x+a_0\in\Q[x]$, tel que $k$ est une
racine de $f_k$. Le polyn™me unitaire ($a_m=1$) de plus petit degrŽ qui possde $k$ comme racine est
appelŽ {\it polyn™me minimal de $k$}, notŽ $min_\Q(k)$ ou $min(k)$. On peut montrer que si
$K=\Q[\alpha]$, alors $min(\alpha)$ est irrŽductible sur $\Q$ et $K$ est isomorphe ˆ l'anneau quotient
$\Q[X]/(min(\alpha))$. On peut voir aussi que si $L\supset K\supset \Q$ sont des corps de nombres,
alors $[L:\Q]=[L:K]\cdot[K:\Q]$.

Un ŽlŽment est dit {\it entier sur $\Z$} s'il existe un polyn™me unitaire $ x^m+
a_{m-1}x^{m-1}+\cdots +a_1x+a_0\in\Z[x]$ qui annule cet ŽlŽment. On peut montrer que l'ensemble des
ŽlŽments entiers d'un corps de nombre est un anneau, qu'on note $O_K$. On peut montrer que les
idŽaux premiers de $O_K$ sont aussi maximaux et que le quotient de $O_K$ par un de ces idŽaux $\cal
P$ est un corps fini, dont le cardinal se note $\N({\cal P})$. On appelle ce nombre {\it norme absolue de
$\P$}. 

Si $\euf a$ et $\euf b$ sont des idŽaux copremiers de $O_K$ (c'est-ˆ-dire que ${\euf a}+{\euf b}=O_K$) alors $O_K/{\euf
ab}$ est isomorphe ˆ $O_K/{\euf a}\times O_K/{\euf b}$ (ce thŽorme est connu sous le nom de {\it thŽorme chinois}, et
il est vrai pour tout anneau commutatif).

On dŽfinit $I(O_K)$ l'ensemble des {\it idŽaux fractionnaires} de $K$, qui est l'ensemble des
sous-$O_K$-modules non-nuls $\euf a$ de $K$ tels qu'il existe $x\in K^*$ tel que $x{\euf a}\subset A$. Si
${\euf a}$ est un idŽal fractionnaire, on dŽfinit ${\euf a}^{-1}=\{ x\in K\mid x{\euf a}\subset O_K\}$. On
a
${\euf a}{\euf a}^{-1}=O_K$, et $I(O_K)$ est un groupe abŽlien engendrŽ librement par les idŽaux premiers de
$O_K$. C'est-ˆ-dire que tout idŽal fractionnaire s'Žcrit de manire unique ${\euf
a}=\prod_{\P\in\bbP}\P^{v_\P({\euf a})}$ avec, pour tout
$\P\in\bbP$, $v_\P({\euf a})\in\Z$, appelŽ {\it valuation $\P$-adique de $\euf a$}, et $\bbP=\bbP(O_K)$ est
l'ensemble des idŽaux premiers de $O_K$. On Žtend alors la dŽfinitions de copremier et disant que deux
idŽaux fractionnaires $\euf a$ et $\euf b$ sont copremiers si $v_\P({\euf a})\cdot v_\P({\euf b})=0 $ pour
tout idŽal premier $\P$.

Un idŽal fractionnaire
${\euf a}$ est dit {\it principal} s'il est de la forme $xO_K$ pour un $x\in K$. Le sous-groupe des idŽaux
fractionnaires principaux se note $P(O_K)$. On peut montrer que le {\it groupe des classes d'idŽaux},
qui est notŽ ${\cal CL}_K$, et qui est le groupe quotient $I(O_K)/P(O_K)$ est fini. Son cardinal se note $h(K)$. 

Si $K=\Q(\zeta_m)$, on le nommera {\it $m$-ime corps cyclotomique}. On peut montrer que dans ce
cas-lˆ, $O_K=\Z[\zeta_m]$ qu'on notera souvent $E_m$. On peut voir aussi que
$[\Q(\zeta_m):\Q]=\varphi(m)$ o $\varphi(m)$ est {\it l'indicateur d'Euler} qui est le cardinal de
$(\Z/m\Z)^*$ et qui se calcule par la formule $\varphi(mm')=\varphi(m)\varphi(m')$ si $(m,m')=1$ et
$\varphi(p^k)=p^k-p^{k-1}$ si $p$ est un nombre premier. Remarquons enfin que si $m$ est impair, alors
$\Q(\zeta_m)=\Q(\zeta_{2m})$ et donc dans ce cas, les racines de l'unitŽ de ce corps sont les $\pm \zeta_m$, $0\leq i\leq
m$. 

}

{\bf Preuve}

Tous ces rŽsultats sont prouvŽs dans [Nar].\qed
\bigskip
On va maintenant donner un lemme qu'on rŽutilisera au moins 7 fois par la suite~:
\bigskip 
{\bf Lemme IMP}

{\sl Soit $m\in \N$ et $K$ un corps de nombre contenant $\zeta_m$. Soit $p\in\gfP$ tel que $p\not
\hskip-0.4pt |\    m$. Soit $\P$ un idŽal de $O_K$ au-dessus de $p$, c'est-ˆ-dire que $\P\cap\Z=p\Z$.
L'application
$\phi \, : \,  O_K\lra O_K/\P:=\F$ envoie Žvidemment $\Z$ sur $\F_p=\Z/p\Z$ et donc $|\F|=\N(\P)=p^f=q$ o
$f=[O_K/\P:\Z/p\Z]$. Alors on a $m|q-1$. Plus prŽcisŽment, le groupe $\{ \zeta_m^i\mid 1\leq i\leq
m\}\subset O_K$ est envoyŽ injectivement par $\phi$. Son image est donc un sous-groupe cyclique
d'ordre de $m$ de $\F^*$.

}

{\bf Preuve}

Il suffit de vŽrifier que $\zeta_m\not\equiv \zeta_m^i\pmod\P$ pour tout $2\leq i\leq m$. Soit
$f=\prod_{i=1}^m(x-\alpha_i)$ un polyn™me. Alors le polyn™me dŽrivŽ ŽvaluŽ en $\alpha_1$ vaut
$f'(\alpha_1)=\prod_{i=2}^m(\alpha_1-\alpha_i)$. On applique cela ˆ
$f=X^m-1=\prod_{i=1}^m(X-\zeta_m^i)$. On trouve alors
$\prod_{i=2}^m(\zeta_m-\zeta_m^i)=m\zeta_m^{m-1}$. Puisque $\P$ est un idŽal premier et que ni $m$,
ni $\zeta_m$ ne sont dans $\P$, on en dŽduit que $\prod_{i=2}^m(\zeta_m-\zeta_m^i)\not\in\P$, donc
aucun $\zeta_m-\zeta_m^i$ n'est dans $\P$.\qed

\bigskip

En corollaire de ce rŽsultat, il nous est possible de poser la dŽfinition suivante

\bigskip

{\bf DŽfinition}

Soit $m\in \N$ et $K$ un corps de nombres contenant $\zeta_m$. Soit $\P$ un idŽal premier de $O_K$ et
$\alpha\in O_K\setminus \P$. Il est clair que $\alpha^{q-1}\equiv 1\pmod \P$ o $q=\N(\P)$. Donc
l'image de $\alpha$ dans $\F=O_K/\P$ notŽe $\overline{\alpha}$ est telle que
$\overline{\alpha}^{q-1}=1$. Ainsi, $\overline{\alpha}^{q-1\over m}$ est une racine $m$-ime de
l'unitŽ dans $\F$. Par le lemme prŽcŽdent, il existe une unique racine $m$-ime de l'unitŽ dans $K$,
notŽe $\left ({\alpha\over \P}\right )_m$ telle que 

$$\left ({\alpha\over \P}\right )_m\equiv \alpha^{q-1\over m}\pmod\P.$$

On complte la dŽfinition par $\left ({\alpha\over \P}\right )_m=0$ si $\alpha\in\P$.

Si on pose $\mu_m$ l'ensemble des racines $m$-ime de l'unitŽ de $\C$, on a ainsi une application 
$\left ({\cdot\over \P}\right )_m :O_K\lra\mu_m\cup\{0\}$ qu'on appellera {\it symbole de puissance
$m$-ime rŽsiduelle}. 
\bigskip

{\bf Lemme 4}

{\sl Sous les mmes hypothses que la dŽfinition prŽcŽdente, on a 

\art{a)}$\left ({\alpha\over \P}\right )_m=\left ({\beta\over \P}\right )_m$ si
$\alpha\equiv\beta\pmod\P$.

\art{b)}$\left ({\alpha\over \P}\right )_m\equiv \alpha^{q-1\over m}\pmod\P$ pour tout $\alpha\in O_K$.

\art{c)}$\left ({\alpha\beta\over \P}\right )_m=\left ({\alpha\over \P}\right )_m\cdot\left
({\beta\over \P}\right )_m$.

\art{d)}$\left ({\alpha\over \P}\right )_m=1$ si et seulement s'il existe $\beta\in O_K\setminus\P$
tel que $\alpha\equiv \beta^{m}\pmod \P$.

\art{e)}Si $m=2$, $\zeta_2=-1$ et $K=\Q$, on retrouve le symbole de Legendre.

}

{\bf preuve}

les partie a), b) et c) dŽcoulent de la dŽfinition.  Pour la partie d), s'il existe $\beta\in
O_k\setminus\P$ tel que $\alpha\equiv \beta^{m}\pmod \P$, alors $\alpha^{q-1\over
m}\equiv(\beta^m)^{q-1\over m}=\beta^{q-1}\equiv 1\pmod{\P}$. RŽciproquement, si $\left ({\alpha\over
\P}\right )_m=1$, alors $\alpha^{q-1\over m}\equiv 1\pmod\P$.  On se souvient que $\F^*$ est un
groupe cyclique engendrŽ par un ŽlŽment disons $\overline{\gamma}$. Donc
$\alpha\equiv\gamma^s\pmod\P$ pour un certain $1\leq s\leq q-1$. Ainsi $\alpha^{q-1\over m}\equiv
\gamma^{s\cdot {q-1\over m}}\equiv 1\pmod\P$. Ainsi, l'ordre de $\gamma$ qui est $q-1$ divise
$s\cdot {q-1\over m}$, c'est ˆ dire que
$m$ divise $s$, disons, $s=km$. Finalement, $\alpha\equiv (\gamma^k)^m=\beta^m\pmod\P$ en posant
$\beta=\gamma^k$. La partie e) est un corollaire immŽdiat de la partie d).\qed

\bigskip

On a montrŽ que ce symbole passe au quotient et dŽfinit ainsi un caractre multiplicatif d'ordre $m$
du corps $\F$.

Voici encore une sŽrie de rŽsultats classiques~:
\bigskip
{\bf DŽfinitions-ThŽormes ÒRappels sur les corps de nombres et la thŽorie de Galois"}

{\sl

Soit $L/K$ une extension de corps de nombres de degrŽ $n$. Soit $\P$ un idŽal premier de $O_K$ et
$\gP$ un idŽal premier de $O_L$. On dit que $\gP$ est {\it au-dessus de $\P$}, et on Žcrit $\gP|\P$,
si $\P=\gP\cap O_K$, ou, ce qui est Žquivalent, $\gP$ appara"t dans la dŽcomposition en idŽaux
premiers de $\P O_L$. On peut alors identifier $O_K/\P$ ˆ un sous-corps de $O_L/\gP$. On notera
$f(\gP|\P)=[O_L/\gP:O_K/\P]$, qu'on appelle {\it degrŽ rŽsiduel de $\gP|\P$}. Si on Žcrit $\P
O_L=\gP_1^{e_1}\cdots\gP_r^{e_r}$, et notant $f_i$ pour $f(\gP_i|\P)$, alors on peut montrer que
$\sum_{i=1}^r e_i f_i=n$. Les $e_i$, souvent notŽs $e(\gP_i|\P)$ s'appellent les {\it degrŽs rŽsiduels
de $\gP_i|\P$}. On dit que $\P$ n'est {\it pas ramifiŽ} dans $L$ si $e_i=1$, pout tout $i$. On peut montrer
que le nombre de $\P$ qui ramifient est fini. Supposons que $L/K$ soit une {\it extension
galoisienne}, c'est-ˆ-dire que l'ensemble ${\rm Aut}_K(L)$ des $K$-automorphismes de $L$ est d'ordre
$n$. Dans ce cas, ${\rm Aut}_K(L)$ se note ${\rm Gal}(L/K)$ ou $G$, s'il n'y a pas d'ambigu•tŽ. On
peut voir que tout $\sigma\in G$ donne un $O_K$-automorphisme de $O_L$ et qu'ainsi $G$ agit
transitivement sur les idŽaux premiers de $O_L$ qui sont au-dessus d'un $\P\subset O_K$ fixŽ. Comme
consŽquence de cela, si $\P O_L=\gP_1^{e_1}\cdots\gP_r^{e_r}$, alors on a $e_1=e_2=\cdots =e_r:=e$ et
$f_1=f_2=\cdots =f_r:=f$. Donc $efr=n$. Posons $Z(\gP_i|\P)=\{\sigma\in G\mid \sigma(\gP_i)=\gP_i\}$. En
gŽnŽral, les $Z(\gP_i|\P)$ sont conjuguŽs entre eux (c'est-ˆ-dire  pour tout $i\ne j$, il existe
$\sigma\in G$ tel que $Z(\gP_i|\P)=\sigma^{-1}Z(\gP_j|\P)\sigma$ et on a
$|Z(\gP_i|\P)|=ef$. Ainsi, si
$G$ est abŽlien (on dit alors que $L/K$ est une extension abŽlienne), alors les $Z(\gP_i|\P)$ sont
Žgaux ˆ un seul sous-groupe de $G$ qu'on note $Z(\P)$. Soit $M\supset L\supset K$ des corps de nombres
tels que $L/K$ et $M/L$ soient galoisiennes. Si $\ggP\supset \gP\supset \P$ sont des idŽaux premiers de $M$,
$L$ et $K$ respectivement, alors $e(\ggP,\P)=e(\ggP,\gP)\cdot e(\gP,\P)$. Il en est de mme avec les $f$ et
les $r$.

Soit $L/K$ est une extension de corps de nombres. Soit $\alpha\in L$ et $\mu_\alpha\, :\, L\rightarrow L$
dŽfinie par $\mu_\alpha(\beta)=\alpha\cdot\beta$. C'est un endomorphisme $K$-linŽaire de $L$. On dŽfinit
$N_{L/K}(\alpha)=\det (\mu_\alpha)$, c'est la {\it norme} de l'extension $L/K$. On a
$N_{L/K}(\alpha\cdot\beta)=N_{L/K}(\alpha)\cdot N_{L/K}(\beta)$. Si $\alpha\in K$,
$N_{L/K}(\alpha)=\alpha^n$. Si l'extension est galoisienne de groupe $G$,
$N_{L/K}(\alpha)=\prod_{\sigma\in G}\sigma(\alpha)$. On a $|N_{L/\Q}(\alpha)|=\N(\alpha\cdot O_K)$, la
norme absolue. Si $K\subset L\subset E$ sont des corps de nombres, et $\alpha\in E$, alors
$N_{E/K}(\alpha)=N_{L/K}(N_{E/L}(\alpha))$. Enfin, il y a une troisime norme diffŽrente, qu'on appellera
{\it Norme relative de $L$ sur $K$} dŽfinie comme suit~: si $\gP|\P$, avec $\P$ idŽal premier de $K$ et
$\gP$ idŽal de $K$, on dŽfinit $N_{L/K}(\gP)=\P^{f(\gP|\P)}$. On prolonge multiplicativement cette norme ˆ
tous les idŽaux fractionnaires. Si $\euf a$ est un idŽal fractionnaire de $K$, $N_{L/K}({\euf
a}\cdot O_L)={\euf a}^n$ o $n=[L:K]$. On a aussi $N_{L/K}(a\cdot O_L)=N_{L/K}(a)\cdot O_K$, o
$N_{L/K}(a)$ est la norme dŽfinie prŽcŽdemment. Enfin, si $L/K$ est galoisienne de groupe $G$ et si $\euf a$
est un idŽal fractionnaire de $L$, alors $N_{L/K}({\euf a})\cdot O_L=\prod_{\sigma\in G}\sigma({\euf a})$. 

Soit $m$ un entiers positif. Alors l'extension $\Q(\zeta_m)/\Q$ est une
extension abŽlienne de groupe de Galois $G$ isomorphe ˆ $(\Z/m\Z)^*$. L'isomorphisme est canonique~:
$\sigma_t\in G$ dŽfini par
$\zeta_m\mapsto \zeta_m^t$ est envoyŽ sur la classe de $t$ modulo $m$. On peut voir que que si $p\not
\hskip-0.4pt |\   m$, alors $p$ (c'est-ˆ-dire $p\Z$) ne ramifie pas dans $\Q(\zeta_m)$. D'autre part, si
$m=p$ est un nombre premier, alors $p\Z[\zeta_p]=(1-\zeta_p)^{p-1}\Z[\zeta_p]$ (nous dŽmontrerons ce
rŽsultat au Chapitre 6).

Soit $\P$ un idŽal de $\Q(\zeta_m)$ au-dessus de $p\not \hskip-0.4pt |\   m$. Alors $f(\P|p)$ est l'ordre
de $p$ modulo $m$ et $Z(\P)$ est le sous-groupe de $G$ engendrŽ par $\sigma_p$. De plus, pour tout
$\alpha\in\Z[\zeta_m]$, on a $\sigma_p(\alpha)\equiv \alpha^p\pmod\P$, on appelle $\sigma_p$ {\it
l'automorphisme de Frobenius de $\P$ sur $p$}.

}

{\bf Preuve}

Tous ces rŽsultats se trouvent dans [Nar].
\bigskip

Nous voici alors fin prs pour Žnoncer et dŽmontrer le premier rŽsultat de Stickelberger qu'on appellera
la {\it congruence de Stickelberger}. On fixe $p$ un nombre premier et $q=p^f$. Dans $\Q(\zeta_{q-1})$,
$p$ ne ramifie pas. Soit $\gP$ un idŽal premier de $E_{q-1}$ au-dessus de $p$. Soit encore
$f_0=f(\gP|p)$. On vient de voir que $f_0$ est l'ordre de $p$ modulo $q-1$, c'est donc le plus petit
entier positif $k$ tel que $p^k-1$ est un multiple de $q-1$. Donc, $f_0=f$, et donc
$\F:={E_{q-1}/\gP}$ est d'ordre $q$. Notons $\omega$ le caractre multiplicatif dŽfini par
$\left({\cdot\over\gP}\right )_{q-1}^{-1}$. On a vu que la somme de Gauss
$G(\omega^a)\in\Q(\zeta_{p(q-1)})$ pour tout $a\in\N$. Soit $\ggP$ un idŽal premier de $E_{p(q-1)}$
au-dessus de~$\gP$.

Soit $0\leq a<q-1$. Alors $a$ s'Žcrit de manire unique $a=a_0+a_1p+\cdots +a_{f-1}p^{f-1}$ avec
$0\leq a_i\leq p-1$. On dŽfinit $s(a)=\sum_{i=0}^{f-1} a_i$ et $\gamma(a)=a_0 !a_1!\cdots a_{f-1}!$.
Et on prolonge $s$ et $\gamma$ ˆ $\N$ tout entier en dŽcrŽtant qu'ils sont de pŽriode $q-1$.
Autrement dit, si $a\equiv \overline{a}\pmod{q-1}$ avec $0\leq a<q-1$. Alors on pose
$s(a)=s(\overline{a})$ et $\gamma(a)=\gamma(\overline{a})$. Finalement, on note $\pi$ pour $\zeta_p-1$ et on se souvient
que
$pE_p=\pi^{p-1}E_p$.

\bigskip

{\bf ThŽorme (la congruence de Stickelberger)}

\bigskip

{\sl Soit $a\in\N$. On a la congruence suivante~:

$${G(\omega^a)\over \pi^{s(a)}}\equiv {1\over \gamma(a)}\pmod\ggP.$$

}

{\bf Preuve}

Il suffit de prouver le thŽorme pour $0\leq a<q-1$. On va faire une rŽcurrence sur $s(a)$.
Remarquons d'abord une chose~: si $0\leq a=bp<q-1$, alors $s(a)=s(b)$ et $\gamma(a)=\gamma(b)$. De
plus, si $t\in\F^*$, alors $Tr(t^p)=t^p+t^{p^2}+\cdots +\underbrace{t^{p^f}}_{=t}=Tr(t)$. Donc, en
posant $\chi=\omega^b$, on a~:

$$G(\omega^{bp})=-\sum_{t\in\F^*}\chi^p(t)\zeta_p^{Tr(t)}=-\sum_{t\in\F^*}\chi(t^p)\zeta_p^{Tr(t^p)}=
-\sum_{t\in\F^*}\chi(t)\zeta_p^{Tr(t)}=G(\omega^b),\eqno{(i)}$$

car si $t$ parcourt $\F^*$, alors $t^p$ aussi ($p$ est premiers ˆ $q-1$).

\art{a)}Si $s(a)=0$, alors $a=0$ et $G({\bf 1})=1$ et $\gamma(0)=1$. Donc, c'est en ordre.

\art{b)}Si $s(a)=1$. La relation $(i)$ nous montre qu'on peut supposer $a=1$. On se souvient (relation
(5)) que $\sum_{t\in\F^*}\omega(t)=0$. Ainsi, 

$$G(\omega)=-\sum_{t\in\F^*}\omega(t)\zeta_p^{Tr(t)}=-\sum_{t\in\F^*}\omega(t)(\zeta_p^{Tr(t)}-1)\eqno({ii)}$$

Si $t\in\F^*$, on notera $t'\in E_{q-1}$ l'unique racine $q-1$-ime de l'unitŽ reprŽsentant $t$
(Lemme IMP), si bien que $\omega(t)=t'^{-1}$. Si $m$ est un entier reprŽsentant $Tr(t)$ modulo $p$.
Alors on a
$\zeta_p^{Tr(t)}-1=\zeta_p^m-1=(\zeta_p-1)(\zeta_p^{m-1}+\cdots + \zeta_p+1)$. Or, Il est
Žvident que pour tout $r$, en utilisant une mme relation, $\zeta_p^r\equiv 1\pmod{\pi E_p}$. Ainsi donc,
on a

$${\zeta_p^{Tr(t)}-1\over \pi}\equiv m\pmod {\pi E_p}.$$

Mais, d'autre part, $m\equiv Tr(t)\equiv t+t^p+t^{p^2}+\cdots +t^{p^{f-1}}\pmod p$. Donc, $m\equiv
t'+t'^p+\cdots +t'^{p^{f-1}}\pmod\gP$. Or, $\pi E_p\subset \ggP$ et $\gP\subset \ggP$, donc ces
congruences sont {\it a fortiori} vraie modulo $\ggP$. Ainsi, 

$${-G(\omega)\over\pi}\equiv\sum_{t\in\F^*}\omega(t)(t'+t'^p+\cdots +t'^{p^{f-1}})\pmod\ggP.$$

Or, on se souvient que $\omega(t)=t'^{-1}$. D'autre part, comme $|\F^*|=q-1$, il existe un
gŽnŽrateur $t\in\F^*$ tel que $t'$ est $\zeta_{q-1}$. On obtient alors~:

$$\eqalign{{-G(\omega)\over\pi}&\equiv\sum_{t\in\F^*}(1+t'^{p-1}+\cdots
+t'^{p^{f-1}-1})\cr &=
\sum_{t\in\F^*}(1)+\sum_{n=0}^{q-2}(\zeta_{q-1}^{p-1})^n+
\cdots+\sum_{n=0}^{q-2}(\zeta_{q-1}^{p^{f-1}-1})^n\cr
&=q-1+{\zeta_{q-1}^{(q-1)(p-1)}-1\over \zeta_{q-1}^{p-1}-1}+\cdots
+{\zeta_{q-1}^{(q-1)(p^{f-1}-1)}-1\over \zeta_{q-1}^{p^{f-1}-1}-1}
=q-1\equiv -1\pmod\ggP.\cr}$$

On a donc montrŽ le cas $s(a)=1$.

\art{c)}On suppose $0\leq a<q-1$, $s(a)>1$ et le thŽorme vrai pour tout $b$ tel que $0\leq b<q-1$
et $s(b)<s(a)$. Posons $a=a_0+a_1p+\cdots +a_{f-1}p^{f-1}$. Par la relation $(i)$, on peut supposer
que $a_0\ne 0$. On a alors $s(a-1)=s(a)-1$ et $a-1\geq 1$. Par le Lemme 2, $G(\omega^{a-1})\cdot
G(\omega)=G(\omega^a)\cdot J(\omega,\omega^{a-1})$. Posons $b=q-a=(q-1)-(a-1)>1$. Si $t\ne 1$ et
$u=1-t$, alors $\omega(1-t)=u'^{-1}\in E_{q-1}$ avec $u'\equiv 1-t'\pmod \gP$. On alors
$\omega^{a-1}(1-t)=u'^{-(a-1)}=u'^b\equiv (1-t')^b\pmod\gP$. Remarquons que la dernire congruence
est aussi valable si $t=1$. On a alors

$$\eqalignno{-J(\omega,\omega^{a-1})&=\sum_{t\in\F}\omega(t)\omega^{a-1}(1-t)\cr
&=\sum_{t\in\F^*}\omega(t)\omega^{a-1}(1-t)\equiv\sum_{t\in\F^*} t'^{-1}(1-t')^b\cr
&=\sum_{t\in\F^*}t'^{-1}\sum_{j=0}^b\pmatrix{b\cr j\cr}(-1)^jt'^j
\equiv\sum_{j=0}^b\pmatrix{b\cr j\cr}(-1)^j\underbrace{\sum_{t\in\F^*}t'^{j-1}}_{=0\ {\rm sauf\ si}\
j=1}\cr
&\equiv -b\cdot (q-1)\equiv -a\equiv -a_0\pmod\ggP.&(iii)\cr}$$

Or, $0<a_0\leq p-1$, donc $p\not \hskip-0.4pt |\   a_0$. Donc $a_0$ est inversible modulo $\ggP$. Par
hypothse de rŽcurrence, et le cas $a=1$, on a ${G(\omega^{a-1})\over \pi^{s(a-1)}}\cdot {G(\omega)\over
\pi}\equiv {1\over \gamma(a-1)}\cdot 1\pmod\ggP$. Et finalement,

$${G(\omega^a)\over\pi^{s(a)}}\equiv{G(\omega^{a-1}\cdot
G(\omega)\over\pi^{s(a-1)}\cdot\pi}\cdot{1\over J(\omega^{a-1},\omega)}\buildrel (iii)\over={1\over
\gamma(a-1)\cdot a_0}\equiv{1\over\gamma(a)}\pmod\ggP$$
\qed
\bigskip\goodbreak

{\bf DŽfinition}

Soit $A$ un {\it anneau de Dedekind} (un anneau de Dedekind est un anneau dont l'ensemble des idŽaux
fractionnaires est un groupe abŽlien, par exemple l'anneau des entiers d'un corps
de nombres). Soit $\euf a$ un ideal fractionnaire de $A$. On sait que ${\euf a}=\P_1^{r_1}\cdots\P_r^{r_s}$
o les idŽaux $\P_i$ sont premiers et les $r_i\in\Z$. En outre, cette Žcriture est unique. Pour $i=1,\ldots
,r$, on pose $v_{\P_i}({\euf a})=r_i$ et si $\P$ est un idŽal premier diffŽrent des $\P_i$, on pose
$v_\P({\euf a})=0$. Si $K$ est le corps des fractions de $A$ et $x\in K$, on pose $v_\P(x)=v_\P((x))$, o
$(x)$ est l'idŽal fractionnaire engendrŽ par $x$.

\bigskip

{\bf Corollaire}

Sous les mmes hypothses et notations que celles du thŽorme prŽcŽdent, on a, pour tout $a\in\N$~:

$$v_\ggP(G(\omega^a))=s(a)\eqno{(6)}$$

{\bf Preuve}

RŽsumons-nous~: on a $\pi=\zeta_p-1$ et $pE_p=\pi^{p-1}E_p$, donc $e(\pi E_p|p\Z)=p-1$.  D'autre part $\gP$
est un idŽal premier de $E_{q-1}$ au-dessus de $p$, comme $p$ ne ramifie pas dans $\Q(\zeta_{q-1})$,
on a $v_\gP(p)=1$, car $e(\gP|p\Z)=1$. Si on rŽunit le tout dans $E_{p(q-1)}$, on trouve $\gP
E_{p(q-1)}=\ggP^{(p-1)}$. En effet, $e(\ggP|p\Z)=e(\ggP|\gP)\cdot e(\gP|p\Z)=e(\ggP|\gP)\leq
[\Q(\zeta_{p(q-1)}):\Q(\zeta_{q-1})]=p-1$; d'autre part, $e(\ggP|p\Z)=e(\ggP|\pi E_p)\cdot e(\pi E_p|p\Z)=e(\ggP|\pi
E_p)\cdot (p-1)\geq p-1$. On en dŽduit que $e(\ggP|p\Z)=p-1$, donc $e(\ggP|\pi
E_p)=1$ et ainsi $v_\ggP(\pi)=1$. Donc, comme $\gamma(a)$ est inversible modulo
$\ggP$ et qu'on a montrŽ que ${G(\omega^a)\over \pi^{s(a)}}\equiv {1\over \gamma(a)}\pmod\ggP$, on en
dŽduit le corollaire.\qed

{\bf DŽfinition}

Soit $p$ un nombre premier, $m>1$, tel que $p\not \hskip-0.4pt |\  m$. Posons $\P$ un idŽal premier de
$E_m$ au-dessus de $p$. On pose $\F=E_m/\P$, $f=f(\P|p)$. On a $|\F|=p^f=:q$. Posons $G={\rm
Gal}(\Q(\zeta_m)/\Q)$ qui est isomorphe ˆ $(\Z/m\Z)^*$, $\sigma_t\leftrightarrow t$ avec
$\sigma_t(\zeta_m)=\zeta_m^t$. Evidemment, il ne faudra pas confondre $G$ avec une somme de Gauss, mais le
contexte permettra de diffŽrencier les deux objets 

Soit $A$ un anneau commutatif, on dŽfinit $A[G]=\{\sum_{t\in(\Z/m\Z)^*}a_t\sigma_t\mid a_t\in A\}$.
C'est un anneau (dont l'addition se fait terme ˆ terme et la multiplication est hŽritŽe de celle de
$A$ et de la loi de composition de $G$). On appelle $A[G]$ {\it l'algbre de 
$G$ sur $A$}. Dans un premier temps, nous Žtudierons $A[G]$ pour $A=\Z$, puis, nous serons obligŽ de
passer ˆ un corps fini pour pouvoir utiliser le fait que $K[G]$ est {\it semi-simple} si $K$ est un
corps et si sa caractŽristique ne divise pas $|G|$. Mais n'anticipons pas, nous verrons cela en temps voulu
!! On fait agir (exponentiellement) $\Z[G]$ sur $\Q(\zeta_m)$~: soit
$\lambda=\sum_{t\in(\Z/m\Z)^*}a_t\sigma_t\in\Z[G]$ et $x\in\Q(\zeta_m)$; on pose

$$x^\lambda=\prod_{t\in(\Z/m\Z)^*}\sigma_t(x)^{a_t}.$$

De mme si $\euf a$ est un idŽal fractionnaire de $\Q(\zeta_m)$ ou mme une classe d'idŽaux, on
pose\goodbreak
${\euf a}^\lambda=\prod_{t\in(\Z/m\Z)^*}\sigma_t({\euf a})^{a_t}$. Et on vŽrifie facilement que
$(x^\lambda)^\mu=x^{\lambda\mu}$, $x^\lambda x^\mu=x^{\lambda+\mu}$ et $(xy)^\lambda=x^\lambda
y^\lambda$. Et de mme pour les idŽaux ou les classes d'idŽaux. Remarquons que si nous avions dŽfini une
action multiplicative du genre $\lambda\cdot x=\sum a_t\sigma_t(x)$, les choses ne se seraient pas
si bien passŽes au niveau des idŽaux... On note $\chi$ le caractre dŽfinit par $\left ({\cdot\over
\P}\right )_m^{-1}$. Notons encore $\P_t=\P^{\sigma_t^{-1}}$ si $t$ parcourt $(\Z/m\Z)^*$. Par
transitivitŽ de l'action du groupe de Galois sur les idŽaux au-dessus d'un idŽal fixŽ, les $\P_t$ parcourent les idŽaux
de $\Q(\zeta_m)$ au-dessus de $p$, chacun apparaissant $f$ fois.

\bigskip
\goodbreak
{\bf Lemme 5}

{\sl Soit $G(\chi)$ la somme de Gauss du caractre dŽfini prŽcŽdemment. Alors
$G(\chi)^m\in\Q(\zeta_m)$ et on a
$$G(\chi)^m E_m=\prod_{t\in(\Z/m\Z)^*}\P_t^{r_t},$$

avec, pour tout $t$,  $r_t={m\over p-1}s\left({t(q-1)\over m}\right)$.

}

{\bf Preuve}

Le fait que $G(\chi)^m\in\Q(\zeta_m)$ provient lemme 3 c) et du fait que $\chi$ est un caractre
d'ordre $m$. Le mme lemme 3, nous apprend que $G(\chi)\cdot \overline{G(\chi)}=q=p^f$. Donc, les
idŽaux qui apparaissent dans la factorisation de l'idŽal engendrŽ par $G(\chi)^m$ sont les $\P_t$.
Il suffit donc de montrer que $v_{\P_t}(G(\chi)^m)={m\over p-1}s\left({t(q-1)\over m}\right)$.

Puisque $m|q-1$(lemme IMP), on a une tour de corps $\Q\subset\Q(\zeta_m)\subset
\Q(\zeta_{q-1})\subset\Q(\zeta_{p(q-1)})$. On choisit un idŽal $\gP$ de $E_{q-1}$ au-dessus de $\P$
et un idŽal $\ggP$ de $E_{p(q-1)}$ au-dessus de $\gP$. On a que $\gP$ n'est pas ramifiŽ au-dessus de
$p$ (car $p$ ne divise pas $q-1$), donc, $\gP$ n'est pas ramifiŽ au-dessus de $\P$. Donc,
$$v_\P=v_\gP,\eqno{(i)}$$ c'est-ˆ-dire $v_\gP({\euf a} E_{q-1})=v_\P({\euf a})$ pour tout idŽal $\euf a$ de
$E_m$. En revanche $\gP E_{p(q-1)}=\ggP^{p-1}$ (voir corollaire prŽcŽdent). Donc,
$$v_{\ggP}=(p-1)v_{\gP},\eqno{(ii)}$$ (c'est-ˆ-dire $(p-1)v_\gP({\euf a})=v_\ggP({\euf a}E_{p(q-1)})$ pour
tout idŽal $\euf a$ de
$\Q(\zeta_{q-1})$. Il est clair que $E_{q-1}/\gP$ contient canoniquement $E_m/\P$, or on a vu (juste
avant la congruence de Stickelberger) que $|E_{q-1}/\gP|=q$, mais on a aussi $|E_m/\P|=q$. Donc
$\F=E_{q-1}/\gP=E_m/\P$. On note ˆ nouveau $\omega$ le caractre multiplicatif de $\F$ $\left
({\cdot\over\gP}\right)_{q-1}^{-1}$. On a alors 
$$\chi=\omega^{q-1\over m},\eqno{(iii)}$$ en effet, soit
$\alpha\in E_m$. on a~: $\chi^{-1}(\alpha)=\left({\alpha\over\P}\right)_{m}\equiv
\alpha^{q-1\over m}\pmod\P$ et $\omega^{-1}(\alpha)=\left ({\alpha\over\gP}\right)_{q-1}\equiv
\alpha\pmod\gP$. Donc,
$\left ({\alpha\over\gP}\right)_{q-1}^{q-1\over m}\equiv \alpha\pmod\gP$. Cette congruence est aussi
valable modulo $\P=\gP\cap E_m$. Donc $\left ({\alpha\over\gP}\right)_{q-1}^{q-1\over m}\equiv \left
({\alpha\over\gP}\right)_{m}\pmod\P$, ce qui veut dire que $\left ({\alpha\over\gP}\right)_{q-1}^{q-1\over
m}= \left ({\alpha\over\gP}\right)_{m} $ car on a vu au Lemme IMP que toutes les racines $m$-ime de
l'unitŽs Žtaient distinctes modulo $\P$. 

Soit $\sigma_t\in G$. On note $\overline{\sigma_t}$ l'ŽlŽment de ${\rm Gal}(\Q(\zeta_{mp})/\Q)$ tel que
$\overline{\sigma_t}\mid_{\Q(\zeta_m)}=\sigma_t$ et $\overline{\sigma_t}\mid_{\Q(\zeta_p)}={\rm
Id}_{\Q(\zeta_p)}$ (ce rŽsultat est aussi un rŽsultat classique de la thŽorie de Galois). Calculons~:

$$G(\chi)^{\overline{\sigma_t}}=\left
(-\sum_{x\in\F^*}\chi(x)\psi(x)\right)^{\overline{\sigma_t}}=-\sum_{x\in\F^*}\chi(x)^t\psi(x)=G(\chi^t).\eqno{(iv)}$$

D'autre part, 

$$\left( G(\chi)^m\right )^{\sigma_t}=\left( G(\chi)^m\right )^{\overline{\sigma_t}}=\left(
G(\chi)^{\overline{\sigma_t}}\right )^{m}\buildrel (iv)\over = G(\chi^t)^m.\eqno{(v)}$$

On trouve enfin, en posant $a={t(q-1)\over m}$~:

$$\eqalign{v_{\P_t}(G(\chi)^m)&=v_\P(G(\chi)^m)^{\sigma_t})\buildrel (v)\over =v_\P(G(\chi^t)^m)
=mv_\P(G(\chi^t))\buildrel (iii)\over =mv_\P(G(\omega^a))
\buildrel (i)\ {\rm et}\ (ii)\over = {m\over p-1} v_\ggP(G(\omega^a))\cr &\buildrel \rm (6)\over ={m\over
p-1}s\left ({t(q-1)\over m}\right).\cr}$$\qed
\bigskip

{\bf Lemme 6}

{\sl Soit $a\in\N$, on rappelle que $s(a)=\sum_{i=0}^{f-1}a_i$ o $a\equiv a_0+a_1p+\cdots +a_{f-1}p^{f-1}\pmod{q-1}$
avec $0\leq a_i\leq p-1$. Alors on a~:
$$s(a)=(p-1)\cdot\sum_{i=0}^{f-1}\left <{p^i a\over q-1}\right >\ \hbox{o $\left <x\right >$ est la partie
fractionnaire de $x$.}$$

}

{\bf Preuve}

Comme $p^f\equiv 1\pmod{q-1}$, on a le systme de congruence 

$$\eqalign{a&\equiv a_0+a_1p+\cdots +a_{f-1}p^{f-1}\pmod{q-1}\cr
a\cdot p&\equiv a_{f-1}+a_0p+\cdots +a_{f-2}p^{f-1}\pmod{q-1}\cr
&\vdots\dotfill\cr
a\cdot p^{f-1}&\equiv a_1+a_2p+\cdots +a_{0}p^{f-1}\pmod{q-1}.\cr}$$

On remarque que le membre de droite de la $i$-ime congruence divisŽ par $q-1$ est $\left<{p^i a\over q-1}\right >$
pour $0\leq i\leq f-1$. En sommant le tout, on obtient 
$$\sum_{i=0}^{f-1}\left <{p^i a\over q-1}\right >={s(a)\over q-1}(1+p+\cdots+p^{f-1})={s(a)\over p-1}.$$\qed

\bigskip\goodbreak

{\bf ThŽorme (la relation de Stickelberger)}

{\sl On se souvient que $\P$ est un idŽal premier de $E_m=\Z[\zeta_m]$ au-dessus de $p\not
\hskip-0.4pt |\  m$ et $\chi$ est le caractre $\left ({\cdot\over
\P}\right )_m^{-1}$. Alors on a

$$G(\chi)^mE_m=\P^{m\cdot \Theta},$$

o 
$$\Theta=\sum_{t\in(\Z/m\Z)^*}\left<{t\over m}\right >\sigma_t^{-1}=\sum_{\matrix{t=0,\ldots , m-1\cr
(t,m)=1\cr}}{t\over m}\sigma_t^{-1}.$$

L'ŽlŽment $\Theta$ s'appelle {\it l'ŽlŽment de Stickelberger}.

}

{\bf Preuve}

Soit $t_1,\ldots ,t_g$ un systme de reprŽsentant de $\left (\Z/m\Z\right)^*$ modulo le sous-groupe engendrŽ par $p$
(remarquons que $t_1^{-1},\ldots ,t_g^{-1}$ est aussi un systme de reprŽsentant modulo le sous-groupe engendrŽ par $p$).
On a vu lors du rappel sur la thŽorie de Galois que $Z(\P)$ Žtait le sous-groupe de $G$ engendrŽ par $\sigma_p$. Donc
$\P_{t_1},\ldots ,\P_{t_g}$ est la liste (sans redondances) des idŽaux premiers de $E_m$ au-dessus de $p$. Donc,
$$G(\chi)^mE_m=\prod_{i=1}^g\P_{t_i}^{v_{\P_{t_i}}(G(\chi)^m)}\buildrel {\rm Lemme\ 5}\over
=\P^{\gamma'}\hbox{ o }\gamma'={m\over p-1}\sum_{i=1}^g s\left ( t_i(q-1)\over m\right
)\sigma_{t_i}^{-1}.$$

On remarque que $p^jt_i$ avec $1\leq i\leq g$ et $0\leq j\leq f-1$ reprŽsentent tous les ŽlŽments de $\left
(\Z/m\Z\right)^*$ et on sait (puisque $Z(\P)$ est le sous-groupe de $G$ engendrŽ par $\sigma_p$) que
$\P^{\sigma_{t_i}^{-1}}=\P^{\sigma_{{p^jt_i}}^{-1}}$. Gr‰ce, au lemme 6, on trouve alors~:

$$\P^{\gamma'}=\P^{{m\over p-1}\sum_{i=1}^g(p-1)\sum_{j=0}^{f-1}\left <{p^j t_i\over m}\right
>\sigma_{t_i}^{-1}}=\P^{m\sum_{i=1}^g\sum_{j=0}^{f-1}\left <{p^j t_i\over m}\right>\sigma_{p^j
t_i}^{-1}}=\P^{m\sum_{t}^{m*}\left<{t\over m}\right >\sigma_t^{-1}}=\P^{m\Theta},$$
avec la convention que $\sum_{t}^{m*}$ veut dire $\dst\sum_{t\in(\Z/m\Z)^*}$.\qed

\bigskip

{\bf Lemme 7}

{\sl Soit $K$ un corps de nombres galoisien, $\euf a$ un idŽal fractionnaire de $K$ et $m\in\N$. ConsidŽrons
${\cal C}$ la classe de
$\euf a$. Alors il existe ${\euf b}\in{\cal C}$ tel que $m$ soit premier ˆ $\euf b$, 
c'est-ˆ-dire $v_\P({\euf b})\cdot v_\P( m O_K)=0$ pour tout idŽal premier $\P$.

}

{\bf Preuve}

C'est un corollaire du thŽorme chinois, on pourrait le mettre en exercice, mais comme c'est court, on le
donne quand mme~: on a ${\euf a}=\prod_{\P\in\bbP}\P^{v_\P({\euf a})}$ et $ m
O_K=\prod_{\P\in\bbP}\P^{v_\P({m O_K})}$. Soit $V$ l'ensemble des idŽaux premiers qui divisent $\euf a$ ou
$mO_K$. Pour tout $\P\in V$, on choisit $x_\P\in
\P^{v_\P({\euf a})}\setminus\P^{v_\P({\euf a})+1}$. Par le thŽorme chinois, il existe $a\in O_K$ tel que
$a\equiv x_\P\pmod{\P^{v_\P({\euf a})+1}}$. L'idŽal fractionnaire ${1\over a}{\euf a}$ rŽpond ˆ la
question, car $v_\P( {1\over a}{\euf a})=v_\P({1\over a})+v_\P({\euf a})=0$ pour tout $\P\in V$; et si
$\P\not\in V$, alors
$v_\P(mO_K)=0$.\qed

\bigskip
\goodbreak
{\bf Corollaire }

{\sl Soit $\euf a$ un idŽal  fractionnaire de $\Q(\zeta_m)$. Alors  ${\euf a}^{m\Theta}$ est principal

}

{\bf preuve}

Par le lemme 7, on peut supposer que $\euf a$ est premier ˆ $m$. Pour tout idŽal premier $\P$ divisant $\euf a$, on a
gr‰ce ˆ la relation de Stickelberger que $\P^{m\Theta}=(G(\chi))^m$ et on conclut gr‰ce ˆ la multiplicativitŽ de l'action
de $\Z[G]$ sur les idŽaux.\qed

\bigskip

{\bf DŽfinition}

On note $I_{st}(\Q(\zeta_m))=\Z[G]\cap \Theta\Z[G]$, appelŽ {\it l'idŽal de Stickelberger}. Soit $b$ est un nombre entier
premier ˆ $m$, on note $\Theta_b=(\sigma_b-b)\Theta\in\Theta\Z[G]$.

\bigskip

{\bf Lemme 8}

{\sl
Sous les mmes hypothses, on a~:

$$\Theta_b\in I_{st}(\Q(\zeta_m)).$$

}

{\bf Preuve}

Il suffit de montrer que $\Theta_b\in\Z[G]$. On a $$\sigma_b\Theta=\sum_{t}^{m*}\left <{t\over
m}\right >\sigma_b\sigma_t^{-1}=\sum_{t}^{m*}\left <{t\over
m}\right>\sigma_{b^{-1}t}^{-1}\buildrelÒt=b^{-1} t"\over=\sum_{t}^{m*}\left <{bt\over m}\right>\sigma_t^{-1}.$$
Attention, quand on note $b^{-1}$, a veut dire qu'on considre un entier $b^{-1}$ tel que
$bb^{-1}\equiv 1\pmod m$. On trouve 

$$\Theta_b=\sum_{\matrix{t=0,\ldots , m-1\cr (t,m)=1\cr}}\left (\left <{bt\over m}\right >-b{t\over m}\right
)\sigma_t^{-1}.$$

Posons $0\leq r<m$ et $s$ tels que $bt=sm+r$ alors on a $\left <{bt\over m}\right >={r\over m}=b{t\over m}-s$ et
$s=\left [ {bt\over m}\right ]$, o $\left[ x\right]$ dŽsigne la partie entire de $x$ . Et ainsi,

$$\Theta_b=- \sum_{\matrix{t=0,\ldots , m-1\cr (t,m)=1\cr}}\left [{bt\over
m}\right]\sigma_t^{-1}\in\Z[G].\eqno{(7)}$$

\qed

\bigskip

{\bf Exemple important}

Si $m=p$ est un nombre premier impair, $b=2$ et $1\leq t\leq p-1$ alors

$$\left [{2t\over p}\right ]=\cases{0&si $t\leq {p-1\over 2}$\cr 1&si $t\geq {p+1\over 2}$\cr}$$

Ainsi,

$$\Theta_2=-\sum_{t=1}^{p-1}\left [{2t\over p}\right]\sigma_t^{-1}=-\sum_{t={p+1\over
2}}^{p-1}\sigma_t^{-1}.\eqno{(8)}$$

\bigskip

{\bf Lemme 9}

{\sl L'idŽal de Stickelberger $I_{st}(\Q(\zeta_m))$ est engendrŽ par les $\Theta_b$ en tant qu'idŽal (et mme en tant que
$\Z$-module)

}

{\bf Preuve}

Soit $\beta\Theta\in I_{st}(\Q(\zeta_m))$. Ce la veut dire que $\beta\in\Z[G]$ et $\beta\Theta\in \Z[G]$.
Posons
$\beta=\sum_t^{m*}b_t\sigma_t$ avec $b_t\in\Z$ pour tout $t$. On a

$$\eqalign{\beta\Theta&=\left (\sum_{t\in(\Z/m\Z)^*}b_t\sigma_t\right )\left (\sum_{a\in(\Z/m\Z)^*}\left <{a\over m}\right
>\sigma_a^{-1}\right )=\sum_{t,a\in(\Z/m\Z)^*}b_t\left<{a\over m}\right >\sigma_{at^{-1}}^{-1}\cr &\buildrel a=ct\over
=\sum_{c\in(\Z/m\Z)^*}\underbrace{\left (\sum_{t\in(\Z/m\Z)^*}b_t\left <{ct\over m}\right >\right )}_{\in\Z{\rm\
par\ hyp.}}\sigma_c^{-1}.}$$

En particulier, si $c=1$  on a $\sum_t^{m*}b_t\left <{t\over m}\right >\in\Z$. Donc 

$$\sum_{\matrix{0\leq t\leq m\cr (t,m)=1\cr}}b_t{t\over m}=:u\in\Z\quad \hbox{ ou encore
}\sum_{\matrix{0\leq t\leq m\cr (t,m)=1\cr}}b_tt=um.$$

Remarquons la petite astuce suivante~: $m\Theta=(m+1-\underbrace{\sigma_{m+1}}_{= \rm Id})\Theta=-\Theta_{m+1}$. Ainsi,
on a 
$$\beta\Theta=\left (\sum_{t\in(\Z/m\Z)^*}b_t\sigma_t\right )\Theta=\sum_{\matrix{0\leq t\leq m\cr
(t,m)=1\cr}}b_t(\sigma_t-t)\Theta+\Bigg (\sum_{\matrix{0\leq t\leq m\cr (t,m)=1\cr}}b_t t\Bigg
)\Theta=\sum_{\matrix{0\leq t\leq m\cr (t,m)=1\cr}}b_t\Theta_t+um\Theta$$
est bel et bien engendrŽ par les $\Theta_b$.
\qed
\bigskip\goodbreak

{\bf ThŽorme de Stickelberger}

{\sl  L'idŽal de Stickelberger $I_{st}(\Q(\zeta_m))$ annule le groupe des classes d'idŽaux de
$\Q(\zeta_m)$.

}

{\bf Preuve}

Soit $\P$ un idŽal premier de $E_m$ premier ˆ $m$. La relation de Stickelberger nous donne~:
$\P^{m\Theta}=G(\chi)^mE_m$, idŽal qu'on note $(G(\chi)^m)$. Elevons ceci ˆ la puissance $\sigma_b-b$, on trouve
$\P^{m\Theta_b}=(G(\chi)^m)^{\sigma_b-b}=(G(\chi)^{\overline{\sigma}_b-b})^m$, o
$\overline{\sigma}_b$ est l'extension de $\sigma_b$ ˆ $\Q(\zeta_{pm})$ telle que $\overline{\sigma}_b(\zeta_p)=\zeta_p$
($p$ Žtant le nombre premier au-dessous de $\P$). 

{\bf Affirmation~:}  $G(\chi)^{\overline{\sigma}_b-b}\in\Q(\zeta_m)$.

L'affirmation montre qu'on peut Òenlever le $m$", c'est-ˆ-dire $\P^{\Theta_b}=(G(\chi)^{\overline{\sigma}_b-b})$. On a
donc montrŽ que ${\euf a}^{\Theta_b}$ est principal pour tout idŽal $\euf a$ premier ˆ $mE_m$. Mais on a vu (Lemme 7) que
dans toute classe d'idŽaux il existe un reprŽsentant premier ˆ $mE_m$. On en dŽduit le thŽorme puisqu'on vient de voir que
les $\Theta_b$ engendraient l'idŽal de Stickelberger.

Pour prouver l'affirmation, il suffit de montrer que $G(\chi)^{\overline{\sigma}_b-b}$ est invariant par $G':={\rm
Gal}(\Q(\zeta_{pm})/\Q(\zeta_m))$. Soit $c\in \Z$ tel que $(c,pm)=1$. Tout ŽlŽment du groupe ${\rm
Gal}(\Q(\zeta_{pm})/ \Q)$ s'Žcrit $\tau_c$ tel que $\tau_c(\zeta_{pm})=\zeta_{pm}^c$. On a alors $G'=\{\tau_c\mid
c\equiv 1\pmod m\}$. On se souvient que, pour tout $x\in\F^*=E_m/\P$, on a que $\chi(x)$ est une racine $m$-ime
de l'unitŽ et que
$\psi(x)$ est une racine $p$-ime de l'unitŽ, de plus si $x$ parcourt $\F^*$, alors $cx$ aussi . On a alors~:

$$\eqalignno{(G(\chi)^{\overline{\sigma}_b})^{\tau_c}&=\left
(-\sum_{x\in\F^*}\chi(x)\psi(x)\right)^{\overline{\sigma}_b\tau_c}=-\sum_{x\in\F^*}\chi(x)^b\psi(x)^c
=-\sum_{x\in\F^*}\chi^b(x)\psi(cx)\cr
&=-\chi^{-b}(c)\sum_{x\in\F^*}\chi^b(cx)\psi(cx)=-\chi^{-b}(c)\sum_{x\in\F^*}\chi^b(x)\psi(x)\cr
&=\chi^{-b}(c)G(\chi)^{\overline{\sigma}_b}.&(i)\cr }$$

Le mme calcul avec $b=1$ donne $(G(\chi))^{\tau_c}=\chi^{-1}(c)G(\chi)$ qui, ŽlevŽ ˆ la puissance $b$ donne
$$(G(\chi)^b)^{\tau_c}=\chi^{-b}(c)G(\chi)^b.\eqno{(ii)}$$ 
Le quotient des deux ŽgalitŽs $(i)$ et $(ii)$ nous donne bien

$$(G(\chi)^{\overline{\sigma}_b-b})^{\tau_c}=G(\chi)^{\overline{\sigma}_b-b}$$

Ce qui prouve l'affirmation et donc le thŽorme.\qed
\bigskip

Nous allons maintenant montrer un rŽsultat qui sera utilisŽ Òen passant" au chapitre 7. Le problme est qu'il
faut faire pas mal de dŽfinition et de Òrappels". Pour la fin de ce paragraphe, on va supposer que $m=p$ est un
nombre premier. Bien sžr, il y a un problme de notation, car avant, $p$ Žtait un nombre premier qui ne divise pas
$m$. Mais nous ne rencontrerons plus ce $p$-lˆ. Donc $G={\rm Gal}(\Q(\zeta_p)/\Q)$ est de cardinal $p-1$. 
\bigskip
\goodbreak
{\bf DŽfinition}

Soit $\chi$ un caractre multiplicatif de $G$ {\it impair}, c'est-ˆ-dire $\chi(-1)=-1$. On sait que $G$
est isomorphe ˆ $(\Z/p\Z)^*$. On peut prolonger
$\chi$ ˆ tous les entiers en posant $\chi(x)=\chi(\overline{x})$  ou $\overline{x}$ est la classe de $x$ si $x$
n'est pas un multiple de $p$, et $\chi(x)=0$ sinon. On dŽfinit la {\it sŽrie $L$ de Dirichlet} 
$$L(1,\chi)=\sum_{n=1}^{\infty}{\chi(n)\over n}.$$

{\bf Lemme 10}

{\sl Sous les mmes hypothses celles de la dŽfinition prŽcŽdente, on a~:

$$0\ne L(1,\chi)={i\cdot\pi\cdot G(\chi)\over p^2}\cdot \sum_{a=1}^{p-1} a\overline{\chi(a)}.$$

}
{\bf Preuve}

Le fait que $L(1,\chi)\ne 0$ est prouvŽ dans [Was, Corollary 4.4, p. 34] et l'autre ŽgalitŽ est donnŽe
dans le mme ouvrage [Was, Theorem 4.9, p.38].\qed

\bigskip\goodbreak

{\bf DŽfinition}

Notons $\iota$ pour $\sigma_{-1}$ (la conjugaison complexe) qui est l'unique ŽlŽment d'ordre 2 dans $G$. Soit
$R$ un anneau commutatif et $M$ un $R[G]$-module. On crŽe deux nouveaux modules $M^{\pm}=\{ x\in M\mid
\iota x=\pm x\}$. Ce sont des sous-$R[G]$-modules de $M$. Si ${1\over 2}\in R$, alors on pose
$\varepsilon^{\pm}={1\pm\iota\over 2}$. On vŽrifie facilement que ${\varepsilon^{\pm }}^2=\varepsilon^{\pm
}$, $\varepsilon^+\varepsilon^-=0$ et $\varepsilon^++\varepsilon^-=1$. Ainsi,
$M=\varepsilon^+M\oplus\varepsilon^-M$ et $M^{\pm}=\varepsilon^{\pm}M$. Si ${1\over 2}\not\in R$, alors
on pose $\varepsilon^{\pm}={1\pm\iota}$, on a encore une somme directe et $\varepsilon^+M\oplus
\varepsilon^-M\subset M$, et l'indice est une puissance de 2. On observe en passant que
$\sigma_a\cdot\iota=\sigma_{-a}$.

Supposons $R=\Z$. Si $x=\sum_{\sigma\in G}n_\sigma\sigma$, on note $\|x\|$ pour $\sum_{\sigma \in
G}|n_\sigma |$. 

On dŽfinit alors $I=I_p=(1-\iota) I_{st}$. On a Žvidemment $I\subset I_{st}^-:=I_{st}\cap\Z[G]^-$ o
dans ce cas, $\Z[G]^-=(1-\iota)\Z[G]$. C'est une vŽrification facile~: il est clair que $(1-\iota)\Z[G]\subset
\Z[G]^-$. D'autre part, si $x=\sum_{j=1}^{p-1} a_j\sigma_j$ est tel que que $\iota x=-x$, alors on voit que pour
tout $j$, $a_{p-j}=-a_j$. Donc, en posant $y=\sum_{j=1}^{p-1\over 2} a_j\sigma_j$, alors on a $x=(1-\iota) y$. 

On sait (Lemme 9) que $I_{st}$ est engendrŽ par
$g_b:=-\Theta_b=(b-\sigma_b)\Theta\buildrel (7)\over=\sum_{a=1}^{p-1}\left [{ab\over p}\right
]\sigma_a^{-1}$, $b=1,\ldots , p-1$
et $g_p:=p\Theta=((p+1)-\sigma_{p+1})\Theta=-\Theta_{p+1}=\sum_{j=1}^{p-1}\left [ {jp\over p}\right
]\sigma_j^{-1}$. Posons, pour $i=1,\ldots ,p-1$~:
$$f_i=g_{i+1}-g_i=\sum_{a=1}^{p-1}\left (\left [{a(i+1)\over p}\right ]-\left [{ai\over p}\right ]\right
)\sigma_a^{-1}.$$

On observe que les coefficients des $f_i$ sont 0 ou 1. Puisque  $g_1,g_2,\ldots ,g_p$ engendrent $I_{st}$ et que
$g_1=0$, alors $f_1,\ldots ,f_{p-1}$ engendrent aussi $I_{st}$. Finalement, on remarque que
$$f_{p-1}=\sum_{a=1}^{p-1}\left (\left [{ap\over p}\right ]-\left [{a(p-1)\over p}\right ]\right
)\sigma_a^{-1}=\sum_{a=1}^{p-1}\left (a-\underbrace{\left [a-{a\over p} \right]}_{=a-1}\right
)\sigma_a^{-1}=\sum_{a=1}^{p-1}\sigma_a=:s(G).$$

DŽfinissons enfin $e_i=(1-\iota)f_i$, $i=1,\ldots ,p-1$. 

\bigskip\goodbreak
{\bf Lemme 11}

{\sl Sous les mmes hypothses, on a 

\art{a)}$\|f_i\|={p-1\over 2}$.

\art{b)}$f_1,\ldots ,f_{p-1\over 2}$ et $s(G)$ engendrent $I_{st}$ qui est, comme $\Z$-module, de rang ${p+1\over
2}$.

\art{c)}$e_1,\ldots ,e_{p-1\over 2}$ engendrent $I$ qui est, comme $\Z$-module, de rang ${p-1\over
2}$.

}

{\bf Preuve}

Soit $1\leq i\leq {p-1\over 2}$ et $1\leq a\leq p-1$. On a $a(i+1)=kp+r$ ($0<r<p$). D'o
$(p-a)(i+1)=(i+1-k)p-r=(i-k)p+(p-r)$. D'autre part, $ai=lp+s$ ($0<s<p$). D'o
$(p-a)i=(i-l)p-s=(i-1-l)p+(p-s)$. On en dŽduit~:
$$\left (\left [{a(i+1)\over p}\right ]-\left [{ai\over p}\right ]\right )+\left (\left [{(p-a)(i+1)\over p}\right
]-\left [{(p-a)i\over p}\right ]\right )=(k-l)+((i-k)-(i-1-l))=1.$$

Donc, si le coefficient de $\sigma_a^{-1}$ dans $f_i$ est 1, celui de $\sigma_{-a}^{-1}$ est 0, et vice versa. On
en dŽduit que $\|f_i\|={p-1\over 2}$, donc, a) est prouvŽ. On a aussi $(1+\iota)f_i=s(G)$.

Puisque $e_i=(1-\iota)f_i$, on a $\|e_i\|=(p-1)$ (car, lˆ o il y avait des 0, on trouve des -1 et les 1 restent).
De manire analogue, si $a(i+1)=kp+r$ ($0<r<p$), alors $a(p-i-1)=(a-k-1)p+(p-r)$. Et si $ai=lp+s$ ($0<s<p$),
alors $a(p-i)=(a-l+1)p+(p-s)$. On en dŽduit~:
$$\left [{a(i+1)\over p}\right ]-\left [{ai\over p}\right ]=k-l\hbox{ et } \left [{a(p-i)\over p}\right ]-\left
[{a(p-1-i)\over p}\right ]=(a-l-1)-(a-k-1)=k-l.$$
Donc, $f_i=f_{p-1-i}$. Cela prouve que $f_1,\ldots ,f_{p-1\over 2}$ et $s(G)$ engendrent $I_{st}$. En
appliquant $(1-\iota)$, on trouve que $e_1,\ldots , e_{p-1\over 2}$ engendrent $I$. Si on montre que $e_1,
\ldots ,e_{p-1\over 2}$ est une base, alors $f_1,\ldots ,f_{p-1\over 2}$ et $s(G)$ formeraient aussi une
base~: si
$\lambda_1 f_1+\cdots +\lambda_{p-1\over 2}f_{p-1\over 2}+\lambda_{p+1\over 2}s(G)=0$, appliquant
$(1-\iota)$, on trouve que
$\lambda_1=\cdots =\lambda_{p-1\over 2}=0$ et donc $\lambda_{p+1\over 2}=0$. Pour cela, il suffirait de montrer
que la matrice $\left(\left [{ij\over p}\right ]\right )_{2\leq j\leq{p+1\over 2}}^{2\leq i\leq{p+1\over 2}}$ est
inversible !! mais vous pouvez toujours essayer, c'est vachement dur !!

On va s'y prendre de manire un peu dŽtournŽe (on utilisera entre autre le lemme prŽcŽdent, ('faut bien
qu'il serve ˆ quelque chose)). Il suffit de montrer que $I$ est de $\Z$-rang ${p-1\over 2}$. Mais,
$2I_{st}^-\subset I\subset I_{st}^-$. La premire inclusion vient du fait que si $x\in I_{st}^-$ alors
$2x=(1-\iota)x\in I$, la deuxime est triviale. Il suffit donc de montrer que $I_{st}^-$ est de $\Z$-rang
${p-1\over 2}$. Mais
$I_{st}=\Z[G]\Theta\cap\Z[G]$ et donc $I_{st}^-=I_{st}\cap\Z[G]^-=\Z[G]\Theta\cap \Z[G]^-$. Ainsi, $J:=
\Z[G]^-\Theta\cap \Z[G]^-\subset I_{st}^-$. Donc, si on arrive ˆ voir que $J$ est de rang ${p-1\over 2}$ ou mme
que $\Z[G]^-\Theta$ est de rang ${p-1\over 2}$, on a gagnŽ. Pour cela, il suffit de montrer que la
multiplication par
$\Theta$ de $\C[G]^-$ dans lui-mme est injective (et mme un automorphisme $\C$-linŽaire). On sait que
$\widehat G$, l'ensemble des caractres de $G$ est isomorphe aux caractres de $(\Z/p\Z)^*$, donc on notera
$\chi(a)$ au lieu de $\chi(\sigma_a)$. Posons, pour $\chi\in\widehat G$ $\varepsilon_\chi={1\over
p-1}\sum_{\sigma\in G}\chi(\sigma)\sigma^{-1}={1\over p-1}\sum_{a=1}^{p-1}\chi(a)\sigma_a^{-1}$. On
vŽrifie que
$\varepsilon_\chi^2=\varepsilon_\chi$, $\varepsilon_\chi\varepsilon_\psi=0$ si $\chi\ne\psi$, et
$\sum_{\chi\in\widehat G}\varepsilon_\chi=1$. On a donc, $\C[G]=\bigoplus_{\chi\in\widehat
G}\C[G]\varepsilon_\chi$. En comparant les dimensions, on voit que chacun des $\C[G]\varepsilon_\chi$ est
de dimension 1 sur $\C$. Donc
$\C[G]\varepsilon_\chi=\C\varepsilon_\chi$. Et, puisque ce sont des idŽaux, on a $\Theta\C\varepsilon_\chi\subset
\C\varepsilon_\chi$. Donc, la base des $\varepsilon_\chi$ diagonalise la multiplication par $\Theta$. On
vŽrifie aussi que pour tout $\sigma\in G$ et $\chi\in\widehat G$, on a
$\sigma\varepsilon_\chi=\chi(\sigma)\varepsilon_\chi$. Donc $\Theta\varepsilon_\chi=\left ({1\over
p}\sum_{a=1}^{p-1} a\overline{\chi(a)}\right )\varepsilon_\chi$. En particulier, $\iota\varepsilon_\chi=
\chi(-1)\varepsilon_\chi$. Donc,
$\C[G]^-=\bigoplus_{\chi\rm\ impair}\C[G]\varepsilon_\chi$. Mais, on sait que (Lemme 10)~:

$$0\ne L(1,\chi)={i\cdot\pi\cdot G(\chi)\over p}\cdot{1\over p} \sum_{a=1}^{p-1} a\overline{\chi(a)}.$$

Donc, la matrice de la multiplication par $\Theta$ vue dans la base des $\varepsilon_\chi$ est diagonale et
chaque ŽlŽment de la diagonale est non nul. Donc la multiplication par $\Theta$ est un automorphisme de
$\C[G]^-$, ce qui prouve notre lemme.\qed

\vfill\eject

\long\def\art#1{{\parindent0pt\item{#1}}\hangindent=7mm\hangafter=-20}
\long\def\artart#1{{\parindent0pt\item{#1}}\hangindent=12mm\hangafter=-20}
\font\para=cmbx12 at 18pt
\def\O{\hbox{$\cal O$}}
\def\U{\hbox{$\cal U$}}
\def\m{\hbox{\rs m\!}}
\def\dst{\displaystyle}
\font\doub=msbm10 at 10pt
\def\lra{\longrightarrow}
\def\qed{\hfill$\square$}
\def\gfP{\relax\ifmmode\bbP\else $\bbP$\fi}
\def\gP{{\euf P}}
\def\P{{\cal P}}
\def\QQ{{\cal Q}}

\def\ggP{{\bf P}}
\newcount\chapnomb \chapnomb=1
\newcount\parnomb \parnomb=1
\pageno =29

\parindent0pt

\centerline{\para CHAPITRE 6}

\bigskip

{\para  Premier ThŽorme de Mih$\taille{18}\breve{\bf a}$ilescu}
\bigskip

On va, dans ce chapitre, montrer un thŽorme qui prolonge les relations de Cassels. Il sera une brique
essentielle pour la preuve finale~:
\bigskip

{\bf ThŽorme 1 (ThŽorme 1 de Mih$\breve{\bf a}$ilescu)}

{\sl Soit $p$ et $q$ des nombres premiers impairs, et $x$, $y$ des entiers non nuls tels que $x^p-y^q=1$.
Alors on a~:

$$p^2|y,\quad q^2|x,\quad p^{q-1}\equiv 1\pmod{q^2}\quad\hbox{et}\quad q^{p-1}\equiv
1\pmod{p^2}\eqno{(7)}$$

}
\bigskip
Evidemment, il va falloir faire quelques lemmes et rappels. Nous aurons en particulier besoin du thŽorme
de Stickelberger et des relations de Cassels. Nous allons commencer par un lemme simple, mais gŽnial~:
\bigskip
{\bf Lemme 1}

{\sl Soit $K$ un corps de nombres et $\P$ un idŽal premier de $O_K$ tel que $\N(\P)=\left | O_K/\P\right
|=p^f$ avec $p$ un nombre premier. Supposons que $\alpha$, $\beta\in O_K$ soient tels que $\alpha^p\equiv
\beta^p\pmod\P$. Alors $\alpha^p\equiv \beta^p\pmod{\P^2}$.

}

{\bf Preuve}

Pour tout $x\in O_K$, on a Žvidemment $x^{p^f}\equiv x\pmod\P$ (car $O_K/\P$ est un corps fini ˆ $p^f$
ŽlŽments). Ainsi, de $\alpha^p\equiv \beta^p\pmod\P$ et Žlevant ˆ la puissance $p^{f-1}$, on obtient
$\alpha\equiv\alpha^{p^f}\equiv\beta^{p^f}\equiv\beta\pmod\P$. Posons $\gamma=\alpha-\beta\in\P$. D'autre
part, $p\in\P$, ainsi~:

$$\alpha^p-\beta^p=(\beta+\gamma)^p-\beta^p=\sum_{j=1}^p\pmatrix{p\cr j\cr}\gamma^j\beta^{p-j}\in
p\gamma O_K+\gamma^p O_K\subset\P^2.$$
\qed
\bigskip
{\bf Lemme 2}
{\sl 
\art{a)}Soit $p$ un nombre premier. Alors pour tout $i=1,\ldots ,p-1$ l'ŽlŽment $1-\zeta_p^i\over
1-\zeta_p$ est un inversible de $\Z[\zeta_p]$ ce qui montre que $p\Z$ est un idŽal qui ramifie totalement
dans $\Z[\zeta_p]$, en outre, il est le seul qui ramifie. 

\art{b)}Si $u$ est un inversible de $\Z[\zeta_p]$, alors $\overline{u}\over u$ est une racine $2p$-ime
de l'unitŽ o ici $\overline{u}$ dŽsigne le conjuguŽ complexe de $u$.

}

{\bf Preuve}

La partie a) est un corollaire du rappel qu'on a fait sur les corps cyclotomiques au chapitre 5, mais
c'est court et joli ˆ dŽmontrer, on va donc le faire~:

$${x^p-1\over x-1}=x^{p-1}+x^{p-2}+\cdots +x+1=\phi_p(x)=\prod_{i=1}^{p-1}(x-\zeta_p^i)$$

En Žvaluant en 1, on trouve $p=\prod_{i=1}^{p-1}(1-\zeta_p^i)$. D'autre part, si $1\leq i\leq p-1$, on a 
$${1-\zeta_p^i\over 1-\zeta_p}=\zeta_p^{i-1}+\cdots
+\zeta_p+1\in\Z[\zeta_p]=E_p\hbox{ et }N(1-\zeta_p^i)=\prod_{j=1}^{p-1}(1-\zeta_p^{ij})\buildrel k=ij\over
=\prod_{k=1}^{p-1}(1-\zeta_p^{k})=N(1-\zeta_p),$$

o $N=N_{\Q(\zeta_p)/\Q}$. Donc, $N((1-\zeta_p^i)/(1-\zeta_p))=1$. Or, pour tout $x$,
$N(x)=x\cdot\prod_{\sigma}\sigma(x)$ o
$\sigma$ parcourt les ŽlŽments du groupe de Galois de $\Q(\zeta_p)/\Q$. Donc, on a montrŽ que
$(1-\zeta_p^i)/(1-\zeta_p)$ est une unitŽ de $E_p$ et donc que
$pE_p=(1-\zeta_p)^{p-1}E_p$. Le fait qu'il soit le seul ˆ ramifier vient du fait que le {\it
discriminant} de $\Q(\zeta_p)$ sur $\Q$ vaut $\pm p^{p-2}$, car seuls
les premiers divisant le discriminant ramifient (c'est un thŽorme qui est montrŽ dans Samuel,
[Sam, ThŽorme 1, ¤ 5.3]).

\art{b)} Cela se dŽmontre gr‰ce au sous-lemme suivant~:

{\bf Sous-Lemme}

{\sl Soit $K$ un corps de nombres de degrŽ $n$ et $\alpha\in O_K$. Si $\left |\sigma(\alpha)\right |=1$ pour
tout plongement complexe $\sigma$ de $K$. Alors $\alpha$ est une racine de l'unitŽ.

}

{\bf Preuve du Sous-Lemme}

Les nombres $\alpha^2,\alpha^3,\ldots$ sont dans $K$ et donc leur polyn™me minimal est de degrŽ infŽrieur
ou Žgal ˆ $n$. De plus, tous les
$\sigma(\alpha^k)$ sont de module 1. Fixons un $k$ et posons $f(x)=x^m+a_{m-1}x^{n-1}+\cdots
+a_1x+a_0=\sum_{\sigma}(x-\sigma(\alpha^k))$ ($m\leq n$), le polyn™me minimal de $\alpha^k$. Puisque les
$\sigma(\alpha^k)$ sont de module 1, on voit que $|a_i|\leq \max_{j=1,\ldots ,m}\left(\pmatrix{j\cr
m\cr}\right )$. Donc les $a_i$ possibles ne sont qu'en nombre fini. Il existe donc $k_1>k_2$ tels que
$\alpha^{k_1}=\alpha^{k_2}$ ou encore $\alpha^{k_1-k_2}=1$ ce qui veut dire que $\alpha$ est une racine de
l'unitŽ.\qed

Pour terminer la partie b) du lemme, il suffit de remarquer que $\left |\sigma({\overline{u}\over
u})\right |=\left |{\sigma(\overline{u})\over
\sigma(u)}\right |=\left |{\overline{\sigma(u)}\over
\sigma(u)}\right |=1$, pour tout $\sigma$ dans le groupe de Galois de $\Q(\zeta_p)/\Q$, car ce groupe est
abŽlien et que la conjugaison complexe en est un ŽlŽment. D'autre part, puisque $u$ est dans
$\Q(\zeta_p)=\Q(\zeta_{2p})$, il doit tre une racine $2p$-ime de l'unitŽ (pas forcŽment primitive).\qed

\bigskip

{\bf Preuve du ThŽorme 1}

Nous aurons besoin de deux ingrŽdients essentiels dans cette preuve~: les identitŽs de Cassels et le
thŽorme de Stickelberger. Supposons donc que $x^p-y^q=1$ avec $p,q$ premiers impairs. Les identitŽs de
Cassels nous donnent que $x\equiv 1\pmod p$. Donc, $x-1\in pE_p\subset (1-\zeta_p^i)E_p$, pour tout
$i=1,\ldots ,p-1$. Cela veut dire que
$x-\zeta_p^i=x-1+1-\zeta_p^i\in (1-\zeta_p^i)E_p$. On a alors $\beta_i:={x-\zeta_p^i\over 1-\zeta_p^i}\in
E_p$. Pour tout $j\ne i$ on observe la relation 
$$(1-\zeta_p^i)\beta_i-(1-\zeta_p^j)\beta_j=\zeta_p^j-\zeta_p^i.\eqno{(i)}$$ Mais,
il est Žvident (en vertu du lemme 2 a)) que
$(1-\zeta_p^i)E_p=(1-\zeta_p^j)E_p=(\zeta_p^j-\zeta_p^i)E_p=(1-\zeta_p)E_p$. En divisant la relation $(i)$
par
$1-\zeta_p$, on trouve que les idŽaux $\beta_iE_p$ et $\beta_jE_p$ sont premiers entre eux.

Maintenant vient une jolie idŽe~: les relations de Cassels nous apprennent que 
${x^p-1\over x-1}=p\cdot v^q$ pour un certain $v\in\Z$ et que $q|x.$ Ainsi,
$$\prod_{i=1}^{p-1}\beta_i=\prod_{i=1}^{p-1}{x-\zeta_p^i\over 1-\zeta_p^i}={x^p-1\over (x-1)\cdot p}=v^q.$$

On vient de voir que les $\beta_i E_p$ Žtaient premiers entre eux. Donc, en vertu de l'unicitŽ de
l'Žcriture des idŽaux en produits d'idŽaux premiers, on en dŽduit que pour chaque $i$, il existe un idŽal
${\cal V}_i$ tel que $$\beta_i E_p={\cal V}_i^q\eqno{(9)}$$ 

Soit $\theta\in I_{st}(\Q(\zeta_p))$ un ŽlŽment de l'idŽal de Sickelberger (en fait, on prendra
$\theta=-\Theta_2$ dŽfini lors du chapitre prŽcŽdent). Posons $\beta=\beta_1$. On a vu que
$(\beta E_p)^\theta=({\cal V}_1^q)^\theta=({\cal V}_1^\theta)^q=\alpha^qE_p$ pour un certain $\alpha\in
\Q(\zeta_p)$. Calculons~:

$$\left ({1-x\zeta_p^{-1} \over 1-\zeta_p^{-1}}\right )^\theta=\left ({\zeta_p-x\over \zeta_p-1}\right
)^\theta=\beta^\theta=\varepsilon\cdot \alpha^q,\eqno{(ii)}$$

pour une une unitŽ $\varepsilon$ de $E_p$; ceci parce que $(\beta
E_p)^\theta=\beta^\theta E_p^\theta=\beta^\theta E_p=\alpha^q E_p$. Posons
$\lambda=(1-\zeta_p^{-1})^\theta$. Alors il existe une racine $2p$-ime de l'unitŽ $\delta$ telle que
$\overline{\lambda}\overline{\varepsilon}=\lambda\varepsilon\delta^q$. En effet~: Žcrivons
$\theta=\sum_{i=1}^{p-1} a_i\sigma_i\in\Z[G]$ (rappelons que $\sigma_i$ est l'application qui envoie
$\zeta_p$ sur $\zeta_p^i$). D'une part, on a

$${\overline{\lambda}\over\lambda}={(1-\zeta_p)^\theta\over
(1-\zeta_p^{-1})^\theta}=\prod_{i=1}^{p-1}{(1-\zeta_p^{i})^{a_i}\over
(1-\zeta_p^{-i})^{a_i}}=\prod_{i=1}^{p-1}-\zeta_p^{ia_i}=\zeta_p^{\sum_{i=1}^{p-1} i a_i}$$

qui est une racine $2p$-ime de l'unitŽ. D'autre part, on a vu au lemme 2 que
${\overline{\varepsilon}\over\varepsilon}$ Žtait aussi une telle racine, donc le produit des deux aussi.
Finalement, toute racine $2p$-ime de l'unitŽ est une puissance $q$-ime d'une autre racine $2p$-ime de
l'unitŽ, car $(2p,q)=1$. De $(ii)$, on a~: $(1-x\zeta_p^{-1})^\theta=\lambda\varepsilon\alpha^q$. En
prenant le conjuguŽ, on a
$(1-x\zeta_p)^\theta=\overline{\lambda}\overline{\varepsilon}\,\overline{\alpha}^q=
\lambda\varepsilon\delta^q\overline{\alpha}^q$. Donc,

$$(1-x\zeta_p^{-1})^\theta-(1-x\zeta_p)^\theta=\lambda\varepsilon (\alpha^q-(\delta\overline{\alpha})^q),\hbox{ o $\delta$ est
une racine de l'unitŽ.}$$

Supposons maintenant (pour simplifier) que les coefficients de $\theta$ sont des entiers positifs (c'est le
cas pour $-\Theta_2$) ainsi, on reste dans les Òvrais idŽaux" plut™t que dans les idŽaux fractionnaires et
donc en particulier  $\alpha\in E_p$. Puisque $q|x$, on en dŽduit que
$(1-x\zeta_p^{-1})^\theta-(1-x\zeta_p)^\theta\equiv 1-1 \equiv 0\pmod {q E_p}$. Soit $\QQ$ un idŽal premier
de
$E_p$ au-dessus de $q$. Puisque $\varepsilon$ est une unitŽ et que $\lambda$ est une puissance de
$1-\zeta_p$ fois une unitŽ, et donc inversible modulo $\QQ$, on trouve que $\alpha^q\equiv
(\delta\overline{\alpha})^q\pmod\QQ$. Et, gr‰ce au Lemme 1, on en dŽduit que $\alpha^q\equiv
(\delta\overline{\alpha})^q\pmod{\QQ^2}$. On a

$$(1-x\zeta_p)^\theta=\prod_{i=1}^{p-1}(1-x\zeta_p^i)^{a_i}\equiv 1-x\sum_{i=1}^{p-1} a_i \zeta_p^i\pmod
{x^2}.$$

De mme,

$$(1-x\zeta_p^{-1})^\theta\equiv 1-x\sum_{i=1}^{p-1} a_i \zeta_p^{-i}\equiv 1-x\sum_{i=1}^{p-1} a_{p-i}
\zeta_p^{i}\pmod {x^2}.$$

On en dŽduit que $x\cdot\sum_{i=1}^{p-1}(a_i-a_{p-i})\zeta_p^i\equiv 0\pmod {\QQ^2}$, car $x^2\in\QQ^2$.
On a besoin d'un petit lemme facile~: si $\QQ$ est un idŽal premier d'un anneau de Dedekind avec
$ab\in\QQ^2$ et $a\in\QQ$, alors $b\in \QQ$ ou $a\in\QQ^2$. Dans notre cas, cela veut dire que
$\sum_{i=1}^{p-1}(a_i-a_{p-i})\zeta_p^i \equiv 0\pmod {\QQ}$ ou alors $x\equiv 0\pmod {\QQ^2}$.

Supposons que $\sum_{i=1}^{p-1}(a_i-a_{p-i})\zeta_p^i \equiv 0\pmod {\QQ}$ pour tout $\QQ$
au-dessus de $q$. Puisqu'il n'y a pas de ramification au-dessus de $q$, cela implique que 
$\sum_{i=1}^{p-1}(a_i-a_{p-i})\zeta_p^i\in qE_p$, ou encore $q|a_i-a_{p-i}$ pour tout $i$. On a dit qu'on
prenait  $\theta=-\Theta_2=\sum_{i={p-1\over 2}}^{p-1}\sigma_i^{-1}=\sum_{i=1}^{p-1}
b_i\sigma_i^{-1}=\sum_{i=1}^{p-1} a_i\sigma_i$. On remarque que l'on a, si
$p\geq 3$, $a_{p-1\over 2}=b_{p-2}=1$ et
$a_{p+1\over 2}=b_2=0$. Cela voudrait dire que $q| a_{p-1\over 2}-a_{p+1\over 2}=1$ ce qui est absurde~!

On suppose donc que $\QQ^2$ divise $x$, pour un certain $\QQ$ au-dessus de $q$, donc
$\sigma(\QQ)^2$ divise $\sigma(x)=x$ pour tout $\sigma\in {\rm Gal}(\Q(\zeta_p)/\Q)$. Or, $\sigma(\QQ)$
parcourt tous les idŽaux au-dessus de $q$, cela implique que $q^2|x$, toujours parce qu'il n'y a
pas de ramification au-dessus de $q$. Donc, par les relations de Cassels, il existe $u\in \Z$ tel que 
$p^{q-1}\cdot u^q=x-1\equiv -1\pmod{q^2}$. Cela entra"ne, par le petit thŽorme de Fermat, que $u^q\equiv
-1\equiv (-1)^{q}\pmod q$, donc, par le Lemme~1, $u^{q}\equiv -1\pmod{q^2}$. Revenant aux mmes relations
de Cassels, on obtient
$p^{q-1} \cdot u^q\equiv -p^{q-1}\equiv -1\pmod {q^2}$, donc, en multipliant par $-1$, on a donc bien prouvŽ
que 
$$p^{q-1}\equiv 1\pmod {q^2}.$$

Le fait que $p^2|y$ et que $q^{p-1}\equiv 1\pmod{p^2}$ se dŽmontre de manire identique.\qed

\vfill\eject

\long\def\art#1{{\parindent0pt\item{#1}}\hangindent=7mm\hangafter=-20}
\long\def\artart#1{{\parindent0pt\item{#1}}\hangindent=12mm\hangafter=-20}
\font\para=cmbx12 at 18pt
\def\O{\hbox{$\cal O$}}
\def\U{\hbox{$\cal U$}}
\def\m{\hbox{\rs m\!}}
\def\dst{\displaystyle}
\font\doub=msbm10 at 10pt
\def\lra{\longrightarrow}
\def\qed{\hfill$\square$}
\def\gfP{\relax\ifmmode\bbP\else $\bbP$\fi}
\def\gP{{\euf P}}
\def\P{{\cal P}}
\def\QQ{{\cal Q}}

\def\ggP{{\bf P}}
\newcount\chapnomb \chapnomb=1
\newcount\parnomb \parnomb=1
\pageno =32

\parindent0pt
\centerline{\para CHAPITRE 7}

\bigskip
{\para  Premiers contacts avec le groupe $\taille {15} H$ et petites valeurs de
$\taille {15}p$ et $\taille {15}q$}
\bigskip

Nous allons prouver que si $p\not \hskip-0.4pt|\  h_q^{-}$ ou $q\not\hskip-0.4pt |\  h_p^{-}$, alors
$x^p-y^q=1$ n'a pas de solution non triviale. Le nombre $h_p^{-}$ est ce qu'on appelle le {\it nombre de
classes relatifs de
$\Q(\zeta_p)/\Q$} que nous allons dŽfinir tout soudain. Cela va en particulier montrer que $p,q> 7$. En
fait des calculs pas trs long mais sous-jacents ˆ une thŽorie un peu plus longues montreraient que
$p,q\geq 43$, mais on n'aura pas besoin d'tre aussi fin pour la suite.

DŽfinissons dŽjˆ ce $h^-$~: il est bien connu des thŽoriciens des nombres, car il intervient dans un pan de
la dŽmonstration du grand thŽorme de Fermat.
\bigskip

{\bf DŽfinition}

Si $n$ est un nombre entier positif, on dŽfinit $\zeta_n^+:=\zeta_n+\zeta_n^{-1}=\zeta_n+\overline{\zeta_n}$
o $\overline{x}$ est le conjuguŽ complexe de $x$. On peut montrer que $\Q(\zeta_p^+)$ est totalement rŽel
(i.e. tout plongement de $\Q(\zeta_n^+)$ dans $\C$ est rŽel), que $\Q(\zeta_n^+)=\{ x\in\Q(\zeta_n)\mid
x=\overline{x}\}$ est d'indice 2 dans $\Q(\zeta_p)$. On note $E_p$ l'anneau des entier de $\Q(\zeta_p)$ et
$E_p^+$ celui de $\Q(\zeta_p^+)$. On peut voir que $E_p^+=\Z[\zeta_n^+]$. On note $U_n$ les inversibles de
$E_n$ et $U_n^+$ ceux de $E_n^+$. On note encore $\pi=\zeta_p-1$ le gŽnŽrateur de l'unique idŽal au-dessus
de $p$, qui ramifie totalement.
\bigskip

{\bf Lemme de Kummer}

{ \sl Soit $p$ un nombre premier impair. Alors on a 
$$U_p=W_pU_p^+,$$ 
o $W_p:=\{\pm \zeta_p^{k}\mid 0\leq k\leq p-1\}$ est l'ensemble des racines de l'unitŽ de $\Q(\zeta_p)$.
Plus prŽcisŽment, soit $u\in U_p$ alors il existe
$0\leq k\leq p-1$ et $u_0\in U_p^+$ tels que $u=\zeta_p^k u_0$.

}

{\bf Preuve}

Par le lemme 2 du Chapitre 6, ${u\over \overline{u}}=\pm \zeta_p^l$. Montrons que le $\pm$ est de
trop~: supposons que $u=-\zeta_p^l\overline{u}$. supposons que $u=b_0+b_1\zeta_p+\cdots
+b_{p-2}\zeta_p^{p-2}$. De l'ŽgalitŽ b-bte (qu'on a dŽjˆ utilisŽ et qu'on rŽutilisera)
$\zeta_p=1+(\zeta_p-1)$, on dŽduit que 
$$u\equiv b_0+b_1+\cdots +b_{p-2}\equiv \overline{u}\pmod {\pi E_p}.$$

Donc, par hypothse, $u\equiv -\overline{u}\equiv -u\pmod {\pi E_p}.$ Cela veut dire que $2u\in\pi E_p$.
Comme $2\not\in \pi E_p$, on en dŽduit que $u\in \pi E_p$ ce qui est absurde, car $u$ est inversible. Ainsi,
on a $u=\zeta^{l}\overline{u}=\zeta^{2k}\overline{u}$ pour un certain $k$, car $(2,p)=1$. Posons
$u_0=u\zeta^{-k}$. On a
$\overline{u_0}=\overline{u}\zeta_p^k=u\zeta^{-l}\zeta^k=u\zeta^{-k}=u_0\in U_p^+$. Et bien sžr $u=\zeta^k
u_0$.\qed
\bigskip

{\bf Lemme 1}

{\sl Soit $K\subset L$ deux corps de nombres. Alors l'application ${\euf a}\mapsto {\euf a}L$ dŽfinit une
application ${\cal CL}_K\rightarrow {\cal CL}_L$. {\it A priori}, cette application n'est pas injective.
Mais si $K=\Q(\zeta_p^+)$ et $L=\Q(\zeta_p)$, alors cette application est injective.

}

{\bf Preuve}

Les idŽaux fractionnaires principaux de $K$ sont
envoyŽs sur les idŽaux principaux de $L$, donc cette application passe aux classes d'idŽaux.

Montrons d'abord en toute gŽnŽralitŽ que l'application ${\euf a}\mapsto {\euf a}L$ est un homomorphisme
injectif de l'ensemble des idŽaux fractionnaire de $K$ dans ceux de $L$.  Soit $\euf a$ un idŽal
fractionnaire de $K$ tel que ${\euf a}O_L=O_L$. Alors Žvidemment ${\euf a}\subset K$ et ${\euf a}={\euf
a}O_K\subset {\euf a}O_L=O_L$. Donc, ${\euf a}\subset K\cap O_L=O_K$. De mme, ${\euf a}^{-1}O_L=O_L$. Donc,
on dŽduit de la mme manire que ${\euf a}^{-1}\subset O_K$ ou encore ${\euf a}\supset O_K$ (car
$O_K^{-1}=O_K$ et le passage ˆ l'inverse change le sens des inclusions). Donc, ${\euf a}=O_K$.  

Montrons maintenant que dans notre cas, l'application reste injective sur le groupe des classes~: soit
${\euf a}\subset E_p^+$ tel que ${\euf a}E_p=\alpha E_p$ (il est possible de supposer que ${\euf a}\subset
E_p^+$, car on travaille dans les classes d'idŽaux). Puisque
$\overline{\euf a}={\euf a}$ et $\overline{E_p}=E_p$, on a $\alpha E_p=\overline{\alpha E_p}$. Donc il existe
une $\varepsilon\in U$ tel que $\alpha=\varepsilon\overline{\alpha}$. Mais puisque $\varepsilon={\alpha\over
\overline{\alpha}}$ on a vu au Lemme 2 du Chapitre 6 que $\epsilon$ est en fait une racine
$2p$-ime de l'unitŽ. Si $\varepsilon=\eta^2$, alors, en posant, comme au Lemme de Kummer
$\beta=\alpha\overline{\eta}$,
on a
$\overline{\beta}=\overline{\alpha}\eta={\alpha\over\eta^2}\eta=\alpha\eta^{-1}=\alpha\overline{\eta}=\beta$.
Donc, $\alpha E_p=\overline{\alpha}E_p=\beta E_p$. On en dŽduit que ${\euf a}=\beta E_p^+$, par injectivitŽ
de ${\euf a}\mapsto {\euf a}L$ sur les idŽaux fractionnaire de $K$ dans ceux de $L$. Mais, {\it a priori},
$\varepsilon=(-\zeta_p)^l$ n'est pas forcŽment un carrŽ, sauf si on arrive ˆ prouver que $l$ est pair. On
remarque d'abord que si $v=v_{\pi E_p}$ est la valuation de l'idŽal au-dessus de $p$, alors $v({\euf
b}E_p)$ est pair pour tout idŽal $\euf b$ de $E_p^+$, ceci parce que $\pi$ ramifie totalement au-dessus de
$p$. En particulier, $v(\alpha)$ est pair, car il provient d'un idŽal de $K$. D'autre part,
${\pi\over\overline{\pi}}=-\zeta_p$. Ainsi, on a
$\alpha=\left ({\pi\over\overline{\pi}}\right )^l\overline{\alpha}$. Donc,
$\alpha\overline{\pi}^l=\overline{\alpha}\pi^l\in E_p^+$. On trouve alors
$l=v(\alpha\overline{\pi}^l)-v(\alpha)$ est pair.\qed

\bigskip\vskip 1cm\goodbreak\vskip -1cm

{\bf Notation}

Le Lemme prŽcŽdent nous montre, entre autre, que le cardinal du groupe des classes de $\Q(\zeta_p^+):=h_p^+$
divise le cardinal du groupe des classes de $\Q(\zeta_p):=h_p$. Le quotient s'appelle le {\it nombre de
classes relatifs de $\Q(\zeta_p)/\Q$} qu'on note $h_p^-$.

\bigskip
\goodbreak

{\bf DŽfinition}

Soit $p$ un nombre premier impair. On pose comme souvent $G={\rm Gal}(\Q(\zeta_p)/\Q)$, de cardinal $p-1$.
Notons $\iota$ pour
$\sigma_{-1}$ (la conjugaison complexe). Soit
$R$ un anneau commutatif et $M$ un $R[G]$-module. On rappelle que $M^{\pm}=\{ x\in M\mid
\iota x=\pm x\}$ qui sont des sous-$R[G]$-modules de $M$ tels que
$M=M^+\oplus M^-$ si 2 est inversible dans $R$. 

On se souvient que l'idŽal de Stickelberger $I_{st}=I_{st}(\Q(\zeta_p))=\Z[G]\cap
\Theta\Z[G]$ avec
$$\Theta=\sum_{a\in(\Z/p\Z)^*}\left<{t\over p}\right >\sigma_a^{-1}={1\over
p}\sum_{a=1}^{p-1}{a}\sigma_a^{-1}.$$

On a vu au Chapitre 5 (Lemme 11) que $I=I_p=(1-\iota) I_{st}$ Žtait un $\Z$-module libre de rang
$p-1\over 2$.

On pose 
$$E=\{u(1-\zeta_p)^k\mid k\in\Z,u\in U_p\}.$$

On voit que $E$ est un sous-groupe de $\Q(\zeta_p)^*$, c'est aussi un $\Z[G]$-module, car
$\sigma_a(1-\zeta_p)=1-\zeta_p^a=u(1-\zeta_p)$ pour une certaine unitŽ $u$. Enfin on voit que
$E=\Z[\zeta_p,{1\over p}]^*$, car $p=u(1-\zeta_p)^{p-1}$ pour une certaine unitŽ $u$.

Soit $q$ un nombre premier diffŽrent de 2 et de $p$. Si $A$ est un groupe abŽlien on pose
$A^q=\{\alpha^q\mid \alpha\in A\}$ et $A[q]=\{\alpha\in A\mid \alpha^q=1\}$. Si le cardinal de $A$ est fini,
$A[q]$ est de cardinal $1$ ou $q$. DŽfinissons maintenant notre fameux groupe $H$~: on notera

$$H=\{\alpha\in\Q(\zeta_p)^*\mid v_{\euf r}(\alpha)\equiv 0\pmod{q}\hbox{ pour tout idŽal premier }{\euf
r}\ne\pi E_p\}/{\Q(\zeta_p)^*}^q.$$

Il est Žvident que $\alpha\in H$ si et seulement s'il existe $k\in\Z$ et $\euf a$ un idŽal fractionnaire de
$\Q(\zeta_p)$  tels que $\alpha E_p={\euf a}^q(1-\zeta_p)^k$. On voit facilement que $H$ est un
$\Z[G]$-module ($\alpha^{\sum a_\sigma\sigma}=\prod\sigma(\alpha)^{a_\sigma}$ et $\sigma_a(\pi E_p)=\pi
E_p$, pour tout $a$). Il est clair que l'action de $q\Z[G]$ sur $H$ envoie tout sur $1_H$, faisant de $H$
un $\F_q[G]$-module. On utilisera souvent cette structure, car 2 est inversible dans $\F_q$.

Posons encore $G^+$ le groupe ${\rm Gal}(\Q(\zeta_p^+)/\Q)=G/<\iota>$; son cardinal est ${p-1\over 2}$. On
note encore $\cal CL$ le groupe des classes de $\Q(\zeta_p)$. Attention~: ${\cal CL}^+$ qui est l'ensemble
des classes d'idŽaux de $E_p$ invariantes par $\iota$ n'est pas l'ensemble des classes d'idŽaux de
$\Q(\zeta_p^+)$, mais le second s'injecte dans le premiers (lemme 1). Donc $h_p^+||{\cal CL}^+|$. D'autre
part, ${\cal CL}^-$ s'injecte dans ${\cal CL}$ sans rencontrer ${\cal CL}^+$. Donc, ${\cal CL}^-$ s'injecte
dans ${\cal CL}/{\cal CL}^+$, ou encore $|{\cal CL}^-|$ divise ${h_p\over |{\cal CL}^+|}$ qui divise $h_p^-$.

\bigskip

{\bf Lemme 2}

{\sl On prend les mmes hypothses que pour les dŽfinitions prŽcŽdentes. 

\art{a)}On a la suite exacte de $\F_q[G]$-modules~:

$$0\lra E/E^q\lra H\lra {\cal CL}[q]\lra 0.$$

\art{b)}$E/E^q$ est invariant par $\iota$, donc est un $\F_q[G^+]$-module.

\art{c)}$H^-\simeq {\cal CL}[q]^-={\cal CL}^-[q]$. Et on a la suite exacte~:

$$0\lra E/E^q\lra H^+\lra {\cal CL}[q]^+\lra 0.$$

\art{d)}$H$ est annihilŽ par $I$.

}
\goodbreak
{\bf preuve}

\art{a)}La premire flche envoie $u(1-\zeta_p)^k$ sur sa classe modulo ${\Q(\zeta_p)^*}^q$. Le noyau est
Žvidemment
$E^q$. Pour la seconde flche, $\alpha\in H$ est tel que $\alpha E_p={\euf a}^q(1-\zeta_p)^k$. On envoie
alors $\alpha$ sur $\euf a$. Puisque ${\euf a}^{q}=\alpha(1-\zeta_p)^{-k} E_p$, la flche est bien dŽfinie.
Supposons que $\alpha$ appartienne au noyau de la seconde flche. Alors il existe $k\in\Z$ et $\beta\in
\Q(\zeta_p)^*$ tel que $\alpha E_p=\beta^q(1-\zeta_p)^kE_p$ donc, il existe une unitŽ $u$ telle que
$\alpha=u\beta^q(1-\zeta)^k\buildrel \rm dans\ H\over\equiv u(1-\zeta_p)^k\in E$. 

\art{b)}Si $\alpha=u(1-\zeta_p)^k\in E$, alors $u^{1-\iota}={u\over \overline{u}}$ est une racine $2p$-ime
de l'unitŽ (Lemme 2 b) du Chapitre 6), donc appartient ˆ $E^q$ car $(2p,q)=1$. Ainsi $u\equiv
u^\iota\pmod{E^q}$. D'autre part, $(1-\zeta_p)^{1-\iota}={1-\zeta_p\over 1-\zeta_p^{-1}}=-\zeta_p\in
E_q$. Donc, $E/E^q$ est un $\F_q[G^+]$-module.

\art{c)}Puisque 2 est inversible dans $\F_q$, tout $\F_q[G]$-module $M=M^+\oplus M^-$, la suite exacte vue
en a) permet de conclure puisqu'on a montrŽ en b) que $\left ( E/E^q\right)^+=E/E^q$ et $\left (
E/E^q\right)^-=0$.

\art{d)}Puisque $I_{st}$ annihile ${\cal CL}[q]$, que $(1-\iota)$ annihile $E/E^q$ et que
$I=(1-\iota)I_{st}$, la suite exacte de a) nous permet de conclure.\qed
\bigskip
{\bf Lemme 3}

{\sl Si $x^p-y^q=1$ avec $p$, $q$ premiers impairs et $x,y$ entiers non nuls alors

\art{a)}${x-\zeta_p\over 1-\zeta_p}E_p={\euf a}^q$ pour un certain idŽal $\euf a$ de $E_p$.

\art{b)}La classe de $(x-\zeta_p)$ modulo les puissances $q$-ime est dans $H$.

}

{\bf Preuve}

\art{a)}c'est la relation $(9)$ vue au chapitre prŽcŽdent.

\art{b)}dŽcoule de a) en multipliant par $1-\zeta_p$.\qed

\bigskip

{\bf Lemme 4}

{\sl Notons $\mu_p$ l'ensemble de toutes les racines $p$-imes de l'unitŽ. Alors on a 

$$\sum_{\matrix{\zeta\in\mu_p\cr \zeta\ne 1\cr}}{1\over
x-\zeta}={\phi_p'(x)\over\phi_p(x)}\quad \hbox{et}\quad\sum_{\matrix{\zeta\in\mu_p\cr
\zeta\ne 1\cr}}{\zeta\over (1-\zeta)^2}={1-p^2\over 12}$$
o $\phi_p(x)=x^{p-1}+x^{p-2}+\cdots +x+1$ est le $p$-ime polyn™me cyclotomique.

}
{\bf Preuve}

On notera $\sum$ tout court pour $\dst\sum_{\matrix{\zeta\in\mu_p\cr \zeta\ne 1\cr}}$ dans cette preuve.
Alors l'intŽgrale  $\int(\sum (x-\zeta)^{-1})dx=\sum\int(x-\zeta)^{-1}dx=\ln (c\cdot \prod (x-\zeta))=\ln (
c\cdot\phi_p(x))$ pour une constante $c$. En dŽrivant, on trouve la premire ŽgalitŽ cherchŽe.

On en dŽduit que $\sum {1\over 1-\zeta}={\phi_p'(1)\over \phi_p(1)}={p(p-1)\over 2\cdot p}={p-1\over 2}$. En
dŽrivant la premire ŽgalitŽ, et en remplaant $x$ par 1, on trouve

$$\sum{1\over (1-\zeta)^2}={\phi_p'^2(1)-\phi_p''(1)\phi(1)\over \phi_p^2(1)}.$$

Or, $\phi_p''(1)=(p-1)(p-2)+(p-2)(p-3)+\cdots +2\cdot 1={(p-2)(p-1)p\over 3}$ car $\sum_{k=1}^n
k(k+1)={n(n+1)(n+2)\over 3}$. Donc,
$\sum{1\over (1-\zeta)^2}={-(p-5)(p-1)\over 12}$. Et finalement

$$\sum{\zeta\over (1-\zeta)^2}=-\sum{1\over 1-\zeta}+\sum{1\over (1-\zeta)^2}={1-p\over
2}+{-(p-5)(p-1)\over 12}={1-p^2\over 12}.$$\qed

A partir de cet endroit, nous sortons du cadre ҎlŽmentaire". Bien sžr, les choses pouvaient tre parfois
techniques et mme un peu "tordu". Nous avons passŽ une Žtape dŽjˆ importante quand nous avons introduit
les corps de nombres et une partie de la thŽorie de Galois. Maintenant nous allons passer un cran
supplŽmentaire dans l'abstraction, on doit dŽfinir les
$p$-adiques, c'est-ˆ-dire l'anneau $\Z_p$ et le corps $\Q_p$. [Se, pp. 23-26]

\bigskip\goodbreak
{\bf DŽfinitions-ThŽormes (Rappels sur les $p$-adiques)}

Soient $n\in\N$, $n\geq 2$, $p$ premier et 
$\varphi_n$ l'homomorphisme naturel 
de $\Z/p^n\Z$ dans $\Z/p^{n-1}\Z$ qui est Žvidemment surjectif.

On dŽfinit alors $$\Z_p=\lim_\leftarrow
(\Z/p^n\Z,\varphi_n)=\{(x_n)_{n=1}^\infty\in\prod_{n=1}^\infty\Z/p^n\Z\mid\varphi_n(x_n)=x_{n-1}\
\forall n\geq 2\}.$$
L'addition, la multiplication et la topologie sur $\Z_p$ sont hŽritŽes de celles induites par
l'anneau topologique produit $\dst \prod_{n=1}^\infty\Z/p^n\Z$. Les anneaux $\Z/p^n\Z$
Žtant munis de la topologie discrte, nous avons donc que $\dst
\prod_{n=1}^\infty\Z/p^n\Z$ est compact (Lemme de Tychonov), donc $\Z_p$ aussi puisqu'il est fermŽ.
\bigskip

$\Z_p$ possde les propriŽtŽs suivantes~:

\art{a)}$\Z_p/p^n\Z_p=\Z/p^n\Z$

\art{b)}$\Z_p$ est un {\it anneau local} d'idŽal maximal $p\Z_p$, donc les seuls
idŽaux de $\Z_p$ sont les $p^n\Z_p$, $n\in\N$\,; il suit que tout $x\in\Z_p$ s'Žcrit de
manire unique sous la forme $p^n\cdot u$ avec $u$ inversible.

\art{c)}La valuation $p$-adique $$\eqalign{v_p\,:\,\Z_p&\lra\N\cup\{\infty\}\cr
x&\longmapsto n\ \hbox{tel que }\ x=p^n\cdot u\cr
0&\longmapsto\infty\cr}$$
 induit une distance : $d(x,y)=p^{-v_p(x-y)}$ 
qui dŽfinit la topologie de $\Z_p$. On a en outre que $\Z$ est dense dans $\Z_p$ qui est complet. 
\art{d)}$\Z_p\cap\Q=\Z_{(p)}=\{{a\over b}\in\Q\mid p\not |\, b\}$
\bigskip

Notons $\Q_p$ le corps des fractions de $\Z_p$. Vu ce qui prŽcde, on a bien sžr que
$\Q_p=\Z_p[p^{-1}]$, donc tout $x\in \Q_p$ s'Žcrit  aussi de manire unique sous la forme
$p^n\cdot u$, o $u$ est un inversible de $\Z_p$ mais maintenant, $n\in \Z$\,; $n$
s'appellera aussi {\it valuation $p$-adique} que l'on notera aussi $v_p(x)$\,; elle
induira de la mme manire la topologie sur $\Q_p$, et on obtient facilement les rŽsultats
suivants.
\bigskip

\art{e)}Le corps $\Q_p$ , muni de la distance $d(x,y)=p^{-v_p(x-y)}$ est
localement compact et complet\,; le corps $\Q$ est dense dans $\Q_p$.

\art{f)}La distance $d$ est ``ultramŽtrique", c'est-ˆ-dire qu'elle vŽrifie
l'inŽgalitŽ suivante~: $$d(x,y)\leq\max(d(x,z),d(z,y)).$$
Nous obtenons gr‰ce ˆ cela le fait agrŽable que toute sŽrie de $\Q_p$ ou de $\Z_p$ est
convergente si et seulement si son terme gŽnŽral tend vers 0.
\bigskip

Nous aurions pu dŽfinir $\Q_p$, de manire tout ˆ fait analytique, comme
le complŽtŽ de $\Q$ pour la distance $d$, en voyant $\Z_p$ comme la boule unitŽ et
$p\Z_p$ comme la boule unitŽ privŽe de la sphre unitŽ.
\bigskip
{\bf ThŽorme (Lemme de Hensel)}
\medskip
Soient $f\in \Z_p[X_1,\ldots,X_m],\ x\in (\Z_p)^m,\ n,k\in \N$ et $j\in\N_m$.

Supposons que $$0\leq 2k<n,\quad f(x)\equiv 0\pmod {p^n}\quad\hbox{et}\quad v_p\left({\partial
f\over\partial X_j}(x)\right)=k.$$
Alors il existe un zŽro $y$ de $f$ dans $(\Z_p)^m$ qui est congru ˆ $x$ modulo $p^{n-k}$.

Il existe une forme de contraposŽe, un peu plus forte dont nous aurons besoin lors de l'appendice 2~: si
$f\in \Z_p[X]$ est irrŽductible, alors $\overline{f}\in \F_p[X]$ est une puissance d'un polyn™me
irrŽductible.

{\bf Preuve}

Ce thŽorme est dŽmontrŽ dans [Serre, pp. 28-30] pour la premire partie, la contraposŽe est prouvŽe
dans [Jac2, p. 573]. \qed

\bigskip
Nous sommes prt ˆ montrer un gros lemme, dont l'ŽnoncŽ ne paye pas de mine, mais qui sera bien
utile~!!

\bigskip\goodbreak

{\bf Proposition 5}

{\sl Si $x^p-y^q=1$ avec $p$, $q$ premiers impairs et $x,y$ entiers non nuls alors l'ŽlŽment
$(x-\zeta)^{1-\iota}$ est non-trivial dans $H$ (donc dans $H^-$).

}

{\bf Preuve}

On va prouver cela par l'absurde. Supposons donc que $(x-\zeta)^{1-\iota}=1$ dans $H$. Cela veut dire qu'il
existe $\alpha\in\Q(\zeta_p)$ tel que ${x-\zeta_p\over x-\zeta_p^{-1}}=\alpha^q$. Les relations de
Cassels nous apprennent que $x\equiv 1\pmod {p^{q-1}}$. Posons $\mu={x-1\over 1-\zeta_p}$. On a
 $x-1\in p^{q-1}\Z\subset\pi^{(p-1)(q-1)}E_p$, donc $\mu\in E_p$ et $v_\pi(\mu)\geq (p-1)(q-1)-1\geq 4$. On
vŽrifie que ${x-\zeta_p\over 1-\zeta_p}=1+\mu$ et

$${1+\mu\over 1+\overline{\mu}}=-\zeta_p^{-1}\alpha^q=\beta^q\quad\hbox{ avec $\beta=-\zeta_p^{-{1\over
q}}\alpha$,}$$

 c'est toujours possible, car $(q,2p)=1$. Choisissons des racines $q$-ime de $1+\mu$ et de
$1+\overline{\mu}$, qu'on notera $\root q\of{1+\mu}$ et $\root q\of{1+\overline{\mu}}$ (on sort peut-tre de
$\Q(\zeta_p)$), telles que 

$${\root q\of{1+\mu}\over \root q\of{1+\overline{\mu}}}=\beta\eqno{(c)}$$

Posons encore 

$$\eta=\left (\root q\of{1+{\mu}}+\zeta_p^{-{1\over q}}\root q\of{1+\overline{\mu}}\right )^q.$$

On prŽtend que $\eta\in U_p$. En effet~: notons $O$ l'ensemble des entiers algŽbrique de $\C$ sur $\Z$. Si
$z_1,z_2 \in O$, alors $\root q\of z\in O$ et $z_1+z_2$ et $z_1z_1\in O$ donc $\eta\in O$, puisque $\mu\in
E_p$. D'autre part, 

$$\eta=\left (\root q\of{1+\overline{\mu}}\right)^q\cdot (\beta+\zeta_p^{-{1\over
q}})^q=(1+\overline{\mu})(-\alpha\zeta_p^{-{1\over q}}+\zeta_p^{-{1\over
q}})^q=(1+\overline{\mu})(1-\alpha)^q\zeta_p^{-1}\in\Q(\zeta_p).$$

Puisque $A+B$ divise $A^q+B^q$, on a $\root q\of{1+{\mu}}+\zeta_p^{-{1\over q}}\root
q\of{1+\overline{\mu}}$ divise $(1+\mu)+\zeta_p^{-1}(1+\overline{\mu})$ dans $O$. Or,
$(1+\mu)+\zeta_p^{-1}(1+\overline{\mu})={x-\zeta_p\over 1-\zeta_p}+\zeta_p^{-1} {x-\zeta_p^{-1}\over
1-\zeta_p^{-1}}={\zeta_p^{-1}-\zeta_p\over 1-\zeta_p}\in U_p\subset O^*$. Ainsi, $\eta\in \Q(\zeta_p)\cap
O^*=U_p$. Tout ŽlŽment de $U_p$ a une norme qui est {\it a priori} de $\pm 1$. Mais comme $\Q(\zeta_p)$ est
totalement complexe, l'ensemble des plongements est par paires de conjuguŽs, donc la norme est positive et
donc $N(\eta)=1$ o $N=N_{\Q(\zeta_p)/\Q}$. Puisque  $[\Q_p(\zeta_p):\Q_p]=p-1$, $N_{\Q_p(\zeta_p)/\Q_p}$
prolonge $N_{\Q(\zeta_p)/\Q}$. 

On va donc chercher des racines de $1+\mu$ et de $1+\overline{\mu}$
satisfaisant $(c)$ dans $\Q_p(\zeta_p)$. En posant (comme dans la preuve du Lemme 2 du Chapitre 4)
$\pmatrix{{1\over q}\cr\noalign{\vskip 4pt} k\cr}=\left ({1\over q}({1\over q}-1)({1\over q}-2)\cdots
({1\over q}-(k+1))
\right )/k!$, et par convention $\pmatrix{{1\over q}\cr\noalign{\vskip 2pt} 0\cr}=1$ on a, tout
d'abord formellement~:

$$(1+\mu)^{1\over q}=\sum_{k=0}^\infty \pmatrix{{1\over q}\cr\noalign{\vskip
4pt} k\cr}\mu^k.$$

D'abord, $\pmatrix{a\cr k\cr}\in\Z_p$ si  $a\in\Z_p$, car c'est vrai si $a\in\Z$ et donc, par continuitŽ et
puisque $\Z$ est dense dans $\Z_p$, et Žvidemment, ${1\over q}\in\Z_p$. Ensuite, la sŽrie converge, car
$v_\pi(\mu)\geq 4$ et donc la norme ($p$-adique) des termes de la sŽrie tend vers 0. Donc, dans
$\Q_p(\zeta_p)$ une racine $q$-ime de $(1+\mu)$ existe. 

Il va falloir maintenant travailler pour affirmer
qu'il est possible de s'arranger pour satisfaire la condition $(c)$. 

Dans le cas o
$p\not\equiv 1\pmod q$, $\Q_p(\zeta_p)$ n'a pas de racine $q$-ime de l'unitŽ autre que 1. En effet, si
$\omega$ est une racine $q$-ime de l'unitŽ diffŽrente de 1, alors $\omega\equiv 1\pmod{\pi
\Z_p[\zeta_p]}$, car 1 est l'unique racine $q$-ime de l'unitŽ dans le corps rŽsiduel
$\Z_p[\zeta_p]/\pi\Z_p[\zeta_p]=\F_p$ (sinon, $q|p-1$). Et si $\omega=1+\pi u$, alors $1=\omega^q=1+q\pi
u+O(\pi^2)$, avec la notation $O(\pi^2)$ qui signifie ici seulement qu'il s'agit d'un multiple de
$\Z_p[\zeta_p]$ de $\pi^2$; donc $\pi|q$, ce qui est impossible. Dans ce cas, la condition $(c)$ est
automatiquement satisfaite, car les racines $q$-ime, si elles existent, sont uniques. Et on peut poser
$\root q\of {1+\mu}=(1+\mu)^{1\over q}$ et $\root q\of {1+\overline{\mu}}=(1+\overline{\mu})^{1\over q}$ et
dŽfinir $\mu$ comme avant.

Dans le cas o $p\equiv 1\pmod q$, le corps rŽsiduel contient $q$ racines distinctes de l'unitŽ (Lemme
IMP, Chapitre 5) qui, par le Lemme de Hensel, se relvent en $q$ racines de l'unitŽ dans $\Z_p[\zeta_p]$.
Soit
$\omega$, une telle racine diffŽrente de 1. {\it A priori}, il existe $k$ tel que 
${(1+\mu)^{1\over q}\over \omega^k(1+\overline{\mu})^{1\over q}}=\beta$. On montre facilement que
$1+\pi\Z_p[\zeta_p]$ est un groupe multiplicatif, donc ${(1+\mu)^{1\over q}\over
(1+\overline{\mu})^{1\over q}}\equiv 1\pmod
{\pi \Z_p[\zeta_p] }$. Puisque $\Z_p[\zeta_p]/\pi\Z_p[\zeta_p]=\F_p$, il existe $u_0\in \Z$ et
$v\in\Z_p[\zeta_p]$ tels que
$\beta =u_0+\pi v$ et $0\leq u_0<p$. On veut montrer que $k=0$ ou, ce qui est Žquivalent, $u_0=1$ (car les
$\omega^k$ sont distincts modulo $\pi$). D'abord, on a $\beta^q={1+\mu\over 1+\overline{\mu}}\equiv
1\pmod\pi$. d'o,
$u_0^q\equiv 1\pmod\pi$, ou encore $u_0^q\equiv 1\pmod p$ (car $u_0\in\Z$). D'autre part, 
$(\beta\overline{\beta})^q={1+\mu\over 1+\overline{\mu}}\cdot{1+\overline{\mu}\over 1+{\mu}}=1=1^q$, donc
$\beta\overline{\beta}=1$, car dans $\Q(\zeta_p)$, il n'y a que 1 comme racine $q$-ime de 1. Mais
$\overline{\beta}=u_0+\overline{ \pi v}$. On en dŽduit que $u_0^2\equiv 1\pmod p$ et on vient de voir
que $u_0^q\equiv 1\pmod p$, donc $u_0\equiv 1\pmod p$, car $(2,q)=1$. 

Donc, dans tous les cas, la condition $(c)$ est satisfaite dans $\Q_p(\zeta_p)$. Et alors, il est possible
de dŽfinir $\eta$ relativement ˆ ce choix. Posons 
$$u=(1+\mu)^{1\over q}+\zeta^{-{1\over q}} (1+\overline{\mu})^{1\over q}.$$

Puisque ${\mu\over\overline{\mu}}=-\overline{\zeta}$ est une unitŽ, en effectuant le dŽveloppement binomial,
on trouve que
$u=(1+{\mu\over q})+\zeta_p^r(1+{\overline{\mu}\over q})+O(\mu^2)$, et $r\in\Z$ tel que $r\equiv -{1\over
q}\pmod {p\Z_p[\zeta_p]}$. Cela veut dire que 
$$u\equiv (1+{\mu\over q})+\zeta_p^r(1+{\overline{\mu}\over q})\pmod{\mu^2\Z_p[\zeta_p]}.$$

Donc, en faisant un petit calcul, on trouve $u=(1+\zeta_p^r)(1+{x-1\over q}{1-\zeta_p^{r+1}\over
(1-\zeta_p)(1+\zeta_p^r)})+O(\mu^2)$.

Or, $1+\zeta_p^r={1-\zeta_p^{2r}\over 1-\zeta_p^r}$ est une unitŽ, sa norme vaut donc $1$. On a ainsi~:

$$\eqalignno{N(u)&=\prod_{\matrix{\zeta\in\mu_p\cr \zeta\ne 1}}(1+\underbrace{{x-1\over
(1-\zeta)}}_{\simeq\mu}{1-\zeta^{r+1}\over q(1+\zeta^r)})+O(\mu^2)\cr
&=1+{x-1\over q}\sum_{\matrix{\zeta\in\mu_p\cr \zeta\ne 1}}\left ( {1-\zeta^{r+1}\over
(1-\zeta)(1+\zeta^r)}\right )+O(\mu^2).&(i)\cr}$$

Posons $\pi'=\zeta-1$, bien sžr $\pi=\pi'$ si $\zeta=\zeta_p$. Alors on a 
$${1-\zeta^{r+1}\over (1-\zeta)(1+\zeta^r)}={1-(1+\pi')^{r+1}\over
-\pi'(1+(1+\pi')^r)}={-(r+1)\pi'+O(\pi^2)\over -\pi'(2+r\pi')+O(\pi^3)}={r+1\over 2}+O(\pi),$$ puisque 2
est inversible dans $\Z_p[\zeta_p]$. Or, $\pi^3 |x-1$, donc $\pi(x-1)|\mu^2={(x-1)^2\over\pi^2}$ et ainsi, 
$$ N(u)=1+{x-1\over q}{(r+1)(p-1)\over 2}+O(\pi(x-1)).$$

En effectuant le dŽveloppement binomial, et en utilisant le fait que $\pi(x-1)|(x-1)^2$,
on trouve que~:

$$1=N(\eta)=N(u)^q=1+\underbrace{{(x-1)(r+1)(p-1)\over 2}+O(\pi(x-1))}_{=0}.$$

Donc, $\pi(x-1)|{(x-1)(r+1)(p-1)\over 2}$. Donc $\pi|r+1$ ce qui implique que $-1\equiv r\equiv -{1\over
q}\pmod p$. Donc $q\equiv 1\pmod p$. On supposera donc, dans la suite que $q\equiv 1\pmod p$, et on reprend
le calcul de $u$ ˆ $O(\mu^3)$ prs. Mais avant cela, on voit que $\zeta_p^{-{1\over q}}=\zeta_p^{-1}$ et
$\mu+\zeta^{-{1\over q}}\overline{\mu}=0$. On obtient alors~:
$$\eqalign{u&=1+\pmatrix{{1\over q}\cr\noalign{\vskip
4pt} 2\cr}\mu^2+\zeta_p^{-1}(1+\pmatrix{{1\over q}\cr\noalign{\vskip
4pt} 2\cr}\overline{\mu}^2)+O(\mu^3)\cr
&=(1+\zeta_p^{-1})\left (1+{1-q\over 2q^2}(x-1)^2{\zeta_p\over (1-\zeta_p)^2}\right )+O(\mu^3),\cr}$$

car $\dst\pmatrix{{1\over q}\cr\noalign{\vskip 4pt} 2\cr}={1-q\over 2q^2}$ et
$(1-\zeta_p^{-1})^2=\zeta_p^{-2}(1-\zeta_p)^2$. En prenant la norme, en remarquant que $(1+\zeta_p^{-1})$
est de norme 1 et en raisonnant comme en $(i)$, on trouve
$$N(u)=1+{1-q\over 2q^2}(x-1)^2\sum_{\matrix{\zeta\in\mu_p\cr \zeta\ne 1}}{\zeta\over (1-\zeta)^2}\
+O(\mu^3).$$

On a vu au Lemme 4 que $\dst\sum_{\matrix{\zeta\in\mu_p\cr
\zeta\ne 1\cr}}{\zeta\over (1-\zeta)^2}={1-p^2\over 12}$.
A nouveau en effectuant le dŽveloppement binomial, et en utilisant le fait que $\mu^3|(x-1)^4$,
on trouve que~:

$$1=N(\eta)=N(u)^q=1+\underbrace{{(1-q)(x-1)^2(1-p^2)\over 2\cdot q\cdot 12}+O(\mu^3)}_{=0}.$$

Or, $\dst {(1-q)(1-p^2)\over 2\cdot1\cdot 12}$ est un ŽlŽment de $\Z_p[\zeta_p]$, mme quand $p=3$, car
$q\equiv 1\pmod p$. Donc, $\dst\mu^3={(x-1)^3\over\pi^3}$ divise $\dst{(1-q)(x-1)^2(1-p^2)\over 2\cdot q\cdot
12}$ et donc, $x-1$ divise ${(q-1)(1-p^2)\pi^3\over  24 q}$. Comme $x\equiv 1\pmod {p^{q-1}}$, on obtient
que
$p^{q-1}|{(q-1)\pi^3\over 3}$, donc $p|\pi^3$ et c'est impossible ou $p^{q-1}| q-1$ et la division est dans
$\Z$, mais c'est aussi impossible, car $p^{q-1}$ est bien plus grand que $q-1$. Donc on a notre
contradiction, ce qui implique que $(x-\zeta)^{1-\iota}$ n'est pas trivial dans $H$.\qed
\bigskip
{\bf Corollaire 6}
 
{\sl Soit $p$ et $q$ des nombres premiers impairs distincts. Si $p\not \hskip-0.5pt|\  h_q^-$ ou $q\not
\hskip-0.5pt|\ h_p^-$, alors l'Žquation $x^p-y^q=1$, n'a pas de solution non nulles $x,y\in\Z$.

}

{\bf Preuve}

Par symŽtrie, supposons que $q\not\hskip-0.8pt|\  h_p^-$. On a vu (Lemme 2 c)) que $H^-\simeq {\cal
CL}^-[q]$. Puisque
${\cal CL}^-[q]=1$ ou $q$ divise $|{\cal CL}^-|$ qui divise $h_p^-$ (cf. DŽfinition entre les Lemmes 1 et
2), on en dŽduit que
$H^-$ est trivial. Or, si $x^p-y^q=1$, la proposition 5 nous dit que
$(x-\zeta)^{1-\iota}$ est un ŽlŽment non-trivial de $H^-$. C'est une contradiction et le  corollaire est
prouvŽ.\qed
\bigskip

{\bf Lemme 7}

{\sl Soit $p$ un nombre premier impair. Soit $\gamma\in\N$ une {\it racine primitive modulo $p$},
c'est-ˆ-dire  que la classe de $\gamma$ engendre cycliquement $\F_p^*$. Posons, pour $i=1,\ldots ,p-1$,
$1\leq\gamma_i\leq p-1$, tel que $\gamma_i\equiv \gamma^i\pmod p$. Posons $F_p(X)=\sum_{i=1}^{p-1}\gamma_i
X^i$. Soit enfin
$\zeta_{p-1}$ une racine primitive $p-1$-ime de l'unitŽ. Alors 

$$h_p^-={\left |\dst\prod_{k=1}^{p-1\over 2}F_p(\zeta_{p-1}^{2k-1})\right |\over (2p)^{p-3\over 2}}.$$

}

{\bf Preuve}

 Ce Lemme est assez long ˆ prouver, mais pas trop dur. Une version semblable ˆ cet ŽnoncŽ se trouve dans
[Edw, p. 225].\qed
\bigskip

Nous allons montrer que la conjecture de Catalan est vraie si $p$ ou $q$ est
infŽrieur ou Žgal ˆ 43. Mais nous n'aurons besoin de ce rŽsultat que pour $p$ ou $q\leq 11$. Alors, on va
calculer ҈ la main" $h_3^-$, $h_5^-$ et $h_7^-$. Ainsi, toute la preuve de Catalan sera fait sans calcul
informatique... mais pour tre franc, j'ai tout de mme vŽrifiŽ mes calculs par ordinateur, pas con le
bourdon ! 
\bigskip
{\bf Lemme 8}

{\sl $h_3^-=h_5^-=h_7^-=1$ }

{\bf Preuve}

Si $p=3$, $\zeta_{p-1}=-1$, $\gamma=2$ et $F_3(X)=X^2+2X$. Donc $h_3^-={|F_3(-1)|\over (6)^0}=|-1|=1$.

Si $p=5$, $\zeta_{p-1}=i$, $\gamma=2$ et $F_5(X)=X^4+3X^3+4X^2+2X$, $F_5(i)=(-3-i)$ et $F_5(i^3)=(-3+i)$. Et
donc $h_5^-={F_5(i)\cdot F_5(i^3)\over (10)^1}=1$.

Enfin, si $p=7$, $\zeta_{p-1}={1+\sqrt{3} i\over 2}$, $\gamma=3$ et $F_p(X)=X^6+5X^5+4X^4+6X^3+2X^2+3X$. Et
donc, $h_7^-={|F_7(\zeta_6)\cdot F_7(\zeta_6^3)\cdot F_7(\zeta_6^5)|\over
(14)^2}={|(-4-2\sqrt{3}i)(-7)(-4+2\sqrt{3}i)|\over 196}=1$.\qed

\bigskip

{\bf ThŽorme 9 (ThŽorme 2 de Mih$\breve{\bf a}$ilescu)}

{\sl Si $p$ et $q$ sont des premiers impairs distincts avec $p$ ou $q\leq 43$, alors l'Žquation $x^p-y^q=1$
n'a pas de solutions entires non nulles
} 

{\bf Preuve}

En faisant des calculs comme au Lemme 8, mais un peu plus long (et ˆ l'ordinateur, cette fois), on trouve
que $h_p^-=1$ si $p\leq 19$. Pour $p=23,29,31,37,41,43$, on a respectivement $h_p^-=3,8,9,37,11^2,211$.
Comme $h_2=h_3^-=h_{11}^-=1$, il n'y a pas de problme pour $p=23,29,31,37,41$. Pour $43$, il faut calculer
$h_{211}^-$. On calcule que

$$h_{211}^-=3^2\cdot 7^2\cdot 41\cdot 71\cdot 181\cdot 281^2\cdot 421\cdot1051\cdot 12251\cdot
113981701\cdot 4343510221$$

et donc, $43\not \hskip-0.4pt|\  h_{211}^-$. En revanche, $h_{47}^-=5\cdot 139$ et 
$$h_{139}^-=3^2\cdot 47^2\cdot 277^2\cdot 967\cdot 1188961909.$$

En revanche, $47^{138}\equiv 7507\pmod {139^2}$ et $139^{46}\equiv 1035\pmod{47^2}$, donc, si on le
voulait, on pourrait encore ajouter $47$ ˆ la liste, mais comme je l'ai dŽjˆ dit, c'est seulement pour $3$,
$5$, et $7$ que c'est vraiment important.\qed 

\vfill\eject

\long\def\art#1{{\parindent0pt\item{#1}}\hangindent=7mm\hangafter=-20}
\long\def\artart#1{{\parindent0pt\item{#1}}\hangindent=12mm\hangafter=-20}
\font\para=cmbx12 at 18pt
\def\O{\hbox{$\cal O$}}
\def\U{\hbox{$\cal U$}}
\def\m{\hbox{\rs m\!}}
\def\dst{\displaystyle}
\font\doub=msbm10 at 10pt
\def\lra{\longrightarrow}
\def\qed{\hfill$\square$}
\def\gfP{\relax\ifmmode\bbP\else $\bbP$\fi}
\def\gP{{\euf P}}
\def\P{{\cal P}}
\def\QQ{{\cal Q}}
\def\Log{{\rm Log}}

\def\ggP{{\bf P}}
\newcount\chapnomb \chapnomb=1
\newcount\parnomb \parnomb=1
\pageno =40

\parindent0pt
\centerline{\para CHAPITRE 8}
\bigskip
{\para Troisime thŽorme de Mih$\taille{18}\breve{\bf a}$ilescu~: 
\bigskip

\centerline {$\taille {15} p<4q^2$ et $\taille {15}
q<4p^2$ }}
\bigskip
Le titre est assez explicite pour voir ce qu'on va prouver. Ce thŽorme permet d'Žviter d'avoir recourt ˆ
l'ordinateur et d'utiliser des rŽsultats difficiles sur les formes logarithmiques.
\bigskip

{\bf Lemme 1}

{\sl Soit $3\leq p<q$ deux nombres premiers. Alors on a les rŽsultats suivants~:

\art{a)}${\log(q)\over q-1}-{1\over q}>{1\over q}$.

\art{b)}${1\over x}>\log (1+{1\over x})$ si $x\geq 1$.

\art{c)}$p^{q-1}> q^{p-1}+q$.

}

{\bf Preuve}

Prouvons a)~:  $q\log(q)-2(q-1)>0$, car $q\geq 5$ et la fonction $q\log(q)-2(q-1)$ est croissante (sa
dŽrivŽe vaut $\log(q)-1$). Donc
${\log(q)\over q-1}>{2\over q}$, ou encore ${\log(q)\over q-1}-{1\over q}>{1\over q}$.

Prouvons b)~: soit $0<y={1\over x}\leq 1$. On a $e^{y}-y>1$, car l'inŽgalitŽ non stricte est vraie pour
$y=0$ et la fonction $e^{y}-y$ est croissante.  Donc,  ${1\over x}>\log (1+{1\over x})$.

Prouvons c)~: posons $f(x)={\log(x)\over x-1}$. Prouver la partie c) revient ˆ prouver que $f(p)-f(q)>
\log(1+{1\over q^{p-2}})\cdot {1\over (p-1)(q-1)}$. Par le thŽorme des accroissements finis, il existe
$c\in ]p,q[$ tel que $f(p)-f(q)=(p-q)\cdot f'(c)$. On voit que $f'(x)={-(x\log(x)-(x-1))\over x(x-1)^2}<0$
si $x>2$. De plus, $f''(x)={2x^2\log(x)-3x^2+4x-1\over x^2(x-1)^3}>0$ si $x>2$. Ainsi, $f'(p)
<f'(c)<f'(q)<0$. Donc $(p-q) f'(c)\geq (p-q)f'(q)$. Donc, il suffit de montrer que 

$$(q-p){q\log (q)-(q-1)\over q(q-1)^2}> \log(1+{1\over q^{p-2}})\cdot {1\over (p-1)(q-1)}.$$

Or, gr‰ce ˆ la partie b),  $\log (1+{1\over q^{p-2}})\leq {1\over q^{p-2}}$. Donc, l'inŽgalitŽ ˆ prouver
est vraie si $(q-p){q\log(q)-(q-1)\over q(q-1)}> {1\over q^{p-2} (p-1)}$, ou encore si
$(q-p)(p-1)q^{p-2}\cdot\left ({\log(q)\over q-1}-{1\over q}\right )> 1$. Mais cela est vrai, car on a vu
en a) que ${\log(q)\over q-1}-{1\over q}>{1\over q}$.\qed

\bigskip

{\bf Lemme 2}

{\sl Si $x^p-y^q=1$ avec $p$, $q$ premiers impairs tels que $p,q\geq 11$, et $x,y$ entiers non nuls,
alors  $$|x|\geq \max (p^{q-1}-1,q^{p-1}+q)\quad\hbox{et}\quad  |y|\geq \max (q^{p-1}-1,p^{q-1}+p).$$

}
{\bf Preuve}

Montrons le lemme pour $|x|$. Le fait que $|x|\geq
p^{q-1}-1$ vient de l'existence, par les relations de Cassels, de $u\in\Z$ tel que $ x-1=p^{q-1} u^q $. Si
$p<q$, on a vu au lemme prŽcŽdent que $p^{q-1}-1\geq q^{p-1}+q$. Donc le Lemme 2 est prouvŽ pour ce cas.
Nous sommes passŽs comme chat sur braise sur le fait que lors de la preuve des relations de Cassels
on avait considŽrŽ des $x$ et $y$ positifs. Mais si ce n'est pas le cas on remplace $x$ par $-x$ et $y$ par
$-y$ et l'affirmation est correcte !

Supposons $p>q$. On a ${y^q+1\over y+1}={((y+1)-1)^q-1\over y+1}=(y+1)^{q-1}-q(y+1)^{q-2}+\cdots +q\equiv
q\pmod{y+1}$. Or, les relations de Cassels nous disent qu'il existe $c\in \Z$ tel que ${y^q+1\over
y+1}=qc^p$; (on remarque que $c>0$). Les mmes relations nous disent que $q^{p-1}$ divise $y+1$. Donc, on a
la congruence~: $q\cdot c^p\equiv q\pmod{q^{p-1}}$ ou encore $c^p\equiv 1\pmod q^{p-2}$. Cela veut dire que
l'ordre de $c$ modulo $q^{p-2}$ est 1 ou $p$. Si c'est $p$, il doit de toute faon diviser
$\varphi(q^{p-2})=(q-1)\cdot q^{p-3}$. Or $p$ ne divise pas $q-1$ (il est supŽrieur ˆ $q$) et ne divise pas
non plus $q^{p-3}$. Donc $c\equiv 1\pmod {q^{p-2}}$. Si $c=1$, alors ${y^q+1\over y+1}=q$. Voyons que c'est
impossible~: si $y>2$, on a ${y^q+1\over y+1}=y^{q-1}-y^{q-2}+y^{q-3}+\cdots
+1=y^{q-2}(y-1)+y^{q-4}(y-1)+\cdots +y(y-1)+1\geq y^{q-2}+1>2^{q-2}+1>q$ car on a vu que $q>7$. Si $y<-2$,
posons $z=-y$.  On a ${y^q+1\over y+1}=z^{q-1}+z^{q-2}+\cdots +z+1\geq z^{q-1}>q$. Si $y=2$ alors $2^q+1=3q$
implique $q=3$ et $q\geq 11$, c'est donc impossible. Si $y=-2$, $q=1-2^q$, a devient loufoque ! Si $y=1$,
$q=1$, aussi. Et enfin si $y=-1$, $x=0$, ˆ nouveau impossible. Tout a implique que $c\geq 1+q^{p-2}$.
Mais alors, les relation de Cassels nous disent encore que $|x|=q|b|c>qc\geq q^{p-1}+q$. 

On dŽmontre le Lemme pour $|y|$ de manire identique.\qed

\bigskip
\goodbreak

{\bf Lemme 3}

{\sl Si $x^p-y^q=1$ avec $p$, $q$ premiers impairs et $x,y$ entiers non nuls. Rappelons que $G={\rm
Gal}(\Q(\zeta_p)/\Q)$. Posons
$X={\rm Ann}_{\Z[G]}([x-\zeta_p]\in H):=\{\theta\in\Z[G]\mid \exists \alpha\in \Q(\zeta_p)^*,\hbox{ avec }
(x-\zeta_p)^\theta=\alpha^q\}$. Il est clair que $X$ est un sous-groupe additif de $\Z[G]$. Alors
l'application

$$\eqalign{X&\longrightarrow \Q(\zeta_p)^*\cr
\theta&\longmapsto\alpha,\cr}$$

o $\alpha$ est l'unique ŽlŽment de $\Q(\zeta_p)^*$ tel que $(x-\zeta_p)^\theta=\alpha^q$, est un
homomorphisme de groupe injectif.

} 
\goodbreak
{\bf Preuve}

L'application est bien dŽfinie, car $\Q(\zeta_p)$ ne contient pas de racine $q$-ime de 1 autre que 1
lui-mme, donc $\alpha$ est unique. Le fait que ce soit un homomorphisme de groupe est une simple
vŽrification. Reste ˆ montrer l'injectivitŽ~: soit $\theta=\sum_{\tau\in G}n_\tau\tau$ tel que
$(x-\zeta_p)^\theta=1$. La relation (9) vu au Chapitre 6 nous montre que les conjuguŽs de
$\beta:={x-\zeta_p\over 1-\zeta_p}$ sont premiers entre eux. De plus, ce ne sont pas des unitŽs, car pour
toute racine primitive $p$-ime de l'unitŽ $\zeta$, on a~:

$$N_{\Q(\zeta_p)/\Q}\left({x-\zeta\over 1-\zeta}\right)=N\left({x-\zeta\over
1-\zeta}\right)={|\prod_{i=1}^{p-1}(x-\zeta^i)|\over p}\geq {(|x|-1)^{p-1}\over p}>1,$$

La dernire inŽgalitŽ venant du fait que $|x|\geq 3$ (Lemme 2). On en dŽduit donc que
les conjuguŽs de $\beta$ sont divisibles par des idŽaux premiers distincts. Soit $\sigma\in G$. Puisque
$(x-\zeta_p)^\theta=1=\sigma(1)$, alors $(x-\sigma(\zeta_p))^\theta=1$. On en dŽduit que
$N(x-\zeta_p)^{\sum_{\tau\in G}n_\tau}=1$, donc $\sum_{\tau\in G}n_\tau=0$, car $N(x-\zeta_p)> 1$. Ainsi, 
${(x-\zeta_p)^\theta\over (1-\zeta_p)^0}= {\prod_{\tau\in G} (x-\tau(\zeta_p))^{n_\tau}\over \prod_{\tau\in
G}(1-\zeta_p)^{n_\tau}}=1$. Or, les
$1-\tau(\zeta_p)$ diffrent de $1-\zeta_p$ d'une unitŽ. Ainsi, $ \prod_{\tau\in G}
{\left(x-\tau(\zeta_p)\over 1-\tau(\zeta_p)\right)}^{n_\tau}$ est une unitŽ. On en dŽduit que chaque
$n_\tau$ est nul, car sinon les ${x-\tau(\zeta_p)\over 1-\tau(\zeta_p)}$ auraient des facteurs en commun,
ce qui n'est pas le cas comme nous venons de le voir. Donc $\theta=0$, ce qui prouve l'injectivitŽ.\qed
\bigskip

{\bf DŽfinition}

Soit $\alpha\in\Q(\zeta_p)$. Alors, il est clair que $\alpha E_p={\euf a}{\euf b}^{-1}$ avec ${\euf a}$ et
$\euf b$ des idŽaux de $E_p$ premiers entre eux. L'idŽal $\euf b$ est appelŽ {\it l'idŽal dŽnominateur de
$\alpha$}.
\bigskip

{\bf Lemme 4}

{\sl Sous les mme hypothses que pour la dŽfinition prŽcŽdentes, on a

\art{a)}${\euf a}{\euf b}^{-1}\cap E_p={\euf a}$

\art{b)}${\euf b}=\{ x\in E_p\mid x\cdot\alpha\in E_p\}$

\art{c)}$\alpha$ et $\alpha-1$ ont mme idŽal dŽnominateur.

}

{\bf preuve}

\art{a)}Il est Žvident que ${\euf a}\subset {\euf a}{\euf b}^{-1}\cap E_p:={\euf a}'$. D'autre part, ${\euf
a}{\euf b}^{-1}\supset {\euf a}'$. Donc ${\euf a}\supset {\euf a}'\cdot {\euf b}$, ce qui veut dire que
${\euf a}\supset {\euf a}'$, car
$\euf a$ et $\euf b$ sont premiers entre eux.

\art{b)}Il est Žvident que ${\euf b}\subset \{ x\in E_p\mid x\cdot\alpha\in E_p\}:={\euf b}'$. Soit
$x\in{\euf b}'$, alors $x\alpha\in E_p\cap {\euf a}{\euf b}^{-1}\buildrel a)\over ={\euf a}$. Donc $x\in
\alpha^{-1}E_p\cdot {\euf a}={\euf b}'$.

\art{c)}DŽcoule de b).\qed  
\bigskip

On avait dŽjˆ dŽfini la branche principale du logarithme au Lemme 10 du Chapitre 5, nous en avons ˆ
nouveau besoin, ainsi que de quelques propriŽtŽs~:
\bigskip

{\bf DŽfinition}

Soit $z\in\C$. Posons $\Log(z)=\log|z|+i\arg (z)$ o $-\pi<\arg(z)\leq\pi$. On appelle $\Log(\cdot)$ la {\it
branche principale du logarithme}. On peut voir que $e^{\Log(z)}=z$, pour tout $z\in\C^*$. Si $w\in\C$ est
tel que $e^{w}=z$, alors il existe $k\in\Z$ tel que $w=\Log(z)+2k\pi i$. 
\bigskip
{\bf Lemme 5}

{\sl 

\art{a)}Si $|z|<1$, alors $\Log(1-z)=-\sum_{n=1}^\infty {z^n\over n}$.

\art{b)}Si $|z|=1$, alors $|\arg (z)|\leq |f(z)|$ pour toute dŽtermination $f$ du logarithme.

\art{c)}Soit $z\in \C$. Si $r\in\N$ est un nombre impair tel que $|\arg (z^r)|<\lambda\leq\pi$, alors il
existe un entier $k$, avec ${1-r\over 2}\leq k\leq {r-1\over 2}$ tel que 

$$\left |\arg(z)-{k\cdot 2\pi\over r}\right |<{\lambda\over r}\leq{\pi\over r}.$$

De plus, si $|\arg (z)|\leq {\pi\over r}$, alors $k=0$.

}
\goodbreak

{\bf Preuve}

\art{a)} On peut lire ce rŽsultat dans tout livre portant sur l'analyse complexe, par exemple
[Mac, pp.77 et 134].

\art{b)}Si $f$ est la branche principale, c'est Žvident. Sinon, il existe $k\in\Z\setminus\{0\}$ tel que 

$$|f(z)|=|\underbrace{\log(|z|)}_{=0}+i\arg(z)+k\cdot 2\pi i|=|\arg(z)+k\cdot 2\pi| \buildrel k\ne
0\over\geq
\min(|2k-1|,|2k+1|)\cdot\pi\geq\pi\geq|\arg(z)|.$$

\art{c)}Il existe $k\in\Z$ tel que $\arg(z^r)=r\cdot\arg(z)+k\cdot 2\pi$. On a donc
$-\pi<r\cdot\arg(z)+k\cdot 2\pi\leq\pi$. Mais puisque $-\pi<\arg (z)\leq \pi$, on a $-\pi<r\cdot\pi+k\cdot
2\pi$, donc $k>{-1-r\over 2}$ ou $k\geq {1-r\over 2}$, car $r$ est impair. On montre de mme que $k\leq
{r-1\over 2}$. De $|\arg (z^r)|<\lambda$, on tire $|r\cdot\arg(z)+k\cdot 2\pi|<\lambda $ et donc $$
|\arg(z)+{k\cdot 2\pi\over r}|<{\lambda\over r},\eqno{(i)}$$ le signe $+$ ne dŽrange pas, car $k$ varie
entre
$-{r-1\over 2}$ et ${r-1\over 2}$.

\art{}Si $-{\pi\over r}\leq \arg(z)\leq{\pi\over r}$, alors $-{\pi\over r}+{k\cdot
2\pi\over r}\leq \arg(z)+{k\cdot 2\pi\over r}\leq{\pi\over r}+{k\cdot 2\pi\over r}$. Donc, si $k\ne 0$, 

$$\left |\arg(z)+{k\cdot 2\pi\over r}\right |\geq \min(|2k-1|{\pi\over r},|2k+1|{\pi\over r})\geq
{\pi\over r}\geq {\lambda\over r}.$$

 cela est une contradiction avec l'Žquation $(i)$. Donc $k=0$.\qed

\bigskip

{\bf Proposition 6}

{\sl Si $x^p-y^q=1$ avec $p,q\geq 11$ premiers impairs et $x,y$ entiers non nuls. Posons $X={\rm
Ann}_{\Z[G]}([x-\zeta_p]\in H)$ comme pour le lemme 3. Soit $0\ne\theta=\sum_{\tau\in
G}n_\tau\tau\in X\cap (1-\iota)\Z[G]$, tel que $\|\theta\|=\sum_{\tau\in G}|n_\tau|\leq {3q\over p-1}$;
et $\alpha\in\Q(\zeta_p)^*$ tel que $(x-\zeta)^\theta=\alpha^q$; ainsi que $\sigma\in G$. Alors

$$|{\arg}(\sigma(\alpha))|> {\pi\over q},$$

avec l'argument tel que $-\pi<{\arg(z)}\leq \pi$.

}

{\bf Preuve}

Remarquons dŽjˆ que puisque $\theta=\sum_{\tau\in G}n_\tau\tau\in (1-\iota)\Z[G]$, alors,
$n_\tau=-n_{\iota\tau}$ pour tout $\tau\in G$, et donc $\sum_{\tau\in G}n_\tau=0$. Ainsi $|\alpha|=1$,
car, si on pose $\overline{\theta}=\sum_{\tau\in G}n_{\iota\tau}\tau$, on a
$$|(x-\zeta_p)^\theta|^2=(x-\zeta_p)^\theta\cdot\overline{(x-\zeta)^\theta}=
(x-\zeta_p)^\theta\cdot(x-\zeta_p)^{\overline{\theta}}=(x-\zeta_p)^0=1,$$

de mme, $|\tau(\alpha)|=1$, pour tout $\tau\in G$. De plus, 

$$(x-\zeta_p)^\theta=\prod_{\tau\in G}(x-\tau(\zeta_p))^{n_\tau}=\underbrace {x^{\sum_{\tau\in
G}n_\tau}}_{=1}\cdot \prod_{\tau\in G}(1-{\tau(\zeta_p)\over x})^{n_\tau}.$$

 On a
$\sigma(\alpha)^q=\prod_{\tau\in G}(1-{\tau(\zeta)\over x})^{n_\tau}$, avec $\zeta=\sigma(\zeta_p)$.
Donc, $\sum_{\tau\in G} n_\tau\Log(1-{\tau(\zeta\over x})=f(\sigma(\alpha)^q)$, o $f(\cdot)$ est une
fonction logarithmique, pas forcŽment la branche principale. D'autre part, pour tout $\tau\in G$, on a~: 

$$|\Log(1-{\tau(\zeta)\over
x})|=\left |\sum_{n=1}^\infty{\left ({\tau(\zeta)\over x }\right )^{n_\tau}\over n}\right |\leq
\sum_{n=1}^\infty\left ({1\over |x|}\right )^n={1\over |x|-1}.$$

Et alors~:

$$|\arg(\sigma(\alpha)^q)|\buildrel {\rm Lemme\ 5\ b)}\over\leq |f(\sigma(\alpha)^q)|\leq \sum_{\tau\in
G}|n_\tau |\cdot |\Log(1-{\tau(\zeta)\over x})|\leq {\|\theta\|\over |x|-1}.\eqno{(11)}$$

Supposons maintenant  par l'absurde que $|{\arg}(\sigma(\alpha))|\leq{\pi\over q}$. Le Lemme 5 c) nous
montre que 

$$|\arg(\sigma(\alpha))|\leq {\|\theta\|\over q(|x|-1)}.$$

\goodbreak
D'autre part,
$|\sigma(\alpha)-1|<|\arg(\sigma(\alpha))|$; on peut voir cela en sachant que pour tout $0<\beta<{\pi\over
q}$, on a $2\sin({\beta\over 2})<\beta$, ou en observant ce petit dessin~:




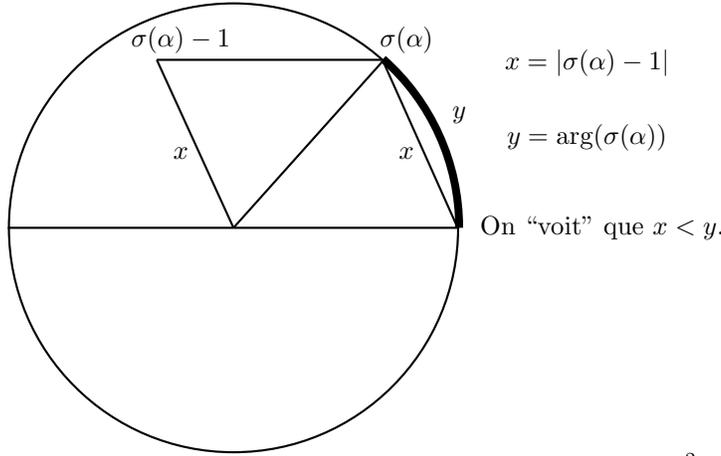
\input pstricks
\input pst-node

\vskip 3cm
\pscircle (5,0){3}
\psline(2,0)(8,0)
\psline(5,0)(7,2.235)
\psline(7.98,0)(6.98,2.235)
\psarc[linewidth=3pt](5,0){3}{0}{48.5} 
\psline(6.98,2.235)(3.98,2.235)
\psline(5,0)(3.98,2.235)
\rput (4.3,1){$x$}
\rput (4.3,2.5){$\sigma(\alpha)-1$}
\rput (7.3,2.5){$\sigma(\alpha)$}
\rput (7.3,1){$x$}
\rput (8,1.5){$y$}
\rput (9.7,2.2){$x=|\sigma(\alpha)-1|$}
\rput (9.7,1.2){$y={\rm arg}(\sigma(\alpha))$}
\rput (9.9,0){On Òvoit'' que  $x<y$.}
\vskip 2.5cm

Donc,
$$|N_{\Q(\zeta_p)/\Q}(\alpha-1)|=|N(\alpha-1)|<\left ({\|\theta\|\over q(|x|-1)}\right )^2\cdot
2^{p-3},\eqno{(i)}$$

car le terme $\left ({\|\theta\|\over q(|x|-1)}\right )^2$ vient de $\sigma$ et $\overline \sigma$, et le
$2^{p-3}$ vient du fait que $|\tau (\alpha)|=1$, et donc $|\tau(\alpha)-1|\leq 2$, pour tout $\tau\in G$. 

D'autre part, 

$$\alpha^q={\dst\prod_{n_\tau\geq 0}(x-\tau(\zeta_p))^{n_\tau}\over \dst\prod_{n_\tau\leq
0}(x-\tau(\zeta_p))^{|n_\tau|}}={\nu_1\over\nu_2}.$$

Comme, $n_\tau=-n_{\iota\tau}$, alors $\nu_1$ et $\nu_2$ sont conjuguŽs l'un de l'autre. Donc,
$N(\nu_2)^2=N(\nu_1)\cdot N(\nu_2)=N(\nu_1\cdot\nu_2)=N(\prod_{\tau\in G}(x-\tau(\zeta_p))^{|n_\tau|})\leq
(1+|x|)^{\|\theta\|(p-1)}$. Donc $|N(\nu_2)|\leq (1+|x|)^{\|\theta\|(p-1)\over 2}$. {\it A fortiori}, si
$\euf b^q$ est l'idŽal dŽnominateur de
$\alpha^q$, il contient l'idŽal engendrŽ par $\nu_2$. On a donc, $\N({\euf b})\leq |N(\nu_2)|^{1\over
q}\leq (1+|x|)^{\|\theta\|(p-1)\over 2q}$. En inversant, on obtient~:

$$(1+|x|)^{-\|\theta\|(p-1)\over 2q}\leq \N({\euf b})^{-1}\buildrel (*)\over \leq |N(\alpha-1)|\buildrel
(i)\over\leq \left ({\|\theta\|\over q(|x|-1)}\right )^2\cdot 2^{p-3}$$

L'inŽgalitŽ $(*)$ vient du fait que l'idŽal dŽnominateur de $\alpha-1$ est le mme que celui d'$\alpha$
(lemme 4) et puisque $\theta\ne 0$, alors, $\alpha\ne 1$ (lemme 3). En multipliant par $2 (|x|-1)^2$, et
puisque $|x|\geq 6$, $(1+|x|)^2\leq 2(|x|-1)^2$, on obtient~:

$$(1+|x|)^{2-{p-1\over 2q}\|\theta\|}\leq 2^{p-1}\left ({\|\theta\|\over q}\right )^2.$$

Mais, $\|\theta\|\leq {3q\over p-1}$, donc

$$q^{p-1\over 2}\buildrel {\rm Lemme\ 2}\over\leq (1+|x|)^{1\over 2}\leq 2^{p-1}\left ({3\over p-1}\right
)^2\leq 2^{p-1}=4^{p-1\over 2}.$$

Comme $q\geq 5$, c'est une contradiction, ce qui prouve que $|\arg(\sigma(\alpha))|> {\pi\over q}$ et
donc, la proposition.\qed

\bigskip
Maintenant, un petit lemme de combinatoire, trs facile~:
\bigskip
{\bf Lemme 7}

$$\left |\left \{(\lambda_1,\ldots ,\lambda_k)\in\N^k\mid \sum_{i=1}^k\lambda_i\leq s\right\}\right
|=\pmatrix{s+k\cr k\cr}$$

{\bf Preuve}

Posons $\Sigma_1=\left \{(\lambda_1,\ldots ,\lambda_k)\in\N^k\mid \sum_{i=1}^k\lambda_i\leq
s\right\}$, $\Sigma_2=\left \{(\lambda_1,\ldots ,\lambda_k)\in\N^k\mid
\lambda_i\leq 1,\ \sum_{i=1}^k\lambda_i\leq s+k\right\}$ et $\Sigma_3=\left \{A \mid
A\subset\{1,2,3,\ldots , s+k\}\ \hbox{et}\ |A|=k\right \}$. Eh bien, ces trois ensembles sont en
bijections~:
$(\lambda_1,\ldots ,\lambda_k)\mapsto (\lambda_1+1,\ldots ,\lambda_k+1) $ est une bijection de $\Sigma_1$
sur $\Sigma_2$; et $(\lambda_1,\ldots ,\lambda_k)\mapsto\{\sum_{i=1}^l\lambda_i \mid l=1,\ldots ,k\}$ est
une bijection de $\Sigma_2$ sur $\Sigma_3$ qui est Žvidemment de cardinal $\pmatrix{s+k\cr k\cr}$.\qed

\bigskip
Le lemme suivant n'est pas trs dur non plus~:
\bigskip
\goodbreak
{\bf Lemme 8}

{\sl Soit $k\geq 2$ et $s\geq 6$ tels que $s+2k\geq 13$. Alors on a

$$\pmatrix{s+k\cr k\cr}=\pmatrix{s+k\cr s\cr}>{4\over 3}(s+1) k^2+1.$$

}
{\bf Preuve}

On vŽrifie d'abord que le lemme est vrai pour les paires $(s,k)=(6,4),(7,3),(9,2)$. Il suffit alors de
voir que si le lemme est vrai pour $(s,k)$, alors il est vrai pour $(s+1,k)$ et pour $(s,k+1)$. Supposons
donc le lemme vrai pour $(s,k)$. On a

$$\eqalign{\pmatrix{s+k+1\cr s+1\cr}&={s+1+k\over s+1}\pmatrix{s+k\cr k\cr}\buildrel \rm par\ hyp.\over
>{s+1+k\over s+1}\cdot ({4\over 3}(s+1)k^2+1)\cr &={4\over 3}(s+1+k)k^2+{s+1+k\over s+1}>{4\over
3}(s+2)k^2+1.\cr}$$
Donc, le cas $(s+1,k)$ est rŽglŽ. On a aussi,

$$\pmatrix{s+k+1\cr s\cr}={s+1+k\over k+1}\pmatrix{s+k\cr s\cr}\buildrel \rm par\ hyp.\over
>{s+1+k\over k+1}\cdot ({4\over 3}(s+1)k^2+1)>{4\over 3}{(s+1)(s+k+1)\over k+1}\cdot k^2+1.$$

Il reste ˆ voir que $(s+k+1)k^2>(k+1)^3$. C'est Žquivalent ˆ $k^2(s-2)>3k+1$. Comme $s-2\geq 4$ et $k\geq
2$, on a $k^2(s-2)>4k^2>4k-1\leq 3k+2-1$.\qed

\bigskip

{\bf Proposition 9}

{\sl Soit $p,q\geq 11$ des nombres premiers tels que $q> 4p^2$, et $x$, $y$ des entiers non nuls tels que 
$x^p-y^q=1$, alors pour tout $\sigma\in G={\rm Gal}(\Q(\zeta_p)/\Q)$, il existe $0\ne\theta \in I$ tel
que $\|\theta\|\leq {3q\over p-1}$ et tel que si $\alpha\in\Q(\zeta_p)^*$ est l'ŽlŽment tel que
$(x-\zeta_p)^\theta=\alpha^q$, alors $|\arg(\sigma(\alpha))|\leq {\pi\over q}$.

}

{\bf Preuve}

Il suffit  de montrer qu'il y au moins $q+1$ ŽlŽments distincts $\theta\in I=(1-\iota)I_{st}$ tel que
$\|\theta\|\leq {3q\over 2(p-1)}$. En effet, ˆ chacun de ces $\theta$ correspond $\alpha\in\Q(\zeta_p)$ tel
que
$(x-\zeta_p)^\theta=\alpha^q$. Puisque $\theta\in I\subset(1-\iota)\Z[G]$ la realtion $(11)$
s'applique~:  $|\arg(\sigma(\alpha)^q)|\leq {\|\theta\|\over |x|-1}$. Par le Lemme 5 c), il existe
$k\in\Z$ tel que ${1-r\over 2}\leq k\leq {r-1\over 2}$ et

$$\left | \arg(\sigma(\alpha))-{2k\pi\over q}\right |\leq {\|\theta\|\over q(|x|-1)}.\eqno{(i)}$$

Par le principe des tiroirs il existe $\theta_1$, $\theta_2$ tels que $\theta_1\ne\theta_2$ qui
correspondent au mme $k$. Prenons $\theta=\theta_1-\theta_2$.
L'$\alpha$ correspondant ˆ $\theta$ vaut ${\alpha_1\over \alpha_2}=:\alpha$. Le fait qu'ils
correspondent au mme $k$ veut dire que $\arg(\sigma(\alpha_i))={\arg(\sigma(\alpha_i)^q)\over
q}+{k\cdot\pi i\over q}$, $i=1,2$. Donc $|\arg (\sigma(\alpha_1))-\arg (\sigma(\alpha_1))|=\left
|{\arg(\sigma(\alpha_1)^q)\over  q}-{\arg(\sigma(\alpha_1)^q)\over  q}\right |\leq {2\pi\over q}<\pi$.
Donc $\arg (\sigma(\alpha))=\arg (\sigma(\alpha_1))-\arg (\sigma(\alpha_1))$. Ainsi, puisque
$\|\theta\|<{3q\over 2(p-1})$ et que gr‰ce au Lemme~2, $|x|-1>q^{p-1}$, on trouve

$$\eqalign{|\arg(\sigma(\alpha))|&=|\arg (\sigma(\alpha_1))-\arg (\sigma(\alpha_1))|\leq |\arg
(\sigma(\alpha_1))-{2k\pi\over q}|+|\arg (\sigma(\alpha_2))-{2k\pi\over q}|\cr &\buildrel(i)\over\leq
{2\|\theta\|\over q(|x|-1)}\leq {3\over(p-1)(|x|-1)}\leq {3\over (p-1)q^{p-1}}<{\pi\over q}.\cr}$$

Il reste donc ˆ trouver ces $q+1$ $Ò\theta"$ distincts. 

Pour cela, on se souvient (Lemme 11 du Chapitre 5) que $I$ admet une $\Z$-base $e_1,\ldots , e_{p-1\over
2}$ telle que $e_1=\sum_{\tau\in G}\pm\tau$ et donc $\|e_i\|=p-1$, pour tout $i$. ConsidŽrons l'ensemble des
$\theta=\sum_{i=1}^{p-1\over 2}\lambda_ie_i$ avec $\lambda_i\in\N$ et $\sum_{i=1}^{p-1\over 2}\lambda_i\leq
s:=\left [{3 q\over 2(p-1)^2}\right ]$.  Pour un tel $\theta$, on a $\|\theta\|\leq s\cdot (p-1)\leq {3
q\over 2(p-1)}$. En vertu du Lemme 7, le nombre de ces $\theta$ est $\pmatrix{s+{p-1\over 2}\cr s\cr}$. En
ajoutant les opposŽs de ces $\theta$, on obtient $M:= 2\cdot \pmatrix{s+{p-1\over 2}\cr s\cr}-1$ (le $-1$
vient du $0$) $\theta\in I$ tels que $\|\theta\|\geq {3q\over 2 (p-1)}$. 

Il reste donc ˆ voir que $M\geq q+1$. Montrons d'abord que le Lemme 8 s'applique pour $s=\left [{3 q\over
2(p-1)^2}\right ]$ et $k={p-1\over 2}$. D'abord, $k\geq 2$, car $p\geq 5$. Comme, par hypothse,
$q>4p^2>4(p-1)^2$, on a 
$$s\geq {3q\over 2(p-1)^2}\geq {12(p-1)^2\over 2(p-1)^2}=6.$$
 Enfin, on a $s+2k\geq 13$. En effet, si
$s+2k<13$, alors $k=2,3,4,5$ ou $6$. Si $k=2$ ou $3$, $p=5$ ou $7$, trop petit; si $k=4,5$ ou $6$, $s<5$,
ce qui n'est pas le cas comme on vient de le voir. Donc le Lemme 8 s'applique et on trouve

$$M\geq 2 ({4\over 3}(s+1) k^2+1)-1={8\over 3}(s+1)k^2+1\geq {8\over 3}\cdot{3\over 2} {q\over (p-1)^2}\cdot
{(p-1)^2\over 4}+1=q+1.$$

Donc la proposition est prouvŽe.\qed

\bigskip

{\bf ThŽorme 10 (ThŽorme 3 de Mih$\breve{\bf a}$ilescu)}

{\sl Si $x$ et $y$ sont des entiers tels que $x^p-y^q=1$, avec $p,q\geq 11$ premiers, alors on $p< 4q^2$ et
$q<4p^2$.
}

{\bf Preuve}

Par symŽtrie, on suppose par l'absurde que $q>4 p^2$. Puisque $H$ est annihilŽ par $I$ (Lemme 2 Chapitre
7), donc, en particulier $x-\zeta_p$ est annihilŽ par $I$. Les propositions 6 et 9 se contredisent alors.
Le thŽorme est ainsi prouvŽ.\qed 

\vfill\eject

\long\def\art#1{{\parindent0pt\item{#1}}\hangindent=7mm\hangafter=-20}
\long\def\artart#1{{\parindent0pt\item{#1}}\hangindent=12mm\hangafter=-20}
\font\para=cmbx12 at 18pt
\font\indi=cmr10 at 7pt
\def\O{\hbox{$\cal O$}}
\def\U{\hbox{$\cal U$}}
\def\m{\hbox{\rs m\!}}
\def\dst{\displaystyle}
\font\doub=msbm10 at 10pt
\def\lra{\longrightarrow}
\def\qed{\hfill$\square$}
\def\gfP{\relax\ifmmode\bbP\else $\bbP$\fi}
\def\gP{{\euf P}}
\def\P{{\cal P}}
\def\QQ{{\cal Q}}
\def\Log{{\rm Log}}

\def\ggP{{\bf P}}
\newcount\chapnomb \chapnomb=1
\newcount\parnomb \parnomb=1
\pageno =46

\parindent0pt
\centerline{\para CHAPITRE 9}
\bigskip
{\para Quatrime thŽorme de Mih$\taille{18}\breve{\bf a}$ilescu~: 
\bigskip
\centerline{$\taille {15} p\equiv 1\pmod q$ ou $\taille {15} q\equiv 1\pmod p$ }}
\bigskip
Comme pour le Chapitre prŽcŽdent, le titre est clair. C'est le dernier morceau avant la fin. Ce sera
aussi le rŽsultat le plus dŽlicat ˆ montrer !! 

D'abord un petit lemme~:

{\bf Lemme 1}

{\sl Soit $R$ un anneau intgre de caractŽristique 0 et $q\in R$. On considre les sŽries formelles
$f=\sum_{k=0}^\infty {a_k\over k!} T^k$ et $g=\sum_{k=0}^\infty {b_k\over k!} T^k$, avec pour tout $k$,
$a_k,b_k\in R$ tels qu'il existe $a$ et $b\in R$ avec $a_k\equiv a^k\pmod{q R}$ et $b_k\equiv b^k\pmod{q
R}$. Alors
$f\cdot g=\sum_{k=0}^\infty{c_k\over k!} T^k$, avec $c_k\in R$ et $c_k\equiv (a+b)^k\pmod{q R}$.

}

{\bf Preuve}

${c_k\over k!}=\sum_{l+m=k}{a_l\over l!}{b_m\over m!}={1\over k!}\cdot
\underbrace{\sum_{l+m=k}\pmatrix{k\cr m\cr} a_l b_m}_{:=c_k\in R}$ et $c_k\equiv \sum_{l+m=k}\pmatrix{k\cr
m\cr} a^l b^m=(a+b)^k\pmod{q R}$.\qed

\bigskip

{\bf Lemme 2}

{\sl Si $f$ est une fonction rŽelle continuement dŽrivable $k+1$ fois sur l'intervalle $[0;t]$, alors il
existe $c\in[0;t]$ tel que

$$\left | f(t)-\sum_{i=0}^k{f^{(i)}(0)\over i!} t^i\right |={1\over (k+1)!} f^{(k+1)}(c)t^{k+1}.$$ 

}

{\bf Preuve}

C'est le thŽorme de Taylor [Ru,ThŽorme 5.15,pp. 110-111]\qed

\bigskip

Maintenant, nous allons fixer quelques notations.

Comme toujours, on va supposer que $p$ et $q$ sont deux nombres premiers impairs distincts et que $G={\rm
Gal}(\Q(\zeta_p)/\Q)$. Soit $\theta=\sum_{\tau\in G}n_\tau\tau\in \Z[G]$. On pose (formellement)

$$(1-\tau(\zeta_p) T)^{n_\tau\over q}=\sum_{k=0}^\infty \pmatrix{{n_\tau\over q}\cr
k\cr}(-\tau(\zeta_p)T)^k\in\Q(\zeta_p)[[T]]\,\hbox{et}\,(1-\zeta_p T)^{\theta\over q}=\prod_{\tau\in
G}(1-\tau(\zeta_p) T)^{n_\tau\over q}=:F(T)=F_\theta(T)$$

Remarquons que si $T=z\in\C$ est tel que $|z|<1$, alors le critre de D'Alembert nous dit que ces sŽries
convergent , car $\left |{{n_\tau\over q}-k\over k+1}\right |\cdot |z|\longrightarrow |z|<1$.

Enfin, si $\sigma\in G$, on note $F^\sigma$ la sŽrie obtenue en appliquant $\sigma$ ˆ chaque coefficient de
$F$.
\bigskip

{\bf Proposition 3}

{\sl Sous les mmes hypothses et notations que ci-dessus, on a

\art{a)}les coefficients de $F$ sont des entiers algŽbriques en dehors de $q$. C'est-ˆ-dire qu'ils sont
de la forme $a\over q^s$ o $a\in E_p=\Z[\zeta_p]$ est un entier cyclotomique.

\art{b)}De plus, $F(T)=\sum_{k=0}^{\infty}{a_k\over k!\cdot q^k} T^k$, avec, pour tout $k$,
$a_k\in\Z[\zeta_p]$ tel que $$a_k\equiv (-\sum_{\tau\in G} n_\tau\tau(\zeta_p))^k\pmod {q\Z[\zeta_p]}.$$

\art{c)}Si $\sigma\in G$ et $t\in\C$ est tel que $|t|<1$, alors la sŽrie $F^\sigma(t)$ converge. De plus,
si $0\leq n_\tau\leq q$ pour tout $\tau\in G$, alors pour tout $k\geq 0$, on a~:

$$\left | F^\sigma(t)-F_k^\sigma(t)\right |\leq \left |\pmatrix{-m\cr k+1\cr}\right |\cdot {|t|^{k+1}\over
(1-|t|)^{m+k+1}},$$

o $m={1\over q}\sum_{\tau\in G} n_\tau$ et $F_k^\sigma$ est la somme des termes de $F^\sigma$ de degrŽ
infŽrieur ou Žgal ˆ $k$.

}

{\bf Preuve}

\art{a)}Soit $\tau\in G$. On se souvient du Lemme 2 du Chapitre 4~: si $l\ne q$, alors
on a

$$v_l(k!)\leq v_l\left ({n_\tau\over q}\left ({n_\tau\over q}-1\right )\cdots  \left ({n_\tau\over
q}-(k-1)\right )\right).$$

Cela nous montre la partie a), car cela veut dire que chaque coefficient binomial
$\pmatrix{{n_\tau\over q}\cr k\cr}$ est un entier
divisŽ par une puissance de $q$. Les autres facteurs sont dans $\Z[\zeta_p]$.

\art{b)}Il est clair que $F(T)$ ˆ la forme voulue, car chaque coefficient binomial possde au plus un
puissance $k$-ime de $q$ au dŽnominateur et par le fait que si $s_1+\cdots +s_r=k$ alors $s_1!\cdots
s_r!$ divise
$k!$ (le quotient est un coefficient multinomial Òclassique"). Donc, $F(T)=\sum_{k=0}^{\infty}{a_k\over
k!\cdot q^k} T^k$, avec, pour tout $k$,
$a_k\in\Z[\zeta_p]$. Reste ˆ voir la congruence. On a

$$(1-\tau(\zeta_p) q T)^{n_\tau\over q}=\sum_{k=0}^\infty{n_\tau(n_\tau-q)\cdots (n_\tau-q(k-1))\over
k!}(-\tau(\zeta_p) T)^k=\sum_{k=0}^\infty{b_k\over k!} T^k,$$

avec 
$$b_k=(-\tau(\zeta_p)^k) n_\tau (n_\tau-q)\cdots (n_\tau-q(k-1))\equiv (-\tau(\zeta_p)n_\tau)^k:=b^k\pmod
{q\Z[\zeta_p]}.\eqno{(i)}$$

Ainsi, $F(qT)=\prod_{\tau\in G}(1-\tau(\zeta_p)qT)^{n_\tau\over q}=\sum_{k=0}^\infty {a_k\over k!} T^k$,
avec $a_k\equiv(-\sum_{\tau\in G}n_\tau \tau(\zeta_p))^k\pmod {q\Z[\zeta_p]}$, en vertu de $(i)$ et du Lemme
1, itŽrŽ $k-1$ fois.

\art{c)}La convergence a dŽjˆ ŽtŽ prouvŽe gr‰ce au critre de D'Alembert. Pour tout $\tau\in G$, on a $0\leq
{n_\tau\over q}$, donc,

$$\left | \pmatrix{{n_\tau\over q}\cr k}\right |={1\over k!}\cdot{n_\tau\over q}\cdot\left |{n_\tau\over
q}-1\right |\cdots\left |{n_\tau\over q}-(k-1)\right |\leq {1\over k!}\cdot \left |-{n_\tau\over q}\right
|\cdot \left |-{n_\tau\over q}-1\right |\cdots \left |-{n_\tau\over q}-(k-1)\right |=\left |
\pmatrix{-{n_\tau\over q}\cr k}\right |.$$

Ainsi, les coefficients, en modules, de la sŽrie $\sum_{k=0}^\infty  \pmatrix{{n_\tau\over q}\cr
k}(-\tau(\zeta_p)T)^k$ sont infŽrieurs ou Žgaux aux coefficients de la sŽrie $\sum_{k=0}^\infty 
\pmatrix{-{n_\tau\over q}\cr k}(-T)^k$ (qui sont tous positifs) $=(1-T)^{-n_\tau\over q}$. De mme, les
coefficients de $F(T)$ sont majorŽs par ceux de $\prod_{\tau\in G}(-T)^{-{n_\tau\over
q}}=(1-T)^{-m}=\sum_{k=0}^\infty\pmatrix{-m\cr k\cr}(1-T)^k=:S(T)$. On trouve alors, pour tout $t\in\C$ tel
que $|t|<1$~:

$$\eqalign{\left | F^\sigma(t)-F_k^\sigma(t)\right |&\leq \left ||(1-|t|)^{-m}-\sum_{j=0}^k\pmatrix{-m\cr
j\cr}(-|t|)^k\right |=|S(|t|)-S_k(|t|)|\cr
&\buildrel\rm Taylor\over\leq {1\over (k+1)!}|S^{(k+1)}(|t|)|\cdot|t|^{k+1}\leq \left |\pmatrix{-m\cr
k+1\cr}\right |\cdot {|t|^{k+1}\over (1-|t|)^{m+k+1}},\cr}$$

car $S(|t|)$ est croissante si $|t|<1$ et $S^{(k+1)}(|t|)=(-1)^{k+1}(-m)(-m-1)\cdots
(-m-k)(1-|t|)^{-m-k-1}$.\qed

\bigskip

{\bf Proposition 4}

{\sl Sous les mmes hypothses, avec en plus $\theta\in (1+\iota)\Z[G]$. Alors on a

\art{a)}$F_\theta=F\in \Q(\zeta_p^+)[[T]]$

\art{b)}Supposons que $t\in\Q$ est tel que $|t|<1$ et qu'il existe $\alpha\in \Q(\zeta_p)$ tel que
$(1-t\zeta_p)^\theta=\alpha^q$, alors $\alpha\in\Q(\zeta_p^+) $ et pour tout $\sigma\in G$, on a
$F^\sigma(t)=\sigma(\alpha)$.

}

{\bf Preuve}

\art{a)}Puisque $\theta=\sum_{\tau\in G}n_\tau\tau\in (1+\iota)\Z[G]$, alors $\iota\theta=\theta$, et donc
$n_\tau=n_{\iota\tau}$, pour tout $\tau\in G$. Donc, il est possible de partager dans le produit qui dŽfinit
$F$ un ŽlŽment et son conjuguŽ, donc les coefficients se trouvent dans $\R\cap \Q(\zeta_p)=\Q(\zeta_p^+)$.

\art{b)} Par un mme raisonnement qu'en a), on voit que $\beta:=(1-t\zeta_p)^\theta\in \R$. Ainsi,
$\overline{\alpha}^q=\overline{\alpha^q}=\overline{\beta}=\beta=\alpha^q$. Donc, $\overline{\alpha}=\alpha$,
car, dans $\Q(\zeta_p)$, les racines $q$-imes sont uniques. Cela implique
qu'$\alpha\in\Q(\zeta_p)\cap\R=\Q(\zeta_p^+)$. Pour la deuxime partie, remarquons dŽjˆ que si $\alpha\in
\Q(\zeta_p^+)$, on a $\sigma(\alpha)\in \Q(\zeta_p^+)$, car
$\overline{\sigma(\alpha))}=\iota\sigma(\alpha)\buildrel G\ \rm ab\hbox{\indi Ž}lien\over =
\sigma\iota(\alpha)=\sigma(\alpha)$. D'autre part, 
$\sigma(\alpha)^q=(1-t\sigma(\zeta_p))^\theta=F^\sigma(t)^q$, donc $\sigma(\alpha)=F^\sigma(t)$, car ils
sont rŽels et $q$ est impair.\qed

\bigskip\goodbreak

{\bf ThŽorme 5}

{\sl Si $x$ et $y$ sont des entiers tels que $x^p-y^q=1$, avec $p,q\geq 11$ premiers. Le
sous-$\F_q[G^+]$-module de $H^+$ engendrŽ par la classe de $(x-\zeta_p)^{1+\iota}$ est libre. Cela veut dire
que
${\rm Ann}_{\F_q[G^+]}([x-\zeta_p]^{1+\iota}\in H)=\{0\}$.
 On rappelle (cf. Chapitre 7) que 
$$H=\{\alpha\in\Q(\zeta_p)^*\mid v_{\euf r}(\alpha)\equiv 0\pmod{q}\hbox{
pour tout idŽal premier }{\euf r}\ne\pi E_p\}/{\Q(\zeta_p)^*}^q,$$ 
et que $H^+=\{[\alpha]\in H\mid
[\alpha]^\iota=[\alpha]\}$.

 }

{\bf Preuve}

Il faut donc dŽmontrer que si $\overline{\psi}\in \F_q[G^+]$ est tel que la classe $\left [
(x-\zeta_p)^{1+\iota}\right ]^{\overline{\psi}}=1_H$, alors $\overline{\psi}=0_{\F_q[G^+]}$. Soit donc un
tel $\overline{\psi}$. Ecrivons $\overline{\psi}=\sum_{\sigma\in G^+}\nu_\sigma\sigma$, avec $\nu_G\in\F_q$.
Soit $P\subset G$ un systme de reprŽsentant de $G^+=G/<\iota >$. Pour chaque $\sigma\in G^+$, on note
encore $\sigma$, l'ŽlŽment de $P$ qui le reprŽsente. Soit $\psi=\sum_{\sigma\in P}\nu_\sigma\sigma\in
\F_q[G]$. On a alors $1_H=\left [(x-\zeta_p)^{1+\iota}\right ]^{\overline{\psi}}=\left
[(x-\zeta_p)^{1+\iota}\right ]^{\psi}=\left [(x-\zeta_p)\right ]^{(1+\iota)\psi}$. Ainsi, par dŽfinition de
$H$, pour tout $\theta\in\Z[G]$ qui relve $\pm (1+\iota)\psi$, on a 
$(x-\zeta_p)^\theta\in{\Q(\zeta_p)^{*}}^q$. Posons pour tout $\sigma\in P$, $\nu_{\iota\sigma}=\nu_\sigma$.
On a donc $(1+\iota)\psi=\sum_{\sigma\in G}\nu_\sigma\sigma\in \F_q[G]$. Ecrivons $\theta=\sum_{\sigma\in G}
n_\sigma\sigma$ qui relve $\pm (1+\iota)\psi$. Evidemment, on peut choisir $\theta$ de manire que $0\leq
n_\sigma<q$. On va voir qu'on peut mme le choisir de sorte que $\|\theta\|\leq{p-1\over 2}\cdot q$. En
effet, si on pose $\theta'=q\sum_{\sigma\in G}\sigma-\theta=\sum_{\sigma_\in G} m_\sigma\sigma$, avec
$0<m_\sigma\leq q$. Puisque $\theta+\theta'=q\sum_{\sigma\in G}\sigma$, on a
$\|\theta\|+\|\theta'\|=(p-1)\cdot q$. On choisit donc celui de $\theta$ ou de $\theta'$ dont la norme
est infŽrieure ou Žgale ˆ ${p-1\over 2}\cdot q$; et si c'est $\theta'$, on remplace les coefficients $q$
Žventuels par $0$, on reste dans la mme classe et a rend la norme encore plus petite. Soit donc $\theta$
possŽdant cette propriŽtŽ. ConsidŽrons $\alpha\in\Q(\zeta_p)^*$ tel que $(x-\zeta_p)^\theta =\alpha^q$.
Notons $\pi$ pour $1-\zeta_p$. On a vu (preuve du ThŽorme 1, Chapitre 6) que ${x-\zeta_p\over
1-\zeta_p}\in E_p$ et que les conjuguŽs de $x-\zeta_p\over 1-\zeta_p$ Žtaient premiers entre eux. Donc,
pour tout $\sigma\in G$, on a $v_\pi(x-\sigma(\zeta_p))=1$. On peut alors remarquer que

$$\eqalign{\|\theta\|&=\sum_{\sigma\in G}n_\sigma=\sum_{\sigma\in G}n_\sigma \cdot
v_\pi(x-\sigma(\zeta_p))=v_\pi\bigg (\prod_{\sigma\in G}(x-\sigma(\zeta_p))^{n_\sigma}\bigg
)\cr &=v_\pi(x-\zeta_p)^\theta)=v_\pi(\alpha^q)=q\cdot v_\pi(\alpha)\equiv 0\pmod q.\cr}$$

Il existe donc $m$ entier avec $0\leq m\leq {p-1\over 2}$ tel que $\|\theta\|=mq$. En outre, puisque les
$n_\sigma$ et $n_{\iota\sigma}$  ont la mme rŽduction modulo $q$ (c'est $\nu_\sigma$), et qu'il sont entre
0 et $q$, on en dŽduit que $n_\sigma=n_{\iota\sigma}$.  Cela implique que $\theta=(1+\iota)\varphi$ o
$\varphi=\sum_{\sigma\in P} n_\sigma\sigma$ est un relevŽ de $\psi$. On en dŽduit que
$(x-\zeta_p)^\theta=\big ((x-\zeta_p)(x-\overline{\zeta_p})\big )^\varphi$ est rŽel ainsi que ses conjuguŽs.
Il en est de mme pour $\alpha$, car les racines $q$-ime, si elles existent sont uniques dans
$\Q(\zeta_p)$. Pour chaque $\tau\in G$, on a $(x-\tau(\zeta_p))^\theta=\tau(\alpha)^q$. Donc, par ce qu'on
vient de voir $(1-{1\over x}\tau(\zeta_p))^\theta=\left ({\tau(\alpha)\over x^m}\right )^q$. Par la
proposition 4 b) (appliquŽe ˆ ${1\over x}$ et ˆ ${\alpha\over x^m}$), on trouve~: 

$$\sigma(\alpha)=x^m F^{\sigma}\big ({1\over x}\big)\hbox{ pour tout }\sigma\in G.\eqno{(i)}$$

On affirme que 

$$q^{m+v_q(m!)}\cdot \left |\sigma(\alpha)-x^mF_m^\sigma\big ({1\over x}\big )\right |<1.\eqno{(ii)}$$

Prouvons l'affirmation $(ii)$. Remarquons d'abord que  $\left |\pmatrix{-m\cr m+1\cr}\right
|=\pmatrix{2m\cr m+1}\leq 4^{m}$, (la dernire inŽgalitŽ se voit en dŽveloppant $(1+1)^{2m}$ par le bin™me
de Newton). D'autre part, le Lemme 6 du Chapitre 4 montre que $v_q(m!)\leq {m\over q-1}$. On a alors 

$$\eqalign{q^{m+v_q(m!)}\cdot \left |\sigma(\alpha)-x^mF_m^\sigma\big ({1\over x}\big )\right
|&=q^{m+v_q(m!)}\cdot |x|^m \left |F^\sigma\big ({1\over x}\big )-F_m^\sigma\big ({1\over x}\big )\right
|\cr &\buildrel \rm Lemme\ 3 c)\over \leq q^{m+v_q(m!)}\cdot{1\over |x|}\cdot \left |\pmatrix{-m\cr
m+1\cr}\right |\cdot \bigg(1-{1\over |x|}\bigg)^{-2m-1}\cr
&\leq q^{m+{m\over q-1}+m\cdot{\log(4)\over\log(q)}}\cdot {1\over |x|}\cdot\bigg(1-{1\over
|x|}\bigg)^{-2m-1}\cr
&\buildrel (*)\over\leq {1\over |x|}\cdot q^{{p-1\over 2}\cdot(1+{1\over q-1}+{\log(4)\over
\log(q)})}\cdot\big (1-{1\over |x|}\big )^{-p}\cr
&\buildrel (**)\over \leq  q^{{p-1\over 2}\cdot(-1+{1\over q-1}+{\log(4)\over
\log(q)})}\cdot\big (1-{1\over q^{p-1}}\big )^{-p}. \cr}$$

L'inŽgalitŽ $(*)$ venant du fait que $m\leq {p-1\over 2}$ et l'inŽgalitŽ $(**)$ venant du fait que $|x|\geq
q^{p-1}$ (Lemme 2 du Chapitre 8). Pour voir que le dernier terme de cette sŽrie d'inŽgalitŽs est infŽrieur
ˆ 1, il suffit de vŽrifier que son $\log_q$ est infŽrieur ˆ 0. Ce $\log_q$ vaut

$${p-1\over 2}\cdot \big (-1+{1\over q-1}+{\log(4)\over\log(q)}\big)-p{\log(1-{1\over
q^{p-1}})\over\log(q)}.$$

Or, $-\log(1-{1\over q^{p-1}})=\log(q^{p-1})-\log(q^{p-1}-1)\buildrel (+)\over\leq{1\over q^{p-1}-1}\leq
{1\over q^2-1}$. L'inŽgalitŽ $(+)$ vient du fait que $\log(x+1)-\log(x)\leq {1\over x}$ qui est une
consŽquence immŽdiate du thŽorme des accroissements finis. D'autre part, on a supposŽ que $q\geq 7$.
Donc, le
$\log_q$ cherchŽ est infŽrieur ˆ

$$\leq {p-1\over 2}\big (-1+{1\over 6}+{\log(4)\over \log(7)}\big )+{p\over 48\cdot\log(7)}\leq
-0.487\cdot p+0.061,$$

qui est nŽgatif, ds que $p\geq 2$. L'affirmation $(ii)$ est donc prouvŽe.

D'autre part, les coefficients du polyn™me $q^{m+v_q(m!)}\cdot F_m^\sigma(T)$ sont des ŽlŽments de
$E_p$, car les coefficients de degrŽ $k$ sont de la forme ${a_k\over k! q^k}$ avec $a_k\in E_p$
(Proposition 3); la puissance de $q$ qui divise le dŽnominateur est au plus $k+v_q(k!)\leq m+ v_q(m!)$ si
$0\leq k\leq m$. 

On a aussi que $\alpha\in E_p$, car $\alpha\in \Q(\zeta_p)$ et parce que l'Žquation
$(x-\zeta_p)^\theta=\alpha^q$ est une Žquation entire, les coefficients de $\theta$ Žtant positifs
(cette intŽgralitŽ vient du fait que si un ŽlŽment satisfait une Žquation du type $x^n+a_{n-1}x^{n-1}+\cdots
+a_1+a_0=0$ avec $a_i$ entiers sur $\Z$ pour tout $i$, alors cet ŽlŽment est aussi entiers $\Z$, voir
[Nar, Theorem 1.2]).

Cela implique que $\gamma:=q^{m+v_q(m!)}\cdot (\alpha-x^m F_m({1\over x}))\in E_p$, donc la norme de
$\gamma$ appartient ˆ $\N$. Mais, l'inŽgalitŽ
$(i)$ montre alors que cette norme est infŽrieure ˆ 1. Donc elle et nulle, et donc $\gamma$ aussi. On a 

$$q^{m+v_q(m!)}\cdot \alpha=\sum_{k=0}^m q^{m+v_q(m!)}{a_k\over q^k\cdot k!} x^{m-k}.$$

En raisonnant modulo $q$, on dŽduit que $a_m\equiv 0\pmod {q E_p}$. Or, la Proposition 3 b) nous apprend
que $a_m\equiv \big (-\sum_{\sigma\in G} n_\sigma\sigma(\zeta_p)\big )^m\pmod {q E_p}$. L'anneau quotient
$E_p/q E_p$ n'a pas d'ŽlŽment nilpotent (thŽorme chinois et $q$ ne ramifie pas). Donc, $\sum_{\sigma\in
G}n_\sigma\sigma(\zeta_p)\equiv 0\pmod {q E_p}$. Comme les $\sigma(\zeta_p)$ forment une $\Z$-base de
$E_p$ ($1,\zeta_p,\ldots ,\zeta_p^{p-2}$ en est une, et en multipliant par $\zeta_p$, s'en est encore
une). On en tire que $n_\sigma\equiv 0\pmod q$, et donc $n_\sigma=0$ pour tout $\sigma\in G$. Donc
$\theta=0$ et {\it a fortiori} $\psi$ et $\overline{\psi}$ aussi.\qed
\bigskip

C'est une jolie preuve, vous ne trouvez pas ? Maintenant, nous allons prouver un thŽorme qui permettra de
terminer la conjecture de Catalan Ce sera le ThŽorme 12. Tout d'abord on va prŽsenter quelques nouveaux
objets et fixer quelques notations~:

\bigskip

{\bf DŽfinition}

On ne rappellera pas les dŽfinitions de $H$ et de $H^+$ on les a revues au ThŽorme 5.  On pose

$E=\{u(1-\zeta_p)^k\mid u\in U_p, k\in \Z\}$ avec $U_p=E_p^*=\Z[\zeta_p]^*$. Il est Žvident que
$E=\Z[\zeta_p,{1\over p}]^*$. On pose aussi

$$H'=\bigg\{[\alpha]\in H\mid \alpha=\beta^q+q^2\gamma, \ \beta,\gamma\in\Z[\zeta_p,{1\over p}]\ \hbox{ et
}\beta\ \hbox{ inversible }\pmod{q^2\Z[\zeta_p,{1\over p}]}\bigg\}$$

et $E'=\{u\in E\mid [u]\in H'\}$. On pose encore $C$, qu'on appelle {\it $p$-unitŽs cyclotomiques},
l'ensemble $\{u\in E\mid u={1-\zeta_p^i\over 1-\zeta_p^j}\omega(1-\zeta_p)^k,\ \omega \hbox{ est une racine
de l'unitŽ et } i,j\in\N\ \hbox{ et } k\in\Z\}$ et $C'=C\cap E'$. On rappelle que ${\cal CL}={\cal
CL}_{\Q(\zeta_p)}$, ${\cal CL}^{pl}={\cal CL}_{\Q(\zeta_p^+)}$ sont les classes d'idŽaux de $\Q(\zeta_p)$
et de $\Q(\zeta_p^+)$. On rappelle (Chapitre 7, dŽfinition p. 33) que ${\cal CL}^{pl}\subset {\cal
CL}^+$. Enfin, ${\cal CL}[q]=\{ [{\euf a}]\in {\cal CL}\mid {\euf a}^q$ est principal $\}$.

\bigskip

Maintenant que les acteurs sont prŽsentŽs, on va faire une sŽrie de lemmes.
\bigskip
{\bf Lemme 6}

{\sl Sous les mmes hypothses que plus haut, on a
$${\cal CL}^{pl}[q]={\cal CL}[q]^+.$$

}

{\bf Preuve}

Le fait que ${\cal CL}^{pl}[q]\subset {\cal CL}[q]^+$ vient de ${\cal CL}^{pl}\subset {\cal CL}^+$.
Prouvons l'autre inclusion. On pourrait dire que puisque ${\cal CL}[q]$ est un $\F_q[G]$-module et
que $2$ est inversible modulo $q$, alors ${\cal CL}[q]^+={\cal CL}[q]^{1+\iota\over 2}={\cal
CL}[q]^{1+\iota}$ et le tour est jouŽ ! C'est juste, mais un peu rapide. Soyons plus terre ˆ terre~:  soit
$[{\euf a}]\in {\cal CL}[q]^+$. Cela veut dire qu'il existe $\alpha,\beta\in \Q(\zeta_p)^*$ tel que
${\euf a}\overline{\euf a}^{-1}=\alpha E_p$ et ${\euf a}^q=\beta E_p$. Puisque $q\ne 2$, il existe $m$ et
$k\in\Z$ tels que $2m-1=kq$. On a, si $[{\euf a}]$ dŽsigne la classe de $\euf a$ dans $\cal CL$,

$$[{\euf a}]=[{\euf a}\cdot\beta^k\cdot\alpha^{-m}]=[{\euf a}\cdot{\euf a}^{kq}\cdot\alpha^{-m}]=[{\euf
a}^{kq+1}\cdot\alpha^{-m}]=[{\euf a}^{2m}\cdot\alpha^{-m}]=[{\euf a}^m\cdot({\euf
a}^m\cdot\alpha^{-m})]=[{\euf a}^m\cdot\overline{\euf a}^m]=[{\euf b}]$$

o ${\euf b}=({\euf a}\overline{\euf a})^m\in {\cal CL}^{pl}[q]$.\qed

\bigskip

On va montrer que $E/E^q$ est un $\F_q[G^+]$-module libre de rang 1. Pour montrer cela, on doit utiliser
le fait que $\F_q[G^+]$ est un anneau {\it semi-simple}, c'est-ˆ-dire dans notre cas que c'est un produit
cartŽsien de corps. On a besoin d'un rŽsultat sur les modules sur un anneau semi-simple. Nous montrerons
ce rŽsultat dans l'appendice 1. Enonons tout de mme ce rŽsultat, mais tout d'abord quelques
dŽfinitions.
\bigskip

{\bf DŽfinitions}

Soit $R$ un anneau commutatif. Un idŽal $\euf b$ de $R$ est dit {\it radical} si $\alpha^n\in {\euf
b}$ implique $\alpha\in{\euf b}$, ou ce qui est Žquivalent, $R/{\euf b}$ n'a pas d'ŽlŽment nilpotent. Si
$M$ est un $R$-module et $S\subset M$, on note ${\rm Ann}_R(M)$ ou ${\rm Ann}(M)$ l'ensemble $\{\alpha\in
R\mid \alpha x=0\hbox{ pour tout } x\in S\}$.
\bigskip

{\bf Lemme 7}

{\sl Soit $R$ un anneau semi-simple commutatif. 

\art{a)}Si $M$ est un $R$-module de type fini, et $\euf b$ un idŽal tel que ${\rm Ann}(M)+{\euf b}$ est
radical. Si $\psi$ est la projection de $R$ sur $R/{\euf b}$. Alors on a $\psi({\rm Ann(M)})={\rm
Ann}_{R/{\euf b}}(M/{\euf b}M)$ (ce rŽsultat est vrai mme si $R$ n'est pas semi-simple).

\art{b)}Soit $M$ un $R$-module. Alors $M$ contient un sous-module isomorphe ˆ $R/{\rm Ann}(M)$, avec
ŽgalitŽ si et seulement si $M$ est cyclique. De plus, dans ce cas, tout sous-module
est aussi cyclique.

\art{c)}Si $R$ est fini et $M$ est un $R$-module tel que le cardinal $|M|=|R/{\rm Ann}(M)|$, alors $M$ est
isomorphe $R/{\rm Ann}(M)$.

}

{\bf Preuve} 

On montre ces rŽsultats dans l'Appendice 1~: Proposition 3 pour a) et Proposition 6 pour b). La partie
c) est une consŽquence de b).\qed

\bigskip

On a tout fait pour s'en passer, mais il faudra quand mme Žnoncer le cŽlbre
\bigskip
{\bf ThŽorme de Dirichlet}

{\sl Soit $K$ un corps de nombre de signature $(r_1,r_2)$ (c'est-ˆ-dire qu'il y a $r_1$ plongements rŽels de
$K$ et $2r_2$ plongements complexes, donc $[K:\Q]=r_1+2r_2$). Soit aussi $U(K)$ les unitŽs de l'anneau des
entiers de $K$ sur $\Q$. Alors $U(K)$ est isomorphe (en tant que $\Z$-module ou en tant que groupe) ˆ
$W\times \Z^{r_1+r_2-1}$ o $W$ est le groupe des racines de l'unitŽs de $K$.

}

{\bf Preuve}

C'est un thŽorme trs classique qui est prouvŽ dans [Sam pp. 72-75, ThŽorme 1 du Chapitre 4].  On
montre (entre autre) la chose suivante~: soit $\sigma_1,\ldots
,\sigma_{r_1},\sigma_{r_1+1},\sigma_{r_1+r_2}$ les plongements de $K$ dans $\C$ (on ne compte qu'une
fois les paires de conjuguŽs). Posons $r=r_1+r_2$. Soit $x\in K$. L'application $x\mapsto
l(x):=(l_1(x),\ldots , l_r(x))\in \R^r$ avec $l_i(x)=\delta_i\log|\sigma_i(x)|$, o
$\delta_i=1$ si $i=1,\ldots ,r_1$ et $\delta_i=2$ si $i=r_1,\ldots ,r$ s'appelle le {\it plongement
logarithmique de $K$}. La preuve du thŽorme de Dirichlet consiste ˆ montrer que $l(U(K))$
engendre l'hyperplan de $\R^r$ d'Žquation $\sum_{i=1}^r x_i=0$.\qed

\bigskip\vskip 1cm\goodbreak\vskip -1cm

{\bf Lemme 8}

{\sl Si $p$ et $q$ sont des nombre premiers impairs distincts tel que $p\ne 1\pmod q$ et $G^+$ est le
groupe de Galois de l'extension $\Q(\zeta_p^+)/\Q$, alors $E/E^q$ est un
$\F_q[G^+]$-module libre de rang 1.

}

{\bf Preuve}

D'abord, $\F_q[G^+]$ est un anneau semi-simple~: $G^+$ est un groupe cyclique d'ordre ${p-1\over 2}$ et
donc $\F_q[G^+]$ est isomorphe ˆ $\F_q[X]/(X^{p-1\over 2}-1)$ et le polyn™me $X^{p-1\over 2}-1$ n'a que
des racines simples, car sa dŽrivŽe qui vaut ${p-1\over 2}X^{p-3\over 2}$ ne s'annule qu'en 0 ($q$ ne
divisant pas ${p-1\over 2}$). Donc la dŽcomposition de $X^{p-1\over 2}-1$ en polyn™me irrŽductibles dans
$\F_q[X]$ n'est fait que de polyn™mes ˆ la puissance un. 

Pour montrer le Lemme, il suffit de montrer (en vertu du Lemme 7) que $$|E/E^q|=q^{p-1\over
2}=|\F_q[G^+]|\hbox{ et que }{\rm Ann}_{\F_q[G^+]}(E/E^q)=\{0\}.$$

Montrons la premire affirmation. L'application $(u\cdot \pi^n)\longmapsto (u,n)$ montre que $E$ est
isomorphe (en tant que groupe) ˆ $U\times\Z$, o $U$ est l'ensemble des unitŽs de $E_p=\Z[\zeta_p]$. Donc
$E/E^q$ est isomorphe $U/U^q\times \Z/q\Z$. Remarquons que $U/U^q$ est un sous $\F_q[G^+]$-module de
$E/E^q$. Le thŽorme de Dirichlet nous dit que
$U$ est isomorphe ˆ $\mu_{2p}\times \Z^{p-3\over 2}$, $\mu_{2p}$ Žtant le groupe des racine $2p$-ime de
l'unitŽ. On a aussi $U^q\simeq \mu_{2p}^q\times (q\Z)^{p-3\over 2}=\mu_{2p}\times (q\Z)^{p-3\over 2}$; la
dernire ŽgalitŽ vient du fait que toute racine $2p$-ime de l'unitŽ est une puissance $q$-ime d'un autre
racine $2p$-ime de l'unitŽ, car $(2p,q)=1$ (on a dŽjˆ utilisŽ ce fait au ThŽorme 1 du Chapitre 6). Donc
le quotient $U/U^q$ est isomorphe ˆ $\{1\}\times (\Z/q\Z)^{p-3\over 2}$, cela prouve la premire
affirmation.

Pour la seconde, on remarque d'abord que $U/U^q$ est isomorphe ˆ $U^+/{U^+}^q$, o $U^+$ est le groupe
des unitŽ de $\Z[\zeta_p^+]$. En effet, l'application composŽe $U^+\hookrightarrow U\rightarrow U/U^q$
est surjective car, par le Lemme de Kummer (vu au Chapitre 6), on sait que tout $u\in U$ peut s'Žcrire
$u=u_0\cdot w$ o $u_0\in U^+$ et $w$ une racine $2p$-ime de l'unitŽ qui appartient, comme on vient de
la voir, ˆ $U^q$. D'autre part, $U^+\cap U^q={U^+}^q$; en effet, $\supset$ est clair, et si $u\in U^+$
s'Žcrit $u=v^q$ avec $v\in U$, alors $u=\overline{v}^q=v^q$, donc $v=\overline{v}\in U^+$ (car dans
$E_p$, si une racine $q$-ime existe, elle est unique). On en dŽduit que 

$$E/E^q\simeq U^+/{U^+}^q\times \Z/q\Z\simeq U^+/\{\pm 1\}/\big (U^+/\{\pm 1\})^q\times \Z/q\Z.$$

Regardons ${\rm Ann}_{\Z[G^+]}(U^+/\{\pm 1\})$~: 

il est clair que $\sum_{\sigma\in G^+}a_\sigma\sigma$
annule $U^+/\{\pm 1\}$ si et seulement si $\prod_{\sigma\in G^+}\sigma(u)^{a_\sigma}=\pm 1$ pour tout $u\in
U^+$. C'est Žquivalent ˆ

$$\sum_{\sigma\in G^+}a_\sigma\log|\sigma(u)|=0\quad\hbox{pour tout }u\in U^+.\eqno{(i)}$$

Or, la preuve du thŽorme de Dirichlet nous apprend que $\sum_{\sigma\in G^+}\log |\sigma
(u)|=0$ pour tout $u\in U^+$ et que si $(x_\sigma)_{\sigma\in G^+}\in \R^{G^+}$ est tel que
$\sum_{\sigma\in G^+}x_\sigma=0$, alors il existe $n\in\N$, $u_1,\ldots ,u_n\in U^+$ et $\lambda_1,\ldots
,\lambda_n\in\R$ tel que $x_\sigma=\sum_{i=1}^n\lambda_i\log|\sigma(u_i)|$. Fixons $\tau,\tau'\in G$,
$\tau\ne \tau'$ et posons $x_\tau=1$, $x_\tau'=-1$ et $x_\sigma=0$ si $\sigma\ne \tau,\tau'$. ConsidŽrons
le $n$, les $u_i$ et les $\lambda_i$ associŽs ˆ ce $(x_\sigma)_{\sigma\in G^+}$-lˆ. On a

$$a_\tau-a_{\tau'}=\sum_{\sigma\in G^+} a_\sigma \cdot x_\sigma=\sum_{\sigma\in
G^+}a_\sigma\cdot\sum_{i=1}^n\lambda_i\log|\sigma(u_i)|=\sum_{i=1}^n\lambda_i\cdot\sum_{\sigma\in
G^+}a_\sigma\log|\sigma(u_i)|\buildrel(i)\over =0.$$

Cela veut dire que tout les $a_\sigma$ sont Žgaux ˆ, disons, $a$. Ainsi, $\sum_{\sigma\in
G^+}a_\sigma\sigma=a\cdot \sum_{\sigma\in G^+}\sigma$. Ce qui veut dire que 

$${\rm Ann}_{\Z[G^+]}(U^+/\{\pm 1\})=\big (\sum_{\sigma\in G^+}\sigma\big )\Z[G^+].$$

Quotientons par l'idŽal engendrŽ par $q$. Pour pouvoir appliquer la partie a) du Lemme 7, il faut vŽrifier
que $I:={\rm Ann}_{\Z[G^+]}(U^+/\{\pm 1\})+q\Z[G^+]$ est radical. En effet,
$\Z[G^+]/I=\F_q[G^+]/(\sum_{\sigma\in G^+}\sigma )\F_q[G^+]\simeq \F_q[X]/(X^{p-3\over
2}+\cdots+X+1)$. Donc, $$\F_q[G^+]\simeq \F_q[X]/(X^{p-1\over 2}-1)\simeq \F_q[X]/(X-1)\times
\F_q[X]/(X^{p-3\over 2}+\cdots+X+1)\simeq \F_q\times \Z[G^+]/I.$$  Donc, $\Z[G^+]/I$ est semi-simple, donc
sans nilpotent. On peut donc appliquer le Lemme 7, et on trouve 

$${\rm Ann}_{\F_q[G^+]}(U^+/\{\pm 1\})/\big (U^+/\{\pm 1\}\big )^q=\big (\sum_{\sigma\in G^+}\sigma\big
)\F_q[G^+].\eqno{(12)}$$

De ceci, on dŽduit que ${\rm Ann}_{\F_q[G^+]}(E/E^q)\subset (\sum_{\sigma\in G^+}\sigma\big
)\F_q[G^+]$. Posons $\alpha:=\sum_{\sigma\in G^+}a_\sigma\sigma$ et  $s:=\sum_{\sigma\in
G^+}\sigma$ \goodbreak avec $\alpha\in
\Z[G^+]$ tel que la classe $\alpha\cdot s$ soit dans ${\rm Ann}_{\F_q[G^+]}(E/E^q)$. En particulier,
$\alpha\cdot s$ annule la classe de~$\pi$, i.e.
$\pi^{\alpha\cdot s}=u\cdot\pi^{\sum_{\sigma\in G^+}a_\sigma\cdot {p-1\over 2}}\in E^q$, pour une certaine
unitŽ $u$. Cela implique que $\sum_{\sigma\in G^+}a_\sigma\cdot {p-1\over 2}\equiv 0\pmod q$. Ainsi
$k:=\sum_{\sigma\in G^+}a_\sigma\equiv 0\pmod q$, car $p\not\equiv 1\pmod q$. Ainsi, 

$$\alpha\cdot s=\sum_{\matrix{\scriptstyle\sigma\in G^+\cr\scriptstyle \tau\in
G^+}}a_\sigma\cdot\sigma\tau\buildrel\mu=\sigma\tau\over =\sum_{\mu\in G^+}\cdot\big (\sum_{\tau\in G^+}
a_{\mu\tau^{-1}}\big )\cdot\mu=k\cdot s\in q\Z\cdot s\subset q\Z[G^+].$$

Cela prouve que ${\rm Ann}_{\F_q[G^+]}(E/E^q)=\{0\}$, et donc le lemme.\qed

\bigskip
On reprendra ce lemme lors de l'appendice 2 qui sera consacrŽ au thŽorme de Thaine (que nous Žnoncerons
bient™t.
\bigskip
{\bf Lemme 9}

{\sl On se souvient que $E'=\{u\in E\mid [u]\in H'\}$ (pour les autres dŽfinitions, voir pp 46 et 48). On a
alors

$$E'=\{u\in E\mid u=\beta^q+q^2\gamma, \beta,\gamma\in\Z[\zeta_p,{1\over p}]\}.$$

}
\goodbreak
{\bf Preuve}

La partie $\supset$ est triviale, il suffit de montrer que $\beta$ est inversible modulo
$q^2\Z[\zeta_p,{1\over p}]$, mais c'est Žvident, car $E=\Z[\zeta_p,{1\over p}]^*$.

Montrons l'autre inclusion. Soit $u\in E$ tel que $[u]\in H'$. Alors il existe $\alpha\in \Q(\zeta_p)^*$
tel que $\alpha^q u=\beta^q+q^2\gamma$, avec $\beta,\gamma\in\Z[\zeta_p,{1\over p}]$ et $\beta$
inversible modulo $q^2\Z[\zeta_p,{1\over p}]$. En regardant cette ŽgalitŽ, on voit que pour tout idŽal
$\gP$ de $E_p$, diffŽrent de celui engendrŽ par $\pi$, on a $v_\gP(\alpha)\geq 0$. Donc $\alpha\in
\Z[\zeta_p,{1\over p}]$. Regardons modulo $q^2\Z[\zeta_p,{1\over p}]$. On a
$\overline{\alpha}^q\cdot\overline{u}=\overline{\beta}^q$. Donc, $\overline{\alpha}$ est inversible et
$\overline{u}=(\overline{\alpha}^{-1}\overline{\beta})^q$. Si $\beta_0\in\Z[\zeta_p,{1\over p}]$
reprŽsente $\overline{\alpha}^{-1}\overline{\beta}$, on a $u=\beta_0^q+q^2\gamma_0$, pour un 
$\gamma_0\in\Z[\zeta_p,{1\over p}]$.\qed
\bigskip

{\bf Lemme 10}

{\sl Soit $p$ et $q$ sont des nombres premiers impairs distincts, $C$ les $p$-unitŽs cyclotomiques de
$\Q(\zeta_p)$, et $C'=C\cap E'$. Si $C=C'$, alors $p<q$.

}

{\bf Preuve}

Soit $\zeta$ une racine primitive $p$-ime de l'unitŽ. Alors $1+\zeta^q={1-\zeta^{2q}\over 1-\zeta^q}\in
C$. Par hypothse, on a alors $1+\zeta^q\in C'$. 

D'autre part, on a $\Z[\zeta_p,{1\over p}]/(q^2)\simeq
\Z[\zeta_p]/(q^2)$. En effet, il suffit de prouver que l'application $x\mapsto x\pmod {q^2}$ est un
homomorphisme surjectif de $\Z[\zeta_p]$ sur $\Z[\zeta_p,{1\over p}]/(q^2)$; c'est-ˆ-dire, si ${y\over
p^n}\in \Z[\zeta_p,{1\over p}]$ avec $y\in\Z[\zeta_p]$, alors on doit voir qu'il existe $x\in\Z[\zeta_p]$
et $\gamma\in \Z[\zeta_p,{1\over p}]$ tels que $x+q^2\gamma={y\over p^n}$; c'est Žvident~: par le thŽorme
de Bezout et puisque $(q^2,p^n)=1$, il existe $\alpha,\beta\in\Z$ tel que $\alpha p^n+\beta q^2=1$. Donc
$y\alpha +q^2{\beta\over p^n} ={y\over p^n}$.

Revenons ˆ notre $1+\zeta^q\in C'$. Gr‰ce ˆ l'isomorphisme prouvŽ ci-dessus, il existe $\beta,\gamma\in
\Z[\zeta_p]$ tels que $1+\zeta^q=\beta^q+q^2\gamma$; c'est-ˆ-dire $1+\zeta^q\equiv
\beta^q\pmod{q^2\Z[\zeta_p]}$. D'autre part, puisque $q$ divise les coefficients binomiaux $\pmatrix{q\cr
j\cr}$ si $j=1,\ldots ,q-1$, on a $(1+\zeta^q)\equiv (1+\zeta)^q\pmod {q\Z[\zeta_p]}$. Donc, on a
$(1+\zeta)^q\equiv \beta^q\pmod{q\Z[\zeta_p]}$. Soit ${\cal Q}$ un idŽal premier divisant
$q\Z[\zeta_p]$. On a  alors $(1+\zeta)^q\equiv \beta^q\pmod {\cal Q}$. Par le Lemme 1 du Chapitre
6, on aussi $(1+\zeta)^q\equiv \beta^q\pmod {{\cal Q}^2}$, donc par le thŽorme chinois, et puisque
$q$ ne ramifie pas, $(1+\zeta)^q\equiv \beta^q\pmod {{q\Z[\zeta_p]}^2}$. Finalement, on trouve
$(1+\zeta)^q\equiv 1+\zeta^q\pmod{{q\Z[\zeta_p]}^2}$, ou encore
$${(1+\zeta)^q-1-\zeta^q\over q\cdot\zeta}\in q\Z[\zeta_p].$$ 
Soit $F(X)={(1+X)^q-1-X^q\over q\cdot X}$. On a $F(X)\in \Z[X]$, est de degrŽ $q-2$ et $F(\zeta)\in
q\Z[\zeta_p]$. Soit 
${\cal Q}$ un idŽal premier divisant $q\Z[\zeta_p]$. Passant au quotient $\Z[\zeta_p]/{\cal Q}$, on a
$\overline{F}(\overline{\zeta})=0$, avec $\overline{F}\in\F_q[X]$. Comme les racine $p$-ime de l'unitŽ
sont distinctes modulo $\cal Q$ (Lemme IMP, Chapitre 5). Donc, $\overline{F}$ a, dans $\Z[\zeta_p]/{\cal
Q}$, au moins $p-1$ racines distinctes. Donc $p-1\leq q-2$, ou encore $p<q$.\qed

\goodbreak

\bigskip

{\bf Lemme 11 (ThŽorme de Thaine)}

{\sl Tout annulateur (dans $\F_q[G^+]$) de $E/CE^q$ annule aussi ${\cal CL}^{pl}[q]$.

}

{\bf Preuve}

Ce thŽorme n'a l'air de rien, mais il est trs long et il sera prouvŽ dans l'appendice 2. Il est aussi
prouvŽ dans [Was, ¤15.2, pp. 334-341]\qed
\bigskip

Maintenant nous sommes en mesure d'Žnoncer le thŽorme central de ce chapitre.

\bigskip

{\bf ThŽorme 12}

{\sl Supposons que $p>q$ et $p\not\equiv 1\pmod q$. Alors le $\F_q[G^+]$-module $H^+\cap H'$ a un
annulateur non nul.

}

{\bf Preuve}

Gr‰ce aux Lemmes 7 et 8, on remarque que tout sous-module de $E/E^q$ est isomorphe ˆ $\F_q[G^+]/{\rm
Ann}(M)\simeq\F_q[X]/(f)$ o $f$ est un diviseur unitaire de $X^{p-1\over 2}-1$, et donc, ${\rm
dim}_{\F_q}(M)={\rm deg}(f)$. On se souvient (lemme 2, Chapitre 7) que
$0\lra E/E^q\lra H^+\lra {\cal CL}[q]^+\lra 0.$ En restreignant cette suite ˆ $H^+\cap H'$ et en se
souvenant de la dŽfinition de $E'$, on a

$$0\lra E'/E^q\lra H^+\cap H'\lra {\cal CL}[q]^+.\eqno{(i)}$$

La dernire flche n'Žtant pas forcŽment surjective. On se souvient aussi que $C$ est le sous-groupe de
$E$ des $p$-unitŽs cyclotomiques et $C'=C\cap E'$. On a les inclusions suivantes~:

$$0\underbrace{\quad \subset\quad }_{E_1}C'E^q/E^q\underbrace{\quad \subset\quad }_{{\rm quotient}=E_2}
CE^q/E^q\underbrace{\quad \subset\quad }_{{\rm quotient}=E_3}E/E^q.$$

On a donc $E_1=C'E^q/E^q$, $E_2=CE^q/C'E^q$ et $E_3=E/CE^q$. On a montrŽ au lemme 8 que $\F_q[G^+]$ Žtait
semi-simple. Donc tout module sur $\F_q[G^+]$ est semi-simple. Donc toute suite exacte est scindŽe (cf.
Corollaire 5 de l'appendice 1). Remarquons que si on a les inclusions $0\subset A\subset B\subset C$ de
$R$-modules semi-simples, alors $C$ est isomorphe (comme $R$-module) ˆ $A\oplus (B/A)\oplus (C/B)$.
Puisque (Lemme 8)
$E/E^q$ est un $\F_q[G^+]$-module libre de rang 1, on a l'isomorphisme

$$E_1\oplus E_2\oplus E_3\simeq \F_q[G^+]\simeq \F_q[X]/(X^{p-1\over 2}-1).\eqno{(ii)}$$ 

Donc, $E_1,E_2$ et $E_3$ sont isomorphes ˆ des sous-module de $\F_q[X]/(X^{p-1\over 2}-1)$, qui sont, en
vertu du lemme 7 b), de la forme $\F_q[G^+]/{\rm Ann}(E_i)$. Les ${\rm Ann}(E_i)$ Žtant des idŽaux de
$\F_q[G^+]$, il existe, pour $i=1,2,3$, $\mu_i$, des facteurs unitaires de $X^{p-1\over 2}-1$ tels que
$E_i\simeq \F_q[X]/(\mu_i)$ et bien sžr ${\rm dim}_{\F_q}(E_i)={\deg \mu_i}$.  Par $(ii)$,  et par
comparaison des dimensions sur $\F_q$, on a $\mu_1\cdot\mu_2\cdot\mu_3=X^{p-1\over 2}-1$. 

Clairement, $C'E^q\subset E'$. Donc, on a la suite exacte 

$$1\lra\underbrace{C' E^q/E^q}_{=E_1}\lra E'/E^q\lra E'/C'E^q\lra 1.$$

Donc,

$$E'/E_q\simeq E_1\oplus E'/C' E^q.\eqno{(iii)} $$ 

D'autre part, l'application composŽe $E'\hookrightarrow E\rightarrow E/CE^q$ a pour noyau $E'\cap CE^q$
($Ò\supset"$ est trivial; pour $Ò\subset"$~: soit $c\cdot e^q\in E'$, avec $c\in C$ et $e\in E$. Puisque
$e^q\in E'$, alors $c\in E'\cap C=C'$, donc $c\cdot e^q\in C' E^q$.) Donc, 

$$E'/C'E^q\hookrightarrow E/CE^q=E_3.\eqno{(iv)}$$ 

On obtient (puisque dans notre cas toute suite exacte est scindŽe)~:

$$H'\cap H^+\buildrel (i)\over \hookrightarrow E'/E^q\oplus {\cal CL}[q]^+\buildrel (iii)\over\simeq
E'/C'E_q\oplus E_1\oplus {\cal CL}[q]^+\buildrel (iv)\over \hookrightarrow E_1\oplus E_3\oplus {\cal
CL}[q]^+.$$

Le thŽorme de Thaine nous apprend que tout annulateur de $E_3=E/CE^q$ annule aussi ${\cal CL}[q]^+$. Donc
$\mu_1\cdot\mu_3$ annule $H'\cap H^+$. Supposons par l'absurde que l'annulateur de $H'\cap H^+$ est nul.
Cela implique que $\mu_1\cdot\mu_3=0$ dans $\F_q[X]/(X^{p-1\over 2}-1)$. Cela implique que $\mu_2=1$, ou
encore que $E_2=0$. Donc $C'E^q=CE^q$. Puisque $C'\subset C$ et $C'\cap E^q=C\cap E^q$ (vŽrification
facile, car $E^q\subset E'$), on trouve que
$C=C'$. En effet, soit $c\in C$. Alors $c=c\cdot 1\buildrel C'\cap E^q=C\cap E^q\over =c'\cdot e^q$,
avec $c\in C'$ et $e\in E$. On a $e^q={c\over c'}\in C\cap E^q=C'\cap E^q\subset C'$; donc $c=c'\cdot
e^q\in C'$.

Le Lemme 10 nous dit qu'alors
$p<q$, ce qui contredit l'hypothse.\qed

\bigskip

{\bf ThŽorme 13  (ThŽorme 4 de Mih$\breve{\bf a}$ilescu) }

{\sl Si $x$ et $y$ sont des entiers tels que $x^p-y^q=1$, avec $p,q\geq 11$ premiers, alors $q\equiv 1\pmod
p$ ou $p\equiv 1\pmod q$.

}

{\bf Preuve}

Le ThŽorme 5 nous dit que l'annulateur ${\rm Ann}_{\F_q[G^+]}([x-\zeta_p]^{1+\iota}\in H)=\{0\}$. D'autre
part, on se souvient que $q^2|x$ (ThŽorme 1, Chapitre 6) et que $-\zeta_p$ est une puissance
$q$-me d'une autre racine de l'unitŽ, car $(q,2p)=1$ . Donc
$x-\zeta_p=(-\zeta_p)+x=\beta^q+q^2\gamma\in H'$. Donc, $[(x-\zeta_p)^{1+\iota}]\in H'\cap H^+$. 
  
Supposons par symŽtrie que $p>q$, donc en particulier, $q\not\equiv 1\pmod p$.  Supposons par l'absurde que
$p\not\equiv 1\pmod q$. Le thŽorme 12 nous dit que $H^+\cap H'$ a un annulateur non trivial, ce qui est
contradictoire.\qed

\vfill\eject

\long\def\art#1{{\parindent0pt\item{#1}}\hangindent=7mm\hangafter=-20}
\long\def\artart#1{{\parindent0pt\item{#1}}\hangindent=12mm\hangafter=-20}
\font\para=cmbx12 at 18pt
\def\O{\hbox{$\cal O$}}
\def\U{\hbox{$\cal U$}}
\def\m{\hbox{\rs m\!}}
\def\dst{\displaystyle}
\font\doub=msbm10 at 10pt
\def\lra{\longrightarrow}
\def\qed{\hfill$\square$}
\def\gfP{\relax\ifmmode\bbP\else $\bbP$\fi}
\def\gP{{\euf P}}
\def\P{{\cal P}}
\def\QQ{{\cal Q}}
\def\Log{{\rm Log}}

\def\ggP{{\bf P}}
\newcount\chapnomb \chapnomb=1
\newcount\parnomb \parnomb=1
\pageno =54

\parindent0pt
\centerline{\para CHAPITRE 10}

\bigskip
{\para   Preuve de la Conjecture de Catalan }
\bigskip

A partir de maintenant, la dŽmonstration est trs courte, mais il s'agit de mettre ensemble tous les
ingrŽdients que nous avons patiemment prŽparŽ jusqu'ˆ maintenant. 
\bigskip
Rassemblons le tout dans le thŽorme suivant~:
\bigskip
{\bf ThŽorme}
{\sl 
\art{a)}Les seules solutions  $(x,y)\in \Z^2$, avec $x,y\ne 0$ de l'Žquation

$$x^2-y^3=1$$

sont donnŽes par $(x,y)=(\pm 3,2)$.

\art{b)} Il n'existe pas de solution $(x,y)\in\Z^2$ avec $x,y\ne 0$, tels que $$x^m-y^2=1$$ avec
$m\in\N$.

\art{c)}Il n'existe pas de solution $(x,y)\in\Z^2$ avec $x,y\ne 0$, tels que $$x^2-y^q=1$$ et $q>3$
premier.

\art{d)}Soit $p$ et $q$ des nombres premiers impairs et $x$, $y$ des entiers non nuls tels que
$x^p-y^q=1$. Alors on a~:

\artart{I.}$ p^{q-1}\equiv 1\pmod{q^2}$ et $q^{p-1}\equiv 1\pmod{p^2}$.

\artart{II.}$p,q\geq 11$.

\artart{III.}$p<4q^2$ et $q<4p^2$.

\artart{IV.}$p\equiv 1\pmod q$ ou $q\equiv 1\pmod p$.

}

{\bf Preuve}

La partie a) est le thŽorme d'Euler du premier chapitre. La partie b) est le thŽorme de Lebesgue du
deuxime chapitre. La partie c) est le thŽorme de Ko-Chao du Chapitre 3. La partie d)\,I. est le
thŽorme 1 de Mihailescu du Chapitre 6. La partie d)\, II. est le le thŽorme 2 de Mihailescu du
Chapitre 7. La partie d)\, III. est le le thŽorme 3 de Mihailescu du Chapitre 8. La partie d)\, IV.
est le le thŽorme 4 de Mihailescu du Chapitre 9.\qed

\bigskip

{\bf ThŽorme}

{\sl  Soit $m$ et $n$ des entiers supŽrieurs ou Žgaux ˆ 2. Alors les seules solutions  $(x,y)\in \Z^2$,
avec $x,y\ne 0$ de l'Žquation

$$x^n-y^m=1$$

sont donnŽes par $(x,y,n,m)=(\pm 3,2,2,3)$.

}

{\bf Preuve}

On peut supposer que $n=p$ et $m=q$ sont des nombres premiers. Le cas $p=2$ et $q=2$ est trivialement
impossible (cf introduction). Le cas
$p=2$ ou $q=2$ sont ŽcartŽs gr‰ce aux parties a), b) et c) du thŽorme prŽcŽdent.

Par symŽtrie, en remplaant Žventuellement $(x,y,p,q)$ par $(-y,-x,q,p)$, on peut supposer, gr‰ce ˆ la
partie d)\, IV. du thŽorme  prŽcŽdent que
$p\equiv 1\pmod q$. On affirme qu'alors $p\equiv 1\pmod {q^2}$. En effet on a la suite exacte
$(\Z/q^2\Z)^*\lra (\Z/q\Z)^*\lra\{1\}.$ 
Donc la classe de $p$ modulo $q^2$ est dans le noyau de cette application. Comme l'ordre de $(\Z/q^2\Z)^*$
vaut
$q(q-1)$ et l'ordre de $(\Z/q\Z)^*$ vaut $q-1$, l'ordre de ce noyau vaut $q$. Mais on sait (partie
d)I. du thŽorme prŽcŽdent) que 
$p^{q-1}\equiv 1\pmod {q^2}$. Donc  l'ordre de $p$ modulo $q^2$ divise $q$ et $q-1$. Il vaut donc 1.
Puisque $p<4q^2$ (partie d)\,III.), il nous reste les cas $p=1+q^2$, $p=1+2q^2$ et $p=1+3q^2$. Les cas
$p=1+q^2$ et $p=1+3q^2$ sont impossibles car cela voudrait dire que $p$ ou $q$ est pair. Le cas
$p=1+2q^2$ est impossible aussi, car en rŽduisant modulo 3, on trouverait $p\equiv 0\pmod 3$, donc $p=3$
et $q=1$ ou alors  $q=3$ et $p=19$, mais c'est impossible, car cela contredit la partie d)\, II du
thŽorme prŽcŽdent.\qed 
\vfill\eject

\long\def\art#1{{\parindent0pt\item{#1}}\hangindent=7mm\hangafter=-20}
\long\def\artart#1{{\parindent0pt\item{#1}}\hangindent=12mm\hangafter=-20}
\font\para=cmbx12 at 18pt
\def\O{\hbox{$\cal O$}}
\def\U{\hbox{$\cal U$}}
\def\m{\hbox{\rs m\!}}
\def\dst{\displaystyle}
\font\doub=msbm10 at 10pt
\def\lra{\longrightarrow}
\def\qed{\hfill$\square$}
\def\gfP{\relax\ifmmode\bbP\else $\bbP$\fi}
\def\gP{{\euf P}}
\def\P{{\cal P}}
\def\QQ{{\cal Q}}
\def\Log{{\rm Log}}

\def\ggP{{\bf P}}
\newcount\chapnomb \chapnomb=1
\newcount\parnomb \parnomb=1
\pageno =55
\parindent0pt
\bigskip
\centerline{\para  Appendice 1}
\bigskip
\centerline{\para  Deux mots sur les anneaux semi-simples }
\bigskip
{\bf DŽfinition}

Soit $R$ un anneau commutatif et ${\euf a}$ un idŽal de $R$. On dit que $\euf a$ est {\it radical} si
$R/{\euf a}$ n'a pas de nilpotent, ou, ce qui est Žquivalent, pour tout $x\in R$, si $x^n\in{\euf a}$,
alors $x\in {\euf a}$.

\bigskip

{\bf Lemme 1}

{\sl Soit $R$ un anneau commutatif et ${\euf a}$ et ${\euf b}$ des idŽaux co-premiers (${\euf a}+{\euf
b}=R$) tels que ${\euf a}\cdot{\euf b}={\euf a}$. Alors ${\euf b}=R$.

}

{\bf Preuve}

$R={\euf a}+{\euf b}={\euf a}\cdot{\euf b}+{\euf b}={\euf b}$.\qed

\bigskip

{\bf Lemme 2}

{\sl Soit $R$ un anneau commutatif, $\euf b$ un idŽal, $M$ un $R$-module de type fini et $\varphi\in {\rm
End}_R(M)$ tels que $\varphi(M)\subset {\euf b}M$. Alors il existe $k\in\N$ et $b_1,\ldots ,b_n$ tels que 

$$\varphi^k+b_1\varphi^{k-1}+\cdots +b_k I_M=0_M,$$

o $I_M$ et $0_M$ sont l'endomorphisme identitŽ, respectivement nul sur $M$, alors qu'on notera $0_R$
pour l'ŽlŽment nul de $R$ et $0$ pour l'ŽlŽment nul de $M$.

}

{\bf Preuve}

On choisit des gŽnŽrateurs $m_1,\ldots m_n$ de $M$ et une matrice $(b_{ij})\in M_n({\euf b})\subset
M_n(R)$ tels que $\varphi(m_i)=\sum_{j=1}^n b_{ij}m_j$, $1\leq i\leq n$. On peut voir $M$ comme un
$R[\varphi]$-module  ($R[\varphi]$ est un sous-anneau commutatif de ${\rm End}_R(M)$), via la rgle
$f(\varphi)\cdot x=f(\varphi)(x)$ pour tout $f\in R[X]$ et $x\in M$. Les relations prŽcŽdentes s'Žcrivent 

$$\left (\delta_{ij}\varphi- b_{ij}\right )\pmatrix{m_1\cr\vdots \cr m_n\cr}=\pmatrix{0\cr\vdots\cr
0\cr},\hbox { avec }\delta_{ij}=\cases{1&si $i=j$\cr 0&sinon.\cr}$$

Puisque $R[\varphi]$ est commutatif, on peut former la co-matrice de $\left (\delta_{ij}\varphi-
b_{ij}\right )$. (pour toute matrice $A$, la co-matrice est l'unique matrice $\widetilde{A}$ telle que
$A\widetilde{A}=\widetilde{A}A=\det(A)I_n$). En multipliant par cette co-matrice, on trouve que $\det \left
(\delta_{ij}\varphi- b_{ij}\right )(m_s)=0$ pour tout $s=1,\ldots ,n$. Ainsi $\det \left
(\delta_{ij}\varphi- b_{ij}\right )=0_M$, et on vŽrifie facilement que ce dŽterminant est un polyn™me
unitaire $\in R[\varphi]$.\qed

\bigskip
{\bf Proposition 3}

{\sl  Soit $R$ un anneau commutatif, $\euf b$ un idŽal et $M$ un $R$-module de type fini. Supposons
que ${\rm Ann}(M)+{\euf b}$ soit radical. Si $\psi\, : \, R\rightarrow$ est la projection canonique,
alors on a 

$$\psi({\rm Ann}(M))={\rm Ann}_{R/{\euf b}}(M/{\euf b}M).$$

Avec bien sžr ${\rm Ann}_R(S)=\{\alpha\in R\mid \alpha x=0$ pour tout $x\in S\}$ pour tout $S\subset M$.

}

{\bf Preuve}

On notera $\overline{\alpha}$ pour $\psi(\alpha)$.
 
L'inclusion $\subset$ est immŽdiate, sans mme supposer que ${\rm Ann}(M)+{\euf b}$ soit radical~:  si
$\alpha\in {\rm Ann}(M)$, alors $\alpha M=\{0\}$, donc {\it a fortiori} $\alpha\cdot M/{\euf b}M=\{0\}$
et donc  $\overline{\alpha}\cdot M/{\euf b}M=\{0\}$. Montrons l'autre inclusion. Supposons que
$\overline{\alpha}\in {\rm Ann}_{R/{\euf b}}(M/{\euf b}M)$. Cela veut dire que $\alpha M\subset {\euf
b}M$. Posons $\widehat{\alpha}$ la multiplication par $\alpha$, et appliquons le lemme prŽcŽdent ˆ cet
endomorphisme, il existe $b_1,\ldots ,b_n\in{\euf b}$ tels que $\widehat{\alpha}^n+
b_1\widehat{\alpha}^{n-1}+\cdots +b_n I_M=0_M$. Cela signifie que $\alpha^n+ b_1\alpha^{n-1}+\cdots
+b_n\in{\rm Ann}_R(M)$. Donc, $\alpha^n\in {\euf b}+{\rm Ann}_R(M)$ qui est supposŽ radical. Donc
$\alpha\in {\euf b}+{\rm Ann}_R(M)$. C'est-ˆ-dire $\alpha=b+x$ avec $b\in {\euf b}$ et $x\in 
{\rm Ann}_R(M)$. On trouve alors $\overline{\alpha}=\psi(\alpha)=\psi(x)\in\psi({\rm Ann}(M))$.\qed

\bigskip

{\bf DŽfinition}

On dira (pour faire simple ha-ha) qu'un anneau commutatif $R$ est {\it semi-simple} si c'est un produit
cartŽsien fini de corps. C'est-ˆ-dire $R=\prod_{\alpha\in A} K_\alpha$, avec $A$ fini et $K_\alpha$
corps (commutatif), pour tout $\alpha$. 

Si $A$ est un anneau (pas forcŽment semi-simple) et si $M\ne\{0\}$ est un $A$-module, alors $M$ est dit
{\it simple} s'il n'a pas d'autre sous-module que $\{0\}$ et lui-mme. Un $A$-module qui est somme directe
de modules simples est appelŽ {\it semi-simple}. Evidemment, un anneau semi-simple est un module semi-simple
comme module sur lui-mme.

Soit $R=\prod_{\alpha\in A} K_\alpha$ un anneau semi-simple. Soit 
$\beta\in A$. On pose
$K'_\beta=\{(x_\alpha)\in
\prod K_\alpha\mid x_\alpha=0$ si $\alpha\ne\beta\}$. Les $K_\alpha'$ sont des idŽaux de $R$, et
$R=\bigoplus_{\alpha\in A} K'_\alpha$, et ${\rm Ann}(K'_\alpha)=\bigoplus_{\beta\ne\alpha} K'_\beta$.
Evidemment, $K'_\alpha$ est isomorphe ˆ $K_\alpha$ additivement et multiplicativement, mais on Òvoit"
qu'ils sont sensiblement diffŽrents. De plus, $K'_\alpha\cdot K'_\beta=\cases{0&si $\alpha\ne\beta$\cr
K'_\alpha& si $\alpha=\beta$\cr}$. De plus, chaque $K'_\alpha$ est un $R$-module simple. On en dŽduit que
les idŽaux de $R$ sont de la forme $\bigoplus_{\alpha\in B} K'_\alpha$, avec $B\subset A$. On en tire
aussi que si ${\euf a}$ et $\euf b$ sont des idŽaux de $R$, alors ${\euf a}\cdot {\euf b}={\euf
a}\cap{\euf b}$. Si $a\in {\euf a}\cdot{\euf b}$, alors il existe $a_1\in{\euf a}$ et $a_2\in{\euf b}$
tels que $a=a_1a_2$. En particulier, ${\euf a}^2={\euf a}$, donc, pour tout $a\in {\euf a}$, il existe
$a_1,a_2\in{\euf a}$ tels que $a=a_1a_2$.
\bigskip

Maintenant on va voir un thŽorme trs connu sur les modules semi-simples. Comme il n'est pas trop long,
on va en donner la preuve.

\bigskip

{\bf ThŽorme 4}

{\sl Soit $R$ un anneau commutatif (pas forcŽment semi-simple) et $M$ un $R$-module. Les conditions
suivantes sont Žquivalentes~:

\art{a)} $M$ est une somme de modules simples.

\art{b)} $M$ est une somme directe de modules simples (i.e. $M$ est semi-simple).

\art{c)} Tout sous-module de $M$ est un facteur direct dans $M$.

}

{\bf Preuve }

Prouvons a) $\Rightarrow$ b). Supposons donc que $M=\sum_{i\in I}P_i$ o $P_i$ est simple. On peut
supposer que $P_i\ne P_j$ si $i\ne j$. Pour tout $J\subset I$, on pose $P_J=\sum_{i\in J}P_i$.
ConsidŽrons la famille de tous les $J$ tels que $P_J=\oplus_{i\in J} P_i$. Cette famille est non vide
(elle contient au moins $J=\emptyset$). Cette famille est partiellement ordonnŽe par l'inclusion. On veut
utiliser le Lemme de Zorn~: soit $(J_\alpha)$ une cha"ne. Posons $J'=\cup J_\alpha$. Alors $P_{J'}$ est
une somme directe. En effet, si $m_1+\cdots +m_s=0$, alors, puisque les $m_i$ sont en nombre fini, il
existe $P_{J_\alpha}$ dans lequel se trouvent tous les $m_i$. C'est une contradiction, car $P_{J_\alpha}$
est une somme directe. Par le Lemme de Zorn, notre famille contient un ŽlŽment maximal, disons $J$. On
affirme que $P_J=M$. Soit $i\in I$. On a $P_J\cap P_i\subset P_i$. Comme $P_i$ est simple, on a $P_J\cap
P_i=P_i$,  ce qui signifie que $P_i\subset P_J$; ou alors $P_J\cap P_i=\{0\}$, ce qui veut dire de
$P_J\oplus P_j$ est une somme directe, ce qui contredit la maximalitŽ de $J$. Donc $P_i\subset P_J$ pour
tout $i\in I$, et ainsi $M=P_J$ est une somme directe de modules simples.

Prouvons b) $\Rightarrow$ c). Soit $N$ un sous-module de $M$ et $f\, :\, M\rightarrow M/N$ la projection
canonique. Alors, si $M=\oplus_i P_i$, on a $M/N=\sum f(P_i)$, et $f(P_i)$ est soit isomorphe ˆ $P_i$, soit
$\{0\}$ (le noyau d'un homomorphisme sur un module simple est un sous-module, qui est soit tout, soit
rien...). Donc, $M/N$ est aussi une somme de module simple. Par la partie a) $\Rightarrow$ b), on en
dŽduit que cette somme peut tre supposŽe directe. Plus prŽcisŽment, il existe $K\subset I$ tel que
les $f(P_i)$ sont distincts et non nuls et tels que

$$M/N=\bigoplus_{i\in K} f(P_i).\eqno{(i)}$$

Posons $N'=\sum_KP_i$. Cette somme est directe, et $f(N')=\sum_kf(P_i)=f(M)$. En prenant les images
inverses, on trouve que $M=N+N'$. Puisque $(i)$ est une somme directe, $f$ est injective sur $N'$, ce
qui veut dire que $N\cap N'=\{0\}$ et donc $M=N'\oplus N$.

Prouvons c) $\Rightarrow$ a). Soit $N$ la somme de tous les sous-modules simples de $M$. Par hypothse,
il existe $N'$ tel que $M=N\bigoplus N'$. nous allons montrer que $N'=\{0\}$. Supposons le contraire, donc
que $N\ne \{0\}$. Soit $P$ un sous-module cyclique de $N'$. On peut montrer facilement gr‰ce au lemme de
Zorn que $P$ contient (au moins) un sous-module propre maximal $Q$ (car il est de type fini). Par
hypothse, il existe $Q'$ tel que $M=Q\oplus Q'$. Clairement, $P=Q\oplus(Q'\cap P)$. Puisque $Q$
est maximal dans $P$, on a $Q'\cap P\simeq P/Q$ est simple. Donc $P$ et par suite $N'$ contient un
sous-module simple et donc $N\cap N'\ne\{0\}$. C'est une contradiction, donc $N'=\{0\}$ et $N=M$.\qed

\bigskip

{\bf Corollaire 5}

{\sl Soit $R=\prod_{\alpha\in A} K_\alpha$ un anneau semi-simple. Alors les affirmations suivantes
sont vraies~:

\art{a)} Tout $R$-module est semi simple.

\art{b)} Tout $R$-module $M$ s'Žcrit $\bigoplus_{\alpha\in A} M_\alpha$ o $M_\alpha$ est la somme des
sous-modules simples de $M$ isomorphes ˆ $K'_\alpha$. 

\art{c)} Toute suite exacte de $R$-module est scindŽe.

}

{\bf Preuve}

Prouvons a) et b). On a $M=\sum_{x\in M} Rx=\sum_{x\in M}\sum_{\alpha\in A} K'_\alpha x$. Pour tout $x\in M$
et $\alpha\in A$, l'application de $K'_\alpha$ sur $K'_\alpha x$, $y\mapsto yx$ est $R$-linŽaire, surjective
et son noyau est $\{0\}$ ou $K'_\alpha$. Donc $K'_\alpha x$ est isomorphe ˆ $K'_\alpha$ ou ˆ $\{0\}$.
Donc, $M$ est la somme de modules simples, il est donc semi-simple en vertu du lemme prŽcŽdent. Par
suite on a bien $M=\sum_{\alpha\in A} M_\alpha$ o $M_\alpha$ est la somme des
sous-modules simples de $M$ isomorphes ˆ $K'_\alpha$. Voyons que cette somme est directe. Soit $x\in
M_\alpha\cap \big (\sum_{\beta\ne\alpha}M_\beta\big )$. Comme $x\in M_\alpha$, il est annulŽ par
$\bigoplus_{\beta\ne\alpha}K'_\beta$; et comme $x\in \sum_{\beta\ne\alpha}M_\beta$, il est annulŽ par
$K'_\alpha$. Ainsi, ${\rm Ann}_R(x)=R$, et donc $x=0$. La somme est donc directe.

Prouvons c). Soit $0\rightarrow M_1\buildrel f\over \rightarrow M \buildrel g\over \rightarrow
M_2\rightarrow 0$ une suite exacte de $R$-modules. Puisque $M$ est semi-simple (partie a)), il existe
$M_3$ tel que $M=f(M_1)\oplus M_3$. Puisque $f(M_1)$ est le noyau de $g$, on a $g(M_3)=g(M)=M_2$ et si
$m_3$ et $m_3'$ sont tels que $g(m_3)=g(m_3')$, alors $m_3-m_3'\in f(M_1)\cap M_3=\{0\}$. On a montrŽ que
$M_3$ et $M_2$ sont isomorphes et donc que
$M\simeq M_1\oplus M_2$, ce qui prouve que le suite est scindŽe.\qed

\bigskip

{\bf Proposition 6}

{\sl Soit $R=\prod_{\alpha\in A} K_\alpha$ un anneau semi-simple et $M$ un $R$-module. Il existe
$a\in M$ tel que ${\rm Ann}_R(a)={\rm Ann}_R(M)$. Donc, $M$ contient un sous-module cyclique (celui
engendrŽ par $a$) isomorphe ˆ $R/{\rm Ann}_R(M)$. En particulier, si $R$ et $M$ sont finis, alors
$|M|\geq |R/{\rm Ann}_R(M)|$, avec ŽgalitŽ si et seulement si $M$ est cyclique.

}

{\bf Preuve}

Posons $B=\{\alpha\in A\mid M_\alpha\ne\{0\}\}$. On a $M=\bigoplus_{\alpha\in B} M_\alpha$ et ${\rm
Ann}_R(M)=\bigoplus_{\beta\in A\setminus B}K'_\beta$. Pour chaque $\alpha\in B$, on choisit $0\ne
x_\alpha\in M_\alpha$. On a ${\rm Ann}_R(x_\alpha)=\bigoplus_{\beta\ne\alpha} K'_\beta$. Posons
$a=\sum_{\alpha\in B}x_\alpha$. On a ${\rm Ann}_R(a)=\bigcap_{\alpha\in B}{\rm Ann}_R(x_\alpha)$, car les
somme des $M_\alpha$ est directe. Ainsi, ${\rm Ann}_R(a)=\bigoplus_{\beta\in A\setminus B} K'_\beta={\rm
Ann}_R(M)$. Le reste en dŽcoule facilement.\qed

\bigskip

{\bf Proposition 7}

{\sl Soit $R=\prod_{\alpha\in A} K_\alpha$ et $M$ un $R$-module cyclique. Soit $M'$ un sous-module de
$M$. Alors $M'$ est cyclique, ${\rm Ann}_R(M')\cdot{\rm Ann}_R(M/M')={\rm Ann}_R(M)$, et les deux idŽaux
${\rm Ann}_R(M')$ et ${\rm Ann}_R(M/M')$ sont co-premiers.

}

{\bf Preuve}

Soit $B\subset A$ tel que ${\rm Ann}_R(M)=\bigoplus_{\beta\in A\setminus B}K'_\beta$. Alors $M\simeq R/{\rm
Ann}_R(M)\simeq \bigoplus_{\beta\in B}K'_\beta$ (idŽal de $R$). Via cet isomorphisme, $M'$ est isomorphe ˆ
un sous-idŽal; c'est-ˆ-dire $M'\simeq \bigoplus_{\beta\in B'} K'_\beta$ pour $B'\subset B$. Et alors,
${\rm Ann}_R(M')=\bigoplus_{\beta\in A\setminus B'}K'_\beta$, et $M'\simeq R/{\rm Ann}_R(M')$. Donc, en
vertu de la proposition prŽcŽdente, $M'$ est cyclique. D'autre part, $M/M'\simeq \bigoplus_{\beta\in
B\setminus B'}K'_\beta$ et donc, ${\rm Ann}_R(M/M')\simeq \bigoplus_{\beta\not\in B\setminus B'}K'_\beta$.
Ce qui prouve que ${\rm Ann}_R(M/M')\cdot{\rm Ann}_R(M')=\bigoplus_{\beta\in A\setminus B}K'_\beta={\rm
Ann}_R(M)$ et que ces deux idŽaux sont premiers entre eux.\qed

\vfill\eject

\long\def\art#1{{\parindent0pt\item{#1}}\hangindent=7mm\hangafter=-20}
\long\def\artart#1{{\parindent0pt\item{#1}}\hangindent=12mm\hangafter=-20}
\font\para=cmbx12 at 18pt
\def\O{\hbox{$\cal O$}}
\def\U{\hbox{$\cal U$}}
\def\m{\hbox{\rs m\!}}
\def\dst{\displaystyle}
\font\doub=msbm10 at 10pt
\def\lra{\longrightarrow}
\def\qed{\hfill$\square$}
\def\gfP{\relax\ifmmode\bbP\else $\bbP$\fi}
\def\gP{{\euf P}}
\def\P{{\cal P}}
\def\QQ{{\cal Q}}
\def\Log{{\rm Log}}

\def\ggP{{\bf P}}
\newcount\chapnomb \chapnomb=1
\newcount\parnomb \parnomb=1
\pageno =58

\parindent0pt
\bigskip
\centerline{\para  Appendice 2}
\bigskip

{\para  Le thŽorme de Thaine }
\bigskip
Tout d'abord quelques lemmes pour Žnoncer le thŽorme de manire sensiblement diffŽrente de celle
qu'on a vue au chapitre 9~:

\bigskip

{\bf Lemme 1}

{\sl Soit $A$ un anneau local commutatif (c'est-ˆ-dire qui ne possde qu'un seul idŽal maximal,
disons $\P$).  Soit $f\, :\, M\lra N$ un homomorphisme injectif de $A$-modules libres de mme rang.
Supposons que $\overline{f}\, :\,M/\P M\lra N/\P N$ soit un isomorphisme. Alors $f$ est un
isomorphisme.

}

{\bf Preuve}

Choisissons un base de $M$ et une base de $N$. Posons $M_f$ la matrice de $f$ relativement ˆ ces
deux bases. Soit $\overline{M_f}$ la matrice de $\overline{f}$. Puisque $\overline{f}$ est un
isomorphisme, son dŽterminant est inversible dans $A/\P$. Cela veut dire que
$\det(M_f)\not\in \P$, donc $\det(M_f)$ est inversible, car si ce n'Žtait pas le cas, l'idŽal
engendrŽ par $\det(M_f)$ et $\P$ contiendrait strictement $\P$ qui est maximal.\qed

\bigskip  

{\bf Lemme 2}

{\sl Soit $G$ un groupe cyclique (attention, on Žcrira additivement ce groupe, mais on appliquera
ce lemme ˆ des groupes notŽs multiplicativement) d'ordre $n$, et $q$ un nombre premier tel que $q\not
\hskip-2pt |\ n$. Soit $M$ un $\Z[G]$-module fini. Alors $M[q]=\{x\in M\mid qx=0\}$ et $M/qM$
sont isomorphes comme
$\Z[G]$ (ou
$\F_q[G]$)-modules.

}

{\bf Preuve}

On peut supposer que $M$ est un $q$-groupe, car si $p\ne q$, les $p$-sous-groupes de Sylow qui sont des
sous-$\Z[G]$-modules de $M$ ne contribuent ni ˆ $M[q]$, ni ˆ $M/qM$. Remarquons de plus que $M[q]$ et
$M/qM$ sont isomorphes comme $\Z$-modules. En effet, si $M=\Z/q^{r_1}\Z\oplus \Z/q^{r_2}\Z\oplus\cdots
\oplus \Z/q^{r_s}\Z$, alors $M[q]=\Z/q\Z\oplus\cdots \oplus \Z/q\Z$ avec $s$ facteurs. On a aussi
$qM=\Z/q^{r_1-1}\Z\oplus\cdots\oplus \Z/q^{r_s-1}\Z$, donc $M/qM\simeq \Z/q\Z\oplus\cdots \oplus \Z/q\Z$.
On peut gŽnŽraliser ce raisonnement ˆ tout anneau principal (et mme pour tout produit d'anneaux
principaux)~: si $A$ est un produit d'anneaux principaux, et si $q\in A$ est premier et si $M$ est de
gŽnŽration finie tel que $q^m\cdot M=\{0\}$ pour $m$ assez grand, alors $M[q]\simeq M/qM$, comme
$A$-modules.

Autre rŽduction~: comme $M$ est de $q$-torsion, il est canoniquement un $\Z_q$ module~: si $x\in M$ est
d'ordre $q^r$, et $a\in \Z_q$, alors on pose $a\cdot x=b\cdot x$ o $b\in \Z$ tel que $b\equiv
a\pmod{q^r\Z_q}$ (souvenons-nous que $\Z_q/q^r\Z_q\simeq \Z/q^r\Z$ (cf. Chapitre 7)). En rŽsumŽ, $M$
est un $\Z_q[G]$-module de $q$-torsion, fini.

On a $\Z_q[G]\simeq \Z_q[X]/(X^n-1)$. Puisque $\Z_q[X]$ est factoriel, Žcrivons $X^n-1=\prod_{i=1}^s
f_i$ la factorisation en polyn™mes irrŽductibles dans $\Z_q[X]$ de $X^n-1$. On peut supposer les
$f_i$ unitaires, car on est en caractŽristique 0. On a donc
un homomorphisme injectif de $\Z_q$-algbres~:

$$\Z_q[G]\simeq \Z_q[X]/(X^n-1) \lhook\joinrel\longrightarrow \prod_{i=1}^s\Z_q[X]/(f_i),\eqno{(i)}$$

car $X^n-1$ n'a pas de racine multiple.

La source et le but sont des $\Z_q$-modules libres de rang $n$. En quotientant par l'idŽal engendrŽ par
$q$, on obtient 

$$\F_q[X]/(X^n-1)\longrightarrow \prod_{i=1}^s\F_q[X]/(\overline{f_i}),$$

o $\overline{f_i}$ est, pour tout $i$, la rŽduction modulo $q$ de $f_i$. Par le thŽorme de
Hensel (cf. Chapitre 7), les $\overline{f_i}$ sont irrŽductibles et premiers entre eux, car
$X^n-1$ n'a que des racines simples ($q\not
\hskip-0.4pt |\ n$). Donc, l'application est encore injective.
Elle est donc surjective car un homomorphisme injectif entre deux espaces vectoriels de mme
dimension est surjectif. Par le Lemme 1, l'homomorphisme $(i)$ est un isomorphisme. Il reste ˆ
voir que pour tout $i$, l'idŽal $\Z_q[X]/(f_i)$ est principal. Il suffit de voir que tout
idŽal premier est principal. Soit donc $\Z_q[X]\supsetneqq\P\supsetneqq (f_i)$. Voyons que
$\P\supset (q,f_i)$. Soit $g\in\P$, $g\not\in (f_i)$. Quitte ˆ prendre le reste de la division
euclidienne de $g$ par $f_i$, au lieu de $g$, on peut supposer que $g$ et $f_i$ sont premiers
entres eux dans $\Q_p[X]$. Par le thŽorme de Bezout, et en remultipliant par les
dŽnominateurs Žventuels, il existe $r,s\in \Z_q[x]$, $m\in \N$, $m\geq 1$ et $u\in\Z_q^*$ tels que
$r\cdot g+s\cdot f_i=u\cdot q^m$. Donc, $q^m\in\P$, donc $q\in\P$, car $\P$ est un idŽal premier.
On a montrŽ que $\P\supset (q,f_i)$. En quotientant $\Z_q[X]$ par $(q,f_i)$, on obtient
$\F_q[X]/(\overline{f_i})$ qui est un corps (donc simple). Ainsi $\P=(q,f_i)$ et ainsi,
$\overline{\P}=(q)$ est principal. Ce qui prouve le lemme.\qed  

\bigskip\goodbreak

On rappelle que $U_p^+=\Z[\zeta_p^+]^*$, que $E=\{u(1-\zeta_p)^k\mid u\in U_p, k\in
\Z\}=\Z[\zeta_p,{1\over p}]^*$, que ${\cal CL}^{pl}={\cal CL}_{\Q(\zeta_p^+)}$ sont les classes d'idŽaux
de $\Q(\zeta_p^+)$ et que $C=\{u\in E\mid u={1-\zeta_p^i\over 1-\zeta_p^j}\omega(1-\zeta_p)^k,\ \omega
\hbox{ est une racine de l'unitŽ et } i,j\in\N\ \hbox{ et } k\in\Z\}$, les $p$-unitŽs cyclotomiques de
$\Q(\zeta_p)$. Si $c\in C$ est telle que $v_\pi(c)=0$, avec $\pi=1-\zeta_p$, on dit que $c$ est
une {\it unitŽ cyclotomique} et on note $C_0$ l'ensemble des unitŽs cyclotomiques.
\bigskip

\goodbreak
{\bf Lemme 3}

{\sl Notons $C_0^+=C_0\cap \R$. Alors les $\F_q[G^+]$-modules $E/CE^q$ et
$U^+/C_0^+{U^+}^q$ sont isomorphes. De plus, si on pose  pour tout $a\in\N$ tel que
$p\not\hskip-0.4pt |\   a$ $\xi_a=\zeta_p^{1-a\over 2}\cdot {\zeta_p^a-1\over\zeta_p-1}$, alors
$C_0^+=<-1,\xi_a>$ est le sous-groupe de $U^+$ engendrŽ par $-1$ et les $\xi_a$.

 }
{\bf Preuve}

Soit $x=u\cdot (1-\zeta_p)^k\in E$ avec $u\in E_p^*=U_p$, et $[x]$ la classe de $x$ modulo $CE^q$.
Par le Lemme de Kummer (cf. Chapitre 7), il existe $u_0\in U^+$ et $s\in\Z$ tel que
$u=u_0\cdot\zeta_p^s$. Ainsi $[x]=[u\cdot (1-\zeta_p)^k]=[u]=[u_0\cdot\zeta_p^s]=[u_0]$. Donc
l'homomorphisme de $\F_q[G^+]$-modules $U_p^+\lra E/CE^q$, $u_0\longmapsto [u_0]$ est surjective.
Voyons le noyau de cette application. Soit $u_0\in U_p^+$, tel que $u_0=c\cdot e^q$ avec $c\in C$ et
$e\in E$. On peut supposer que $v_\pi(c)=v_\pi(e^q)=0$. Donc $c\in C_0$ et $e\in U_p$. Toujours gr‰ce au
Lemme de Kummer, $e=\zeta_p^s\cdot e_0$, avec $e_0\in U_p^+$. Donc, $u_0=(c\cdot\zeta_p^{sq})\cdot
e_0^q$. On vŽrifie que $(c\cdot\zeta_p^{sq})\in C_0\cap \R$, et donc $u_0\in C_0^+{U_p^+}^q$.

Pour la seconde partie de la preuve, soit $\alpha\in C_0^+$. Alors, puisque $\alpha\in C_0$, il
existe $d\in \Z$ tel que $\alpha=\pm\zeta_p^d\prod_{a=1}^{p-1}(1-\zeta_p^a)^{c_a}$ avec
$c_a\in\Z$ pour tout $a$. Puisque $\alpha$ est inversible, on a l'ŽgalitŽ $E_p=\alpha
E_p=\prod_{a=1}^{p-1}(1-\zeta^a)^{c_a}E_p=(\pi E_p)^{\sum c_a}$, car les idŽaux engendrŽ par
les $1-\zeta_p^a$ sont tous Žgaux ˆ $\pi E_p$ (cf. Lemme 2, Chapitre 6). Puisque $\pi E_p$
est un idŽal premier, on en dŽduit que $\sum_{a=1}^{p-1} c_a=0$. Ainsi,

$$\alpha=\pm \zeta^{d}\prod_{a=1}^{p-1}\left ({1-\zeta_p^a\over 1-\zeta_p}\right
)^{c_a}=\pm\zeta_p^{d-\sum_a c_a({1-a\over 2})}\prod_{a=1}^{p-1}\xi_a^{c_a}.$$

On vŽrifie facilement que, pour tout $a$, $\xi_a\in \R$, et  donc $\zeta_p^{d-\sum c_a({1-a\over
2})}=\pm 1$, ce qui prouve le lemme.\qed

\bigskip

On peut maintenant Žnoncer un thŽorme Žquivalent au thŽorme de Thaine tel qu'on l'a ŽnoncŽ au
Chapitre 9. 
\bigskip
{\bf ThŽorme 4 (ThŽorme de Thaine seconde version)}

{\sl Soit $p,q$ des nombres premiers impairs distincts tels que $p\not\equiv 1\pmod q$. Posons
$F=\Q(\zeta_p^+)$, $G^+={\rm Gal}(\Q(\zeta_p^+)/\Q)$, $U_p^+=\Z[\zeta_p^+]^*$, $C_0^+$, les
unitŽs cyclotomiques de
$\Q(\zeta_p^+)$, qui est, comme on vient de le voir, le sous-groupe de $U_p^+$ engendrŽ par $-1$ les
$\xi_a=\zeta_p^{1-a\over 2}\cdot{\zeta_p^a-1\over\zeta_p-1}$. Notons encore ${\cal CL}^{pl}$ le
groupe des classes d'idŽaux de $F$. 
\medskip
\centerline{Si $\theta\in\Z[G^+]$ annule  $U_p^+/C_0^+{U_p^+}^q$ alors, il annule ${\cal
CL}^{pl}/{{\cal CL}^{pl}}^q$.}

}
\bigskip

Il est Žvident que ce thŽorme est Žquivalent ˆ celui ŽnoncŽ au Chapitre 9. Le thŽorme
ŽnoncŽ au chapitre 9 dit que tout annulateur de $E/CE^q$ annule aussi ${\cal CL}^{pl}[q]$. Or,
on a vu au Lemme 3 que $E/CE^q$ et $U^+/C_0^+{U^+}^q$ sont isomorphes et au lemme 2 que ${\cal
CL}^{pl}[q]$ est isomorphe ˆ ${\cal CL}^{pl}/{{\cal CL}^{pl}}^q$, ceci parce que $p\not\equiv
1\pmod q$. C'est cette version que nous allons dŽmontrer.
\medskip
Maintenant, un petit peu de thŽorie de Galois supplŽmentaire.

\bigskip
{\bf Notation }

Si $K$ et $L$ sont des corps de nombres, on note $KL$ le plus petit corps contenant
$K$ et $L$. 
\bigskip
{\bf Lemme 5}

{\sl 
\art{a)}Soit $K\subset L\subset E$ des corps de nombres. Supposons que $E/K$ soit une extension
galoisienne de groupe $G$. Alors $E/L$ est galoisienne et ${\rm Gal}(E/L)=\{g\in G\mid g|_L={\rm
Id}_L\}:=H$. De plus $L/K$ est galoisienne si et seulement si $H$ est un sous groupe normal de
$G$ et dans ce cas, ${\rm Gal}(L/K)=G/H$.

\art{b)}Soit $K\subset L$ et $K\subset E$ deux extensions de corps de nombres. On suppose que
$L/K$ est galoisienne. Alors $EL/L$ est aussi galoisienne et ${\rm Gal}(EL/E)\simeq {\rm
Gal}(L/L\cap E)$.

\art{c)}Soit $K\subset L$ et $K\subset E$ deux extensions galoisiennes de corps de nombres
telles que $L\cap E=K$. Alors $EL/K$ est une extension galoisienne et ${\rm Gal}(EL/K)\simeq
{\rm Gal}(L/K)\times {\rm Gal}(E/K) $.

} 
{\bf Preuve}

Ces rŽsultats se trouvent dans [Lang1], au Chapitre VI, ¤1. Il s'agit du ThŽorme 1.10 et du
Corollaire 1.9, p. 265, pour la partie a); du ThŽorme 1.12, p. 266, pour la partie b); et du
ThŽorme 1.14, p. 267, pour la partie c).\qed

\bigskip\goodbreak

{\bf Notations}

 Posons $\Delta={\rm Gal}(\Q(\zeta_q)/\Q)={\rm Gal}(F(\zeta_q)/F)$ (cf. Lemme
prŽcŽdent partie b)). On a de mme $G^+={\rm Gal}(F/\Q)={\rm Gal}(F(\zeta_q)/\Q(\zeta_q))$.
Puisque
$F$ et $\Q(\zeta_q)$ sont linŽairement disjoints (i.e. $F\cap\Q(\zeta_q)=\Q$), on a que
$F(\zeta_q)/\Q$ est une extension galoisienne et
${\rm Gal}(F(\zeta_q)/\Q)=G^+\times \Delta$ (cf. Lemme prŽcŽdent partie c)). IntŽressons-nous
maintenant au corps
$F(\zeta_q,\root q\of {U_p^+})$. Puisque $U_p^+$ est un $\Z$-module de gŽnŽration finie, alors
$U_p^+/{U_p^+}^q$ est fini. Donc $F(\zeta_q,\root q\of {U_p^+})$ est un corps de nombres. De
plus, chaque fois qu'une racine $q$-ime d'un ŽlŽment est dans ce corps, toutes les autres y
sont aussi, car $\zeta_q$ en fait partie. Ainsi, l'extension $F(\zeta_q,\root q\of {U_p^+})/F(\zeta_q)$
est galoisienne. Notons $G_0$ le groupe de Galois de cette extension.

En rŽsumŽ, on a la situation suivante~:



\vskip 3cm

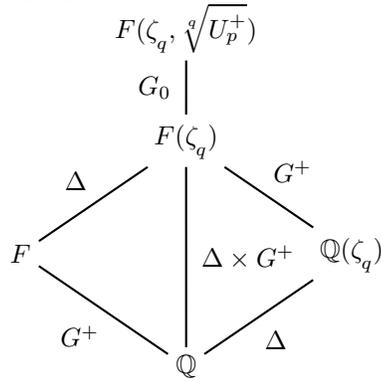
\rput(5,0){\rnode{F}{$F $}} 
\rput(9.4,0){\rnode{Qzeta}{$\Q(\zeta_{q})$}}
\rput(7.2,-1.5){\rnode{Q}{$\Q$}}
\rput(7.2,1.5){\rnode{Fzeta}{$F(\zeta_{q})$}}
\rput(7.2,3){\rnode{FF}{$F(\zeta_{q},\root q\of {U_{p}^+})$}}
\ncline[nodesep=3pt]{Q}{F}
\Aput{$G^+$} 
\ncline[nodesep=3pt]{Q}{Qzeta}
\Bput{$\Delta $} 
\ncline[nodesep=3pt]{Q}{Fzeta}
\Bput{$\Delta\times G^+$} 
\ncline[nodesep=3pt]{F}{Fzeta}
\Aput{$\Delta$} 
\ncline[nodesep=3pt]{Qzeta}{Fzeta}
\Bput{$G^+$} 
\ncline[nodesep=3pt]{Fzeta}{FF}
\Aput{$G_{0}$} 

\vskip 2cm

On montrera mme au Lemme 9 que $F(\zeta_q,\root q\of {U_p^+})/\Q$ est galoisienne.\goodbreak
\bigskip
{\bf Lemme 6 (second lemme de Kummer)} 

{\sl Sous les mmes hypothses que prŽcŽdemment, en posant $\mu_q$ l'ensemble des racines
$q$-ime de l'unitŽ, on a 

$$G_0={\rm Gal}(F(\zeta_q,\root q\of {U_p^+})/F(\zeta_q))\simeq {\rm
Hom}(U_p^+/{U_p^+}^q,\mu_q).$$

Cet isomorphisme est un isomorphisme de groupe et mme de $\F_q$-espace vectoriel. 

}

{\bf Preuve}

ConsidŽrons l'application 
$$\eqalign{\phi\, :\, G_0\times U_p^+/{U_p^+}^q&\lra\mu_q\cr
(g,[u])&\longmapsto {g(\root q\of u)\over\root q\of u}\cr}$$

o $\root q\of u$ est n'importe quelle racine $q$-ime de $u$. L'application est bien dŽfinie,
car, d'une part toute racine $q$-ime diffre d'une autre d'une racine $q$-ime de l'unitŽ, sur
laquelle $g$ est l'identitŽ et d'autre part, si $u'=u\cdot v^q$, avec $u',v\in U_p^+$ alors
${g(\root q\of u)\over\root q\of u}={g(\root q\of u')\over\root q\of u'}$, car $g$ est l'identitŽ
sur $U_p^+$. L'application $\phi$ (qu'on appelle parfois ÒKummer pairing") est
bi-multiplicative. La bi-multiplicativitŽ ˆ droite est Žvidente. Pour celle de gauche, on a

$$\phi(gg',[u])={gg'(\root q\of u)\over\root q\of u}={gg'(\root q\of u)\over g'(\root q\of
u)}\cdot {g'(\root q\of u)\over\root q\of u}=\phi(g,[u])\cdot\phi(g',[u]),$$

car $g'(\root q\of u)^q=g'(u)=u$, donc $g'(\root q\of u)$ est une racine $q$-ime de $u$.

D'autre part, supposons que $\phi(g,[u])=1$ pour tout $u\in U_p^+$. Cela veut dire que $g(\root
q\of u)=\root q\of u$ pour tout $u\in U_p^+$, ce qui veut dire que $g$ est l'identitŽ. Donc
$\phi$ est non dŽgŽnŽrŽe ˆ gauche.\hfill$(i)$

Montrons celle ˆ droite. Remarquons tout d'abord que
l'on a ${F(\zeta)^*}^q\cap U_p^+={U_p^+}^q$. En effet, l'inclusion $\supset$ est triviale.
Pour l'autre inclusion, supposons que $u=u_0^q$, avec $u\in U_p^+$ et $u_0\in F(\zeta_p)$.
Posons $N=N_{F(\zeta_p)/F}$ la norme de l'extension $F(\zeta_p)/F$. Appliquant $N$ ˆ $u=u_0^q$,
on obtient $u^{q-1}=(N(u_0))^q$, ou encore $u=\big ({u\over N(u_0)}\big )^q$. Or, $u_0$ et
$u_0^{-1}$ sont des entiers de $F(\zeta_p)$, car ils sont solutions des Žquations $X^q-u=0$
respectivement $X^q-u^{-1}=0$, donc, $N(u_0^{-1})={1\over N(u_0)}\in U_p^+$. Et donc $ {u\over
N(u_0)}\in U_p^+$. L'autre inclusion est ainsi prouvŽe.

Donc si $\phi(g,[u])=1$ pour tout $g\in G_0$, cela veut dire que $g(\root q\of u)=\root q\of u$,
pour tout $g$ et pour tout racine $q$-ime de $u$, donc $u\in U_p^+\cap {F(\zeta)^*}^q={U_p^+}^q$, ce qui
veut dire que $[u]=1$.\hfill$(ii)$.

On a montrŽ que l'application $\phi$ est non-dŽgŽnŽrŽe ˆ gauche et ˆ droite. Cela prouve le
lemme. En effet, soit l'application 

$$\eqalign{\check{}\, :\, G_0 &\lra{\rm Hom }(U_p^+/{U_p^+}^q,\mu_q)\cr
g&\longmapsto \check g=\phi(g,\cdot)\cr}$$

La relation $(ii)$ montre que $\check{}$ est injective, donc, $G_0$ est un $q$-groupe dont tous
les ŽlŽments sont d'ordre $q$ (ou $1$). Ainsi, tous ces groupes sont des $\F_q$-espaces
vectoriels. Et on a $|G_0|\leq |{\rm Hom}(U_p^+/{U_p^+}^q,\mu_q)=|\widehat
{U_p^+/{U_p^+}^q}|=|U_p^+/{U_p^+}^q|$. Finalement, l'application $[u]\longmapsto
\phi(\cdot,[u])$ de $U_p^+/{U_p^+}^q$ dans ${\rm Hom }(G_0,\mu_q)$ est injective gr‰ce ˆ la
relation $(i)$, ainsi $|U_p^+/{U_p^+}^q|\leq {\rm Hom }(G_0,\mu_q)\buildrel \F_q-\rm
e.v.\over =|G_0|$. Donc
$|U_p^+/{U_p^+}^q|=|G_0|$ et ainsi $\check{}$ est un isomorphisme, qu'on appellera {\it isomorphisme de
Kummer} ce qui prouve le lemme.\qed

\bigskip

On verra mieux au Lemme 9~:  l'isomorphisme de Kummer est en fait un isomorphisme de $\F_q[G^+]$-modules.
\bigskip

{\bf DŽfinition}

Soit $M$ est $\F_q[G^+]$-module. On munit ${\rm Hom}(M,\F_q)={\rm Hom}_{\F_q}(M,\F_q)$
d'une structure de
$\F_q[G^+]$-module gr‰ce ˆ l'action suivante~: si $\sigma\in G^+$ et $\varphi\in {\rm
Hom}(M,\F_q)$ , alors on pose

$$(\varphi^{\sigma})(x)=\varphi(\sigma^{-1}(x))$$ 

On dit alors que ${\rm Hom}(M,\F_p)$ est muni de la {\it structure de $\F_q[G^+]$-module duale}.
Si $\theta=\sum_{\sigma\in G^+}a_\sigma\sigma\in \F_q[G^+]$, on dŽfinit l'application

$$\theta\longmapsto\hat\theta:=\sum_{\sigma\in G^+}a_\sigma\sigma^{-1}.$$

Puisque $G^+$ est abŽlien, cette application est un automorphisme de $\F_q[G^+]$.

Posons encore $s(G^+)=\sum_{\sigma\in G^+}\sigma$ (on l'avait dŽjˆ dŽfinit au Chapitre 5, mais
il n'est pas inutile de le rappeler). Posons encore ${\euf A}=\{\sum_{\sigma\in G^+}a_\sigma\sigma\in
\F_q[G^+]\mid \sum_{\sigma\in G^+}a_\sigma=0\}$. C'est un idŽal de $\F_q[G^+]$ appelŽ {\it l'idŽal
d'augmentation}.
\bigskip
\goodbreak
{\bf Lemme 7}

{\sl Sous les mme hypothses que la dŽfinition prŽcŽdente, on a les isomorphismes de
$\F_q[G^+]$-modules~:

$${\rm Hom}(U_p^+/{U_p^+}^q,\mu_q)\simeq U_p^+/{U_p^+}^q\simeq \F_q[G^+]/(s(G^+))\simeq {\euf
A}.$$
De plus, ils sont des $\F_q[G^+]$-modules libres de rang 1. En outre, $\F_q[G^+]=(s(G^+))\oplus {\euf
A}$. }

{\bf Preuve}

On a vu au Lemme 8 du Chapitre 9 (relation (12)) que ${\rm
Ann}_{\F_q[G^+]}(U_p^+/{U_p^+}^q)=(s(G^+))$ (ˆ l'Žpoque, on avait quotientŽ par $\{\pm 1\}$,
mais a ne change rien). Si on munit ${\rm Hom}(U_p^+/{U_p^+}^q,\mu_q)$ de la structure de
$\F_q[G^+]$-module duale, alors on a aussi que ${\rm Ann}_{\F_q[G^+]}({\rm
Hom}(U_p^+/{U_p^+}^q,\mu_q))=(s(G^+))$. En effet, $\theta\in{\rm Ann}_{\F_q[G^+]}({\rm
Hom}(U_p^+/{U_p^+}^q,\mu_q))$ si et seulement si $\varphi^\theta=1$ pour tout
$\varphi\in{\rm Hom}(U_p^+/{U_p^+}^q,\mu_q)$ si et seulement si $\varphi(x^{\hat\theta})=1$
pour tout $\varphi$ et tout $x\in U_p^+/{U_p^+}^q$ (car $\varphi^\theta(x)=\varphi^{\sum_\sigma
a_\sigma\sigma}(x)=\prod_\sigma(\varphi^\sigma)^{a_\sigma}(x)=\prod_\sigma(\varphi^\sigma(x))^{a_\sigma}=
\prod_\sigma (\varphi(x^{\sigma^{-1}}))^{a_\sigma}=\varphi(x^{\hat\theta})$); si et seulement si
$x^{\hat\theta} =1$ pour tout  $x\in U_p^+/{U_p^+}^q\iff \hat\theta\in {\rm
Ann}(U_p^+/{U_p^+}^q)=(s(G^+))\iff\theta \in (s(G^+))$ (la dernire Žquivalence vient du fait
que dans l'idŽal $(s(G^+))$, tous les $a_\sigma$ sont Žgaux; on avait dŽjˆ montrŽ ce fait au
lemme 8 du Chapitre 9, mais on redit l'argument, car c'Žtait un peu cachŽ~: l'idŽal
$s(G^+)\F_q[G^+]=s(G^+)\F_q$, car si $\alpha=\sum_{\sigma\in G^+}a_\sigma\sigma\in\F_q[G^+]$,
alors

$$\alpha\cdot s(G^+)=\sum_{\matrix{\scriptstyle\sigma\in G^+\cr\scriptstyle \tau\in
G^+}}a_\sigma\cdot\sigma\tau\buildrel\mu=\sigma\tau\over =\sum_{\mu\in G^+}\cdot\big (\sum_{\tau\in G^+}
a_{\mu\tau^{-1}}\big )\cdot\mu=k\cdot s(G^+)\eqno{(i)}$$

o $k=\sum_{\sigma\in G^+}a_\sigma$.

D'autre part, cette mme preuve montre que 
$$q^{p-3\over 2}=\left |\F_q[G^+]/(s(G^+))\right |=\left
|U_p^+/{U_p^+}^q\right |=\left |{\rm Hom}(U_p^+/{U_p^+}^q,\mu_q)\right |.$$

On en dŽduit les deux premiers isomorphismes en vertu du Lemme 7 b) et c) du Chapitre 9.

Reste ˆ voir le dernier isomorphisme. ConsidŽrons l'application

$$\eqalign{f\ :\ \F_q[G^+]&\lra (s(G^+))\cr \sum_{\sigma\in G^+}a_\sigma\sigma&\longmapsto
\left (\sum_{\sigma\in G^+}a_\sigma\right)\cdot s(G^+).\cr}$$

Par $(i)$, voit que $f(x)=x\cdot s(G^+)$, pour tout $x\in\F_q[G^+]$. Ainsi, pour tout
$x,y\in \F_q[G^+]$, on a $f(x+y)=f(x)+f(y)$ et $f(x\cdot y)=x\cdot f(y)$. Ainsi, 
$f$ est un homomorphisme de $\F_q[G^+]$-module. Toujours par $(i)$, elle est surjective et son noyau
est ${\euf A}$. On se souvient que $\F_q[G^+]$ est un anneau semi-simple. Donc, on a
que $\F_q[G^+]=(s(G^+))\oplus {\euf A}$, et ainsi, $\F_q[G^+]/(s(G^+))\simeq{\euf A}$.\qed

\bigskip

{\bf Lemme 8}

{\sl Si $H$ est un sous-groupe abŽlien normal d'un groupe $G$, alors le groupe $Q=G/H$ agit
canoniquement sur $H$ par $\alpha\in Q$ et $h\in H$~:

$$h^\alpha=ghg^{-1}$$

o $g$ est n'importe quel reprŽsentant de $\alpha$.

}

{\bf Preuve}

C'est un vŽrification facile~: si $g'=gh_0$ avec $h_0\in H$, alors

$$g'hg'^{-1}=gh_0hh_0^{-1}g^{-1}=ghg^{-1}$$

car $H$ est abŽlien.\qed
 
\bigskip

{\bf Lemme 9}

{\sl Sous les mmes hypothses que la dŽfinition prŽcŽdente, on a que

\art{a)}L'isomorphisme de Kummer entre $G_0$ et ${\rm Hom}(U_p^+/{U_p^+}^q,\mu_q)$ est en fait
un isomorphisme de $\F_q[G^+]$-modules.

\art{b)}$\Delta$ agit sur $G_0$ par $\F_q[G^+]$-automorphisme.

\art{c)}Si $g\in G_0$, la classe de conjugaison dans ${\cal G}:={\rm Gal}(F(\zeta_q,\root q\of
{U_p^+})/F)$ est l'orbite de $g$ sous l'action de $\Delta={\cal G}/G_0$.

}

{\bf Preuve}

Tout d'abord, l'extension $F(\zeta_q,\root q\of {U_p^+})/\Q$ est une extension galoisienne. En
effet, les conjuguŽs de $\zeta_q$ et de $\zeta_p^+$ sont dans $F(\zeta_q,\root q\of {U_p^+})$. 
Concentrons-nous sur les ŽlŽments de $\root q\of {U_p^+}$. Soit $u\in U_p^+$ et $\root q\of u$
une racine $q$-ime de $u$. Pour montrer que l'extension est galoisienne, il suffit de trouver un
polyn™me annulateur de $\root q\of u$ qui soit dans $\Q[X]$ et dont toutes les racines soient
dans $F(\zeta_q,\root q\of {U_p^+})$. Soit $\sigma\in G^+$. Puisque $\sigma(u)\in U_p^+$, on a
toute $\root q\of{\sigma(u)}\in F(\zeta_q,\root q\of {U_p^+})$. Enfin, le polyn™me
$\prod_{\sigma\in G^+}(X-\sigma(u))\in \Q[X]$ et donc $\prod_{\sigma\in G^+}(X^q-\sigma(u))\in
\Q[X]$ et toutes les racines de ce polyn™mes sont dans $F(\zeta_q,\root q\of {U_p^+})$. Posons
${\cal G}_0={\rm Gal}(F(\zeta_q,\root q\of {U_p^+})/\Q)$ le groupe de Galois de cette extension. Par le
Lemme 8, en posant $G={\cal G}_0$ et $H=G_0={\rm Gal}(F(\zeta_q,\root q\of {U_p^+})/F(\zeta_q))$, on a
que
$\Delta\times G^+\simeq {\cal G}_0/G_0$ agit canoniquement sur $G_0$. En particulier, $G^+$ et $\Delta$
agissent sur $G_0$. 

Voyons dŽjˆ comment $\Delta={\rm Gal}(F(\zeta_q)/F)$ agit sur $G_0$~: soit $\varphi\in G_0$ et
$h\in\Delta$. On cherche $\varphi^h$, ou plut™t si $\check\varphi\in {\rm
Hom}(U_p^+/{U_p^+}^q,\mu_q)$ est l'homomorphisme correspondant ˆ $\varphi$ par l'isomorphisme
de Kummer, on cherche l'ŽlŽment qui correspond ˆ $(\varphi^h)^{\check{}}$. Par dŽfinition de
la  ÒKummer pairing", soit $[u]\in U_p^+/{U_p^+}^q$ et $\root q\of u$ une racine $q$-ime de $u$, on a
$\varphi^h(\root q\of u)=(\varphi^h)^{\check{}}([u])\cdot\root q\of u$. Notons cette ŽgalitŽ $(i)$.
Or, par dŽfinition de l'action canonique du Lemme~7, on a $\varphi^h(\root q\of u)=\tilde
h\varphi\tilde h^{-1}(\root q\of u)$ o $\tilde h$ est un automorphisme de $F(\zeta_q,\root q\of
{U_p^+})$ qui prolonge $h$. Or, $\tilde h^{-1}(\root q\of u)$ est une racine $q$-ime de
$h^{-1}(u)=u$ (car $h$ vaut l'identitŽ sur $F$). Donc, par dŽfinition de l'isomorphisme de Kummer,
on a $\varphi(\tilde h^{-1}(\root q\of u))=\check\varphi([u])\cdot \tilde h^{-1}(\root q\of u)$.
Donc, $\varphi^h(\root q\of u)=\tilde h(\check\varphi([u])\cdot\tilde h^{-1}(\root q\of
u))=h(\tilde\varphi (u))\cdot\root q\of u$. Ceci combinŽ avec $(i)$, en simplifiant par $\root
q\of u$, donne $(\varphi^h)^{\check{}}([u])=h(\check\varphi([u])).$ Pour tre plus prŽcis, si $h=h_a$
est l'automorphisme de $F(\zeta_q)$ tel que
$h_a(\zeta_p^+)=\zeta_p^+$ et $h_a(\zeta_q)=\zeta_q^a$ ($q\not \hskip-0.4pt |\   a$), on a
$$(\varphi^{h_a})^{\check{}}([u])=\check\varphi ([u])^a.\eqno{(ii)}$$

En particulier, si $\iota$ est la conjugaison
complexe dans $F(\zeta_q)$, alors $(\varphi^{\iota})^{\check{}}([u])=\check\varphi ([u])^{-1}$, pour
tout $[u]\in U_p^+/{U_p^+}^q$. Donc,

$$(\varphi^{\iota})^{\check{}}=(\check\varphi)^{-1}.\eqno{(13)}$$

Regardons maintenant l'action de $g\in G^+$. Soit $\varphi$, $\check\varphi$ et 
$(\varphi^g)^{\check{}}$, comme avant.  Soit $\tilde g$ un prolongement de $g$ en un
automorphisme de $F(\zeta_q,\root q\of{U_p^+})$. On a comme avant $\varphi^g(\root q\of
u)=(\varphi^g)^{\check{}}([u])\cdot\root q\of u$ et $\varphi^g(\root q\of u)=\tilde g\varphi\tilde
g^{-1}(\root q\of u)$. Lˆ, une petite diffŽrence~: $\tilde g^{-1}(\root q\of u)$ et une racine
$q$-ime de $g^{-1}(u)$. Donc, $\varphi(\tilde g g^{-1}(\root q\of
u))=\check\varphi([g^{-1}(u)])\cdot\tilde g^{-1}(\root q\of u)$. Ainsi, $\varphi^g(\root q\of
u)=\tilde g(\check\varphi([g^{-1}(u)])\cdot\tilde g^{-1}(\root q\of
u))=g(\check\varphi([g^{-1}(u)]))\cdot\root q\of u=\check\varphi(g^{-1}(u))\cdot\root q\of u$. La
dernire ŽgalitŽ vient du fait que $g(\zeta)=\zeta$ pour tout $\zeta\in\mu_q$. En simplifiant
par $\root q\of u$, on obtient
$(\varphi^g)^{\check{}}([u])=\check\varphi([g^{-1}(u)])=(\check\varphi)^g([u])$ (structure duale),
ceci pour tout $u\in U_p^+$ et parce que $[g^{-1}(u)]=g^{-1}[u]$. Donc,

$$(\varphi^g)^{\check{}}=(\check\varphi)^g.\eqno{(iii)}$$
Cela montre que l'isomorphisme de Kummer $G_0\lra {\rm Hom}(U_p^+/{U_p^+}^q,\mu_q)$ est 
un isomorphisme de $\F_q[G^+]$-modules, donc a) est prouvŽ.

Pour prouver b), il faut voir que si $\varphi\in G_0$, $h=h_a\in\Delta$ et $g\in G^+$, alors
$(\varphi^g)^h=(\varphi^h)^g$. Soit $u\in U_p^+$, on a~:

$$((\varphi^g)^h)^{\check{}}(u)\buildrel (ii)\over =(\varphi^g)^{\check{}}(u^a)\buildrel
(iii)\over =
(\check\varphi)^g(u^a)=\check\varphi(g^{-1}(u^a))=\check\varphi(g^{-1}(u))^a\buildrel
(ii)\over =(\varphi^h)^{\check{}}(g^{-1}(u))\buildrel (iii)\over =
((\varphi^h)^g)^{\check{}}(u).$$ 

Et on conclut pour b), car l'application $\check{}$ de Kummer est un isomorphisme de $\F[G^+]$-modules.

Pour la partie c), c'est ˆ peu prs Žvident, car si $h\in {\cal G}={\rm Gal}(F(\zeta_q,\root q\of
{U_p^+})/F)$, par le Lemme 8, l'action canonique de $\Delta={\cal G}/G_0$ est prŽcisŽment
$hgh^{-1}=g^{h_0}$ o $h_0=h|_{F(\zeta_q)}\in\Delta$.\qed
\bigskip
Maintenant nous allons faire quelques rappels sur l'automorphisme de Frobenius. On en avait dŽjˆ
parlŽ au Chapitre 5, mais uniquement pour le corps $\Q(\zeta_m)$. 

\bigskip

{\bf DŽfinitions-ThŽormes Òl'automorphisme de Frobenius"}

{\sl

Soit $L/K$ une extension galoisienne de corps de nombres de degrŽ $n$ de groupe de Galois $G$. Soit
$\P$ un idŽal premier de $O_K$ et $\gP$ un idŽal premier de $O_L$ au dessus de $\P$. On se souvient que
$Z(\gP|\P)=\{\sigma\in G\mid \sigma(\gP)=\gP\}$ et que si $\sigma\in G$, on a
$Z(\sigma(\gP)|\P)=\sigma Z(\gP|\P)\sigma^{-1}$, donc si $G$ est abŽlien, ce groupe est toujours le
mme et ne dŽpend que de $L$ et $\P$, on le note alors $Z(L/\P)$. On a
$\P O_L=(\gP_1\cdots\gP_r)^e$, $\gP$ est l'un des $\gP_i$. On a $n=efr$, o $f=[O_L/\gP_i:O_K/\P]$
est le {\it degrŽ de $\gP$ sur $K$}. On a aussi $[G:Z(\gP|\P)]=r$ et $|Z(\gP|\P)|=ef$. Si $\sigma\in
Z(\gP|\P)$, il determine un
$\overline{\sigma}\in {\rm Gal}((O_L/\gP)/(O_K/\P))$ et l'application $\sigma\lra\overline{\sigma}$
est un homomorphisme surjectif de $Z(\gP|\P)$ sur ${\rm Gal}((O_L/\gP)/(O_K/\P))$. Son noyau se note
$I(\gP/\P)$ ou $I(\gP/K)$ et s'appelle le {\it groupe d'inertie} de $\gP/\P$ ou de $\gP/K$. On a aussi
$I(\sigma(\gP/\P)=\sigma I(\gP/\P)\sigma^{-1}$ pour tout $\sigma\in G$. On a donc $|I(\gP/\P)|=e$ et
$Z(\gP/\P)/I(\gP/\P)\simeq {\rm Gal}((O_L/\gP)/(O_K/\P))$, de cardinal $f$. Pour tout $\sigma\in
I(\gP/\P)$, on a $\sigma(x)\equiv x\pmod \gP$ pour tout $x\in O_L$. Supposons que $\P$ ne ramifie
pas, c'est-ˆ-dire $e=1$ et donc le groupe d'inertie est trivial et donc l'application $\sigma\lra\overline{\sigma}$
est un isomorphisme de $Z(\gP|\P)$ sur ${\rm Gal}((O_L/\gP)/(O_K/\P))$. Ce groupe de Galois est un
groupe cyclique avec un gŽnŽrateur privilŽgiŽ qui est l'application $\nu\mapsto \nu^{\N(\P)}$ pour
tout $\nu\in O_L/\gP$ et l'unique ŽlŽment de $Z(\gP|\P)$ qui correspond ˆ cet automorphisme s'appelle
{\it l'automorphisme de Frobenius de $\gP/\P$.} On le note ${\rm Frob}(\gP/\P)$. Il est caractŽrisŽ
comme l'ŽlŽment de $G$ qui satisfait~:

$${\rm Frob}(\gP/\P)(x)\equiv x^{\N(\P)}\pmod \gP\quad \hbox{ pour tout }x\in O_L.$$

On a aussi ${\rm Frob}(\sigma(\gP)/\P)=\sigma {\rm Frob}(\gP/\P)\sigma^{-1}$ et donc l'ensemble $\{
{\rm Frob}(\gP/\P)\mid
\gP|\P\}$, qu'on note ${\rm Fr}_{L/K}(\P)$ est une classe de conjugaison dans $G$. Si $G$ est
abŽlien, alors ${\rm Frob}(\gP|\P)$ ne dŽpend que de $\P$, on le notera ${\rm Frob}_{L/K}(\P)$, et
on a

$${\rm Frob}_{L/K}(\P)(x)\equiv x^{\N(\P)}\pmod {\P O_L}\quad \hbox{ pour tout }x\in O_L.$$

Finalement, si $K\subset M\subset L$ sont des corps de nombres et $\P\subset\P_0\subset \gP$ sont
des idŽaux premiers de $K$, $M$ et $L$ respectivement. Alors ${\rm Frob}(\P_0/\P)={\rm
Frob}(\gP/\P)|_{M}$.

} 

\smallskip
{\bf Preuve}

Tout cela se trouve dans [Nar, p. 180].\qed
\bigskip

Maintenant, on va Žnoncer un thŽorme important qui ne sera utilisŽ qu'une fois mais est crucial.
C'est le thŽorme de $\check {\rm C}$ebotarev.

\bigskip\goodbreak

{\bf DŽfinition}

Soit $K$ un corps de nombres et $A$ un ensemble d'idŽaux premiers de $O_K$. On dit que $A$ est {\it
rŽgulier} s'il existe $0\leq a\leq 1$ et $g_A$ une fonction holomorphe sur un voisinage ouvert de
$\{s\in \C\mid \Re (s)> 1\}\cup\{1\}$ telle que pour tout $s\in \C$ avec $\Re (s)>1$ on ait 

$$\sum_{\P\in A}\N(\P)^{-s}=a\log\big ({1\over s-1}\big)+g_A(s).$$

Le nombre $a$ s'appelle la {\it densitŽ de Dirichlet de $A$}. Si $A$ est l'ensemble de tous les
idŽaux premiers de $O_K$, $a=1$.

\bigskip\goodbreak

{\bf ThŽorme 10 ($\check {\bf C}$ebotarev, 1923)}

{\sl Soit $L/K$ une extension galoisienne de corps de nombres de groupe de Galois $G$. Soit $C$ une
classe de conjugaison de $G$. Alors l'ensemble 

$$A=\{\P\mid \P\hbox{ est un idŽal premier de $K$ non ramifiŽ dans $L$
 avec }{\rm Fr}_{L/K}(\P)=C\}$$

est rŽgulier et sa densitŽ de Dirichlet vaut ${|C|\over |G|}$.

}
\smallskip
{\bf Preuve}

[Lang2, Thm 10, p. 169].\qed
\bigskip

{\bf Remarque}

Les idŽaux premiers de degrŽ sur $\Q$ strictement supŽrieur ˆ 1 sont de densitŽ de Dirichlet nuls,
donc dans le thŽorme de $\check {\rm C}$ebotarev, on peut se restreindre aux $\P$ qui sont de
degrŽ 1 sur $\Q$.

\bigskip

{\bf Lemme 11}

{\sl Un nombre premier $l$ ne divisant pas $n\in\N$ se dŽcompose totalement dans $\Q(\zeta_n^+)$ si et
seulement si $$l\equiv \pm 1\pmod{n}.$$

}
{\bf Preuve}

Le nombre $l$ se dŽcompose totalement dans $\Q(\zeta_n)$ si et
seulement si $f=e=1$ si et seulement si l'automorphisme de Frobenius est l'identitŽ si et seulement
si $l\equiv 1\pmod n$ (car l'ordre de l'automorphisme de Frobenius vaut $f$). Or, le
Frobenius ${\rm Frob}_{\Q(\zeta_n^+)/\Q}(l)$ est la restriction ˆ $\Q(\zeta_n^+)$ de ${\rm
Frob}_{\Q(\zeta_n)/\Q}(l)$. De plus,
${\rm Gal}(\Q(\zeta_n)/\Q)=(\Z/n\Z)^*$, donc, ${\rm Gal}(\Q(\zeta_n^+)/\Q)=(\Z/n\Z)^*/<-1>$. On en
dŽduit donc que $l$ se dŽcompose totalement dans $\Q(\zeta_n^+)$ si et
seulement si $l\equiv \pm 1\pmod{n}.$\qed

\bigskip

Maintenant un autre monument de la thŽorie de nombres, liŽ au prŽcŽdent~: le corps de classe de
Hilbert.

\bigskip

{\bf DŽfinition}

Soit $L/K$ une extension abŽlienne de corps de nombres. Soit $S$ un ensemble fini d'idŽaux premiers 
de $K$ contenant ceux qui ramifient dans $L$. Posons $I_S$ le groupe abŽlien libre engendrŽ
par les idŽaux premiers qui ne sont pas dans $S$. Si $\P\not\in S$ est un idŽal premier, on dŽfinit
l'application $\P\longmapsto {\rm Frob}_{L/K}(\P)$ se prolonge en un homomorphisme de groupe 

$$\Phi_{L/K}^S\ :\ I_S\lra G$$

appelŽ l'{\it homomorphisme d'Artin}. Remarquons en passant que le thŽorme de $\check {\rm C}$ebotarev
montre que cette application est surjective (mais on n'en aura pas besoin). Toute la thŽorie du corps
de classe est en fait de dŽcrire le noyau de cette application. 
\bigskip

{\bf ThŽorme 12 (existence du corps de Hilbert)}

{\sl Soit $K$ un corps de nombres. Il existe une (unique) extension $H/K$ telle que

\art{a)}$H/K$ est une extension abŽlienne finie (de groupe disons $G$).

\art{b)}Aucun idŽal premier de $K$ ne ramifie dans $H$.

\art{c)}En posant $S=\emptyset$ et $I_\emptyset=:I_K$, l'homomorphisme $\Phi_{H/K}^\emptyset\ :\ I_K\lra
G$ a pour noyau 
 les idŽaux (fractionnaires) principaux et dŽfinit donc un isomorphisme $\Phi_{H/K}\, :\, {\cal
CL}_K\lra G $.

On appelle $H$ le {\it corps de Hilbert} de $K$.

}

\bigskip \goodbreak

{\bf Lemme 13 (ThŽorme de Hilbert 90)}

{\sl Soit $L/K$ une extension cyclique de corps (i.e. galoisienne de groupe de Galois cyclique). Soit
$\tau$ un gŽnŽrateur de ${\rm Gal}(L/K)$. Soit
$x\in L$ et $N=N_{L/K}$ la norme associŽe ˆ cette extension. Alors on a

$$N(x)=1\iff x={\tau (\alpha)\over\alpha}\hbox{ pour un certain }\alpha\in L^*$$

}

{\bf Preuve}

On peut trouver une preuve de ce thŽorme classique dans [Lang1, Thm 6.1, p.288].\qed

\bigskip
Nous aurons besoin (vers la fin de la preuve) de quelques notions rudimentaires sur le localisŽ d'un
anneau~: 
\bigskip

{\bf DŽfinitions-ThŽorme}

Si $A$ est un anneau intgre commutatif et $Q$ son corps des fractions. Soit $\P$ un idŽal premier de
$A$. On dŽfinit $A_\P:=\{{a\over b}\in Q\mid b\not\in\P\}$. On appelle $A_\P$ le localisŽ de $A$ en
$\P$. C'est un anneau local d'idŽal maximal qu'on note encore $\P$. Voir [Ati, pp. 36-43] pour
plus de dŽtails.
\bigskip\goodbreak

 {\bf Preuve de la deuxime
version du thŽorme de Thaine}

On rappelle les notations~:  $p,q$ sont des nombres premiers impairs distincts tels que $p\not\equiv
1\pmod q$. On a $F=\Q(\zeta_p^+)$, $G^+={\rm Gal}(\Q(\zeta_p^+)/\Q)$, $U_p^+=\Z[\zeta_p^+]^*$, $C_0^+$,
sont les unitŽs cyclotomiques de $\Q(\zeta_p^+)$, qui est comme nous le savons (Lemme 3) le sous-groupe de
$U_p^+$ engendrŽ par $-1$ les $\xi_a=\zeta_p^{1-a\over 2}\cdot{\zeta_p^a-1\over\zeta_p-1}$. Enfin,
${\cal CL}^{pl}$ est le groupe des classes d'idŽaux de $F$.

Soit $\theta\in {\rm Ann}_{\Z[G^+]}( U_p^+/C_0^+{U_p^+}^q)$ et ${\cal C}\in {\cal CL}^{pl}$. On doit
montrer que ${\cal C}^\theta\in {{\cal CL}^{pl}}^q$. Allons-y, dans la joie et la bonne humeur~!

Soit $H$ le corps de Hilbert de $F$. Posons $G_H={\rm Gal}(H/F)$. On considre aussi les corps
$H(\zeta_q)$ et $H(\zeta_q,\root q\of {U_p^+})$ vus comme extensions sur $F$. Tout d'abord, $H\cap
F(\zeta_q)=F$. En effet, si ${\cal Q}$ est un idŽal premier de $F$ au-dessus de $q$, alors ${\cal Q}$
ramifie totalement dans $F(\zeta_q)$, mais il ne ramifie pas dans $H$. Donc, s'il y avait eu un corps
intercalŽ, dans ce corps, ${\cal Q}$ ramifierait et ne ramifierait pas, ce qui est absurde !

Maintenant, plus dŽlicat~: $H(\zeta_q)\cap F( \zeta_q,\root q\of {U_p^+})=F(\zeta_q)$. Rappelons que
$\Delta={\rm Gal}(F(\zeta_q)/F)$ ($\simeq {\rm Gal}(\Q(\zeta_q)/\Q)$) agit sur $G_0={\rm Gal}(
F( \zeta_q,\root q\of {U_p^+})/F(\zeta_q))\simeq {\rm Hom}(U_p^+/{U_p^+}^q,\mu_q)$ (Lemme 9 b) et
a)). De plus, puisque $H\cap F(\zeta_q)=F$, alors l'extension $H(\zeta_q)/F$ est galoisienne de
groupe $G_H\times \Delta$ (Lemme 5 c)), donc est abŽlienne, car $\Delta$ et $G_H$ le sont. Par le Lemme 8
et le Lemme 5 a), $\Delta$ agit canoniquement sur $G_H={\rm Gal}(H/F)\simeq {\rm
Gal}(H(\zeta_q)/F(\zeta_q))$, ceci de manire triviale car $H(\zeta_q)/F$ est abŽlienne. Mais si on prend
dans $\Delta$ la conjugaison complexe $\iota$, et $\varphi\in G_0$, on a

$$(\varphi^{\iota})^{\check{}}\buildrel (13)\over=(\check\varphi)^{-1}\buildrel \check{}\ {\rm
homomorphisme}\over = (\varphi^{-1})^{\check{}}.\eqno{(i)}$$

On en dŽduit que $\varphi^{\iota}=\varphi^{-1}$. Donc, il n'y a pas de point fixe dans cette action,
car $|G_0|=q^{p-3\over 2}$ est impair. Supposons par l'absurde que $H(\zeta_q)\cap
F(\zeta_q,\root q\of {U_p^+})\ne F(\zeta_q)$, c'est-ˆ-dire qu'il existe un corps, disons $K$, intercalŽ.
Nous sommes donc dans la situation suivante~:

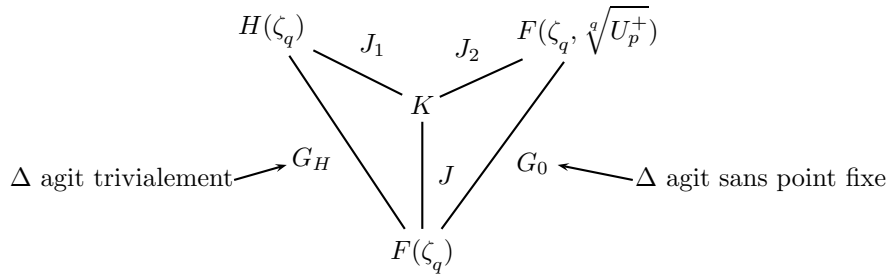
\vglue 0.5cm



\vskip 3cm

\rput(8,0){\rnode{Fzeta}{$F(\zeta_{q}) $}} 
\rput(8,2){\rnode{K}{$K$}}
\rput(6,3){\rnode{Hzeta}{$H(\zeta_{q})$}}
\rput(10.2,3){\rnode{FF}{$F(\zeta_{q},\root q\of {U_{p}^+})$}}
\ncline[nodesep=3pt]{Fzeta}{K}
\Bput{$J$} 
\ncline[nodesep=3pt]{Fzeta}{Hzeta}
\Aput{$G_{H} $} 
\ncline[nodesep=3pt]{Fzeta}{FF}
\Bput{$G_{0}$} 
\ncline[nodesep=3pt]{K}{Hzeta}
\Bput{$J_{1}$} 
\ncline[nodesep=2pt]{K}{FF}
\Aput{$J_{2}$} 
\rput (4,1){$\Delta$ agit trivialement}
\psline{->}(5.5,1)(6.2,1.2)
\rput (12.5,1){$\Delta$ agit sans point fixe}
\psline{->}(10.8,1)(9.8,1.2)
\vskip 1cm

\goodbreak
En effet, puisque $G_H$ est abŽlien, l'extension $K/F(\zeta_q)$ est aussi abŽlienne (tout sous-groupe
 d'un groupe abŽlien est normal), de groupe, disons $J\simeq G_H/J_1\simeq G_0/J_2$. Puisque
$\Delta$ agit trivialement sur $G_H$, il agit trivialement sur $G_H/J_1$. Mais, il agit sans point
fixe sur $G_0$. Puisqu'il agit par conjugaison sur $G_0$ (Lemme 9 c)) et que $J_2$ est un groupe
normal, alors l'action de $\Delta$ sur $J_2$ laisse fixe $J_2$ dans son ensemble. Donc, l'action de
$\Delta$ passe au quotient, car si $\varphi,\varphi_0\in G_0$ sont tels qu'il existe $\psi\in J_2$ avec
$\varphi_0=\varphi\psi$; soit $h\in\Delta$ et $h_0\in {\rm Gal}(F(\zeta_q,\root q\of {U_p^+})/F)$ tel que
$h_0|_{F(\zeta)}=h$, on a $\varphi_0^h=(\varphi\psi)^h=h_0\varphi\psi
h_0^{-1}=h_0\varphi h_0^{-1}\underbrace{h_0\psi h_0^{-1}}_{:=\psi'\in J_2}=\varphi^h\psi'$. Donc
$[\varphi_0^h]=[\varphi^h]=:[\varphi]^h$. Mais cette action n'est pas triviale~: si elle l'Žtait, on
aurait $[\varphi]=[\varphi]^\iota=[\varphi^{\iota}]\buildrel (i)\over=[\varphi^{-1}]$ pour tout
$\varphi\in G_0$. Donc $G_0^2\subset J_2$. Mais puisque l'ordre de $G_0$ est impair, on a
$G_0=G_0^2$ (thŽorme de Bezout), donc, on aurait $J_2=G_0$, ce qui est absurde. On a ainsi prouvŽ
que $H(\zeta_q)\cap F(\zeta_q,\root q\of {U_p^+})= F(\zeta_q)$. On en dŽduit que $H\cap
F(\zeta_q,\root q\of {U_p^+})= F$, car

$$H\cap F(\zeta_q,\root q\of {U_p^+})=H\cap (H(\zeta_q)\cap  F(\zeta_q,\root q\of {U_p^+}))=H\cap
F(\zeta_q)= F.$$

Donc, $H$ et $F(\zeta_q,\root q\of {U_p^+})$ sont linŽairement disjoints. Ainsi, l'extension
$H(\zeta_q,\root q\of {U_p^+})/F$ est galoisienne de groupe $G_H\times {\cal G}\supset G_H\times
G_0\simeq G_H\times {\rm Hom}(U_p^+/{U_p^+}^q,\mu_q)$, o ${\cal G}={\rm Gal}(F(\zeta_q,\root q\of
{U_p^+})/F)$.

Revenons ˆ la classe $\cal C$ fixŽe au dŽbut de notre preuve. L'homomorphisme d'Artin fait
correspondre ${\cal C}$ ˆ $\Phi_{H/F}({\cal C})\in G_H$. ConsidŽrons la classe de conjugaison rŽduite
ˆ $\Phi_{H/F}({\cal C})$. ConsidŽrons aussi ${\cal C}'$ une classe de conjugaison contenant un
gŽnŽrateur de ${\rm Hom}(U_p^+/{U_p^+}^q,\mu_q)\simeq G_0$. Par ce qu'on vient de voir, on a
$\{\Phi_{H/F}({\cal C})\}\times {\cal C}'$ est une classe de conjugaison de ${\rm Gal}(H(\zeta_q,\root
q\of {U_p^+})/F)$.  Le thŽorme de $\check {\rm C}$ebotarev nous apprend qu'il existe un (en fait une
infinitŽ) idŽal premier  de $F$,
$\lambda$, de degrŽ 1 sur $\Q$, tel que ${\rm Fr}_{H(\zeta_q,\root q\of
{U_p^+})/F}(\lambda)=\{\Phi_{H/F}({\cal C})\}\times {\cal C}'$. Cela signifie que

\art{a)}$\lambda$ est de degrŽ 1 sur $\Q$.

\art{b)}${\rm Frob}_{H/F}(\lambda)=\Phi_{H/F}({\cal C})$.

\art{c)}${\rm Fr}_{F(\zeta_q,\root q\of {U_p^+})/F}(\lambda)={\cal C'}$.

Soit $l\in \gfP$ le nombre premier sous $\lambda$ ($\lambda\cap\Z=l\Z$). La partie a) et le
lemme 11 montrent alors que $l\equiv \pm 1\pmod p$. La partie b) implique que $\lambda$ est dans la
classe ${\cal C}$. La partie c) dit plusieurs choses~: la premire est qu'il existe $\Lambda$ un
idŽal premier de $F(\zeta_q,\root q\of {U_p^+})$ au-dessus de $\lambda$ tel que ${\rm
Frob}(\Lambda/\lambda)$ est un gŽnŽrateur de $G_0$. Gr‰ce au Lemme 5 a), on sait que $G_0=\{h\in {\cal G}\mid
h|_{F(\zeta_q)}={\rm Id}_{F(\zeta_q)}\}$. Soit
$\Lambda_0=\Lambda\cap F(\zeta_q)$. On a alors
${\rm Frob}(\Lambda_0/\lambda)={\rm Frob}(\Lambda/\lambda)|_{F(\zeta_q)}={\rm Id}_{F(\zeta_q)}=1_\Delta$. Donc,
$\Lambda_0$ est de degrŽ 1 sur $F$ (car ${\rm Frob}(\Lambda_0/\lambda)$ engendre $Z(\Lambda_0/\lambda)$ qui
est de cardinal $f(\Lambda_0/\lambda)$). Cela implique que
$l$ est totalement dŽcomposŽ dans $\Q(\zeta_q)$ (car si $\cal L$ est un idŽal de $\Q(\zeta_q)$ au-dessus de
$l$, on a  $f({\cal L}/l)=f(\Lambda_0/\lambda)$)  et donc que
$l\equiv 1\pmod q$ (cf. remarque avant le Lemme 11). RŽsumons encore tout cela

$$\lambda\in {\cal C}\quad l\equiv\pm 1\pmod p \quad l\equiv 1\pmod q\quad {\rm
Frob}(\Lambda/\lambda)\hbox{ est un $\Z[G^+]$-gŽnŽrateur de $G_0$.}\eqno{(ii)}$$
\bigskip

Le corps de Hilbert et le thŽorme de $\check{\rm C}$ebotarev ont fait leur office, nous n'en
n'aurons dŽsormais plus besoin.

\bigskip\goodbreak

{\bf Sous-Lemme 14}

{\sl Sous les mmes notations et hypothses, posons encore ${\euf L}=\left (\Z[\zeta_p^+]/(l)\right
)^*$. Alors on a un homomorphisme injectif de $\Z[G^+]$-module~:

$$U_p^+/{U_p^+}^q\lhook\joinrel\lra {\euf L}/{\euf L}^q\simeq \F_q[G^+]$$

et l'image de $U_p^+/{U_p^+}^q$ dans $\F_q[G^+]$ est l'idŽal d'augmentation ${\euf A}$.

}
\goodbreak

{\bf Preuve du sous-lemme 14}

On a $(l)=l\cdot\Z[\zeta_p^+]=\prod_{\sigma\in G^+}\sigma(\lambda)$ et les
$\sigma(\lambda)$ sont tous disjoints, car $l$ se dŽcompose totalement dans $\Z[\zeta_p^+]$. Par le
thŽorme chinois, on a
$$ \Z[\zeta_p^+]/(l)\simeq \prod_{\sigma\in G^+}\Z[\zeta_p^+]/\sigma(\lambda)\buildrel
f=1\over\simeq \prod_{\sigma\in G^+} \F_l\simeq \F_l^{G^+}$$

o $\F_l^{G^+}$ est l'ensemble des applications de $G^+$ dans $\F_l$.

Le groupe $G^+$ agit lˆ-dessus par permutations des coordonnŽes~: tout $\sigma\in G^+$ induit un
isomorphisme de $\Z[\zeta_p^+]/\lambda$ sur $\Z[\zeta_p^+]/\sigma(\lambda)$ qui donne l'identitŽ si
on identifie ces quotients ˆ $\F_l$. Donc ce sont des $\Z[G^+]$-modules. Toujours canoniquement, on a
$\left (\Z[\zeta_p^+]/(l)\right)^*={\euf L}\simeq \prod_{\sigma\in G^+} \F_l^*\simeq {{\F_l}^*}^{G^+}$
et donc ${\euf L}/{\euf L}^q\simeq \prod_{\sigma\in G^+}\F_l^*/{\F_l^*}^q$. D'autre part, comme
groupes, $\F_l^*/{\F_l^*}^q\simeq \F_q$, le premier est notŽ multiplicativement et le second 
additivement, voyons pourquoi~: si $s$ est une racine primitive modulo $l$ et $\overline{s}$ est sa
classe modulo ${\F_l^*}^q$, alors on envoie $\overline{s}$ sur 1 modulo $q$. L'application est bien
dŽfinie, car $q$ divise $l-1$ $(cf. (ii))$, et est un isomorphisme. On a donc ${\euf L}/{\euf L}^q\simeq
\F_q^{G^+}$ D'autre part, $\F_q^{G^+}$ et $\F_q[G^+]$ sont isomorphes comme $\F_q[G^+]$-modules~:
l'application $x\mapsto \sum_{\sigma\in G^+}x(\sigma)\sigma$ en est un isomorphisme (attention,
l'action de $\F_q[G^+]$ sur $\F_q^{G^+}$ se fait comme suit~: $x^\sigma(\tau):=x(\sigma^{-1}\tau)$).

Soit  $\psi$ l'homomorphisme canonique (de $\Z[G^+]$-module) $U_p^+\lra {\euf L}/{\euf L}^q$. Il est
clair que ${U_p^+}^q\subset \ker(\psi)$. On va voir que c'est en fait Žgal. Soit $u\in\ker(\psi)$;  $u$
est donc une puissance $q$-ime dans ${\euf L}$. Cela veut dire que $u^{l-1\over q}\equiv 1\pmod
{l\cdot\Z[\zeta_p^+]}$, ou encore

$$\sigma^{-1}(u)^{l-1\over q}\equiv 1\pmod {l\cdot\Z[\zeta_p^+]}\quad\hbox{pour tout }\sigma\in
G^+.\eqno{(*)}$$

Posons $f={\rm Frob}(\Lambda/\lambda)$ et $\check f$ l'ŽlŽment de ${\rm
Hom}(U_p^+/{U_p^+}^q,\mu_q)$ associŽ ˆ $f$. On se souvient que $f(x)\equiv x^l\pmod\Lambda$ pour tout
$x\in O_{F(\zeta_q,\root q\of {U_p^+})}$ ($l=\N(\lambda)$) et, pour tout $[u']\in U_p^+/{U_p^+}^q$, on a
$\root q\of {u'}\cdot \check f([u'])=f(\root q\of {u'})$ o $\root q\of {u'}$ est n'importe quelle racine
$q$-ime de $u'$. Donc, $\check f([u'])\root q\of {u'}\equiv \left (\root q\of {u'}\right
)^l={u'}^{l-1\over q}\cdot \root q\of {u'} \pmod\Lambda$. Donc, $\check f([u'])\equiv {u'}^{l-1\over q}
\pmod\Lambda$. Revenons ˆ notre $u\in\ker(\psi)$. Par $(*)$, on a donc montrŽ que $\check
f([\sigma^{-1}(u)])\equiv 1\pmod\Lambda$, pour tout $\sigma\in G^+$. Or $\check f([\sigma^{-1}(u)])$ est une
racine $q$-ime de l'unitŽ. Le Lemme IMP du Chapitre 5 nous dit que les racines $q$-ime de l'unitŽ
sont distinctes modulo $\Lambda$, ainsi, $\check f([\sigma^{-1}(u)])=1$ pour tout $\sigma\in G^+$, en d'autre
terme, $(\check f)^\sigma([u])=1$ pour tout $\sigma\in G^+$. Donc, $(\check
f)^{\sum_{\sigma}a_\sigma\sigma}([u])=1$ pour tout $\sum_{\sigma}a_\sigma\sigma\in \Z[G^+]$. Or, par
$(ii)$, on sait que $\check f$ est un gŽnŽrateur de ${\rm Hom}(U_p^+/{U_p^+}^q,\mu_q)$, donc
$\alpha([u])=1$ pour tout $\alpha\in {\rm Hom}(U_p^+/{U_p^+}^q,\mu_q)$. Donc, par dualitŽ, $[u]=1$,
c'est-ˆ-dire $u\in {U_p^+}^q$. On a donc prouvŽ que 

$$U_p^+/{U_p^+}^q\lhook\joinrel\lra {\euf L}/{\euf L}^q\simeq \F_q[G^+].$$

Pour montrer que l'image de $U_p^+/{U_p^+}^q$ est l'idŽal d'augmentation, il suffit de montrer que 
$${\rm Ann}_{\F_q[G^+]}(U^+/{U_p^+}^q)=s(G^+)\cdot\F_q[G^+],$$

car dans ce cas, $U^+/{U_p^+}^q\buildrel {\rm Prop\, 6,\ App\, 1}\over \simeq
\F_q[G^+]/(s(G^+))\buildrel {\rm Lemme\ 7}\over\simeq {\euf A}$.  Mais on l'a dŽjˆ prouvŽ au lemme 8
du Chapitre 9 (relation (12)).\qed
\medskip

Revenons ˆ notre preuve. Nous savons que $U^+/{U_p^+}^q$ est cyclique (Proposition 6, Appendice 1 et le
sous-lemme ci-dessus). Soit $u\in U_p^+$ tel que $\overline{u}$ soit un gŽnŽrateur de
$U_p^+/{U_p^+}^q$. Soit le $\theta$ du dŽbut de notre preuve. Puisque $\theta\in {\rm Ann}_{\Z[G^+]}(
U_p^+/C_0^+{U_p^+}^q)$, alors on a
$u^\theta\in C_0^+{U_p^+}^q$, c'est-ˆ-dire $u^\theta=c\cdot v^q$, avec $c\in C_0^+$ et $v\in
{U_p^+}^q$. On peut supposer que $c=c_0^2$ est un carrŽ. En effet~: par le thŽorme de Bezout, il existe
des entiers $k$ et $k'$ tels que $2k+qk'=1$, donc $u^\theta=({c^k})^2\cdot (c^{k'}\cdot v)^q$. D'autre
part, puisque $c_0\in C_0^+$, on a que $C_0$ est un produit de $\xi_a=\zeta_p^{1-a\over 2}\cdot
{\zeta_p^a-1\over \zeta_p-1}$ avec $a\in\N\setminus\{0,1\}$ et $a\leq p-1$ (Lemme 3). Posons, pour de
tels $a$~:

$$\varepsilon_a=\zeta_p^{1-a\over 2}\cdot {\zeta_p^a-\zeta_l\over
\zeta_p-\zeta_l}\in\Q(\zeta_p,\zeta_l)=\Q(\zeta_{pl})\hbox{ et }\varepsilon'_a=\zeta_p^{1-a\over 2}\cdot
{\zeta_p^a-\zeta_l^{-1}\over
\zeta_p-\zeta_l^{-1}}\in\Q(\zeta_{pl}).$$

Posons encore $\eta_a=\varepsilon_a\cdot\varepsilon'_a $.
\bigskip\goodbreak
{\bf Sous-Lemme 15}

{\sl Sous les mmes notations et hypothses, on a~:

\art{a)}$\varepsilon_a$ et $\varepsilon'_a$ sont des unitŽs cyclotomiques de $\Q(\zeta_{pl})$

\art{b)}$\varepsilon_a\equiv\varepsilon'_a\equiv\xi_a\pmod {(\zeta_l-1)\Z[\zeta_{pl}]}$.

\art{c)}$N_{\Q(\zeta_{pl})/\Q(\zeta_{p})}(\varepsilon_a)=
N_{\Q(\zeta_{pl})/\Q(\zeta_{p})}(\varepsilon'_a)=1$.

\art{d)}$\eta_a$ est une unitŽ de $F(\zeta_l)$

\art{e)}$N_{F(\zeta_{l})/F}(\eta_a)=1$

\art{f)}$\eta_a\equiv \xi_a^2\pmod {(\zeta_l-1)O_{F(\zeta_{l})}}$.

}

{\bf Preuve du Sous-Lemme 15}

\art{a)}On a
${\zeta_p^a-\zeta_l\over\zeta_p-\zeta_l}={\zeta_l^{-1}\zeta_p^a-1\over\zeta_l^{-1}\zeta_p-1}=
{(\zeta_l^{-1}\zeta_p)^b-1\over\zeta_l^{-1}\zeta_p-1}$, o $b\in \Z$ tel que $b\equiv a\pmod p$ et
$b\equiv 1\pmod l$, est une unitŽ cyclotomique. La preuve pour $\varepsilon'_a$ est semblable.

\art{b)}On pourrait na•vement dire que $\zeta_l\equiv 1\pmod {(\zeta_l-1)\Z[\zeta_{pl}]}$, et le tour
est jouŽ. Mais c'est un argument fallacieux! On vŽrifie que
$\xi_a-\varepsilon_a=\zeta_p^{1-a\over 2}\cdot \left ({\zeta_p^a-1\over
\zeta_p-1}-{\zeta_p^a-\zeta_l\over\zeta_p-\zeta_l}\right )=\zeta_p^{1-a\over 2}\cdot \left
({\zeta_p^a-\zeta_p\over
\zeta_p-1}\right )\cdot\left ({1-\zeta_l\over\zeta_p-\zeta_l}\right )$. D'une part, $\zeta_p^{1-a\over 2}\cdot \left
({\zeta_p^a-\zeta_p\over \zeta_p-1}\right )$ est clairement une unitŽ (cf. Lemme 2 a) du Chapitre 6).
D'autre part,
$${1-\zeta_l\over\zeta_p-\zeta_l}={\zeta_l^{-1}-1\over
\zeta_l^{-1}\zeta_p-1}\buildrel (*)\over={(\zeta_l^{-1}\zeta_p)^b-1\over
(\zeta_l^{-1}\zeta_p)-1}=\sum_{k=0}^{b-1}(\zeta_l^{-1}\zeta_p)^k\buildrel (**)\over \equiv
\sum_{k=0}^{b-1}\zeta_p^k={\zeta_p^b-1\over\zeta_p-1}=0\pmod{(\zeta_l-1)\Z[\zeta_{pl}]}.$$
L'ŽgalitŽ $(*)$ est vraie en posant $b\equiv 0\pmod p$, $b\equiv 1\pmod l$ et $b>0$. L'Žquivalence
$(**)$ vient du fait que $\zeta_l^{-1}\equiv 1\pmod{(\zeta_l-1)\Z[\zeta_{pl}]}$, car l'idŽal engendrŽ
par $\zeta_l-1$ est le mme que celui engendrŽ par $\zeta_l^{-1}-1$. Cela prouve la partie b).

\art{c)}On se souvient que $\prod_{k=0}^{l-1}(X-\zeta_l^k)=X^l-1$. Posons
$N=N_{\Q(\zeta_{pl})/\Q(\zeta_{p})}$. On rappelle d'autre part que $l\equiv\pm 1\pmod p$ (relation
$(ii)$). On a

$$\xi_a\cdot N(\varepsilon_a)\buildrel (***)\over =\zeta_p^{l\cdot \left ({1-a\over 2}\right
)}\cdot\prod_{k=0}^{l-1}{\zeta_p^a-\zeta_l^k\over\zeta_p-\zeta_l^k}=\zeta_p^{l\cdot \left ({1-a\over 2}\right
)}\cdot{\zeta_p^{la}-1\over\zeta_p^l-1}=\cases{\xi_a&si $l\equiv 1\pmod p$\cr \overline{\xi_a}=\xi_a &
($\xi_a$ est rŽel) si $l\equiv -1\pmod p$\cr}$$
\art{}Dans l'ŽgalitŽ $(***)$, le $l$ de $\zeta_p^{l\cdot \left ({1-a\over 2}\right )}$ vient du fait que
$l=l-1+1$, le
$l-1$ vient de la norme et le $1$ vient du $\xi_a$, d'autre part le produit qui part de $k=0$ vient aussi du
$\xi_a$. On en dŽduit que $\xi_a \cdot N(\varepsilon_a))=\xi_a$. Donc $N(\varepsilon_a))=1$. La preuve
de $N(\varepsilon'_a))=1$ est identique.

Les parties d), e) et f) se dŽmontrent simultanŽment. Le groupe de Galois ${\rm
Gal}(\Q(\zeta_p,\zeta_l)/F(\zeta_l))$ est d'ordre 2, engendrŽ par l'automorphisme $\sigma_0$ de
$\Q(\zeta_p,\zeta_l)$ qui change $\zeta_p$ en $\overline{\zeta_p}=\zeta_p^{-1}$ et laisse fixe $\zeta_l$. On va
montrer que $\sigma_0(\varepsilon_a)=\varepsilon'_a$~:

$$\sigma_0(\varepsilon_a)=\sigma_0\left (\zeta_p^{1-a\over 2}\cdot {\zeta_p^a-\zeta_l\over
\zeta_p-\zeta_l}\right )=
\zeta_p^{a-1\over 2}\cdot {\zeta_p^{-a}-\zeta_l\over \zeta_p^{-1}-\zeta_l}=\zeta_p^{1-a\over 2}\cdot
{1-\zeta_p^a\zeta_l\over 1-\zeta_p\zeta_l}=\zeta_p^{1-a\over 2}\cdot{\zeta_l^{-1}-\zeta_p^a\over
\zeta_l^{-1}-\zeta_p}=\varepsilon'_a.$$

On en dŽduit que $\eta_a=\varepsilon_a\varepsilon'_a\in F(\zeta_l)$ est une unitŽ, que 
$N_{F(\zeta_{l})/F}(\eta_a)=1$, car $N_{\Q(\zeta_p,\zeta_l)/\Q(\zeta_p)}$ restreint ˆ $F(\zeta_l)$ est 
$N_{F(\zeta_{l})/F}$ (ces deux extensions ont mme groupe de Galois) et que $\eta_a\equiv
\xi_a^2\pmod {(\zeta_l-1)O_{F(\zeta_{l})}}$, par la partie b).\qed
\bigskip
Revenons ˆ nouveau ˆ notre preuve. Nous avons que $u^\theta=c\cdot v^q$ avec
$c=c_0^2=\prod_{i=1}^r\xi_{a_i}^2$. Posons $\varepsilon=\prod_{i=1}^r\eta_{a_i}$. Par le sous-lemme
prŽcŽdent, $\varepsilon$ est une unitŽ de $F(\zeta_l)$, $N(\varepsilon)=1$ et
$c\equiv\varepsilon\pmod{(\zeta_l-1)O_{F(\zeta_{l})}}$. Notons $\tau$ le gŽnŽrateur de ${\rm
Gal}(F(\zeta_l)/F)\simeq {\rm Gal}(\Q(\zeta_l)/\Q)$ qui envoie $\zeta_l$ sur $\zeta_l^s$ ou $s$ est la
racine primitive vue lors de la preuve du Sous-Lemme 14. Par le thŽorme de Hilbert 90, il existe
$\alpha\in F(\zeta_l)^*$ tel que 
$$\varepsilon={\tau(\alpha)\over\alpha}.$$

On en dŽduit en particulier que l'idŽal fractionnaire $(\alpha)$ est invariant par $\tau$. Puisque
$l$ ramifie totalement dans $\Q(\zeta_l)$, l'idŽal $\lambda$ ramifie totalement dans 
$F(\zeta_l)$. Notons $\euf l$ l'unique idŽal premier de $F(\zeta_l)$ au-dessus de $\lambda$. De mme, pour
chaque $\sigma\in G^+$, l'unique idŽal premier au-dessus de $\sigma(\lambda)$ est $\sigma({\euf l})$ (on
identifie Žvidemment $G^+={\rm Gal}(F/\Q)\simeq{\rm Gal}(F(\zeta_l)/\Q(\zeta_l))$). Les $\sigma({\euf
l})$ sont aussi $\tau$-invariants, car il y a ramification totale. Ecrivons

$$\alpha\cdot O_{F(\zeta_l)}=(\alpha)={\euf J}\cdot\prod_{\sigma\in G^+}\sigma({\euf
l})^{r_\sigma},\quad (r_\sigma\in\Z,\ {\euf J}\ \hbox{idŽal fractionnaire premier ˆ $l$}).\eqno{(iii)}$$

Puisque $(\alpha)$ et les $\sigma({\euf l})$ sont invariants par $\tau$, alors $\euf J$ l'est aussi.
Puisque $\euf J$ est $\tau$-invariant et qu'il n'a que des diviseurs premiers non ramifiŽ sur $F$,
$\euf J$ doit provenir d'un idŽal ${\euf J}_0$, c'est-ˆ-dire ${\euf J}_0\cdot O_{F(\zeta_l)}={\euf J}$.
On a donc, ${\euf J}_0^{l-1}\cdot O_{F(\zeta_l)}={\euf J}^{l-1}=N_{F(\zeta_l)/F}({\euf J})\cdot
O_{F(\zeta_l)}$. Donc, $N_{F(\zeta_l)/F}({\euf J})={\euf J}_0^{l-1}$ (cf. Chapitre 5, Rappels sur les corps de
nombres et la thŽorie de Galois). D'autre part, puisque
$\sigma(\lambda)$ ramifie totalement dans $F(\zeta_l)$, alors $N_{F(\zeta_l)/F}(\sigma({\euf
l}))=\sigma(\lambda)$. En prenant la norme de l'Žquation $(iii)$, on trouve~:

$$N_{F(\zeta_l)/F}(\alpha)\cdot O_F={\euf J}_0^{l-1}\cdot\prod_{\sigma\in
G^+}\sigma(\lambda)^{r_\sigma}={\euf J}_0^{l-1}\cdot\lambda^{\Sigma r_\sigma\sigma}.$$

Cette ŽgalitŽ montre que l'ŽlŽment $\sum_{\sigma\in G^+}r_\sigma\sigma\in \Z[G^+]$ est tel que 
$${\cal C}^{\Sigma r_\sigma\sigma}\in {{\cal CL}^{pl}}^q,\eqno{(iv)}$$
ceci parce que $\lambda\in {\cal C}$,
 donc ${\cal C}$ est la classe dŽfinie par $\lambda$, $l\equiv 1\pmod q$ (cf. $(ii)$), donc ${\euf
J}_0^{l-1}\in {{\cal CL}^{pl}}^q$ et parce que $N_{F(\zeta_l)/F}(\alpha)\cdot O_F$ est principal. Donc,
il faut jeter un pont entre $\theta$ et $\sum_{\sigma\in G^+}r_\sigma\sigma$.

On a que l'idŽal $(\zeta_l-1)= \prod_{\sigma\in G^+}\sigma({\euf l})$. Donc, l'ŽlŽment
${\alpha\over (\zeta_l-1)^{r_\sigma}}$ est premier ˆ $\sigma({\euf l})$, il est donc inversible dans
l'anneau localisŽ $(O_{F(\zeta_l)})_{\sigma({\euf l})}$. Le groupe de Galois ${\rm
Gal}(F(\zeta_l)/F)=<\tau>$ est le groupe d'inertie de $\sigma({\euf l})$ sur $\sigma(\lambda)$ car la
ramification est totale. On a donc, $\tau(x)\equiv x\pmod {\sigma({\euf l})}$ pour tout $x\in
O_{F(\zeta_l)}$ (cf. DŽfinitions-ThŽorme sur l'automorphisme de Frobenius), et donc aussi
$\tau(x)\equiv x\pmod {\sigma({\euf l})}$ pour tout $x\in (O_{F(\zeta_l)})_{\sigma({\euf l})}$. En
particulier,

$${\alpha\over (\zeta_l-1)^{r_\sigma}}\equiv \tau\left ({\alpha\over
(\zeta_l-1)^{r_\sigma}}\right )={\varepsilon \cdot\alpha\over (\zeta_l^s-1)^{r_\sigma}}\pmod
{\sigma({\euf l})\cdot (O_{F(\zeta_l)})_{\sigma({\euf l})}}.$$ 

On notera dorŽnavant $\pmod{\sigma({\euf l})'}$ au lieu de $\pmod
{\sigma({\euf l})\cdot (O_{F(\zeta_l)})_{\sigma({\euf l})}}$. Or, ${\varepsilon \cdot\alpha\over
(\zeta_l^s-1)^{r_\sigma}}={\varepsilon
\cdot\alpha\over (\zeta_l-1)^{r_\sigma}}\cdot \left ( {\zeta_l^s-1\over \zeta_l-1}\right )^{-r_\sigma}$
et
${\zeta_l^s-1\over \zeta_l-1}=1+\zeta_l+\cdots +\zeta_l^{s-1}\equiv s\pmod {\zeta_l-1}$, donc {\it a
fortiori} modulo $\sigma({\euf l})'$. Donc $\left ( {\zeta_l^s-1\over \zeta_l-1}\right
)^{-r_\sigma}\equiv s^{-r_\sigma}\pmod {\sigma({\euf l})'}$. Ainsi, ${\alpha\over
(\zeta_l-1)^{r_\sigma}}\equiv {\varepsilon \cdot\alpha\over (\zeta_l^s-1)^{r_\sigma}}\equiv
{\varepsilon\alpha\over (\zeta_l-1)^{r_\sigma}}\cdot s^{-r_\sigma}\pmod{\sigma({\euf
l})'}$. En simplifiant, par ${\alpha\over (\zeta_l-1)^{r_\sigma}}$, qui est, rappelons-le, inversible, on
trouve $\varepsilon\equiv s^{r_\sigma}\pmod{\sigma({\euf l})'}$. Or, on se souvient que $c\equiv
\varepsilon\pmod{(\zeta_l-1)}$, on a donc {\it a fortiori} $c\equiv \varepsilon\pmod{\sigma({\euf
l})'}$. Donc, on a $c\equiv s^{r_\sigma}\pmod{\sigma({\euf l})'}$. Comme $c$ et $s^{r_\sigma}$ sont
dans $(\Z[\zeta_p^+])':= (\Z[\zeta_p^+])_{\sigma(\lambda)}$, on a

$$c\equiv s^{r_\sigma}\pmod{\sigma(\lambda)'},\eqno{(v)}$$ 

car $\sigma({\euf l})'\cap (\Z[\zeta_p^+])'=\sigma(\lambda)'$ qui est l'idŽal maximal de
$(\Z[\zeta_p^+])'$. Puisque $\Z[\zeta_p^+]'/\sigma(\lambda)'\simeq \Z[\zeta_p^+]/\sigma(\lambda)\simeq
\Z/l\Z=\F_l$, la congruence $(v)$ montre que l'image de $c$ dans ${\euf L}=\left
(\Z[\zeta_p^+]/(l)\right )^*\simeq \prod_{\sigma\in G^+}\left (\Z[\zeta_p^+]/\sigma(\lambda)\right
)^*\simeq (\F_l^*)^{G^+}$ est
$(s'^{r_\sigma})_{\sigma\in G^+}$, o $s'\in \F_l^*$ est la classe de $s$ modulo
$l$. Passant encore aux classes modulo les puissances $q$-imes, on obtient que l'image de $u^{\theta}$
qui est la mme que celle de celle de $c$ modulo $U_p^+$ dans ${\euf L}/{\euf L}^q\simeq
(\F_l^*/{\F_l^*}^q)^{G^+}$ est
$(\overline{s}^{r_\sigma})_{\sigma\in G^+}$, o ${\overline{s}}\in \F_l^*/{\F_l^*}^q$ est
l'image canonique de $s'$. Par le Sous-Lemme 14, ${\euf L}/{\euf L}^q$ est isomorphe ˆ  $\F_q[G^+]$ et
cet isomorphisme envoie $(\overline{s}^{r_\sigma})_{\sigma\in G^+}$ sur 
$\sum_{\sigma\in G^+}\overline{r}_\sigma\sigma$.
\goodbreak
D'autre part, $u$ est tel que $\overline{u}$ est un gŽnŽrateur de $U_p^+/{U_p^+}^q$ et nous savons
que l'image de $U_p^+/{U_p^+}^q$ dans $\Z[G^+]$ est l'idŽal d'augmentation ${\euf A}$ (cf. Sous-Lemme
14), donc l'image de $u$ dans  $\Z[G^+]$ correspond ˆ un gŽnŽrateur $\varphi$ de ${\euf A}$. On a ainsi

$$\sum_{\sigma\in G^+}\overline{r}_\sigma\sigma=\varphi\cdot\overline{\theta},$$

o $\overline{\theta}$ est l'image de $\theta$ dans  $\F_q[G^+]$. Ou encore
$$\varphi\cdot \theta\equiv \sum_{\sigma\in G^+}r_\sigma\sigma\pmod{q\cdot\Z[G^+]}.\eqno{(vi)}$$
 Or nous savons que $\F_q[G^+]={\euf
A}\oplus(s(G^+))$, avec rappelons-le $s(G^+)=\sum_{\sigma\in G^+}\sigma$ (cf. Lemme 7) et que
$(s(G^+))=\F_q\cdot s(G^+)$ (relation (i) du Lemme 7). Donc il existe $\gamma\in\Z[G^+]$ et $m\in \Z$
tels que
$$\gamma\cdot \varphi+m\cdot\sum_{\sigma\in G^+}\sigma\equiv 1\pmod {q\cdot\Z[G^+]}.$$ 

En multipliant par $\theta$, on en dŽduit que 

$$\theta\equiv \gamma\cdot\varphi\cdot \theta+m\cdot\theta\cdot\sum_{\sigma\in G^+}\sigma\buildrel
(vi)\over
\equiv
\gamma\cdot \sum_{\sigma\in G^+}r_\sigma\sigma+m\cdot\theta\cdot\sum_{\sigma\in G^+}\sigma\pmod
{q\cdot\Z[G^+]}.$$

Il existe donc $\rho\in \Z[G^+]$ tel que $\theta=\gamma\cdot \sum_{\sigma\in
G^+}r_\sigma\sigma+m\cdot\theta\cdot\sum_{\sigma\in G^+}\sigma+q\rho$. Ainsi,

$${\cal C}^\theta=({\cal C}^{\Sigma r_\sigma\sigma})^\gamma\cdot ({\cal C}^{
\Sigma \sigma})^{m\theta}\cdot ({\cal C}^{\rho})^q$$

Mais, $({\cal C}^{\Sigma r_\sigma\sigma})^\gamma\in {{\cal CL}^{pl}}^q $ (relation $(iv)$). D'autre
part,  $({\cal C}^{ \Sigma \sigma})^{m\theta}$ est une classe principale, car si ${\euf a}$ est un
idŽal de $\cal C$ (par exemple $\lambda$), alors ${\euf a}^{\Sigma \sigma}=N_{F/\Q}({\euf a})\cdot O_F$
(cf. Chapitre 5 (rappel sur les corps de nombres et la thŽorie de Galois)) et  $N_{F/\Q}({\euf a})$ est
un idŽal de $\Z$ qui est principal. Et finalement, $({\cal C}^{\rho})^q\in {{\cal CL}^{pl}}^q$. On a
prouvŽ que ${\cal C}^\theta\in {{\cal CL}^{pl}}^q$, ce qu'il fallait dŽmontrer.\qed

\bigskip

\vfill\eject

\font\para=cmbx12 at 18pt
\pageno=70
\parindent0pt
\centerline{\para Bibliographie}
\vskip2cm
Les livres suivants contiennent des dŽmonstrations de rŽsultats acceptŽ sans preuve dans notre texte.
\medskip
[Ati] : M. F. ATIYAH and I.G. MACDONALD, {\it Introduction to Commutative Algebra},  Addison-Wesley Publishing
Company,  1969.
\medskip
[Edw] : H.M. EDWARDS, {\it Fermat's Last Theorem: A Genetic Introduction to Algebraic Number Theory} 
(Graduate Texts in Mathematics), Springer, 1977.
\medskip
[Jac 1] : N. JACOBSON, {\it Basic Algebra 1}, Second Edition. New York, W.H. Freeman, 1989.
\medskip
[Jac 2] : N. JACOBSON, {\it Basic Algebra 2}, Second Edition. New York, W.H. Freeman, 1989.
\medskip
[Lang1] : S. LANG, {\it Algebra}, Addison-Wesley Publishing Company, 1993.
\medskip
[Lang2] : S. LANG, {\it Algebraic Number Theory}, Addison-Wesley Publishing Company, 1970. 
\medskip
[Mac] : G. W. MACKEY, {\it Lectures of the theory of functions of a complex variable}, D. van Nostrand,
1967.
\medskip
[Mar] : D. MARCUS, {\it Number Fields}, Springer, 1977.
\medskip
[Nar] : W. NARKIEWICZ, {\it Elementary and Analitic Theory of Algebraic Numbers}, Springer, 1990.
\medskip
[Ru] : W. RUDIN,    {\it Principles of Mathematical Analysis}, International Series
in Pure \& Applied Mathematics, McGraw-Hill, 1964.
\medskip
[Sam] : P. SAMUEL, {\it ThŽorie algŽbrique des nombres}, Hermann, 1971.
\medskip
[Se] : J.-P. SERRE,  {\it Cours d'arithmŽtique}, Collection SUP No.2, Presses Universitaires de France,
Paris 1970.
\medskip
[Sier] : W. SIERPINSKI, {\it Elementary Theory of Numbers}, North-Holland, 1988.
\medskip
[Was] : L.C. WASHINGTON, {\it Introduction to cyclotomic fields}, Springer, second edition, 1997.

\medskip

Les livres (ou articles) suivants ont ŽtŽ consultŽs pour l'Žlaboration de notre texte 
\medskip
[Bil] : Y. BILU, {\it Catalan's conjecture (after Mih$\breve{\rm a}$ilescu)}, SŽm. Bourbaki, 2002.
\medskip
[Coh] : P.M. COHN, {\it Basic Algebra}, Springer, 1989.
\medskip

[Lem] : F. LEMMERMEYER,  {\it Reciprocity Laws: Their Evolution from Euler to Artin}. Springer, 2000.
\medskip
[Mih] : P. MIH$\check {\rm A}$ILESCU {\it A class number free criterion for catalan's conjecture}, Journal
of Number theory 99, 2003.
\medskip
[Rib] : P. RIBENBOIM, {\it Catalan's Conjecture : Are 8 and 9 the only Consecutive Powers?}, Academic Press,
1994.
\medskip
[Sch] : R. SCHOOF, {\it Catalan's Conjecture}, www.mat.uniroma2.it/\~{}schoof/catalan.pdf, 2003.

\end